\documentclass[11pt, openright]{report}
\usepackage{amsfonts, amsmath, amssymb, graphics, color}
\usepackage{arev}
\usepackage[a4paper]{geometry}
\usepackage{nattoparxiv}
\usepackage{latexsym}
\usepackage{graphicx}
\usepackage{url}
\usepackage{adjustbox}
\usepackage{authblk}
\usepackage{fancyhdr}
\usepackage{enumitem}
\usepackage{footmisc}
\makeatletter
\renewcommand{\@pnumwidth}{2.25em}
\renewcommand{\@tocrmarg}{3.25em}
\makeatother

\renewcommand{\chaptermark}[1]%
                 {\markboth{#1}{}}
\renewcommand{\sectionmark}[1]%
                  {\markright{#1}}
\begin{document}

\vspace*{-.5cm}
\thispagestyle{empty}
\enlargethispage{-.5cm}
\begin{center}
\scalebox{.3}{{\fontsize{144}{180}\selectfont Natural Topology}}
\newlyne
\vspace*{1.5cm}
\scalebox{.65}{\includegraphics{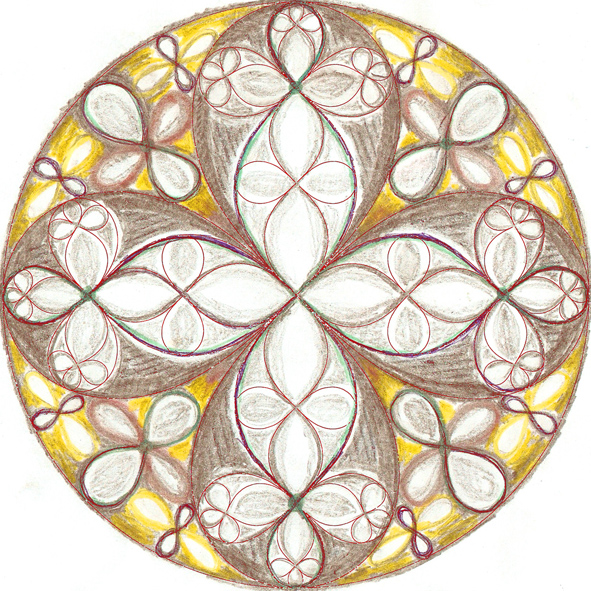}}

\alin\newlyne\vspace*{.8cm}
\scalebox{.95}{{\LARGE Frank Waaldijk}}\raisebox{.2em}{\makebox[0pt][l]{$^{\mbox{\scriptsize\,\maltese}}$}}

\alin

\scalebox{.85}{\Large ---with great support from Wim Couwenberg\,$^{\mbox{\fontfamily{cmr}\selectfont\large\ddag}}$} 

\vfill

\scriptsize \maltese\normalsize\ www.fwaaldijk.nl/mathematics.html
\newlyne
{\fontfamily{cmr}\selectfont\ddag}\ http://members.chello.nl/{\ensuremath{\sim}}w.couwenberg
\end{center}

\clearpage
\thispagestyle{empty}

\normalsize
\vspace*{-1.3cm}\section*{Preface to the second edition}
\vspace*{-.5\baselineskip}

\subsection*{} In the second edition, we have rectified some omissions and minor errors from the first edition. Notably the composition of natural morphisms has now been properly detailed, as well as the definition of (in)finite-product spaces. The bibliography has been updated (but remains quite incomplete). 
We changed the names `path morphism' and `path space' to `trail morphism' and `trail space', because the term `path space' already has a well-used meaning in general topology.
\parr
Also, we have strengthened the part of applied mathematics (the \appl\ perspective). We give more detailed representations of complete metric spaces, and show that natural morphisms are efficient and ubiquitous. We link the theory of star-finite metric developments to efficient computing with morphisms. We hope that this second edition thus provides a unified framework for a smooth transition from theoretical (constructive) topology to applied mathematics.
\parr
For better readability we have changed the typography. The Computer Modern fonts have been replaced by the Arev Sans fonts. This was no small operation (since most of the symbol-with-sub/superscript configurations had to be redesigned) but worthwhile, we believe. It would be nice if more fonts become available for \LaTeX, the choice at this moment is still very limited.

\parr

\texttt{(the author, 14 October 2012)}

\vfill
\parbox{\textwidth}{%
Copyright {\fontfamily{cmr}\selectfont\copyright}\ Frank Arjan Waaldijk, 2011, 2012
\parr
Published by the Brouwer Society, Nijmegen, the Netherlands
\newlyne
All rights reserved

\parr
First edition, July 2011
\newlyne
Second edition, October 2012
\newlyne
\parr
Cover drawing \emph{ocho infinito xxiii} by the author
\newlyne
(`ocho infinito' in co-design with Wim Couwenberg)

}

\clearpage
\thispagestyle{empty}
\normalsize
\vspace*{-1.3cm}\section*{Summary}
\vspace*{-.5\baselineskip}
\subsection*{} We develop a simple framework called `natural topology', which can serve as a theoretical and applicable basis for dealing with real-world phenomena. Natural topology is tailored to make pointwise and pointfree notions go together naturally. As a constructive theory in \bish, it gives a classical mathematician a faithful idea of  important concepts and results in intuitionism.

\parr
Natural topology is well-suited for practical and computational purposes. We give several examples relevant for applied mathematics, such as the decision-support system \hwkin, and various real-number representations.

\parr 
We compare classical mathematics (\class), intuitionistic mathematics (\intu), recursive 
mathematics (\russ), Bishop-style mathematics (\bish) and formal topology,
aiming to reduce the mutual differences to their essence. To do so, our mathematical foundation must be precise and simple. There are links with physics, regarding the topological character of our physical universe.

\parr
Any natural space is isomorphic to a quotient space of Baire space, which therefore is universal. We develop an elegant and concise `genetic induction' scheme, and prove its equivalence on natural spaces to a formal-topological induction style. The inductive Heine-Borel property holds for `compact' or `fanlike' natural subspaces, including the real interval $[\alpha,\beta]$. Inductive morphisms respect this Heine-Borel property, inversely. This partly solves the continuous-function problem for \bish, yet pointwise problems persist in the absence of Brouwer's Thesis.  

\parr
By inductivizing the definitions, a direct correspondence with \intu\ is obtained which allows for a translation of many intuitionistic results into \bish. We thus prove a constructive star-finitary metrization theorem which parallels the classical metrization theorem for strongly paracompact spaces. We also obtain non-metrizable Silva spaces, in infinite-dimensional topology. Natural topology gives a solid basis, we think, for further constructive study of topological lattice theory, algebraic topology and infinite-dimensional topology. 
\parr
The final section reconsiders the question of which mathematics to choose for physics. Compactness issues also play a role here, since the question `can Nature produce a non-recursive sequence?' finds a negative answer in \CTphys. \CTphys, if true, would seem at first glance to point to \russ\ as the mathematics of choice for physics.  To discuss this issue, we wax more philosophical. We also present a simple model of \intu\ in \russ, in the two-player game \LIfE.

\clearpage

\markboth{Contents}{}

\setcounter{tocdepth}{1}
\thispagestyle{empty}
\textheight=23.5cm

\vspace*{-1.3cm}
\section*{Contents}

\vspace*{-1cm}

\tableofcontents
\textheight=23cm

\setcounter{chapter}{-1}

\chepter{Chapter zero}{Introduction}{In this chapter we introduce the subject of natural topology, starting from the natural sciences. The basis for observations and measurements in science seems inescapably to be one of ever-increasing refinement. We do not obtain finished real numbers, but only finite approximations of ever-increasing exactitude.
\parr
This type of `constructive' considerations has become increasingly important with the advent of the computer. The translation of theoretical classical mathematics to applied mathematics is often problematic. Constructive mathematics incorporates finite approximations in its theory, thus providing a smooth theoretical framework for applied mathematics. 
\parr
Also important: natural topology gives a mathematical model of the real-world topology that we encounter when measuring. This we believe to be relevant for physics.
\parr
Last but not least, we believe natural topology to be relevant for the foundations of (constructive) mathematics. 
}

\sectionnb{Introduction}

\sbsc{Background and motivation of this paper}\label{introback}
For a historical background and motivation of this paper, we refer the reader to the appendix, section \ref{background}. Section \ref{examples}\ of the appendix holds some nice mathematical examples which should help clarify our approach. For readability, most proofs are given in the appendix in section \ref{proofs}. Constructive axioms and concepts are given and discussed in section \ref{axioms}, additional remarks can be found in section \ref{addrem}, and the bibliography is in section \ref{bibliography}.

\parr

Summarizing: don't skip the appendix!

\sbsc{Introduction to natural topology}\label{introintro}
Imagine an engineer taking measurements of some natural physical phenomenon. With ever-increasing precision steps she arrives at ever more precise approximations of certain real numbers. In this process she may come across two measurements which at the outset could still indicate the same real number, yet, when more precision is attained, are seen to be really apart.

\parr
The interesting thing about this description lies in the hidden meaning of the word `real number'. Usually this meaning is taken for granted, with an intuitive image of the real line as `foundation'. But in mathematics ---the precision language of science--- the real numbers are commonly defined as equivalence classes of Cauchy-sequences of rational numbers. Then, later, a metric topology can be defined where the basic open sets are the open rational intervals. This topology is optional, as a system the real numbers are usually viewed to exist `on their own'. 

\parr
Compared to the situation that the engineer finds herself in, the above mathematical approach is exactly the other way round. For the engineer first and only encounters a finite number of shrinking rational intervals (the measurements), and then regards these finite measurements as an approximation to a real number (getting better all the time if one can apply more and more precision). 
\enlargethispage{1mm}

\parr
From a topological view, the engineer comes across the topology of rational intervals before ever seeing a real number. It turns out that many problems of translating theoretical mathematics into practical applications hinge around this reversal of approach to real numbers and their natural topology. So let us go into this matter a bit more. 

\sbsc{Classical mathematics}\label{topinapp}
In theoretical mathematics, theorems about the real numbers are often proved in a way which `disregards' the topology. For instance, consider a theorem asserting $\all x\inn\R\driss y\inn \R[D(x,y)]$. Consider a reasonably algorithmic proof even which has the following form: for $x\smlr 0$ do A and for $x\geqq 0$ do B. Looking at it from a topological standpoint, one could say such a proof is `discontinuous' in $0$. When an engineer looks to implement the theorem say on a computer fed by real-time data, the problem arises immediately that for real-time data $x$ the distinction $x\smlr 0$ or $x\geqq 0$ cannot always be made. On data hovering around $0$, the `method' supplied by the theoretical proof might lead to a non-terminating program. 

\sbsc{Constructive mathematics}\label{constructive1}
These types of problems have partly motivated the development of \deff{constructive mathematics}, especially of course since the advent of the computer. Constructive mathematics is a branch of mathematics in which theorems are proved in such a way that the translation of the proof to a working program should be immediate (however in general no claim to efficiency is made). 

\parr 
The surprise for mathematicians and engineers alike is that many theorems from `classical' mathematics can be shown to be `non-constructive', meaning that they can never be translated into a working program. In essence, then, these theorems are simply untrue in constructive mathematics. This questioning of hitherto `solid' theory has led to quite some reluctance amongst theoretical mathematicians to adopt a `constructive' outlook. 

\parr
The downside of this reluctance has been that applicability issues are left to ad-hoc solutioneering, whereas a clearer and more efficient mathematical approach is possible. We hope to give part of such an approach in this paper. Amongst other results, this approach yields a simple topological foundation to some practical issues arising from different representation methods of the real numbers (\mbox{e.g.}\ Cauchy-sequences, decimal notation floating points, binary notation floating points, interval arithmetic and others). 

\sbsc{Our framework: pointwise as well as pointfree}\label{topincon}
The topological framework which we present should also be of theoretical interest, both classically and constructively. For instance, in our (classically valid) framework we present a pathwise connected metric space, which is not arcwise connected. (How this can be? Read and see.) For constructive mathematics, the development of topology has been tackled in different ways. Our approach seems a simple and elegant alternative to the avenues explored so far. 

\parr

For large parts of mathematics the concept of `points' seems natural and elegant. Therefore we focus on this concept. Still, we will show that points arise naturally from 'pointfree' topological constructions.  We believe that our `pointfree' machinery is simple, compared to the framework of formal topology and pointfree topology. This simplicity might look restrictive in the sense that our machinery leads us only a little further than `separable topological T$_1$-spaces'. On the other hand, that is a vast class of spaces. Bishop (the founder of \bish) even went as far as saying that non-separable spaces are a form of pseudogenerality which is to be avoided in constructive math. We hope that by keeping things simple, we can explain the relevance of constructive topology to the `working class' classical/applied mathematician.

\parr

From the foundational perspective, another advantage of keeping things simple is that axiomatic and conceptual assumptions become clear. These assumptions also reflect on physics. As an example we like to state already here that topological compactness of the unit real interval $[0,1]$ turns out to be an independent axiom. There is a perfectly acceptable and beautiful model of the real numbers in which $[0,1]$ is not topologically compact (only very trivial spaces are topologically compact in this model). To our knowledge no one has put forth a convincing argument why reality is not better modeled by this non-compact real model than by the `standard' real model. But the present monograph gives a handhold for the discussion, we believe. In fact we hope that the monograph can serve the foundations of (constructive) mathematics in general.

\parr

Topological spaces always exist in conjunction with continuous mappings between them. 
We define different types of such mappings, to deal with lattice and tree structures which arise naturally from topological investigations. We think that our choices are suited for both theoretical and practical (computational) purposes.

\parr
In order to achieve a constructively valid framework, our logic and our proof methods are constructive. This will not impact greatly on our presentation, which seems quite natural if one keeps in mind that all results should be implementable on a computer. For a precise axiomatic account the reader can read the appendix \ref{axioms}.

\sbsc{Whom it may concern}\label{perspectives}
In the light of the above, we think this paper is of interest from four different perspectives: applied mathematics \& computer science (\appl), general mathematics (\gnrl), constructive foundations of mathematics (\cnst) and foundations of physics (\phys).

\chepter{Chapter one}{Natural Topology}{In this chapter we give the definition of `natural space', starting with the topology and obtaining the points in the process. The natural real numbers are a prime example. Natural morphisms between natural spaces are defined, and shown to be continuous. Conversely, continuous functions going to `basic neighborhood spaces' can be represented by morphisms. Still, for \class\ the equivalence structure determined by isomorphisms is finer than the equivalence structure determined by homeomorphisms. Natural topology is seen to resemble intuitionistic topology.
\parr
Natural Baire space and Cantor space are defined. We show that the class of natural spaces is large, containing (representations of) every complete separable metric space.
\parr
From the \appl\ perspective we discuss the natural topology of binary, ternary and decimal reals. We also look at the well-known Cantor function, and examine the line-calling decision-support system \hwki\ which is used in professional tennis.
}

\sectionnb{Basic definitions and the natural reals}\label{natbasdef}

\sbsc{Topology first, points later}\label{mathintro1}
From the previous introduction our mathematical challenge becomes clear. Namely how to define the real numbers - and more generally a (separable) topological space -- \emph{starting}\ with the topology, and obtaining the points of the topological space in the process. 

\sbsc{Dots and points}\label{mathintro2}
We turn to an intuitive picture of the situation that our previously defined engineer finds herself in. This picture briefly runs as follows. Our engineer in fact encounters only `dots' or `specks', which for the sake of our mathematical argument we think of as being arbitrarily refinable. Two `dots' defined by different processes might at some approximation be seen to definitely lie apart, in which case they represent different (real) numbers. But if the dots are still overlapping at some approximative state, then our engineer cannot tell whether the dots represent different numbers or not.
\parr
So the information of dots lying apart gives more tractability than the information of dots overlapping. The first situation allows a definite conclusion at a finite state (two different real numbers) whereas the second situation still hovers around two possible conclusions. Therefore we will define our spaces using the apartness properties of dots, more than overlap properties. 
   
\parr

In the above intuitive picture, the dots play the central role, and the real numbers arise only as an idealization. Namely, a real number arises as the intersection of an infinite `ever-shrinking' sequence of dots. `Ever-shrinking' can be defined in terms of the apartness of dots, but this turns out to be less convenient than introducing a second notion regarding dots, namely: `being a refinement of', which behaves like a partial order \leqct\ on the countable collection of dots. In accordance with our intuitive picture, when $a\leqc b$, then $a$ is a `refinement' of $b$ and represents a `smaller' dot, contained in $b$.\footnote{Still, we should keep in mind that basic dots arise in the course of a process of measuring points. We will see that when creating pointwise mappings -or in other words, when looking at transformations- it can be meaningful to distinguish between the unit interval $[0,1]$ as a refinement of $[0,2]$ and the `same' unit interval $[0,1]$ as a refinement of $[-1,1]$.}

\sbsc{Pre-natural spaces}\label{nattopdef1}
From the previous introduction we distill the basic mathematical setting: we have a countable set \Vt\ of \deff{basic dots}\ of a \deff{natural topological space}\ \Vnatt\ which we build with a number of definitions in this section. Along with the definitions we give some explanations and examples.

\defi
A pre-natural space is a triple \Vprenatt\  where \Vt\ is a countable\footnote{A set $S$ is \deff{countable} iff there is a bijection from \Nt\ to $S$, and \deff{enumerable} iff there is a surjection from \Nt\ to $S$.}\  set of \deff{basic dots}\ and \aprt\ and \leqct\ are binary relations on \Vt, satisfying the properties following below. Here \aprt\ is a \deff{pre-apartness}\ relation (expressing that two dots lie apart) and \leqct\ is a \deff{refinement}\ relation (expressing that one dot is a refinement of the other, and therefore contained in the other). 

\be

\item[(i)] The relations \aprt\ and \leqct\ are \deff{decidable}\ on the basic dots.\footnote{This means we have a finite procedure to decide the relation. Consider e.g. two rational intervals $[a,b]$ and $[c,d]$, we can decide whether these intervals lie apart or not. We can also decide whether one is a refinement of the other, which in this case is the same as being contained in the other, or not.}
\item[(ii)] For all $a, b \inn V$: $a \aprt b$ (`\aat\ is apart from \bbt') if and only if $b \aprt a$. Pre-apartness is symmetric.
\item[(iii)] For all $a \inn V$: $\neg (a \aprt a)$.  Pre-apartness is antireflexive.
\item[(iv)] For all $a, b, c \inn V$: if $a \leqc b$ (`\aat\ refines \bbt') then $c \aprt b$ implies $c \aprt a$. Pre-apartness is \geqc-monotone.
\item[(v)] The relation \leqct\ is a partial order, so for all $ a, b, c\inn V$: $a\leqc a$ and if $a\leqc b\leqc c$ then $a\leqc c$, and if $a\leqc b \leqc a$ then $a=b$. Refinement is reflexive, transitive and antisymmetric. 

\ee
For basic dots we write $a\touch b$ (`$a$ touches $b$') iff $\neg(a\aprt b)$. Then $\touch$ is the decidable complement of $\aprt$.
\edefi

\rem For the motivating example of the real numbers, the basic dots can be thought of as the rational intervals.\footnote{Using open or closed intervals both yields the same structure, but working with closed intervals fits the intuitive picture better.}\ Two rational intervals $[a,b]$ and $[c,d]$ are said to be \deff{apart}, notation $[a,b]\aprt[c,d]$, iff either $d\smlr a$ or $b\smlr c$. $[a,b]$ \deff{refines} $[c,d]$, notation $[a,b]\leqc[c,d]$, iff $c\leqq a$ and $b\leqq d$. (Also see \ref{defrealnat}\ and appendix \ref{containrefine})
\erem

\sbsc{Points arise from shrinking sequences}\label{defpoints}
Of course one idea is to turn to infinite shrinking sequences (of dots), in order to arrive at points. Looking at our example of rational intervals we see that we need to impose a `sufficient shrinking' condition, otherwise the infinite intersection may contain a whole interval rather than just a point. For an infinite shrinking sequence $\alpha=r_0, r_1, \ldots $ of closed rational intervals ($r_{m+1} \leqc r_{m}$ for all indices $m$) to represent a real number, \alfat\ must `choose' between each pair of apart rational intervals $[a,b]\aprt[c,d]$. By which we mean: for each such pair $[a,b]\aprt[c,d]$, there is an index $m$ such that $r_m \aprt [a,b]$ or $r_m \aprt [c,d]$. (We leave it to the reader to verify that this is indeed equivalent to saying that the infinite intersection of $(r_m)_{\minn}$ contains just one real number.). 

\parr 

The elegance of this approach is that for an infinite shrinking sequence of dots, the property of `being a point' can be expressed by an enumerable condition of pre-apartness. There is no need to talk of `convergence rate' or `Cauchy-sequence', which both presuppose some metric concept. 
To define points, we can simply study the real numbers and transfer certain of their nice properties to our general setting.

\defi
A \deff{point}\ on the pre-natural space \Vprenatt\ is an infinite sequence $p\iz p_0, p_1, p_2 \ldots $ of elements of $V$ that satisfies:

\be
\item[(i)] for all indices $n$ we have: $p_{n+1} \leqc p_{n}$ and there is an index $m$ with $p_m\precc p_n$\hspace*{-0.15pt}.
\item[(ii)] If $a, b \inn V$ and $a \aprt b$ then there is an index $m$ such
that $p_m \aprt a$ or $p_m \aprt b$.
\ee
 
Note that any infinite subsequence of $p$ is itself a point (equivalent to $p$ in the natural sense to be defined).  The set of all points on \Vprenatt\  is denoted by \Vcalt.  
\edefi

Since points are infinite sequences, the set \Vcalt\ is generally not enumerable (but all points in \Vcalt\ could be equivalent).

\sbsc{Apartness on points}\label{defaprtpoints}
The points of our pre-natural space \Vprenatt\ are defined, but clearly we obtain many points which are in some sense equivalent (see our example of rational intervals). The constructive approach to an equivalence relation is to look at its strong opposite, namely an \deff{apartness}\ (see below for the standard properties of an apartness). 

Therefore it is convenient to extend \aprt\ to points in \Vcalt, and also define when points `belong' to dots, in the obvious way:

\defi
For $p=p_0, p_1,\ldots,  q=q_0, q_1,\ldots \inn \Vcal $ and $a \inn V$:

\be
\item[(i)] $a \aprt p$ and $p \aprt a$ iff $a \aprt p_m$ for some index $m$.
\item[(ii)] $p \aprt q$ iff $p_n \aprt q_n$ for some index $n$.
\item[(iii)] $p\equivv q$ iff $\neg(p\aprt q)$.
\item[(iv)] $p \precc a$ iff $p_m\precc a$ for some index $m$.  This relation is
also referred to as `$a$ is a beginning of $p$' or `$p$ begins with $a$' or `$p$ belongs to $a$'. 
\item[(v)] We write \ahatt\ for the set of all points $p$ such that $a$ is a beginning of $p$. Notice that \ahatt\ is not necessarily closed under \equivv. We write \aclost\ for the \equivv-closure of \ahat. 
\ee
\edefi

In terms of complexity, \aprt\ is a $\Sigma_0^1$-property, whereas \equivvt\ is a $\Pi_0^1$-property. This reflects that apartness of two sequences of dots can be seen at some finite stage, but equivalence of two such sequences is an infinite property. Therefore apartness is better suited for constructive and computational purposes (also see appendix \ref{classequiv}).

\parr
To see that \aprt\ is indeed an apartness (and that therefore \equivvt\ is an 
equivalence relation), notice that for all $p=p_0, p_1\ldots$ and $q=q_0, q_1\ldots$ and $r=r_0, r_1\ldots$ in \Vcalt:

\be
\item[(1)]$\neg (p\aprt p)$ (\deff{anti-reflexivity}).
\item[(2)]$p\aprt q$ implies $q\aprt p$ (\deff{symmetry})
\item[(3)] if $p\aprt q$ then $p\aprt r$ or $q\aprt r$ (\deff{co-transitivity}). (If $p\aprt q$ then there is $n$ with $p_n\aprt q_n$, therefore by definition of points (\ref{defpoints}(ii)) we can find an index $m$ with: $r_m\aprt p_m$ or $r_m\aprt q_m$, so $p\aprt r$ or $q\aprt r$).   
\ee

\sbsc{Apartness topology is the natural topology}\label{defnattop}
There is a natural topology on the set of points \Vcalt\ of a pre-natural
space \Vprenat. This topology is expressed in terms of apartness and refinement, we call it the \deff{natural topology}\ and also the \deff{apartness topology}\footnote{See appendix \ref{apartnesshistory}\ for some historical remarks.}, denoted as \Topaprt. \Topaprt\ is the collection of \aprt-open subsets of \Vcalt\ where \aprt-open is defined thus: 

\defi
A set $U\subseteqq \Vcal $ is \deff{\aprt-open} iff for each $x \inn U$ and each $y \inn\Vcal$ we can determine at least one of the following two conditions (they need not be mutually exclusive):

\be
\item[(1)] $y \aprt x$
\item[(2)] there is an index $m$ such that $\hattr{y_m}=\{z\inn\Vcal\midd z \precc y_m\}$ is contained in $U$.
\ee
When the context is clear we simply say `open' instead of `\aprt-open'. 
\edefi

It follows from this definition that an open set is saturated for the equivalence on points (meaning if $U$ is open, $x\inn U$ and $x\equivv y$ then $y\inn U$). We leave this to the reader for easy verification. (This also means that we could replace $\hattr{y_m}$ with $\closr{y_m}$ in (2) above, but in practice this leads to slightly more elaborate proofs). Let us first show that the above indeed defines a topology on \Vcalt:  

\bei
\itemmz{\Topo} Clearly the empty set \emptyyt\ and the entire set \Vcalt\ are open.
\itemmz{\Topto} Let $U, W \subseteqq \Vcal$ be open sets, we wish to show that $U \capp W$ is open. For this suppose $x \inn U \capp W$ and $y \inn \Vcal$, we must show: $y \aprt x$ or there is an index $m$ such that all points beginning with $y_m$ are contained in $U\capp W$. However, since $U$ is open, we can choose \kase{$U$\!(1)} $y \aprt x$ or \kase{$U$\!(2)} there is an index $s$ such that all points beginning with $y_s$ are contained in $U$. Since $W$ is open, we can also choose \kase{$W$\!(1)} $y \aprt x$ or \kase{$W$\!(2)} there is an index $t$ such that all points beginning with $y_t$ are contained in $W$. Combining these two choices, we find: $y \aprt x$ or for $m=\max(s, t)$ all points beginning with $y_m$ are contained in $U\capp W$.
\itemmz{\Topth} Suppose that $U \subseteqq \Vcal$ is a set and each $x \inn U$ has an open neighbourhood $W_x$ such that $x \inn W_x \subseteqq U$. We must show that $U$ is open (this is the constructive formulation of `an arbitrary union of open sets is open'). For this suppose $x \inn U$ and $y \inn \Vcal$, we must show: $y \aprt x$ or there is an index $m$ such that all points beginning with $y_m$ are contained in $U$. Determine an open neighbourhood $W_x$ such that $x \inn W_x \subseteqq U$. Since $W_x$ is open, we find: $y \aprt x$ or there is an index $m$ such that all points beginning with $y_m$ are contained in $W_x$ and therefore in $U$. We see that $U$ is open.
\eei

\sbsc{Natural spaces}\label{defnatspace}
All the ingredients for our main definition have been prepared. Notice that we did not yet stipulate that each dot should at least contain a point. Also it turns out to be necessary to have a maximal dot, which contains the entire space. These then become the final requirements:

\defi
Let \Vprenat\ be a pre-natural space, with corresponding set of points \Vcalt\ and apartness topology \Topaprt. An element $d$ of $V$ is called a \deff{maximal dot} iff $a\leqc d$ for all $a\in V$. Notice that $V$ has at most one maximal dot\footnote{Actually, it also makes sense to reverse the \leqct-notation, and to consider our maximal dot as being the minimal element, which carries the least information. Then each refinement is `larger' because it carries more information than its predecessor.}, which if existent is denoted \maxdotVt\ or simply \maxdott.  \Vnatt\ is a \deff{natural space}\ iff $V$ has a maximal dot and every $a \inn V$ contains a point.
\edefi

\lem
Let \Vnatt\ be a natural space, with corresponding pre-natural space \Vprenat. Let $a\inn V$ be a basic dot. Then the set $\aprt(\hattr{a})=\{z\inn\Vcal\midd z\aprt a\}$ is open in the natural topology. 
\elem

\crl 
For  $x$ in \Vnatt, the set $\{w\inn\Vcal\midd w\aprt x\}$ is open in the natural topology. So a set containing one point (up to equivalence) is closed, showing that every natural space is T$_1$.
\ecrl

\prf let $x$ be in $\aprt(\hattr{a})$, and let $y$ be in \Vcalt. We need to show one of the following two conditions:

\be
\item[(1)] $y \aprt x$
\item[(2)] there is an index $m$ such that $\hattr{y_m}=\{z\inn\Vcal\midd z \precc y_m\}$ is contained in $\aprt(\hattr{a})$.
\ee
   
Since $x \aprt a$, by definition there is an index $s$ with $x_s \aprt a$. Therefore by definition of points there is an index $m$ with $y_m\aprt x_s$ or $y_m\aprt a$. In the first case we find that $y\aprt x$, in the second we find that $\hattr{y_m}$ is contained in $\aprt(\hattr{a})$.
For the corollary, notice that $\{w\inn\Vcal\midd w\aprt x\}=\bigcup_{\ninn}\{w\inn\Vcal\midd w\aprt x_n\}$.
\eprf

\rem
One might think that the lemma shows that in our approach the dots correspond to `closed' subsets in the topology, but this is not always the case. We can also construct the real numbers as a natural space where the basic dots correspond to open intervals, see the alternative definition in paragraph \ref{defrealnat}. In this monograph, different representations of the `same' natural space are studied also to arrive at computational efficiency (for example, consider floating-point versus interval arithmetic). This is relevant for the \appl\ perspective.
\erem

\sbsc{Other ways to define natural spaces}\label{defnattop2}
There are other ways (than the definitions above) to introduce points and spaces. One such way is to simply look at sequences of basic dots which are not necessarily successive refinements, but which all touch (so none are apart) and which fulfill the same condition of `choosing between each pair of apart dots'. However, the resulting point-space can be easily transformed in an `equivalent' space in which the points are again given by a refinement condition. 

\alin

To see this, note that we can always move to new basic dots which are made up of a finite sequence $a_0, \ldots, a_n$ of `old' basic dots where $a_i$ touches $a_j$ for all $i,j\leqq n$. Then a new dot $b\iz b_0, \ldots, b_m$ refines a new dot $a\iz a_0, \ldots, a_n$ whenever $m\geqq n$ and $b_i\iz a_i$ for all $i\leqq n$.

\parr

The above transition from basic dots to finite sequences of basic dots plays a part in our discussion of continuous mappings later on. Two different points may start out differently yet pass through the same basic dot at some later point in time. As we said earlier, we think that for a nice theory of such mappings one should be able to distinguish between $[0,1]$ as a refinement of $[0,2]$ and $[0,1]$ as a refinement of $[-1,1]$. This can be easily realized by looking at the different finite sequences $a\iz ([0,2], [0,1])$ and $b\iz([-1,1], [0,1])$.

\sbsc{The natural real numbers}\label{defrealnat}
After using the rational intervals as a running example for $V$, we can now formally define the \deff{natural real numbers}\ \Rnatt\ as follows:

\defi
Let $\Rrat\isdef\{[p,q]| p,q \inn \Q | p<q\}\cupp\{(-\infty,\infty)\}$. For two rational intervals $[a,b]$ and $[c,d]$ put $[a,b]\aprtR[c,d]$ iff ($d\smlr a$ or $b\smlr c$) and put $[a,b]\leqcR[c,d]$ iff ($c\leq a$ and $b\leq d$). The maximal dot \maxdotRt\ is obviously $(-\infty,\infty)$. The points on the pre-natural space \Rprenatt\ are called the \deff{natural real numbers}\ (also `natural reals'), the set of natural reals is denoted by \Rnatt. The corresponding natural topology is denoted by \TopRaprt. (Also see the remark later in this paragraph).
\parr
Next, let $\zorrat\isdef \{[p,q]| p,q \inn \Q | 0\leqq p\smlr q\leqq 1\}$, then \zorprenatt\ is a pre-natural space with corresponding natural space $(\zornat, \TopRaprt)$ and maximal dot $\maxdot_{\zor}\iz\zor$. 
\edefi 

\thm
\Rnattopt\ is a natural space which is homeomorphic to the topological space of the real numbers \Rt\ equipped with the usual metric topology.
\ethm  

\prf

We prove this in the appendix (\ref{realshomeo}), it is not difficult. Notice that by `homeomorphism' we mean the usual definition (a continuous function from one space to the other which has a continuous inverse; `continuous' meaning that the inverse image of an open set is itself open). Also notice that we are a bit free here, since for a classical theorist we should first move to the quotient space of equivalence classes.
\eprf

\rem An interesting alternative definition of \Rt\ as a natural space is obtained by changing just very little in the definition. For two rational intervals $[a,b]$ and $[c,d]$ in $\Rrat$ put $[a,b]\aprtRo[c,d]$ iff ($d\leqq a$ or $b\leqq c$) and put $[a,b]\leqcRo[c,d]$ iff ($c\smlr a$ and $b\smlr d$). Then\Rprenatot\ is a pre-natural space, and the corresponding natural space is again homeomorphic to \Rt. But one sees that the basic dots $[a, b]$ now correspond to the \emph{open} real intervals $(a, b)$.\,\footnote{It would therefore be better to denote the basic dots as open rational intervals $(a, b)$ under this definition.}\ \Rprenatot\ corresponds to the definition of the formal reals in formal topology (we believe). However, we think compactness is less wieldy in \Rprenatot, which is one reason to stick with \Rprenatt.\footnote{The two spaces are isomorphic in a sense yet to be defined.}
\erem

From a classical point of view the theorem above however also begs the question: `If what we get are the same old real numbers, then what did we gain?'. To answer this question we turn to a final important element of natural spaces: morphisms between them.

\sectionn{Natural morphisms}\label{natmorf}

\sbsc{Why study morphisms?}\label{whymorphisms}
Insight into natural spaces is gleaned from mappings from one space to another which are structure-preserving to some extent. We define different types of such mappings, calling all of them unimaginatively \emph{natural morphisms}. Each natural morphism defines a continuous function with respect to the natural topology. More surprising, from a classical point of view, is that the structure of natural morphisms gives a finer distinction between natural spaces, than the structure of continuous functions between their corresponding topological spaces. This means that the class of natural morphisms forms an interesting subclass of the class of continuous functions. In \ref{basicneighbor}\ we present a nice condition on a natural space \Vnatt, which when satisfied guarantees in \class, \intu\ and \russ\ that a continuous function from another natural space to \Vnatt\ can be represented by a natural morphism, see \ref{basicneighbor}.

\parr

We will show that there is no \emph{isomorphism} between the natural real numbers and the `natural decimal real numbers', whereas classically these spaces are topologically identical.\footnote{This is interesting for representation issues in computer science.} Another interesting result: the space of `natural decimal real numbers' turns out to be `\pathnat\ connected' but not `\arcnat\ connected'. This is fact the translation of an intuitionistic result. Similarly, many intuitionistic results can be translated to the setting of natural spaces and natural morphisms, providing an alternative classical way to view important parts of intuitionism.

\sbsc{Different representations of the `same' space}\label{diffrep}\hspace*{-0.85pt}
In topology, homeomorphisms play a central role. When two spaces are homeomorphic, one can see them as two different representations of the `same' topological space. Yet there is often an intrinsic interest in these different representations. Consider for example \Rt\ and \Rplust. These are two homeomorphic spaces ($(\R, +)$ and $(\Rplus, \cdot)$ are even isomorphic topological groups), but we often have use for one or the other representation, depending on context.

\parr
  
In order to build an elegant theory and prove its correctness, we will need to look at many different representations of `same' natural spaces. As can be expected, in natural topology `sameness' is induced by a special class of natural morphisms called `isomorphisms'. Every isomorphism induces a homeomorphism, but the converse is not true in \class\ (see the above example of the natural decimal real numbers).

\sbsc{Natural morphisms 1: refinement morphisms}\label{naturalmorphisms}
We will distinguish two types of natural morphisms: \emph{refinement morphisms} (denoted \leqc-morphisms) and \emph{trail morphisms} (denoted \pthh-morphisms). The definition of refinement morphisms should pose no surprises in the light of our previous narrative. 

\parr

When going from one natural space to another, a refinement morphism sends basic dots to basic dots, respecting the apartness and refinement relations, in such a way that `points go to points'. This means that any \leqc-morphism\ is an order morphism with respect to the partial order \leqct.\footnote{Not all order morphisms are refinement morphisms though. Our notation ` \leqc-morphism'  can be slightly misleading in this respect.}\ 

\defi
Let \Vnatot\ and \Wnattot\ be two natural spaces, with corresponding pre-natural spaces \Vprenatot\ and \Wprenattot. Let $f$ be a function from $V$ to $W$. Then $f$ is called a \deff{refinement morphism} (notation: \leqc-morphism) from \Vnatot\ to \Wnattot\ iff for all $a,b\in V$ and all $p=p_0,p_1,\ldots\inn\Vcal$:

\be
\item[(i)]$f(a)\aprtto f(b)$ implies $a\aprto b$.
\item[(ii)]$a\leqco b$ implies $f(a)\leqcto f(b)$
\item[(iii)]$f(p)\isdef\ f(p_0), f(p_1), \ldots $ is in \Wcalt.
\ee

As indicated in (iii) above we will write $f$ also for the induced function from \Vcalt\ to \Wcalt. The reader may check that (ii) follows from (iii). By (i), a \leqc-morphism $f$ from \Vnatot\ to \Wnattot\ respects the apartness/equivalence relations on points, since $f(p)\aprtto f(q)$ implies $p\aprto q$ for $p,q\inn\Vcal$. 
\edefi

\thm
Let $f$ be a \leqc-morphism from \Vnatot\ to \Wnattot. Then $f$ is continuous.
\ethm

\prf
Let $U$ be open in \Wnattot. We must show that $T\iz f^{-1}(U)$ is open in \Vnatot. For this let $x\inn T$, and $y\inn \Vcal$. We must show: $x\aprto y$ or there is an index $m$ such that $\hattr{y_m}\subseteqq T$. Since $f(x)$ is in $U$, we can choose \kase{$U$\!(1)}\ $f(x)\aprtto f(y)$, then $x\aprto y$ by (i) above; or \kase{$U$\!(2)}\ there is an index $m$ such that $\hattr{f(y)_m}\subseteqq U$, then by (iii) and (ii) above: $\hattr{y_m}\subseteqq T$.
\eprf

\sbsc{Refinement morphisms are computationally efficient}
The concept of \leqc-morphisms is simple. They have the added advantage of computational efficiency. With a suited `lean' representation \sigRt\ of the natural real numbers, \leqc-morphisms from \sigRt\ to \sigRt\ resemble interval arithmetic, and match the recommendations in \BauKav\ for efficient exact real arithmetic.

\parr

More generally, a continuous function between two `lean' natural spaces can usually be represented by a \leqc-morphism (see thm.\,\ref{basicneighbor}\ and prp.\,\ref{refvstrail2}). Therefore we are interested, already from the \appl\ perspective, in constructing `lean' representations of natural spaces (called `spraids'). Spraids turn out to be fundamental for the \gnrl, \cnst\ and \phys\ perspectives as well. To understand the complexities and to prove our framework correct, we need to define trail morphisms and trail spaces.

\sbsc{Natural morphisms 2: trail morphisms}\label{pathmorphisms}
For the most general setting of natural spaces and pointwise topology, \leqc-morphisms turn out to be too restrictive. This explains our use for the more involved concept of `trail morphism' (denoted \pthh-morphism). Trail morphisms play a necessary role in establishing nice properties of natural spaces. But once these properties have been established, we will primarily use \leqc-morphisms\ (see the previous paragraph). Where \leqc-morphisms are defined naturally on basic dots, one can see \pthh-morphisms as mappings which are naturally defined on points. 

\parr

Actually, a trail morphism from a natural space \Vnatt\ to another space \Wnattot\ is given by a refinement morphism from the `trail space' associated with \Vnatt, to \Wnattot. To define this trail space, we form new basic dots from finite sequences of `old' basic dots (as described in \ref{defnattop2}). 

\defi
Let \Vnatt\ be a natural space derived from \Vprenatt. Let \ninnt, and let $a\iz a_0\geqc\ldots\geqc a_{n-1}$ be a shrinking sequence of basic dots in $V$. The \deff{\precc-trail of $a$}, notation $\asstr$, is the longest subsequence $a_0\succ\ldots\succ a_s$ of $a$.

Let $p\iz p_0, p_1,\ldots$ be a point in \Vcalt, then we write $\pstr(n)$ for the finite sequence $p_0, \ldots, p_{n-1}$ of basic dots in $V$. Notice that $p_0\geqc\ldots\geqc p_{n-1}$, by definition of points. Write $\psstr(n)$ for the \precc-trail 
of $\pstr(n)$.  A finite sequence $a\iz a_0\succ\ldots\succ a_{n-1}$ of basic dots in $V$ is called a \deff{\precc-trail} from $a_0$ to $a_{n-1}$ of length $n$, or simply a trail from $a_0$ to $a_{n-1}$ in \Vprect. The empty sequence is the unique trail of length $0$, and denoted \maxdotpt. The countable set of trails in \Vprect\ is denoted \Vpath, notice that $\Vpath\iz\{\psstr(n)\midd\ninn, p\inn\Vcal\}$. 

\parr
Let $a\iz a_0, \ldots, a_{n-1}$ and  $b\iz b_0, \ldots, b_{m-1}$ be trails in \Vprect\ such that $a_{n-1}\succ b_0$, then we write $a\starr b$ for the concatenation $a_0,\ldots, a_{n-1}, b_0\ldots b_{m-1}$ which is again a trail and so in \Vpath. (Hereby $a\starr\maxdotp$ is defined and equals $a$.).
\parr
The basic dots of our trail space are the trails in \Vprect. For trails $a\iz a_0, \ldots,$ $a_{n-1}$ and  $b\iz b_0, \ldots, b_{m-1}$ 
we put: $a\leqcstar b$ iff there is a trail $c\inn\Vpath$ in such that $a\iz b\starr c$. We also put $a\aprtstar b$ iff $a_{n-1}\aprt b_{m-1}$. 
The natural space \Vpathnatt\ defined by the pre-natural space \Vpathpret\ is called the \deff{trail space}\ of \Vnatt.
\parr
Finally, a \leqc-morphism $f$ from \Vpathnatt\ to another natural space \Wnattot\ is called a \deff{trail morphism}\ (notation $\pthh$-morphism) from \Vnatt\ to \Wnattot. For a point $p\inn\Vcal$ we write $f(p)$ for the point of \Wcalt\ given by $f(\psstr(0)), f(\psstr(1)), \ldots$.

\edefi

\rem
From the pointwise perspective, one readily sees that \Vpathnatt\ is `just another representation' of \Vnatt. Differences in representation should be filtered out by the concept of `isomorphism'. This is the main reason for introducing trail morphisms, since \Vpathnatt\ is not always \leqc-isomorphic to \Vnatt\ (for an example consider the natural real numbers). In fact refinement morphisms preserve the lattice-order properties of the basic neighborhood system which is chosen for a specific representation. Due to the presence of an apartness/equivalence relation, these order properties are not always relevant since we can freely add or distract equivalent basic dots to our system with different lattice properties, without essentially changing the point space. Also see \ref{disclazycon}\ for extra comments.
\erem

\thm
Let \Vnatt\ and \Vpathnatt, \Vpathpret\ be as in the above definition. Then 
\be
\item[(i)] \Vpathpret\ is a pre-natural space and \Vpathnatt\ is a natural space.
\item[(ii)]
\Vpathnatt\ is homeomorphic to \Vnatt\ as a topological space. A homeomorphism is induced by the \pthh-morphism \idptht\ from \Vnatt\ to \Vpathnatt\ given by $\idpth(p)\iz \psstr(0), \psstr(1), \psstr(2),\ldots\inn\Vcalpth$ for $p\inn\Vcal$ (as a refinement morphism \idptht\ is the identity on \Vpatht, with $\idpth(a)\iz a$ for $a\inn\Vpath$). Its inverse homeomorphism is induced by the \leqc-morphism \idstrt\ from \Vpathnatt\ to \Vnatt\ which is defined by putting $\idstr(\maxdotp)\iz\maxdot$ and $\idstr(a)\iz a_{n}$ for a trail $a\iz a_0, \ldots, a_{n}$ in \Vpatht.
\item[(iii)]
Let $f$ be a \pthh-morphism from \Vnatt\ to \Wnattot. Then $f$ is continuous.
\ee
\ethm

\prf
Ad (i): this is a straightforward checking of the definitions, which we leave to the reader.
\parr
Ad (ii): First we show that \idptht\ is continuous. 

Let $U$ be open in \Vpathnat. We must show that $T={\idpth}^{\!-1}(U)$ is open in \Vnatt. For this let $x\inn T$, and $y\inn \Vcal$. We must show: $x\aprt y$ or there is an index $m$ such that $\hattr{y_m}\subseteqq T$. Since $\idpth(x)$ is in $U$, we can choose \kase{$U$\!(1)}\ $\idpth(x)\aprtstar\,\idpth(y)$, then there is \ninnt\ such that $\xsstr(n)\aprtstar\,\ysstr(n)$ which implies $x_{n-1}\aprt y_{n-1}$ and so $x\aprt y$; or \kase{$U$\!(2)}\ there is an index $m$ such that $\hattr{\ysstr(m+1)}\subseteqq U$. Then one easily sees that $\hattr{y_{m}}\subseteqq T$. For let $z\inn\hattr{y_m}$, then let $z_{y_m}\iz y_0, \ldots, y_m, z_{M}, z_{M+1},\ldots$ be the canonical $z$-equivalent point such that $\overline{z_{y_m}}(m+1)\iz\ystr(m+1)$, where $M$ is the first index for which $z_M\precc y_m$. Trivially $\idpth(z_{y_m})$ is in $U$. Since $U$ is open, $U$ is closed under equivalence, and therefore $\idpth(z)$ is in $U$, which shows that $z$ is in $T$.
\parr
The \leqc-morphism \idstrt\ is continuous by theorem \ref{whymorphisms}. That both morphisms are injective\footnote{Meaning: $x\aprt y$ implies $f(x)\aprt f(y)$.}, and therefore homeomorphisms is trivial.

\parr

Ad (iii): Strictly speaking, $f$ is defined as a \leqc-morphism $f^{\pthh}$ from \Vpathnatt\ to \Wnattot. So as a function from \Vnatt\ to \Wnattot, simply notice that $f\iz f^{\pthh}\circ\idpth$, where $f^{\pthh}$ is continuous by theorem \ref{whymorphisms}\ and \idptht\ is continuous by (ii) above, and so $f$ is continuous as well.
\eprf

If $f$ is a \leqc-morphism from \Vnatot\ to \Wnattot, then $f\circ\idstr$ is by definition a \pthh-morphism from \Vnatot\ to \Wnattot, which is clearly equivalent to $f$ on \Vcalt. Therefore we will consider each \leqc-morphism to be a \pthh-morphism as well. 

\sbsc{Natural morphisms' convention}
The difference between the notions \leqc-morphism and \pthh-morphism is often not relevant, which justifies the following:

\conv
If \Vnatot\ and \Wnattot\ are two natural spaces, and $f$ is a \leqc-morphism or a \pthh-morphism from \Vnatot\ to \Wnattot, where the difference is irrelevant, then we simply say: $f$ is a \deff{natural morphism}\ from \Vnatot\ to \Wnattot, or even more simply: a morphism from \Vnatot\ to \Wnattot. 
Only when the difference is relevant will we specify `refinement morphism' and/ or `trail morphism'. This happens mostly in technical proofs or in the context of computation, since refinement morphisms are generally more efficient.
\econv

\sbsc{Composition of natural morphisms}\label{morfcomp}
Given two \leqc-morphisms $f,g$ from natural spaces \Vcalt\ to \Wcalt\ and \Wcalt\ to \Zcalt\ respectively, to form their composition is unproblematic. We leave it to the reader to verify that putting $h(a)\iz g(f(a))$ for all $a\inn V$ defines a \leqc-morphism $h$ from \Vcalt\ to \Zcalt.

\parr

But if $f, g$ are \pthh-morphisms from natural spaces \Vcalt\ to \Wcalt\ and \Wcalt\ to \Zcalt\ respectively, then how do we form the composition? Here \Vcalt, \Wcalt\ and \Zcalt\ are derived from the pre-natural spaces $V, W, Z$, with (pre-natural) trail spaces $\Vpath, \Wpath, \Zpath$ respectively.

\parr

One should notice that $f$ is defined as a \leqc-morphism from $V^{\wr}$ to $W$, but can be uniquely lifted to a \leqc-morphism $f^{\wr}$ from $V^{\wr}$ to $W^{\wr}$. This is straightforward, for a basic dot $a = a_0, a_1,\ldots, a_{n-1}$ in $V^{\wr}$ we look at $b = f(a_0), f(a_0, a_1), {\ldots},$ $f(a_0, a_1, \ldots, a_{n-1})$ and put $f^{\wr}(a)\iz\bsstr$ (the \precc-trail of $b$) which is a basic dot in $W^{\wr}$ since $f$ is a \leqc-morphism.

\parr

Therefore the composition of $f$ and $g$ is defined to be the composition $g\circ f^{\wr}$, which is a \pthh-morphism from \Vcalt\ to \Zcalt. Since any \leqc-morphism can be thought of as a \pthh-morphism (trivially), the composition of a \leqc-morphism with a \pthh-morphism can be similarly dealt with. The composition of a \pthh-morphism with a \leqc-morphism directly yields a new \pthh-morphism.

\sbsc{Isomorphisms}
We can now define a natural parallel to the topological idea of `homeomorphism'. We will call this parallel `isomorphism'. Isomorphisms between natural spaces will automatically be homeomorphisms, but classically we can find homeomorphic natural spaces which are non-isomorphic. This shows that our theory enriches \class\ as well.

\defi
Let \Vnatot\ and \Wnattot\ be two natural spaces, let $f$ be a natural morphism from \Vnatot\ to \Wnattot. Then $f$ is called an \deff{isomorphism}\ iff there is a morphism $g$ from \Wnattot\ to \Vnatot\ such that $g(f(x))\equivvo x$ for all $x$ in \Vcalt\ and $f(g(y))\equivvto y$ for all $y$ in \Wcalt. An isomorphism $f$ from \Vnatot\ to \Vnatot\ is called an \deff{automorphism}\ of \Vnatot, and an \deff{identical automorphism}\ iff $f(x)\equivvo x$ for every $x\inn\Vcal$.
\edefi 

\parr

To see whether certain properties of natural spaces are truly `natural', we check if they are preserved under isomorphisms.

\sectionn{Fundamental natural spaces}

\sbsc{Baire space and Cantor space}
Baire space (\NN) is fundamental because it is a universal natural space (meaning that every natural space can be thought of as a quotient space of Baire space). From chapter two on, we exploit this to simplify the theory considerably. Cantor space (\zoN) is likewise a universal `fan' by which we mean a space generated by a partial order \leqct\ which is a finitely branching tree. Cantor space can be seen as a universal compact space. 

\sbsc{The class of natural spaces is large}\label{natlarge}
Many spaces can be represented by a natural space. In other words, the class of natural spaces is large. A non-exhaustive and also repetitive list of spaces which can be represented as a natural space:
\be
\item[$\bullet$] every complete separable metric space
\item[$\bullet$] the (in)finite product of natural spaces
\item[$\bullet$] \Nt, \Rt, \Ct, the complex p-adic numbers \Cpt, \RNt, Baire space, Cantor space, Hilbert space \Hit, every Banach space, the space of locally uniformly continuous functions from \Rt\ to \Rt, many other continuous-function spaces, and Silva spaces (see chapter four).
\ee

Sometimes, classically defined non-separable spaces (for instance function spaces equipped with the sup-norm) can be constructed under a different metric to become separable. Although the topology is then not equivalent, one can still work with the space constructively as well. For this, one sometimes needs to construct a completion first, to refind the original space as a subset of the completion. Thinking things through, we do not really see a constructive way to define `workable' spaces other than by going through some enumerably converging process. In this sense we concur with Brouwer. Brouwer's definition of spreads in essence parallels  the definition of natural spaces. But unlike Brouwer, we are also engaged in achieving computational efficiency (\appl\ perspective), as well as establishing links between \class, \intu, \russ\ and \bish\ (and formal topology).

\parr

An example of a continuous function space which cannot be represented as a natural space is the space of continuous functions from Baire space to itself. We prove this result (copied from \VelThe) using natural morphisms,  in \ref{Funnotnat}. Still, we will see that there is a subset \Fun\ of Baire space \NNt\ such that every $\alpha\inn\Fun$ codes a natural morphism from Baire space to itself, and every natural morphism from Baire space to itself is coded by some $\alpha\inn\Fun$.
 
\sbsc{Basic-open spaces and basic neighborhood spaces}\label{basicneighbor}
Basic dots do not always represent an open set, or even a neighborhood in the apartness topology.\footnote{In contrast to formal topology, where one only works with opens. We believe this to be unwieldy, both computationally and w.r.t.\ compactness issues.}\ Still, so-called `basic-neighborhood spaces' are fundamental, especially in the context of metric spaces. In \class, \intu\ and \russ\ every continuous function from a natural space to a basic neighborhood space \Vnatt\ can be represented by a natural morphism. The idea is to look at basic dots $a$ which are neighborhoods, meaning $\closr{a}$ contains an inhabited open $U$.

\defi
Let \Vnatt\ be a natural space, with corresponding pre-natural space \Vprenatt. Let $a$ be a basic dot, and let $x\inn\hattr{a}$. Then $a$ is called a \deff{basic (open) neighborhood}\ of $x$ iff $\closr{a}$ is a neighborhood of $x$ (resp. $\closr{a}$ is itself open). Now \Vnatt\ is called a \deff{basic-open space}\ iff $\closr{a}$ is open for every $a\inn V$. \Vnatt\ is called a \deff{basic neighborhood space}\ iff \Vnatt\ is isomorphic to a basic-open space. 
\edefi


\rem
`Basic-open space' is not a `natural' property, meaning that it is not necessarily preserved under isomorphisms (see \ref{defrealnat}, where \Rprenatot\ and $(\Rrat,\aprtR,$ $\leqcR)$ are isomorphic, yet only \Rprenatot\ is basic-open). So we `naturalize' the concept `basic-open space' to `basic neighborhood space', which then trivially is preserved under isomorphisms. This `naturalizing' of desirable properties might seem a bit cheap, but is simply practical. Perhaps basic neighborhood spaces are also definable by the following property, but we leave this as a challenge to the reader.
\erem

\prp
If \Vnatt\ is a basic neighborhood space then there is an identical automorphism $f$ of \Vnatt\ such that $f(x)_n$ is a basic neighborhood of $f(x)$ for every $x\inn\Vcal, \ninn$\ (`there is an identical automorphism which constructs every point as a shrinking sequence of basic neighborhoods of that point').
\eprp

\prf
Let (w.l.o.g., see \ref{basneiprf}) $g, h$ respectively be  \leqc-isomorphisms between \Vnatt\ and the basic-open space \Wnattot. Put $f\iz h\circ g$, then $f$ is an identical automorphism. Let $x\inn \Vcal, \ninn$, then $\closr{g(x_n)}$ is open in \Wnattot\ since \Wnattot\ is basic-open. Since $h$ is also a homeomorphism of the topologies, we see that $h(g(x_n))\iz f(x)_n$ is a basic neighborhood of $f(x)$.
\eprf

If \Vnatt\ is a basic neighborhood space derived from \Vprenat, then $V$ contains a neighborhood basis for the natural topology. The converse does not hold in \class: see \ref{excontmorf}\ for a space where in \class\ every point has an equivalent which arises as a shrinking sequence of basic neighborhoods of that point, and yet there is no identical automorphism sending each point to such an equivalent shrinking sequence. This illustrates that for natural spaces the information `$\all x\driss y\ [P(x,y)]$' only becomes effective if we know that there is a morphism $f$ such that `$\all x\ [P(x,f(x))]$'.\footnote{In \intuf, the statement `$\all x\driss y\ [P(x,y)]$' is generally equivalent to `there is a morphism $f$ such that $\all x\ [P(x,f(x))]$'. This is precisely the content of the \intuf-axiom of continuous choice \acoo\ (see \ref{axiomscont}, and $^x$27.1 in \KleVes).}
\parr
The prime example of a basic neighborhood space is a basic-open space where the basic dots represent open sets (then the identity on $V$ is an identical automorphism which constructs every point as a shrinking sequence of basic neighborhoods of that point). We put forward the main theorem, that in \class, \intu\ and \russ\ continuous functions from a natural space to a basic neighborhood space can be represented by a natural morphism. Later we show that every complete metric space has a basic-open representation (and therefore in \class, \intu\ and \russ\ by the corollary below a unique representation (up to isomorphism) as a basic neighborhood space, see \ref{sepmetnat}).

\thm (in \class, \intu\ and \russ)

Let $f$ be a continuous function from a natural space \Vnatot\ to a basic neighborhood space \Wnattot.  Then there is a natural morphism $g$ from \Vnatot\ to \Wnattot\ such that for all $x$ in \Vcalt: $f(x)\equivvto g(x)$.

\ethm 
\enlargethispage{1mm}
\prf
The not so easy proof is given in the appendix (\ref{prfcontmorf}). In \intu\ the existence of $g$ follows already from the information `$\all x\driss y [y\iz f(x)]$' and the fact that every natural space is `spreadlike' (see \ref{Baireuni}).
\eprf

\crl (in \class, \intu\ and \russ)
If \Vnatot\ and \Wnattot\ are two homeomorphic basic neighborhood spaces, then they are isomorphic.
\ecrl

\rem
The theorem suggests that from a \bish\ point of view, the concept of `natural morphism' adequately captures the notion of continuous function (under the usual topological definition). To (partly) capture the metric property `uniformly continuous on compact subspaces' we will define `inductive morphisms' later on. 
\erem

\sbsc{Complete separable metric spaces are natural}\label{sepmetnat}
To see that the class of natural spaces is large enough to merit interest, we point out with a theorem below that every complete separable metric space is homeomorphic to a natural space. Therefore every separable metric space is homeomorphic to a subspace of a natural space. Some key examples of spaces which can be constructed as a natural space are \Nt, \Rt, \Ct, the complex p-adic numbers \Cpt, \RNt, Baire space, Cantor space, Hilbert space \Hit, and every Banach space. We prove slightly more, because of our interest in different representations of complete metric spaces: 

\thm
Every complete separable metric space \xdt\ is homeomorphic to a basic-open space \Vnatt. 
\ethm

\prf
The rough idea is simple: for a separable metric space \xdt\ with dense subset \anninnt, let for each $n,\sinn$ a basic dot be the open sphere $B(a_n,2^{-s})=\{x\inn X\midd d(x, a_n)\smlr 2^{-s}\}$.  Then we have an enumerable set of dots $V$ by taking $V=\{B(a_n, 2^{-s})\midd n,\sinn\}$. The only trouble now is to define \aprt\ and \leqct\ constructively, since in general for $n, m$ and $s, t$ the containment relation $B(a_n, 2^{-s})\subseteqq B(a_m,2^{-t})$ is not decidable.  We leave this technical trouble, which can be resolved using \aczo\ (countable choice), to the appendix \ref{sepmetnatprf}. 
\eprf

\crl 
In \class, \intu\ and \russ\ the following holds:
\be
\item[(i)]
A continuous function $f$ from a natural space \Wnatt\ to a complete metric space \xdt\ can be represented by a morphism from \Wnatt\ to a basic neighborhood space \Vnatt\ homeomorphic to \xdt, by theorem \ref{basicneighbor}. 
\item[(ii)] A representation of a complete metric space as a basic neighborhood space is unique up to isomorphism.
\ee

\clearpage

In \bish\ the following holds:
\be
\item[(iii)] If \xdt\ and \Vnatt\ are as above in the theorem, then we can define a metric $d'$ on \Vnatt\ (see def.\ \ref{metnat}) by defining $d'(x, y)\iz d(h(x), h(y))$ for $x,y\inn\Vcal$ and $h$ a homeomorphism from \Vnatt\ to \xdt. This metric can be obtained as a morphism from $(\Vcal\timez\Vcal, \Topaprt)$ to \Rnatt\ by the construction of \Vnatt. We then see that the apartness topology and the metric $d'$-topology coincide, in other words \Vnatt\ is metrizable.
We conclude: on a well-chosen basic-neighborhood natural representation of a complete metric space, the metric topology coincides with the apartness topology.
\ee
\ecrl

\rem
\be
\item[(i)]
The construction in the proof sketch above merits a closer look, since we do not simply choose each `rational sphere' $B(a_n, q), q\inn\Q$ to be a basic dot. Yet for \Rt\ and its corresponding natural space \Rnatt, choosing all closed rational intervals works fine. We cannot guarantee in the general case \xdt\ however, that by taking $V=\{B(a_n, q)\midd q\inn\Q,\ninn\}$ we end up with a natural space \Vnatt\ which is homeomorphic to \xdt. We do know that for $X=\Cp$, taking $V=\{B(a_n, q)\midd q\inn\Q,\ninn\}$ gives us a \Vnatt\ which contains `more' points than \Cpt. It might prove a nice challenging exercise to the reader to see why this is the case. In the appendix \ref{nonarc}\ we detail this nice example of a non-archimedean metric natural space.

\item[(ii)]
For most applied-computational purposes, a basic neighborhood representation of a complete metric space seems the best option. We believe that for \Rt, the representation \sigRt\ which we define in the following chapter is a good choice for computational purposes also. Our definition of \sigRt\ and \leqc-morphisms matches the recommendations in \BauKav\ for efficient exact real arithmetic.

\item[(iii)] That the metric topology coincides with the apartness topology on (a well-chosen basic-neighborhood representation of) a complete metric space, allows for simplification later on. Our definition of `direct limit' leading to \mbox{e.g.}\ Silva spaces uses only the apartness topology. 
\ee
\erem

\sbsc{Metrizability of natural spaces}\label{metriz1}
From intuitionistic topology, we can retrieve results on the metrizability of natural spaces. With a definition of the notion `star-finitary' which closely resembles the notion `strongly paracompact', we obtain the constructive theorem that every star-finitary natural space is metrizable. 

\parr

Also, we can easily define natural spaces which are non-metrizable. Comparable to ideas from Urysohn (\Urya), in intuitionistic topology one finds spaces with separation properties `T$_1$ but not T$_2$' and `T$_2$ but not T$_3$' (see \WaaThe).  These spaces can be transposed directly to our setting. 

\parr

However, a different class of non-metrizable natural spaces arises when we look at direct limits in infinite-dimensional topology. As an example we show that the space of `eventually vanishing real sequences' (which is the direct limit of the Euclidean spaces $(\R^n)_{n\in\N}$) can be formed as a non-metrizable natural space.

\parr

We postpone the definitions and theorems to chapter four.

\sbsc{Infinite products are natural}\label{infprod}
Another way to see that the class of natural spaces is large is to look at (in)finite products of natural spaces. We postpone the definitions to paragraph \ref{defprod}, because of some technical concerns and extra issues. The basic idea to arrive at the natural product of \texttt{a)} a finite sequence \texttt{b)} an infinite sequence of natural spaces is simple however. For (weak) basic neighborhood spaces the (in)finite natural product space is homeomorphic to the Tychonoff product-topology space.\footnote{See \ref{projtych}. Products of (weak) basic neighborhood spaces are `faithful'. In \classf\ and \intuf\ also products of `star-finitary' spaces are faithful. We know of no unfaithful products.}

\sectionn{Applied math intermezzo:\hspace*{-.47pt} Hawk-Eye,\hspace*{-.47pt} binary \hspace*{-.24pt}and\hspace*{-.24pt} decimal reals}\label{appmath}

\sbsc{Hawk-Eye}\label{hawkeye}
Already we have introduced enough material to discuss an interesting application of mathematics, in the world of professional tennis. In 2006 the multicamera-fed decision-support system \hwki\ was first officially used to give players an opportunity to correct erroneous in/out calls. \hwki\ uses ball-trajectory data from several precision cameras to calculate whether a given ball was \textscc{in}: `inside the line or touching the line' or \textscc{out}: `outside the line'. \hwki\ is now widely accepted, for decisions which can value at over \$100,000. 

\parr

The measurements of the cameras can be seen as the `dots' or `specks' that we used for illustration in our introduction. Software of \hwki\ must in some way run on these dots. The interesting thing is that \hwki\ does not have the feature of a \textscc{let}: `perhaps the ball was in, perhaps the ball was out, so replay the point'. From this and our work so far we immediately derive: 

\clmm
\hwki, irrespective of the precision of the cameras, will systematically call \textscc{out}\ certain balls which are measurably \textscc{in}\ or vice versa.
\eclmm

\parr

The claim is not per se important for tennis. \hwki\ admits to an inaccuracy of 2-3 mm, and under this carpet the above claim can be conveniently swept (still, one sees `sure' decisions where the margin is smaller). \hwki's inaccuracy is usually blamed on inaccuracy of the camera system. But regardless of camera precision we cannot expect to solve the topological problem that there is no natural morphism from the real numbers to a two-point natural space $\{\mbox{\textscc{in}}, \mbox{\textscc{out}}\}$ which takes both values \textscc{in}\ and \textscc{out}. And our recommendation to \hwki\ is to introduce a \textscc{let} feature, see appendix \ref{exhawkeye}\ for a more detailed description.

\sbsc{Binary, ternary and decimal real numbers}\label{binarydecimal}
Mathematically more challenging than morphisms from \Rnat\ to a two-point space are morphisms from \Rnat\ to the (natural) binary real numbers \Rbint\ and decimal real numbers \Rdect. These morphisms reveal the topology behind different representations of the real numbers on a computer, and transitions between these representations. For simplicity we discuss mainly \Rbint, since the situation with \Rdect\ is completely similar. For some purposes also the ternary real numbers \Rtert\ come in handy. 

\defi
$\Rratbin\!\isdef \{\maxdotR\}\cupp\{[\frac{n}{2^{m}}, \frac{n+1}{2^{m}}]\midd n\inn\Z, \minn\}$, $\Rratter\!\isdef \{\maxdotR\}\cupp\{[\frac{n}{3^{m}}, \frac{n+1}{3^{m}}]\!\mid n\inn\Z,\minn\}$ and $\Rratdec\isdef \{\maxdotR\}\cupp[\frac{n}{10^{m}}, \frac{n+1}{10^{m}}]\midd n\inn\Z, \minn\}$. 

Then 
$\Rbin=(\Rratbin, \aprtR, \leqcR)$ is the natural space of the \deff{binary real numbers}. Similarly we form the corresponding natural spaces \Rtert\ and \Rdect\ of the \deff{ternary}\ and
\deff{decimal real numbers}.

\parr
Put $\zorratbin\isdef \{[\frac{n}{2^{m}}, \frac{n+1}{2^{m}}]\midd n,\minn \midd n\smlr 2^{m}\}$, $\zorratter\isdef \{[\frac{n}{3^{m}}, \frac{n+1}{3^{m}}]\mid n,\minn \mid n < 3^{m}\}$ and $\zorratdec\isdef \{[\frac{n}{10^{m}}, \frac{n+1}{10^{m}}]\midd n,\minn \midd n\smlr 10^{m}\}$ to form the corresponding natural spaces \zorbint, \zortert\ and \zordect, each with the same maximal dot $[0,1]$ denoted by $\maxdot_{\zor}$. 

\parr

Notice that as a partial order, $(\Rratbin, \leqcR)$ is a tree. The reader can think of the natural binary reals as the set of those real numbers $x$ that can also be given as a \deff{binary expansion} $x=(\mino)^s\cdott\Sigma_{\ninn} a_n\cdott 2^{-n+m}$, where $s\inn\{0, 1\}$, $m\inn\N$ and $a_n\inn\{0, 1\}$ for all \ninn, such that $m\bygr 0$ implies $a_0\notiz 0$. We call $s$ the \deff{sign}\ and write $s\iz +, -$ for $s\iz 0, 1$ respectively. We call $m$ the \deff{binary point place}. Then the $(a_n)_{\ninn}$ are the \deff{binary digits} in this binary expansion of $x$, and we write $x\iz (s)\,a_0\, a_1 \ldots a_{m}{\rm \bol{.}} a_{m+1} \ldots$. Notice the \deff{binary point}\ that we write between $a_m$ and $a_{m+1}$ to denote the binary point place.

Replacing `binary, 2' with `ternary, 3' and `decimal, 10' respectively, we obtain the similar definitions for \Rtert\ and \Rdect. 
\edefi

Classically every real number $y$ has an equivalent binary expansion, but in computational practice and in constructive mathematics this is not the case (see \mbox{e.g.}\ \GNSWcom\ for a thorough discussion). So with \Rbint, \Rtert\ and \Rdect\ we in practice obtain \emph{different} representations of the real numbers. We wish to shed some light on the natural topology involved in the (im)possible transition from one such representation to another.

\sbsc{Morphisms to and from the binary reals}\label{morfbinreals}
It turns out that a morphism $f$ from \Rnat\ to \Rbint\ which is order preserving ($x\leqqR y$ implies $f(x)\leqqR f(y)$) has to be `locally constant' around the $f$-originals of the rationals $\{\frac{k}{2^{m}}\midd k,\minn\}$. These are the points where the binary expansion has two alternatives (\mbox{e.g.}\ for $1$ both $0.111\ldots\equivv 0+1\cdott 2^{-1}+1\cdott 2^{-2}+1\cdott 2^{-3}+\ldots$ and $1.000\ldots\equivv 1+0\cdott 2^{-1}+0\cdott 2^{-2}+0\cdott 2^{-3}+\ldots$ are binary representations). Since these binary rational numbers lie dense in \Rt, there can be no injective morphism from \Rnatt\ to \Rbint\ (notice that any injective morphism $f$ from \Rnatt\ to \Rnatt\ is either order preserving, or order reversing in which case a similar argument for local constancy obtains). But this does not mean that all morphisms from \Rnatt\ to \Rbint\ are constant. We will derive a non-constant morphism from \Rnatt\ to \Rtert\ from our metrization theorem in \ref{starfinmet}. We can easily turn this into a non-constant morphism from \zort\ to \zorbint, when we realize that \zorbint\ and \zortert\ are isomorphic.

\parr

The well-known \deff{Cantor function}\ \fcant\ (also known as `the devil's staircase') is another example of a non-constant natural morphism from \zort\ to \zorbint. The Cantor function is most easily described as a refinement morphism from \zortert\ to \zorbint, but also can be given as a trail morphism on \zort.\footnote{We leave this latter statement as a non-trivial exercise to the reader though, see \ref{cantorfun}.}

\parr

We now have an example in \class\ of a continuous function between natural spaces which cannot be represented by a morphism. In \class, the identity is a homeomorphism from \Rnatt\ to \Rbint\ (remember that in \class\ we work with the equivalence classes, and that every real number has an equivalent binary representation). But this identity cannot be represented by a natural morphism, as we pointed out above. In the light of theorem \ref{basicneighbor}, the `reason' for this is that \Rbint\ is not a basic neighborhood space, which we can easily verify by looking at the real number $\half$. In fact, in \Rbint, of the basic dots only the maximal dot is a neighborhood of $\half$.

\parr 

We will make the above statements and definitions precise in the appendix \ref{cantorfun}, also showing the equivalence between the reals allowing a binary expansion and the binary reals. We then use the ternary reals to construct the Cantor set \Canzort, and the ContraCantor set, which is a compact subspace \contrCant\ of \zort\ such that: $\dR(\contrCan, \Canzor)\iz 0$ and yet $d(x, y)\bygr 0$ for all recursive $x\inn \contrCan, y\inn\Canzor$. So in \russ\ we have $\dR(x, \Canzor)\bygr 0$ for $\all x\inn\contrCan$, giving us a \russ\ example of two disjoint complete compact spaces with distance $0$.

\rem
One can show with little effort that for $n, \minn$ the $n$-ary and $m$-ary reals are \leqc-isomorphic. However, we believe the $n$-ary reals can only be identically embedded in the $m$-ary reals if there is a $b\geqq 1$ in \Nt\ such that $m$ divides $n^b$ (for an identical embedding $f$ we have $f(x)\equivvR x$ for all $x$). This gives a natural-topological classification of the different $n$-ary representations of real numbers. 
We leave this as an exercise. \erem

\sectionn{Intuitionistic phenomena in natural topology}\label{intphenat}
\sbsc{Pathwise and arcwise connectedness}\label{patharc}
From the previous paragraph we deduce an interesting property of \Rbint\ (and \Rtert, \Rdect): it is a \pathnat\ connected space which is not \arcnat\ connected. For this we must define:

\defi
A natural space \Vnatt\ is called \deff{\pathnat}\ (resp. \deff{\arcnat}) \deff{connected} iff for all $x, y\inn V$ there is a morphism (resp. an injective morphism) $f$ from \zornatt\ to \Vnatt\ such that $f(0)\equivv x$ and $f(1)\equivv y$.
\edefi

\thm
\Rbint\ (as well as \Rtert, \Rdect) is a \pathnat\ connected space which is not \arcnat\ connected.
\ethm

\prf
A detailed constructive proof for \Rtert\ is given in \WaaThe\ in an intuitionistic setting. We indicate the translation to our setting in the appendix \ref{binrealpatharc}. The reader should have no trouble giving a proof using the above subsections.
\eprf

\sbsc{Intuitionistic phenomena arise naturally}\label{intunat}
The previous examples and theorems are a quite faithful mirror of the phenomena studied in intuitionism (\intu). In \intu, an elegant class of natural spaces (called `spreads') is studied. Hereby the main intuitionistic axiom which is not classically valid is the continuity principle \CP\ (and by implication its stronger versions \acoz\ and \acoo, see the appendix \ref{axioms}). A main consequence of \CP\ is that total functions on spreads are always given by morphisms. In natural topology, by directly considering morphisms we create a simple classical mirror of many intuitionistic results. Not surprisingly, these results have direct computational meaning, and are therefore also of relevance for applied mathematics.

\parr

Let's introduce \CP\ here, to see what sets \intu\ apart from \class\ axiomatically.\footnote{Of course, apart from axioms there is also a fundamental conceptual difference between constructive mathematics and classical mathematics, regarding infinity and omniscience, see also \ref{finvsinf}.}\ \CP\ is formulated for Baire space \NN, we show that Baire space is a universal natural space in the next section. For an element $\alpha$ of Baire space \NNt\ we write $\alfstr(n)$ to denote the finite sequence formed by the first $n$ values of $\alpha$.


\xiom{\CP} Let $A$ be a subset of $\NN\times\N$ such that $\all\alpha\inn\NN\driss\ninn\ [(\alpha, n)\inn A]$. Then $\all\alpha\inn\NN\driss m,n\inn\N$ $\all\beta\inn\NN\,[ \alfstr(m)\iz\betastr(m)\rightarrow (\beta, n)\inn A]$.
\exiom  

\alin

The motivation for this axiom is that in \intu\ infinite sequences arise step-by-step, and that among these sequences are also those about which we know -at any given time $m$- nothing more than the  
first $m$ values. We can also form recursive determinate sequences, no problem, but asserting that for \emph{all} $\alpha$ one can find an $n$ such that $(\alpha, n)\inn A$ means that the indeterminate sequences are included. For any such indeterminate sequence $\alpha$ one has to produce the favourable $n$ with $(\alpha, n)\inn A$ at some point in time, say $m$. At this point, we know only the first $m$ values of $\alpha$, which implies that for $\beta\inn\NN$ with $\betastr(m)\iz\alfstr(m)$, we also have $(\beta, n)\inn A$. 

\parr

The other intuitionistic axioms (apart from the strengthenings of \CP) are all valid classically, 
and the above axiom also even makes sense classically in the right setting.\footnote{The author considers the axiom to be simply true in its intended context. He thinks intuitionistic mathematics deserves a prime role in mathematical investigations.}\ We think that such a setting is obtained naturally when considering a two-player game with limited information, see also section \ref{intuinruss}). In this setting we can even prove \CP. Notice that this setting strongly resembles that of our engineer taking measurements from nature. This gives a philosophical explanation for the aptness of \intu\ for physics. 

\parr
However, by considering morphisms on natural spaces a classical mathematician can skip most of this issue and still see intuitionistic phenomena arising naturally and quite faithfully. 

\parr

We have come to the end of this chapter. Its main purpose, apart from giving the basic definitions and properties, is to give some inkling of the relation between constructive topology and applied mathematics (our \appl\ perspective). In the next chapters we will turn predominantly to the other perspectives \gnrl, \cnst\ and \phys.

\chepter{Chapter two}{Baire space is universal}{Natural Baire space is a universal natural space, meaning that every natural space is the image of Baire space under a natural morphism. 
\parr
Its set of basic opens can be pictured as a tree, and this partial-order property provides a fundamental simplification. For many natural spaces, the related concept of `trea' is similarly useful. 
\parr
Through Baire space, a direct link with intuitionism (\intu) can be established. We define natural spreads and spraids (corresponding to trees and treas) as well as fans and fanns. 
This enables us to show, for the \appl\ perspective, that continuous $\R$-to-$\R$-functions can be represented by computationally efficient morphisms.
\parr
Cantor space is a universal fan. Compactness of Cantor space is seen to depend on the axiom \FT\ (Fan Theorem), derived from \BT\ (Brouwer's Thesis), which holds in \class\ and \intu\ but not in \russ. In \russ, Baire space is isomorphic to Cantor space. 
\parr
Since we wish to work in \bish, we do not adopt \BT, and turn to inductive definitions instead.} 

\sectionnb{Introduction to Baire space}\label{introBaire0}

\sbsc{Classical Baire space}\label{introBaire1}
In \class, Baire space is $\NN$ with the usual product topology. Baire space is a universal Polish space (`Polish' meaning `second countable completely metrizable'), that is: every Polish space is the continuous image of Baire space.\footnote{This can often be used to streamline proofs for separable complete metric spaces.} However, in fact Baire space is universal for a certain larger class of spaces. We can see this classically by looking at any equivalence relation \equivvt\ on \NNt, and forming the corresponding quotient topology and quotient space. In this way we trivially see that Baire space is universal for the class of topological spaces homeomorphic to a quotient space of Baire space. In this paper we study this larger class from our natural perspective. For constructive reasons we limit ourselves to quotient topologies derived from a $\Pi^1_0$-equivalence relation. To see that this class is larger than the class of Polish spaces, it suffices to see that some of these quotient spaces are non-metrizable. 

\sbsc{Introduction to natural Baire space}\label{introBaire2}
One goal of this paper is to give a classical mathematician (our perspective \gnrl) an understanding of what constructive and intuitionistic topology `is all about', in such a way that little foundational terminology is necessary. We believe that the (classically valid) setting of natural spaces and natural morphisms between them, is a faithful mirror of Brouwer's basic intuitionistic setting of `spreads' and `spread-functions'. 

\parr

To see that our class of natural spaces is in fact quite large, we will show that it encompasses (representations of) every Polish space. Also, (in)finite products of natural spaces are representable as a natural space. Another result already noted by Brouwer is that natural Baire space $\Bnat=\NNnat$ is a universal natural space, meaning that every natural space \Vnatt\ is the image of natural Baire space under some natural morphism from \Bnatt\ to \Vnatt. 

\parr

To put it differently, from the \gnrl\ perspective we study (subspaces of) Baire space endowed with a quotient topology derived from a $\Pi^1_0$-equivalence relation. BUT: our view is finer than the usual classical topological view. Using morphisms we can distinguish between interesting spaces that are classically homeomorphic (indistinguishable with topological methods). AND: our view and constructive methods can be fruitful in applied math and physics.

\parr

Quite some work has already been done in intuitionistic topology, with Baire space as fundament.\footnote{Brouwer already gave an intuitionistic proof of the Jordan curve theorem, to mention just one historic result.}\ We will be able to retrieve some nice results (\mbox{e.g.}\ regarding metrizability of natural spaces) in a simple way, once we have shown that our setting of natural spaces is indeed a faithful mirror of Brouwer's setting. Since we take a neutral constructive approach, we will not use any specific classical or intuitionistic axioms. However we will freely use the axioms of countable choice \aczo\ and dependent choice \dco, which are generally accepted as constructive (also see \ref{axioms}).

\sbsc{Natural Baire space}\label{natBaire} 
There is much mathematical beauty in Baire space, and its definition as natural space is likewise elegant. Its set of basic dots is \Nstart, the set of all finite sequences of natural numbers (representing the basic clopen sets of Baire space). We use the definition also to relate natural Baire space to usual Baire space.

\defi 
Let \Nstart\ be the set of all finite sequences of natural numbers. For $a\iz\ a_0, \ldots, a_i$, $b\iz\ b_0, \ldots, b_j\inn\Nstar$ the concatenation $a_0,\ldots, a_i, b_0\ldots, b_j$ is denoted by $a\star b$. Define: $b\leqcom a$ iff there is $c$ such that $b=a\star c$. Define: $a\aprtom b$ iff $a\nleqcom b$ and $b\nleqcom a$. 

\parr

Then \Bprenatt\ is a pre-natural space, with the empty sequence as maximal dot, which we also denote \maxdotBt\ or simply \maxdott. Its corresponding natural space \Bnattopt\ we call \deff{natural Baire space}. We also write \NNnatt\ for \Bnatt.

\parr

In addition, let $\alpha\inn\NN$ and \minnt, then we write $\alfstr(m)$ for the finite sequence $\alpha(0), \ldots, \alpha(m-1)$ consisting of the first $m$ values of $\alpha$. Notice that $\alfstr(m)$ is an element of \Nstart, so the sequence $\alfstr\iz \alfstr(0), \alfstr(1), \ldots $ is a point in \Bnat.

Conversely, for a point $p\inn\Bnat$, there is a unique sequence $\alpha\inn\NN$ such that $p\equivvom\alfstr$. We write $p^*$ for this unique $\alpha$, giving that $p\equivvom\overline{p^*}$ for $p\inn\Bnat$ and $\alpha\iz\alfstr^*$ for $\alpha\inn\NN$.\,\footnote{This also somewhat relates to the axiom of extensionality, see \ref{axiomext}.}
\edefi

\thm
\Bnattop\ is homeomorphic with \Bairut.
\ethm

\prf
We leave it to the reader to verify that the function $\alpha\rightarrow\alfstr$ from $\NN$ to \Bnatt\ defined above is a homeomorphism, with inverse $p\rightarrow p^*$ (also defined above).
\eprf 

\sbsc{Natural Cantor space}\label{defcantor}
We first define the notion `natural subspace', since in natural Cantor space we have a prime example.

\defi Let \Vnatt\ be a natural space derived from \Vprenatt. Let \Wt\ be a countable subset of \Vt, then \Wprenatt\ is a pre-natural space, with corresponding set of points \Wcalt. If \Wnatt\ is a natural space (see def.\,\ref{defnatspace}), then we call \Wnatt\ a \deff{natural subspace}\ of \Vnatt\ iff in addition \Wnatt\ as a natural space coincides with \Wnatt\ as a topological subspace of \Vnatt\ (in the subspace topology `$U\subseteq\Wcal$ is open' is defined thus: there is an open $U'\subseteq\Vcal$ such that $U\iz U'\capp\Wcal$).

Let \zostart\ be the set of finite sequences of elements of \zot. Now \deff{natural Cantor space}\ is the natural subspace \Cnattopt\ of natural Baire space formed by the pre-natural space \Cprenatt\ and its set of points \Cnatt.
\edefi

\rem We will see that any natural subspace of natural Baire space is the image of Baire space under a (continuous) morphism. From descriptive set theory, it follows that many topological subspaces of a natural space cannot be represented as a natural subspace. We will give an example of such a space in \ref{Funnotnat}. The notion of `natural subspace' is weaker than the intuitionistic notion of `subspread'. We will define this notion in our context also. Natural Cantor space is homeomorphic to usual Cantor space, and corresponds directly to Brouwer's fan \sigtot. 
\erem

From now on, when the context is clear we will simply say `Baire space' and `Cantor space' and omit the extra word `natural'.

\sectionn{Lattices, trees and spreads}\label{lattrespre}

\sbsc{Lattices and posets of basic dots}\label{latopen}
In topology, the open sets form a lattice structure under the inclusion relation. This structure is often exploited in various ways. One way is to (somewhat) disregard meet and join operations and focus simply on the partial-order properties (of the `poset' of opens). Our basic dots in general need not form a lattice, but their partial-order properties play an important role. This role could even be too restrictive, without trail morphisms. 

\alin

We now go into these partial-order properties in more detail. This requires some attention from the reader, but the rewards are great. Eventually we will derive from these properties some important simplifications which enable us to forge a direct link with intuitionistic topology. Brouwer's simple and elegant concept of a spread then becomes accessible to us as well.

\sbsc{Trees and treas}\label{treas}
The elegance of Baire space can be seen as stemming from the fact that its poset of basic dots \Nleqcomt\ forms a countable tree. That is: for any dot $a=a_0,\ldots, a_{n-1}\inn\Nstar$, there is a \deff{unique} finite trail of immediate successors/predecessors from \maxdotBt\ to $a$. (Therefore any \precc-trail between dots is finite, and also the successor/predecessor relationship is decidable.). We cannot hope to achieve this elegance for any natural space, but we can show that any natural space \Vnatot\ is isomorphic to a natural space \Wnattot\ where \Wleqctot\ equals \Nleqcomt. Or more practical: where \Wleqctot\ is a full subtree of \Nleqcomt, definition follows.

\parr

This means that we could limit ourselves to natural spaces \Vnatt\ where \Vleqct\ is (a full subtree of) \Nleqcomt. This has a strong simplifying effect, which gives much beauty to Brouwer's intuitionism. The only downside is that for many natural spaces, we have to replace our original basic dots with elements of \Nstart, which can in practice be a tedious encoding. Therefore we propose the compromise notion of a `trea'. One can think of a trea as being a tree wherein certain of the branches have been neatly glued together in a number of places. A more precise characterization of a trea: a countable \precc-directed acyclic graph with a maximal element, where for each node there are finitely many immediate-predecessor trails to the maximal element, all of the same length. Another characterization: a countable \leqct-poset with a maximal element where each point has finitely many immediate-predecessor trails to the maximal element, all of the same length.\footnote{We apologize for the lengthy definitions. The concepts are not too difficult, and will provide elegant simplification later on.}

\defi   
Let \Vnatt\ be a natural space, with corresponding \Vprenat,  and let \Wnatt\ with corresponding \Wprenatt\ be a natural subspace of \Vnatt\ (so $W\subseteq V$). Let $a\precc c$ in $V$.  
\be

\item[(i)]
We say that $a$ is a \deff{successor} of $c$ in \Vleqct\ (notation $a\suczV c$, or simply $a\sucz c$ if the context is clear) iff for all $b\inn V$, if $a\precc b\leqc c$ then $b\iz c$. A sequence  $b_0\predcz\ldots\predcz b_n$ in $V$ is called a
\deff{\sucz-trail of length $n$ from $b_0$ to $b_n$ in \Vleqct}. For $b\inn V$ we put $\suczV(b)\isdef\{d\inn V\midd d\sucz b\}$, and simply write $\suczz(b)$ when the context is clear.
\item[(ii)] 
\Vleqct\ is called a \deff{trea} iff for every $a\inn V$ the set $\{b\inn V\midd a\leqc b\}$ of predecessors of $a$ is finite (then the successor relation \suczt\ is decidable, and for every $a\inn V$ there are finitely many \sucz-trails from \maxdott\ to $a$) and in addition there is an integer $\grdd(a)\inn\N$ such that every \sucz-trail from \maxdott\ to $a$ has length $\grdd(a)$.\footnote{This last condition more or less follows from the first. If we only stipulate that for every $a\inn V$ the set $\{b\inn V\midd a\leqc b\}$ is finite, then we can add extra basic dots to $V$ to end up with an isomorphic space in which all \suczV-trails from \maxdott\ to a given $a$ have the same length.} 
\item[(iii)]
Now let \Vleqct\ be arbitrary, where \Wleqct\ is a tree (trea), then we say that \Wleqct\ is a \deff{subtree (subtrea)} of \Vleqct.
\item[(iv)]
Let \Vleqct\ be a tree (trea), and \Wleqct\ a subtree (subtrea). We then call \Wleqct\ a \deff{full subtree (subtrea)} of \Vleqct\ iff $b\suczW d$ implies $b\suczV d$ for all $b, d\inn W$. (Then each \suczW-trail in \Wleqct\ is a \suczV-trail in \Vleqct).

\ee
\edefi

As stated above, we can show that any natural space \Vnatot\ is isomorphic to a natural space \Wnattot\ where \Wleqctot\ is a full subtree of \Nleqcomt. In intuitionism therefore, most notions are defined for full subtrees of \Nleqcomt. But since we wish to incorporate the possibility to avoid encoding schemes, we will define our notions for treas and full subtreas. Treas behave just like trees (and any tree is a trea). Most of the spaces of interest that we mentioned so far (see \ref{natlarge}) have an intuitive representation as a natural space \Vnatt\ where \Vleqct\ is a trea. One reason for this is that the infinite product of a sequence of treas can be naturally built as a trea.\footnote{Using finite products of basic dots from the treas involved, see \ref{defprod}.}\  One might find another reason in the existence of star-finite refinements of per-enumerable open covers of metric spaces, see \WaaThe.

\xam
For the natural real numbers \Rnatt\ we can easily indicate an isomorphic subspace \signatR\ with corresponding pre-natural space \sigprenatRt, where \sigleqcRt\ is a trea: 

\parr

$\sigprenatR\isdef(\{\maxdotR\}\cupp\{[\frac{n}{2^{m}}, \frac{n+2}{2^{m}}]\midd n\inn\Z, \minn\}, \aprtR, \leqcR)$.

\parr

Our examples in \ref{binarydecimal}\ should show why we cannot hope to find an isomorphic subspace $(\Vcal,\TopRaprt)$ where \Vleqct\ is a tree (!).
\exam

\sbsc{Spreads and spraids}\label{spraids}
The previous example illuminates a bridge towards intuitionistic terminology, which we give in the following definition:

\defi
Let \Vnatt\ be a natural space, with corresponding \Vprenat,  and let \Wnatt\ with corresponding \Wprenatt\ be a \deff{decidable} natural subspace of \Vnatt\ (meaning $W$ is a decidable subset of $V$). 
\be
\item[(i)]
We call \Vnatt\ a \deff{spread (spraid)} iff \Vleqct\ is a tree (trea) and each infinite \precc-trail defines a point (see also \ref{sprddefcom}). Then we call \Wnatt\ a \deff{subspread (subspraid)} of \Vnatt\ iff \Wleqct\ is a full subtree (subtrea) of \Vleqct.
\item[(ii)]
We call \Vnatt\ a \deff{Baire spread} iff \Vnatt\ is a subspread of Baire space.
\item[(iii)]
For any space \Vnatt, if \Wnatt\ is a spread (spraid), then we simply call \Wnatt\ a weak subspread (subspraid) of \Vnatt. (If \Vleqct\ contains an infinite \precc-trail between two dots, then we could drop the prefix `weak'.\footnote{If \Vleqct\ contains such an infinite trail, then \Vleqct\ is not a trea and \Vnatt\ is not a spraid. So the condition for spraids and subspraids, that \Wleqct\ is a full subtrea of \Vleqct, becomes void. However we cannot always decide whether \Vleqct\ is a trea or not.}).  
\ee
 
By extension, \Vnatt\ is \deff{spreadlike}\ iff there is an isomorphism between \Vnatt\ and a spread. 
\edefi

\xam
Important basic examples of subspraids are obtained as follows. For \Vnatt\ a spraid and $a$ in $V$, one easily sees that $\Vsuba\iz\{b\inn V\midd b\leqc a\}\iz\{a\}_{\!\leqc}$ determines a subspraid of $V$ if we put its maximal dot as $\maxdot_a\iz a$. These basic subspraids are important later on in defining `genetic induction'.
\exam

\sectionn{Universal natural spaces}\label{uninat}

\sbsc{Baire space is universal}\label{Baireuni}
Baire space is a universal natural space, by which we mean that each natural space can be seen as the image of Baire space under a natural morphism. Brouwer already realized this, and simplified his concepts accordingly.\footnote{Often in difficult language...} In our setting, we adopt a similar simplification, for esthetic reasons and to save paper and energy. We show that every natural space is spreadlike, which on a meta-level gives us a one-on-one correspondence with many important intuitionistic results. 

\thm
Every natural space is spreadlike. In fact, every natural space \Vnatt\ is isomorphic to a spread \Wnatt\ whose tree is $(\Nstar, \leqcom)$.
\ethm

\crl
\be
\item[(i)] Let \Vnatt\ be a natural space, then there is a surjective \leqc-morphism from Baire space to \Vnatt. (`Baire space is a universal spread', `every natural space is the natural image of Baire space', `every natural space is a quotient topology of Baire space'). 
\item[(ii)] If \Vnatt\ is a basic-open space (see definition \ref{basicneighbor}) then \Vnatt\ is isomorphic to a basic-open spread \Wnatt\ whose tree is $(\Nstar, \leqcom)$. 
\ee
\ecrl

\prf
Not trivial, see \ref{Baireuniprf}\ where we give a self-contained proof.
\eprf

From the theorem we reobtain Brouwer's simplification: \emph{without any loss of generality we may assume that a given natural space \Vnatt\ is a spread}. Points in \Vcalt\ can be constructed step-by-step, as infinite trails in the countable tree \Vleqct. Any such tree can of course be embedded as a full subtree in $(\Nstar, \leqcom)$. 
\parr
The corollary gives the equivalent picture that each natural space \Vnatt\ with corresponding pre-natural space \Vprenatt\ is in fact nothing but a pre-apartness \aprtVt\ and a refinement relation \leqcVt\ on \Nstart\ which respect \aprtomt\ and \leqcomt. To define these decidable relations we only have to `pull back' the decidable relations \aprt\ and \leqct\ using the given surjective morphism $f$ thus: for $a, b\inn\Nstar$ put $a\aprtV b$ resp. $a\leqcV b$ iff $f(a)\aprt f(b)$ resp. $f(a)\leqc f(b)$. Then $a\aprtV b$ resp. $a\leqcom b$ implies $a\aprtom b$ resp. $a\leqcV b$.

\parr

An equivalent situation which avoids encoding arises whenever a natural space \Vnatt\ contains a subspraid on which the identity is an isomorphism. Then from the often vast partial-order universe of \Vleqct\ we can restrict ourselves to a subtrea. We give an important example below where the isomorphic subspace is a spraid. We believe this to be the most common setting for natural spaces. In the uncommon case that we cannot find an isomorphic subspace which is a spraid, we can always invoke Brouwer's encoding to find an isomorphic spread.   

\xam
Looking at the natural real numbers \Rnatt, we can easily indicate an isomorphic subspace which is a spraid as in example \ref{treas}. Put 

\parr

$\sigR\isdef \{\maxdotR\}\cupp\{[\frac{n}{2^{m}}, \frac{n+2}{2^{m}}]\midd n\inn\Z, \minn\}$.

\parr

Then \sigprenatRt\ is a spraid which is an isomorphic subspace of \Rnatt. To turn this spraid into an isomorphic spread, we only have to `unglue'. The best way to do so is to look at the trail space \sigRptht\ of \sigRt\ (see def.\,\ref{pathmorphisms}). Specifically, we look at the \sucz-trails\ in \sigRptht\ which (if not equal to the empty sequence $\maxdot^*$) start with a basic interval in $\suczz(\maxdotR)\iz\{[m, m\pluz 2]\midd m\inn\Z\}$. So put:

\parr

$\sigRungl\isdef\{a\iz a_0,\ldots a_{n-1}\inn \sigRpth\midd \ninn\midd a\ \mbox{is a}\ \mbox{\sucz-trail and}\ n\geqq 1\rightarrow\grdd(a_0)\iz 1 \}$

\parr

Then \sigRunglt\ has as maximal dot $\maxdot^*$, and an example of a basic dot in \sigRunglt\ is the sequence $[0, 2], [1,2]$, which has as unglued twin the basic dot $[1, 3], [1,2]$.
For simplicity, we also write \sigRunglt\ for the spread derived from the pre-natural space $(\sigRungl, \leqcstar, \aprtstar)$, which is the unglued version of \sigRt. Similarly we define:

\parr

$\sigzor\isdef \{[\frac{n}{2^{m}}, \frac{n+2}{2^{m}}]\midd n, \minn\midd n\pluz 2\leqq 2^{m}, m\geqq 1\}$,

\parr
so that taking $\maxdotzor\iz [0,1]$ we get a subfann \sigzort\ of \sigRt\ which is isomorphic to \zort. We can also unglue this subfann, as a subfan \sigzorunglt\ of \sigRunglt. We leave the details to the reader.
\exam

\parr
Another more involved example of a spraid arises when building the natural space $C^{\rm unif}(\zor, \R)_{\rm nat}$ of uniformly continuous functions from \zort\ to \Rt. We will sketch this in the appendix, see \ref{contfunasspraid}, referring for details to earlier work of Brouwer.

\sbsc{Cantor space is a universal fan}\label{Cantuni}
Similar to Baire space being a universal spread, Cantor space is a universal fan, by which we mean that each `finitely branching' spraid can be seen as the image of Cantor space under a natural morphism:\footnote{This is also a result due to Brouwer.} 

\defi
Let \Vnatt\ be a spread (spraid) with corresponding \Vprenat, so \Vleqct\ is a tree (trea). We call \Vleqct\ \deff{finitely branching} iff for all $c\inn V$ the set $\suczz(c)\iz\{a\inn V\midd a\sucz c\}$ is finite. We call \Vnatt\ a \deff{fan (fann)} iff \Vleqct\ is a finitely branching tree (trea). By extension, \Vnatt\ is \deff{fanlike}\ iff \Vnatt\ is isomorphic to a fan.\footnote{Any fann is fanlike so we need not define `fannlike'.}
\edefi

\thm
Let \Vnatt\ be a fann, then there is a surjective morphism from Cantor space to \Vnatt. (`Cantor space is a universal fan').
\ethm

\crl 
Every fann is fanlike. Every fanlike space is the natural image of Cantor space.
\ecrl

\prf
See \ref{Cantuniprf}\ where we give a self-contained proof. Compare the corollary to the `unglueing' that we did in example \ref{Baireuni}.
\eprf

\sbsc{Every compact metric space is homeomorphic to a fan}\label{metcompfan}
If we define a separable metric space to be compact whenever it is totally bounded and complete (as is standard in \bish), then it is a well-known result that every compact metric space is homeomorphic to a metric fan (a fan with a metric respecting the apartness, and endowed with the metric topology, which a fortiori is refined by the apartness topology).

Conversely, in \WaaThe\ it is shown in \intu\ that on a metric fan the apartness topology coincides with the metric topology, and that every apartness fan is metrizable. More general, in \intu\ every star-finitary apartness spread is metrizable. We retrieve this result for \bish\ in section \ref{starfinmet}.

\sectionn{Morphisms on spreads and spraids}\label{spreadmorph}

\sbsc{Refinement versus trail morphisms 1}\label{refvstrail1}
We return briefly to our discussion of refinement morphisms versus trail morphisms. With spreads (which derive from a tree) there is no need for trail morphisms. In fact a spread \Vnatt\ is \leqc-isomorphic to its trail space \Vpathnatt.
Since Baire space is universal (\ref{Baireuni}), we could develop a fruitful theory using only spreads and refinement morphisms (as is done in \intu). Studying refinement morphisms also on spraids is therefore an enlargement of this fruitful theory. But when working with spraids, we sometimes need trail morphisms as well.

\parr

Already for practical ease we wish to work with spraids such as \sigRt\ and refinement morphisms. This is relevant also, we believe, for computational purposes.\footnote{We refer once more to the recommendations in \BauKav, which we follow.} We do not straightaway know any theoretical applications, but the possibility to also use the \leqc-lattice properties seems to good to pass up. A relevant discipline in this respect could possibly be constructive (topological) lattice theory, or perhaps even algebraic topology.    

\parr

Therefore we take some time to show that for many important spraids a trail morphism can already be directly represented by a refinement morphism. This is especially relevant from the \appl\ perspective, we believe. There is the usual drawback: it requires a bit more patience from our readers. 

\sbsc{Unglueing of spraids}\label{unglspraids}
To simplify things, we show that any spraid \Vnatt\ can be unglued in exactly the same manner as \sigRt\ in paragraph \ref{Baireuni}. The idea is to turn to the subspread of \Vpathnatt\ formed by the \sucz-trails in \Vleqct\ (instead of looking at the spread of all trails).

\defi
Let \Vnatt\ be a spraid derived from \Vprenatt. The \deff{unglueing}\ of $(\Vcal,\Topaprt\!)$ is the spread \Vunglnatt\ derived from the pre-natural space \Vunglpre, where $\Vungl=\{a\iz a_0, \ldots a_{n-1}\inn\Vpath\midd \ninn\midd a\ \mbox{is a}\ \mbox{\sucz-trail}$ $\mbox{and}\ n\geqq 1\rightarrow\grdd(a_0)\iz 1\}$.
\edefi

We leave it to the reader to verify that unglueing a spraid \Vnatt\ amounts to adding, for each $a\inn V$, a finite number of copies of $a$ such that each \sucz-trail from \maxdott\ to $a$ is represented by one of the copies. These copies all have $\grdd(a)$ as length in $(\Vungl, \leqcstar)$.

For spraids, working with \Vunglnatt\ is more elegant than working with \Vpathnatt, which again seems relevant for computational practice. Notice also that if we start with a spread \Vnatt, then there is a trivial bijection between $V$ and \Vunglt, showing that spreads are already unglued.

\sbsc{Refinement versus trail morphisms 2}\label{refvstrail2}
Now we can show that for many important spraids a trail morphism can already be directly represented by a refinement morphism. We illustrate this first with \sigRt, our preferred representation of \Rt.

\prp
Let $f$ be a \pthh-morphism from \sigRt\ to \sigRt. Then there is a \leqc-morphism $g$ from \sigRt\ to \sigRt\ such that $f(x)\equivvR g(x)$ for all $x\inn\sigR$.
\eprp

\prf
We see $f$ as a \leqc-morphism from \sigRunglt\ to \sigRt. 
For $a\inn\sigR$ of the form $[\frac{4s+i}{2^{t+2}}, \frac{4s+i+2}{2^{t+2}}]$ where $1\leqq i\leqq 4$ and $t\inn\N$, put $\widehat{a}\iz [\frac{s}{2^{t}}, \frac{s+2}{2^{t}}]$. For all other
$a\inn\sigR$ let $\widehat{a}\iz\maxdotR$. 
Now for $a\inn\sigR$ there are finitely many \sucz-trails from \maxdotRt\ to $a$, say $b_0, \ldots, b_n$ where each $b_i$ is in \sigRunglt. Since the $f(b_i)$'s all touch, $\bigcapp_i\,\widehat{f(b_i)}$ is in \sigRt. We  
put $g(a)\isdef\bigcapp_i\,\widehat{f(b_i)}$. Then $g$ thus defined is a \leqc-morphism from \sigRt\ to \sigRt\ such that $f(x)\equivvR g(x)$ for all $x\inn\sigR$.
\eprf

\crl Let $f$ be a \pthh-morphism from \Rnatt\ to \Rnatt. Then there is a \leqc-morphism $g$ from \Rnatt\ to \Rnatt\ (in fact \sigRt) such that $f(x)\equivvR g(x)$ for all $x\inn\Rnat$.
\ecrl

\prf
For $a\inn\Rrat$ let $h(a)$ be the (unique) smallest interval in \sigRt\ such that $a\leqcR h(a)$. This determines a \leqc-isomorphism $h$ from \Rnatt\ to \sigRt. Again, for $a\inn\sigRungl$ put $\fslang(a)\iz h(f(a))$, which yields a \pthh-morphism \fslangt\ from \sigRt\ to \sigRt. By the proposition, \fslangt\ can be represented by a \leqc-morphism $g'$. Now for $a\inn\Rrat$ simply put $g(a)\iz g'(h(a))$ to obtain the required
\leqc-morphism $g$ representing $f$.
\eprf

The proposition combined with paragraphs \ref{compmetind}\ and \ref{starfinsprd}\ illustrates that for many complete metric spaces, we can find efficient spraid representations such that we can always work with refinement morphisms.

\rem
One relevant property here is that
for a finite intersection of basic dots we can find a basic dot of `small enough diameter' which contains the intersection in its interior. For our standard basic-open complete metric spreads this is obvious, but these spreads are themselves not an efficient representation. See \ref{starfinsprd}.
\erem

\sectionn{A direct link with intuitionism}\label{dirlinkint}

\sbsc{Natural spaces mirror Brouwer's spreads}\label{intmirror}
Our two goals in this section are to establish fundamental properties of natural spaces, and to give a classical mathematician a simple picture of important intuitionistic results. For the latter goal we establish that there is a constructive one-on-one correspondence between Brouwer's notion of spreads, and our notion of spreads as natural spaces given by a tree (for sake of comparison let us call this notion `natural spread'). Brouwer's \deff{spread-functions}\ correspond precisely to natural morphisms between our natural spreads. The practical advantage lies in the retrieval of many intuitionistic results for natural spaces.

\parr

We formulate the fundamental theorem in the language of \WaaThe. To avoid cumbersome mathematical structures, we do so on the meta-level.

\mthm
Brouwer's universal spread \sigomt\ corresponds precisely to natural Baire space. Other intuitionistic spreads can all be seen as an apartness spread \sigaprt\ which can be translated \deff{directly}\ to a corresponding natural spread \Vnatsig, which is homeomorphic in \intu\ to \sigaprtt\ equipped with the apartness topology \Topaprtsigt. Spread-functions between two spreads correspond precisely to \leqc-morphisms between the two corresponding natural spreads. Intuitionistic fans correspond precisely to our natural fans.
\emthm

\prf
The theorem is self-evident for anyone familiar with intuitionistic apartness topology, which was developed in \WaaThe. For reasons of space efficiency, we do not repeat all the relevant definitions here.
\eprf
  
For \class, the theorem offers a quite direct translation of many results in intuitionistic topology to the context of natural spaces (see the remark below). In this monograph we also translate some intuitionistic results to \bish, which is usually a bit more work. At the end of this section we translate a result of Veldman, which has its place in the intuitionistic study of the Borel hierarchy (see \VelBora\ and \VelBorb).\footnote{Other results of Veldman on this hierarchy might also be translatable to our classically valid setting, but we have not studied this.}\ In chapter four, we translate an intuitionistic metrization theorem.

\rem
The `direct' translations for \class\ often co-depend on a classically valid axiom of \intu\ called \BT\ (`Brouwer's Thesis, or equivalently decidable-bar induction \BID) which we did not mention earlier, for simplicity. The validity of \BT\ has been questioned in other branches of constructive mathematics, not in the least because Kleene showed that it fails in the branch of recursive mathematics called \russ. In \russ\ the principal axiom is \CT\ (derived from `Church's Thesis') which states that every infinite sequence $\alpha\inn\NN$ is computed by a known Turing-algorithm. Notice that this does not resemble our description of the engineer/scientist taking ever-more refinable measurements from nature.
\parr
Classical mathematicians are often unaware that compactness (more specific, the Heine-Borel property) of Cantor space fails in \russ. This compactness does however follow from \BT, therefore compactness can be dealt with elegantly in \intu. In the past decades constructive mathematicians have put much effort in developing a variety of constructive theories which respect \russ\ but still allow for a working theory of compactness.
\erem

\sbsc{Brouwer's Thesis, compactness and induction}\label{brothecom}
We discuss Brouwer's Thesis already here, although it more properly belongs in the next sections, because it influences our presentation from here on. In fact one can see \BT\ as a transfinite induction scheme (countable-ordinal induction) combined with the meta-insight that we can construct open covers of Baire space only by such a transfinite induction procedure. That is, if Baire space also contains sequences about which we know only finite initial segments at any given point in time. This explains why \BT\ fails in \russ, since in \russ\ we have a lot of additional information about sequences (namely the algorithms computing them) and we can use this information to construct non-inductive covers of recursive Baire space. This follows from Kleene's construction of an infinite decidable subset of \zostart\ which contains no infinite recursive path (the Kleene Tree, see \ref{canbairus}\ and \BauKle). We will need some definitions to phrase \BT.

\defi
We use natural Baire space \Bnattop, derived from \Bprenat. Let $A, B\subseteqq\Nstar$, then $B$ is called a \deff{bar} on $A$ iff $\all x\inn \hattr{A}\,\dris\ninn\,\dris b\inn B\,[x_n\leqc b]$, where $\hattr{A}\iz\bigcupp_{a\in A}\hattr{a}$. Notice that a bar on $A$ is the same as an open cover of $\hattr{A}$ consisting of basic open sets. 

\parr

Next we introduce `genetic bars' on \Nstar, using a form of countable-ordinal induction. 

\be
\item[\gindzb] The set $\{\maxdotB\}$ is a genetic bar on \Nstar. 
\item[\gindomb] If for each \ninnt\ we have a genetic bar $B_n$, then $B\iz\{n\starr a\midd a\inn B_n, \ninn\}$ is also a genetic bar on \Nstar.
\ee
Repeated application of the rules \gindzb\ and \gindomb\ yields all genetic bars on \Nstar.
\edefi

We believe this form of induction is constructively acceptable, and formulate the appropriate axiom:

\xiom{\PGI}
The definition of genetic bars is valid. Moreover, let $P$ be a property of bars on \Nstart\ such that:
\be
\item[\gindz] The genetic bar $\{\maxdotB\}$ has property $P$.
\item[\gindom] If for each \ninnt\ we have a genetic bar $B_n$ with property $P$ then the genetic bar $B\iz\{n\starr a\midd a\inn B_n, \ninn\}$ also has property P.
\ee
Then all genetic bars on \Nstart\ have property $P$.
\exiom

\defi
Let $C, D\subseteqq \Nstar$ be bars on \Nstar. We say that $C$ \deff{descends} from $D$ iff for each $d\inn D$ there is a $c\inn C$ with $d\leqc c$.
\edefi

\rem
This terminology makes sense if we picture Baire space as an infinite tree, which branches \emph{upward}\ from its maximal element \maxdotBt. Now `$C$ descends from $D$' describes the picture of a bar $D$ on \Nstart, such that below each element of $D$ there is already an element of $C$ (therefore $C$ covers at least what $D$ covers). But we acknowledge that we could have called \maxdotBt\ the minimal element, reversing the \leqct-notation. 
\erem

We can now formulate our version of \BT: \footnote{This version is easily seen to be equivalent to the version given in \WaaArt.}

\xiom{\BT}
\PGI\ holds, and every bar on \Nstart\ descends from a genetic bar on \Nstart.
\exiom
\parr

An intuitionistic plea for \BT\ can be give in the following way. The universe of Baire space that we have in mind is inhabited by choice sequences arising step by step in the course of time. This means that in general only the minimum of information about an element is known. For the author, the axiom expresses that --- given such  a universe --- the only one way to convince ourselves that a subset $B$ of \Nstart\ is indeed a bar, is to show that it descends from something that we can intuitively grasp as a bar, namely a genetic bar.
\parr
In fact the more precise analysis is that there are two basic methods that can be employed in ascertaining the `bar' status of a given subset $B$ of \Nstar. The direct method is to have $B$ given as a genetic bar. The other method is to see that $B$ descends from a previously ascertained bar. These two methods put together yield the method of checking whether $B$ descends from a genetic bar.
\parr
Kleene calls this aptly `reversing the arrows'. The genetic definition in fact mirrors what our intuition tries to do when visualising an arbitrary bar. The author has no trouble accepting this definition, and in this acceptance lies his intuitive justification of \PGI. 
Brouwer's justification looks rather more complex even when explained by Veldman, but we believe it to be essentially the same as our presentation above. In \KleVes, Kleene gives some nice examples of bars descending from genetic bars which cannot be shown to be genetic unless we prove \mbox{e.g.}\ that there are 99 consecutive 9's in the decimal expansion of $\pi$. 
\parr
\BT\ can be proved classically by contradiction as follows: suppose $B$ is a bar on \Nstart\ which does not contain a genetic bar. For $a\inn\Nstar$ write $B^a$ for $\{d\inn \Nstar\midd a\starr d\inn B\}$. Since $B$ does not contain a genetic bar, \maxdotBt\ is not in $B$. Put $b_0\iz\maxdotB$. Then for at least one $b\iz b_1\inn\suczz(b_0)$ the bar $B^b$ on \Nstart\ does not contain a genetic bar, meaning $\maxdotB\notinn B^b$ which implies $b_1\notinn B$. Repetition of this argument yields a shrinking sequence $\beta\iz\bnninn$ such that $\beta$ evades $B$. Contradiction, since $B$ is a bar. Therefore $B$ must contain a genetic bar, and so trivially descends from a genetic bar.    

\rem
If one is willing to accept \BT, then an elegant constructive theory of compactness is possible for natural spaces from the ingredients presented up to now. As stated earlier, this theory closely resembles intuitionistic results. The author believes that in physics, the (mathematics underlying the) universe which we study resembles the setting of sequences arising step by step without much further information. This is why he would especially invite physicists to take notice of intuitionistic mathematics. 
\erem

For the rest of this monograph, we will not assume \BT. Instead, we will introduce inductive definitions which capture most of \BT. The prize we gain is that with these inductive definitions we can build a (limited) theory of compactness also in \russ. The price we pay is that this inductive machinery is not easy, and tends to lead to somewhat involved proofs. 
\parr
Another gain is perhaps more `bridging' in character. Intuitionism has never been very popular, in contrast maybe to \bish\ and constructive formal topology. By presenting natural spaces as spraids and then developing inductive definitions, we can hopefully illuminate intuitionistic concepts and results, while still remaining in \bish. 

\sbsc{Baire morphisms together do not form a natural space}\label{Funnotnat}
In our discussion so far, we have highlighted that the class of natural spaces is large, not to say vast. In fairness let us also state that there are important spaces which cannot be represented as a natural space. In the appendix (see \ref{notnat}) we will comment on this shortly.
\parr
A revealing example in this respect is the space of all Baire morphisms, by which we mean natural morphisms from Baire space to itself. Equivalently (by thm.\,\ref{basicneighbor}) we can see this as the space of all continuous functions from Baire space to itself. The key for our development\footnote{Which is indebted to Veldman's expositions on intuitionism, see the bibliography.}, is that every Baire morphism can be coded by an element of Baire space itself. We detail this explicitly:

\defi
We fix a bijection $\lcod\rcod$ from \Nstart\ to \Nt, with inverse $\rcod\lcod$ such that for all $a, b \inn\Nstar$ we have $\lcod a\rcod \leqq \lcod a\star b\rcod$.  Using $\rcod\lcod$ we can see each $\alpha\inn\NN$ as a sequence of pairs of basic dots $(\rcod n\lcod, \rcod \alpha(n)\lcod)_{\ninn}$ in $\Nstar\timez\Nstar$. We say that $\alpha$ is a \deff{coded Baire morphism}\ iff the so obtained set $\{(\rcod n\lcod, \rcod \alpha(n)\lcod)\midd\ninn\}$ is a Baire morphism (a natural morphism from Baire space to itself).

\parr 

To transfer this idea to $\Bnat\iz\NN_{\rm nat}$,
remember that for $p\inn\Bnat$ we write $p^*$ to denote the unique sequence $\alpha\inn\NN$ corresponding to $p$ (see def.\,\ref{natBaire}). So $p^*(m)$ then denotes $\alpha(m)$ for this $\alpha\inn\NN$ corresponding to $p$. Thus we look on $p$ in \Bnatt\ also as a sequence of pairs of basic dots $({^{\!}}\rcod n\lcod, \rcod p^*(n)\lcod{^{\!}})_{\ninn}$ in $\Nstar\!\timez\Nstar$. We say that $p$ is a \deff{star-coded Baire morphism}\ iff the so obtained set $\{(\rcod n\lcod, \rcod p^*(n)\lcod)\midd\ninn\}$ is a Baire morphism.\footnote{For each basic dot $\rcod n\lcod$, the image under the morphism is the basic dot $\rcod p^*(n)\lcod$.}

\parr

We define: $\Fun\isdef\{p\inn\Bnat\midd\,\mbox{$p$ is a star-coded Baire morphism}\}$. For $f$ in \Fun, write $\tilde{f}$ for the $f$-induced Baire morphism. Now if $f\aprtom g$ for $f, g\inn\Fun$ then the induced morphisms $\tilde{f}, \tilde{g}$ define different continuous functions from Baire space to itself. (Therefore, by slight abuse of notation, we also write: `let $\tilde{f}$ be a morphism, determine the corresponding $f\inn\Fun$'.)
\edefi

The next theorem, slightly modified, is from \VelThe\ (also read \VelBora, \VelBorb), as is its elegant proof (based on a Cantorian diagonal argument):

\thm \hspace*{-.4pt}(Wim Veldman)
Let $F$ be a Baire morphism for which $\Ran(F)\subseteqq\Fun$. Then there is an $f\inn\Fun$ such that $f\aprt F(p)$ for every $p$ in \Bnat.
\ethm

\crl
\Fun\ cannot be represented by a natural space; the space of continuous functions from Baire space to itself cannot be represented by a natural space (`the space of continuous functions from Baire space to itself is not spreadlike').
\ecrl

\prf
Let $p$ be in \Bnat, then $F(p)\inn\Fun$ determines a Baire morphism $\widetilde{F(p)}$. Therefore for each $p\inn\Bnat$, we can send $p$ to $\widetilde{F(p)}(p)$, and this function is given by a morphism $z$. We now construct a Baire morphism $\tilde{f}$ such that $\tilde{f}(p)\aprt z(p)\iz\widetilde{F(p)}(p)$ for all $p\inn\Bnat$.

Since $z$ is a morphism, the set $B=\{b\inn\Nstar\midd z(b)\precc\maxdotB \weddge \all c\succ b [z(c)\iz\maxdotB]\}$ (the set of basic dots on which $z$ determines the first basic dot of the image) is decidable. We take care to ensure that $\tilde{f}(b)\aprt z(b)$ for all $b\inn B$. Then for $a\notinn B$: if there is $b\inn B$ with $b\precc a$, we put $\tilde{f}(a)\iz\maxdotB$, and otherwise there is a $c$ and a $b\inn B$ with $a\iz b\star c$, and we put $\tilde{f}(a)\iz \tilde{f}(b)\star c$.

Now for $p$ in \Bnat, we see that the morphism $\widetilde{F(p)}$ is apart from 
$\tilde{f}$, because $\widetilde{F(p)}(p)= z(p)\aprt\tilde{f}(p)$ by construction of $\tilde{f}$. For $\tilde{f}$, determine the corresponding $f\inn\Fun$, then $F(p)\aprtom{_{\,}}f$ for all $p\inn\Bnat$.
\eprf

\sectionnb{Compactness, Brouwer's Thesis and inductive definitions}\label{compactness}

\sbsc{Compactness in the absence of Brouwer's Thesis}\hspace*{-.91pt}
From this section on we carry out the idea mentioned in \ref{brothecom}\ (where \BT\ was defined): to develop a theory of compactness using inductive definitions, in the absence of \BT. We discuss the relationship between natural topology and some other approaches to constructive mathematics, especially regarding compactness. 
We do not in any way claim final wisdom, since these are complex issues.\footnote{Also there are many varieties of constructive and semi-constructive topology and the author's  knowledge of and insight in these varieties is very limited.} 

\sbsc{Pointwise, pointfree or both?}
We take up some extra space and time to show that natural topology suits both the pointwise and pointfree perspective. In the previous chapters the pointwise approach was our central focus. But in the absence of Brouwer's Thesis \BT, the pointfree approach is necessary for inductive definitions if one wishes to reduplicate compactness-like results. We concentrate on `basic dots' but at the end of the day we also have the usual separable topological spaces and the pointwise perspective at our disposal.  

\parr

For clarity we repeat that the results in this paper are derived within \bish, and that these results are mostly translations of existing intuitionistic results. Intuitionistic results have been called non-effective by some authors, but we believe that the difference lies mainly in the acceptance of induction axioms versus the incorporation of the induction in all the definitions.  We consider this more a question of style than a paramount distinction.

\sbsc{Introduction to constructive compactness issues}\label{introcantor}
To understand compactness issues we turn to Cantor space (\Cant). In \class\ and \intu, Cantor space is a universal separable compact space (a result due to Brouwer). By this we mean that every separable compact space is the continuous image of \Cant. This also holds in \bish, but in \bish\ the definition of compactness is a metrical one, see below. 

\parr

Cantor space can be pictured as a dually branching subtree of Baire space \NNt, the infinitely branching tree. Brouwer's notation for Cantor space was \sigtot, and for Baire space \sigomt. In Brouwer's terminology, every finitely branching subtree of \sigomt\ is called a \deff{fan}. Brouwer's compactness axiom \FT\ (the Fan Theorem) states that every cover of a fan has a finite subcover (a form of the Heine-Borel property, classically equivalent to K\"{o}nig's lemma).

\parr

However, in recursive mathematics (\russ) there is no immediate topological definition of `compact'. In \russ, Baire space is in fact homeomorphic to Cantor space, due to the axiom \CT\ (see \ref{axiomsrussct})  Therefore in our theory of natural topology, the compactness of Cantor space\footnote{The property `every open cover has a finite subcover' (the Heine-Borel property).} cannot be shown unless we adopt Brouwer's \FT, or strengthen the definition of `cover' to `inductive cover' (thereby excluding many interesting non-inductive covers in \russ). \FT\ can be derived from \BT, and so also holds in \class. Therefore our description of natural spaces for a classical mathematician leads to results closely resembling intuitionism. \BT\ fails in \russ, from which the picture often is drawn that \russ\ and \intu\ are incompatible. 

\parr

Bishop, who developed a neutral constructive stance with Bishop-style mathematics (\bish), thought that Brouwer's motivation for \BT\ was mystical (not\-withstanding the simplifications offered by Kleene in \KleVes).  He tried to work around the difficulties associated with topological compactness by only defining `metrical compactness' (meaning `totally bounded and complete', for metric spaces). In the light of the situation in \russ, this seems a good neutral solution and it has a large number of nice applications (see \BisBri, \BriVit). But Bishop also wanted to obtain that continuous functions are uniformly continuous on compact spaces, and being unable to prove this by lack of a compactness axiom, he added this property to the definition of continuous function. In \WaaArt\ it is shown that this is practically equivalent to the adoption of \FT\ (thereby \bish\ somewhat loses its neutral stance, and veers towards intuitionism). In reaction to \WaaArt, in \BriVit\ the definition of `continuous function' has been modified back to the usual epsilon-delta one.

\parr
\enlargethispage{1mm}
In formal topology a way of dealing with compactness is to look at  what we will call `inductive covers' and `inductive morphisms', using a form of transfinite countable-ordinal induction. It is a nice solution which in a sense incorporates \BT\ already in the definitions. In this way, compactness results of formal topology also have a recursive interpretation in \russ, which is an attractive feature shared with \bish. But as in \bish, the non-inductive non-compactness of Cantor space in \russ\ then may be ignored, and formal topology in this way also seems to veer towards intuitionism.\footnote{The author finds the literature on formal topology hard to read, and to avoid mistakes refrains from a more precise mathematical comparison.}

\sbsc{Natural Cantor and Baire space are isomorphic in RUSS}\label{canbairus}
From a well-known basic result from \KleVes\ we derive the equally well-known result that in \russ\ natural Cantor space is isomorphic to natural Baire space. 

\thm
In \russ\ natural Cantor space is isomorphic to natural Baire space.
\ethm

\prf
From \KleVes\ we can directly define a decidable countable subset\enlargethispage{1mm}\footnote{derived from what Andrej Bauer in \BauKle\ aptly calls the Kleene Tree.}\ $K_{\rm bar}=\{k_n\midd\ninn\}\subsett\zostar$ such that $k_n\leqcom k_m$ implies $n\iz m$ for all $n,\minn$ and such that in \russ\ $\{\hattr{k}\midd k\inn K_{\rm bar}\}$ is an open cover of \Cnattop\ which has no finite subcover. We use $K_{\rm bar}$ to define an isomorphism $f$ in \russ\ from \Bnattopt\ to \Cnattopt\ as follows.
\parr
Put $f(\maxdotB)=\maxdotB$. Then let $a=a_0, \ldots, a_m \inn\Nstar$. Put $f(a)=k_{a_0}\star k_{a_1}\star\ldots\star k_{a_m}$.
\parr
The verification that $f$ is an isomorphism is not difficult and left to the reader.
\eprf

\parr

The theorem shows that with our definitions so far we cannot hope to define a natural-topological notion of compactness in \russ. In formal topology, this is partly resolved by using a form of transfinite induction.\enlargethispage{1mm}\footnote{Which for natural spaces seems equivalent to \BT, see prp.\,\ref{indcov}.}\ However, pointwise problems in \bish\ related to these compactness issues persist (see \ref{metaruss}). Yet we will adopt this transfinite-induction strategy as well, for elegance and for purposes of informal comparison.

\sbsc{A model of \intu\ as part of \russ} 
Usually, \intu\ is seen as being at odds with \russ, because of the compactness troubles in \russ. However, it is also possible to informally model \intu\ as a two-player game in \russ. In this model, one can see \intu\ as the part of \russ\ where all covers are inductive. \BT\ then becomes an elegant way of saying that we restrict our recursive world to all things inductive. In this model we can prove the intuitionistic continuity principle \CP. This model has similarities to Weihrauch's \TTE\ (type-two effectivity), but we are no expert and refrain from making direct comparisons.  

\parr
We present the model of \intu\ as part of \russ\ in section \ref{intuinruss}. This model might give an argument for physics why \FT\ is more valid than `not \FT', \deff{even} if \CTphys\ is seen to hold. Here \CTphys\ stands for the statement that nature can only produce recursive sequences. Since \CTphys\ is as of yet undecided, \russ\ might be a relevant mathematical model for physics, and these compactness issues seem worthy of physicists' attention as well. 

\parr

In \bish\ and in constructive formal topology, a preference for \FT\ sometimes seems cloaked in definitions. This leads to an exclusion of (parts of) \russ\ which is not always easy to spot. We believe it more fruitful for our foundational discussion here to give \russ\ a more equal place, and then study the arising topological structures in the light of different axiom systems. This is one reason for developing the concepts `neutrally' in the previous chapters.     

\sbsc{Compactness and inductivity}\label{compind}
So now we need to find a perspective on compactness and inductivity. We adopt from formal and pointfree topology the notion of (what we call) `inductive covers' and `inductive morphisms'. For natural spaces, Brouwer's induction scheme in our eyes is more elegant and wieldy. Therefore we derive our inductive covers from `genetic bars' as defined in \ref{brothecom}, only slightly generalized to arbitrary spraids. Thus emulating formal topology in a simplified way, we can interpret compactness in \russ\ as well. One may however ask oneself, after all the work has been done, whether any real benefit has been gained over intuitionism. For one answer, one could look at the hopefully better understanding of the relevant issues and the relationship between the various branches of constructive mathematics. As a second answer, we also phrase a strong support of intuitionistic mathematics in paragraph \ref{recap}.

\chepter{Chapter three}{Inductivizing our definitions}{We develop genetic induction based on a simplification of bar induction. Inductive covers thus defined are equivalent to `formal inductive covers' coupled with a formal-topology style of induction. 
\parr
Every subfann (including Cantor space and the real interval $[\alpha,\beta]$) has the Heine-Borel property for inductive covers (of that fann in its mother spraid). 
\parr
Inductive morphisms respect inductive covers; they are uniformously continuous on metric (sub)fanns. \bisCo\ functions from \Rt\ to \Rt\ are representable by an inductive morphism. The statement that \bisco\ functions from \Rt\ to \Rplust\ are representable by an inductive morphism from \Rnatt\ to \Rplusnatt\ is equivalent to \FT. Pointwise problems for \bish\ persist, related to the reciprocal function and compactness (the relevant example of ContraCantor space is given in the appendix).
\parr
We discuss Kleene's realizability and other ways to define inductive morphisms.
\parr
Finally we define (in)finite-product spaces, and prove a \bish\ version of Tychonoff's theorem.
}  

\sectionnb{Induction in formal-topology style}\label{indformtop}

\sbsc{Bootstrap method}\label{introinddef}
Our development strategy has a drawback: we already gave the basic definitions, and now we want to build an inductive theory. This involves revisiting the earlier basic definitions, to `inductivize' them. We again ask some patience from the reader.

\sbsc{Basic (open) covers and per-enumerable open covers}\label{basiccovers}
We begin our development of inductive covers by defining a basic covering relation \covbt\ on sets of basic dots.

\defi
Let \Vnatt\ be a natural space with corresponding \Vprenatt. Let $A, B \subseteqq V$. We say that $B$ is a \deff{basic cover}\ of $A$, notation $A\covb B$ iff for all $x\inn\hattr{A}=\{y\inn\Vcal\midd \driss a\inn A[y\leqc a]\}$ there is a $b\inn B$ with $x\leqc b$. In other words, iff $\hattr{A}\subseteqq\hattr{B}$. Notice that \covbt\ is transitive. 
By extension we say that $B$ is a basic cover of \Vcalt\ iff $\{\maxdot\}\covb B$.
Basic covers need not correspond with open sets in the topology (but for Baire space the distinction is moot), so for $A\covb B$ we say that $B$ is a \deff{basic open cover}\ of $A$ if in addition $\closr{B}$ is open in \Vnatt, and then we write $A\ocovb B$.
\edefi

\parr
Notice that this definition is not `pointfree', it relies essentially on the points in \Vcalt. The way in which we have acquired the insight `for all $x\inn\hattr{A}$ there is a $b\inn B$ with $x\leqc b$' is left unspecified. This means that in \russ\ covers derived from the Kleene Tree are also basic covers, and so compactness of Cantor space cannot be derived with regard to basic covers without extra axioms. What does hold in \russ\ as well as in \intu\ and \class\ is that every open cover of Baire space is refined by an enumerable basic open cover of Baire space, which entails a form of the Lindel\"of property (`every open cover has an enumerable refinement').\footnote{See axiom \BDD\ in \ref{axiombdd}.}\ To be able to work with this important Lindel\"of property in \bish\ as well, we translate some definitions of \WaaThe\ to our setting here.

\defi
Let \Vnatt\ as above, and let $\Ucal\subsett\Vcal$ be open in \Vnatt. We say that \Ucalt\ is \deff{enumerably open} iff there is an enumerable subset $U$ of $V$ such that $\Ucal\iz\hattr{U}$.  Let $\ucov\iz\{\Ucal_n\midd\ninn\}\iz\{\hattr{U_n}\midd\ninn\}$ be an enumerable collection of enumerably open subsets $\Ucal_n\iz\hattr{U_n}$ of \Vcalt, then we say that \ucovt\ is a \deff{per-enumerable open collection}. If in addition \ucovt\ is an open cover of a subset \Wcalt\ of \Vcalt, then we say that \ucovt\ is a per-enumerable open cover of \Wcalt\ in \Vnatt.
\edefi  

\parr

Per-enumerable covers have nice properties. Per-enumerable covers of metric spaces for instance allow a subordinate partition of unity, as well as a star-finite refinement (see \WaaThe). These are powerful topological tools, which use the paracompact properties of metric spaces. Also, per-enumerable covers form a connection between basic dots in natural spaces and basic opens of formal topology. For example in \Rnatt, the open real interval $(0,1)$, which is represented by a basic open in formal topology, is enumerably open since it is represented by the enumerable set of basic dots $\{[p,q]\inn\Rrat\midd 0\smlr p\smlr q\smlr 1\}$.

\parr  

The next steps in our development concern inductive basic covers. If we specify covers inductively, we capture a form of compactness, namely the Heine-Borel property that inductive covers of a fann have a finite subcover.

\sbsc{Formal inductive covers}\label{formindcov}
To facilitate a connection with formal topology, we first define a `formal' inductive covering relation \icovbt, using a form of transfinite countable-ordinal induction which is generally considered constructive. We later give a more concise induction scheme which is in essence due to Brouwer.

\defi
Let $b,c\inn V$ and $A, B, C\subseteqq V$ as previously, we define:
\bei
\itemmz{\indob} $b\leqc c$ implies $\{b\}\icovbv\{c\}$.\footnote{One can replace this with the seemingly stronger rule: if $A\covb B$ and $B$ is finite, then $A\icovb B$. We believe it equivalent, but this equivalence is a bit circuitous, depending on Baire space as universal space, so one may prefer this stronger rule.} 
\itemmz{\indtob} if for all $a\inn A$ we have $\{a\}\icovbv B$, then $A\icovbv B$.
\itemmz{\indthrb} if $A\icovbv B\subseteqq C$ then $A\icovbv C$. 
\itemmz{\indfob} if $A\icovbv B\icovbv C$ then $A\icovbv C$.
\itemmz{\indfib} $\{b\}\icovbv\{d\midd d\precc b\}$.
\eei

Repeated application of the rules \indob\ through \indfib\ yields all sets $D, E\subseteqq V$ for which $D\icovbv E$.\footnote{This is transfinite countable-ordinal induction}\ We say that $B$ is a \deff{formal inductive cover} of $A$ w.r.t. \Vnatt, iff $A\icovbv B$ (as an exercise the reader may prove by induction that $A\icovb B$ implies $A\covb B$). If the context space is clear, we omit the subscript and simply write $A\icovb B$. By extension we say that $B$ is a formal inductive cover of \Vcalt\ iff $\{\maxdot\}\icovbv B$.
\edefi

\parr

We state the appropriate axiom which from now on we take to hold:  

\xiom{\PFI} (Principle of Formal Induction): For any natural space \Vnatt, the definition of formal inductive covers is valid. Moreover, let $P$ be a property of pairs of subsets $A, B$ of $V$, such that:
\bei
\itemmz{\indo} $b\leqc c$ implies $P(\{b\},\{c\})$ for $b, c\inn V$.
\itemmz{\indto} if for all $a\inn A$ we have $P(\{a\},B)$, then $P(A, B)$.
\itemmz{\indthr} if $P(A, B)$ and $B\subseteqq C$ then $P(A, C)$. 
\itemmz{\indfo} if $P(A, B)$ and $P(B, C)$ then $P(A, C)$.
\itemmz{\indfi} $P(\{b\}, \{d\midd d\precc b\})$ for all $b\inn V$.
\eei
Then $A\icovb B$ implies $P(A, B)$ (for  $A, B\subseteqq V$).
\exiom

\alin

To show how things work in this formal-topology style of induction, we deduce the Heine-Borel property for formal inductive covers of Cantor space -on its own- as a corollary to the following proposition: 

\prp
In \Cnat, let $A, B\subseteqq\zostar$ such that $A\icovb B$. Then for all $a\inn A$ there is a finite $C\subseteqq B$ such that $\{a\}\icovb C$.
\eprp

\prf
By formal induction, using as property $P(A,B)$: `for all $a\inn A$ there is a finite $C\subseteqq B$ such that $\{a\}\icovb C$'. We check that \indo\ through \indfi\ hold for $P$:
\bei
\itemmz{\indo} trivially, by \indob, $b\leqc c$ implies $P(\{b\},\{c\})$ for $b, c\inn\zostar$.
\itemmz{\indto} if for all $a\inn A$ we have $P(\{a\},B)$, then for $a\inn A$ we know: for all $a'\inn\{a\}$ there is a finite $C'\subseteqq B$ such that $\{a'\}\icovb C'$. It trivially follows that for all $a\inn A$ there is a finite $C\subseteqq B$ such that $\{a\}\icovb B$. 
\itemmz{\indthr} if $P(A, B)$ and $B\subseteqq C$ then for all $a\inn A$ there is a finite $D\subseteqq B$ such that $\{a\}\icovb D$. Since $B\subseteqq C$ we see that $P(A,C)$.
\itemmz{\indfo} if $P(A, B)$ and $P(B, C)$ then for arbitrarily given $a\inn A$ there is a finite $E\subseteqq B$ such that $\{a\}\icovb E$. Since $P(B, C)$, we also have that for $e \inn E$ there is a finite $D_e\subseteqq C$ such that $\{e\}\icovb D_e$. But then, taking $D\iz\bigcupp_{e\in E}D_e$, we see that $D\subseteqq C$ is finite and $\{a\}\icovb E\icovb D$, so by \indfob\ also $\{a\}\icovb D$. Since $a$ is arbitrary we conclude $P(A, C)$. 
\itemmz{\indfi} let $b\inn\zostar$, we wish to show $P(\{b\}, \{d\inn \zostar\midd d\precc b\})$. For this we take the finite subset $C\iz\{b\starr 0, b\starr 1\}$ of $B\iz\{d\midd d\precc b\}$, and we see that $\{b\}\icovb C$, since $B\icovb C$ by \indtob\ and \indob, and $\{b\}\icovb B$ by \indfib, and so $\{b\}\icovb C$ by \indfob.  
\eei
Therefore by \PFI, for any $A, B\subseteqq\zostar$ such that $A\icovb B$ we have $P(A, B)$. In other words: for all $a\inn A$ there is a finite $D\subseteqq C$ such that $\{a\}\icovb D$, so we are done. 
\eprf

\crl
\be
\item[(i)] In \Cnat, let $A\subseteqq\zostar$ be finite and $B\subseteqq\zostar$ such that $A\icovb B$. Then there is a finite $B'\subseteqq B$ such that $A\icovb B'$.
\item[(ii)] In \Cnat, let $B\subseteqq\zostar$ be a formal inductive cover of \Cnatt, then $B$ has a finite subcover. (Heine-Borel property for formal inductive covers of Cantor space).
\ee 
\ecrl

\prf
Ad (i): In \Cnat, by the theorem, for all $a\inn A$ there is a finite $B_a\subseteqq B$ such that $\{a\}\icovb B_a$. Since $A$ is finite, we can take $B'\iz\bigcupp_{a\in A}B_a$. 

Ad (ii): take $A\iz\{\maxdotB\}$ and use (i) above.
\eprf

\sectionn{Genetic induction in Brouwer's style}\label{genindbro}

\sbsc{Inductive covers following Brouwer's Thesis}\label{indcov}
We could continue developing the theory in this formal-topology style, but we feel that for the setting of natural spaces, Brouwer's approach is more precise and concise. Therefore we develop an alternative notion of `inductive cover' for spraids (and prove that it amounts to the same as `formal inductive cover'). To facilitate the connection with intuitionism, we adopt (and adapt) intuitionistic terminology. To be foundationally clear, we formalize our countable-ordinal induction scheme as an axiom (\PGIg). The advantage to this alternative approach is that definitions and proofs become shorter, and that we can more easily adopt results from intuitionism.
\parr
We generalize the Baire space definition \ref{brothecom}\ of genetic bars, genetic induction and \BT\ to arbitrary spraids. (Remember from \ref{Baireuni}\ that without loss of generality a natural space is given by a spraid.).

\defi
Let \Vnatt\ be a spraid with corresponding \Vprenatt\ (so \Vleqct\ is a trea). Let $B\subseteqq V$. If $B$ is a basic cover of \Vnatt\ (see def.\,\ref{formindcov}), then we say that $B$ is a \deff{bar} on $V$ (in \Vnatt, equivalently $\{\maxdot\}\covb B$. We remind the reader that for any $a$ in $V$, the basic subspraid $\Vsuba\iz\{b\inn V\midd b\leqc a\}\iz\{a\}_{\!\leqc}$ is formed by putting its maximal dot as $\maxdot_a\iz a$.\footnote{Taking $a\iz \maxdotV$ gives the entire spraid $V$ (\Vnatt) itself.}\ We extend this notation by putting $C_{\leqc}\iz\bigcupp_{c\in C}\Vsubc$ for any subset $C\subseteqq V$. Also kindly remember that for $a$ in $V$, we write $\suczz(a)$ for $\{b\inn V\midd b\sucz a\}$. 

\parr

Now we inductively define genetic bars on basic subspraids $\Vsuba$ of \Vnatt\ as follows:
\be
\item[\gindzb] For $a\inn V$ the set $\{\maxdot_a\}$ is a genetic bar on $\Vsuba$. 
\item[\gindomb] If for $a\inn V$ and all $b\inn\suczz(a)$, $\Bsubb$ is a genetic bar on $\Vsubb$, then $\bigcupp_{b\in\suczs(a)}\Bsubb$ is a genetic bar on $\Vsuba$.
\ee
Repeated application of the rules \gindzb\ and \gindomb\ yields all genetic bars on basic subspraids of \Vnatt.\footnote{Again this is countable-ordinal transfinite induction. Also note that the version for Baire space is even more elegant. But using this more elegant form requires the encoding of all entities as natural numbers, which we have waived.}
\edefi

\parr
We state the appropriate axiom which from now on we take to hold:  

\xiom{\PGIg} (generalized Principle of Genetic Induction): For any spraid \Vnatt, the definition of genetic bars is valid. Moreover, let $P$ be a property of bars on basic subspraids of \Vnatt\ such that:
\be
\item[\gindz] For $a\inn V$, the genetic bar $\{\maxdot_a\}$ on $\Vsuba$ has property $P$.
\item[\gindom] If for $a\inn V$ and all $b\inn\suczz(a)$, $\Bsubb$ is a genetic bar on $\Vsubb$ with property P, then the genetic bar $\bigcupp_{b\in\suczs(a)}\Bsubb$ on $\Vsuba$ has property P.
\ee
Then all genetic bars on basic subspraids of \Vnatt\ have property $P$.
\exiom

\defi
Let \Vnatt\ as above, with $a\inn V$, and let $B, C, D\subseteqq V$ be bars on $\Vsuba$. Then $C$ \deff{descends} from $D$ iff for all $d\inn D$ there is $c\inn C$ with $d\leqc c$ (iff $D\subseteqq C_{\!\leqc})$. 

We say that $B$ is an \deff{inductive bar} on $\Vsuba$ (and an \deff{inductive cover} of $\Vsuba$) iff there is a genetic bar $G$ on $\Vsuba$ such that $B$ descends from $G$. By extension we say that $B$ is an inductive cover of $\{a\}$, notation $\{a\}\oicovbv B$ or simply $\{a\}\oicovb B$ when the context is clear. Notice that $B$ need not be a subset of $\Vsuba$.

Next, let $E, F\subseteqq V$. We say that $E$ is an inductive cover of $F$, notation $F\oicovbv E$, iff for all $a\inn F$ we have $\{a\}\oicovbv E$. We also write $F\oicovb E$ when the context is clear. 
\parr
From a pointwise perspective (see def. \ref{basiccovers}) a per-enumerable open cover $\ucov\iz\{\hattr{U_n}\midd\ninn\}$ of a subspraid $\Wcal\subseteqq\Vcal$ derived from \Wprenatt\ is called an \deff{inductive open cover} of \Wcalt\ iff $W\oicovbw\bigcupp_{\ninn} U_n$.
Then, for an arbitrary subset $\Acal\subseteqq\Wcal$, we also say \ucovt\ is an inductive open cover of \Acalt. 
See \ref{addremindcov}\ for additional comments.
\edefi

\sbsc{Genetic and formal covers coincide on spraids}\label{geneticisformal}
We can show that for spraids there is no distinction between inductive covers (obtained through genetic bars in Brouwer's style) and formal inductive covers (defined in \ref{formindcov}\ in formal-topology style).

\thm
Let \Vnatt\ as above, and let $E, F\subseteqq V$. Then $F\icovb E$ iff $F\oicovb E$.
\ethm

\prf See appendix \ref{prfindcov}, we use \PFI\ and \PGIg. 
\eprf

\rem
This will allow us to use \indob\ through \indfib\ as properties of \oicovbt\ as well. The theorem hopefully also partly clarifies the relation between \intu\ and formal topology. As explained earlier, we prefer genetic induction on spraids. 
\erem

\sbsc{Brouwer's Thesis generalized to spraids}\label{brother}
We do not adopt \BT\ in our narrative, yet we cannot escape generalizing this axiom to arbitrary spraids:

\xiom{\BTg}
\PGIg\ holds, and every bar on a spraid \Vnatt\ descends from a genetic bar on \Vnat.
\exiom 
\parr
The defense of \BTg\ derives straightforwardly from the defense of \BT\ given in \ref{brothecom}. \BTg\ also follows from \BT, and therefore holds in \class\ and \intu. 

\parr
We see \BTg\ as a deep insight in the nature of how we can construct bars at all, if we have to deal with infinite sequences (of basic dots) about which at any given time we know only initial finite segments. 
\BTg\ fails in \russ\ since in \russ\ we have a finite algorithm for each infinite sequence, giving us far more knowledge of such sequences than in the limited-information setting that Brouwer had in mind. We transpose Brouwer's setting to \russ\ in section \ref{intuinruss}, to show what we mean. 

\prp
Let \Vnatt\ be a spraid derived from \Vprenatt. Then \BTg\ implies that every cover
of $V$ is inductive.
\eprp

\rem
The proof is trivial. The proposition underlines that if we accept \BTg\ we can simply work with basic covers and skip the inductivizing of basic definitions. The resulting theory is common ground of \class\ and \intu, and elegant. The genetic induction scheme remains essential in this theory though, so not all the work done here is superfluous if one accepts \BTg. Using \PFI\ and theorem \ref{Baireuni} (`every natural space is spreadlike') we obtain the following corollary.
\erem 

\crl
Let \Vnatt\ be a natural space derived from \Vprenatt. Then \BTg\ (with \PFI) implies that every cover
of $V$ is formal-inductive.
\ecrl

\sbsc{Genetic bars are decidable}\label{genbardec}
As a first exercise in genetic induction, we prove that a genetic bar on a spraid \Vnatt\ is a decidable subset of the set of basic dots $V$. \footnote{This reflects our preference for genetic bars as a vehicle for countable-ordinal induction, since we feel that genetic bars have an intuitively manageable complexity.}   

\prp
Let \Vnatt\ be a spraid, with corresponding \Vprenatt. If $B$ is a genetic bar on a basic subspraid $\Vsuba$ of \Vnatt, then $B$ is a decidable subset of $\Vsuba$.
\eprp

\prf
By genetic induction, using \PGIg:
\be
\item[\gindz] For $a\inn V$, the genetic bar $\{\maxdot_a\}$ on $\Vsuba$ is a decidable subset of $\Vsuba$.
\item[\gindom] For $a\inn V$, let for all $b\inn\suczz(a)$, $\Bsubb$ be a genetic bar on $\Vsubb$ which is also a decidable subset of $\Vsubb$. Then $\bigcupp_{b\in\suczs(a)}\Bsubb$ is a genetic bar on $\Vsuba$ which is a decidable subset of $\Vsuba$ since for $c$ in $V$  the set $D\iz\{b\inn\suczz(a)\midd c\leqc b\}=\{b\inn\suczz(a)\midd c\inn \Vsubb\}$ is finite, so we can determine if $c\inn\bigcupp_{b\in D}\Bsubb$ or not. And $c\inn\bigcupp_{b\in\suczs(a)}\Bsubb$ is equivalent to $c\inn\bigcupp_{b\in D}\Bsubb$.
\ee
Therefore all genetic bars on basic subspraids $\Vsuba$ of \Vnatt\ are decidable subsets of $\Vsuba$.
\eprf

\rem
Genetic bars on spreads are decidable thin bars, where a thin bar is a bar $B$ for which if $a\inn B$ and $b\precc a$ then $b\notinn B$. With \BT\ one can also show the converse, that every decidable thin bar is genetic.\footnote{See \WaaArt, or do it yourself (nice exercise).}\ On spraids, due to the glue, genetic bars need not be thin. We will show that we can unglue genetic bars also. 
\erem

\sectionnb{Inductive Heine-Borel for (sub)fanns}\label{HeiBorind}

\sbsc{Inductive Heine-Borel for fanns}\label{HBforfanns}
As a second exercise in genetic induction, we first show an inductive Heine-Borel property for fanns, considered as natural space on their own. This can also be considered an inductive version of the fan theorem \FT.

\thm (\FTind)
Every genetic bar on (a basic subspraid of) a fann \Vnatt\ derived from \Vprenatt\ is finite.
\ethm

\prf
Using genetic induction:
\be
\item[\gindz] For $a\inn V$, the genetic bar $\{\maxdot_a\}$ on $\Vsuba$ is finite.
\item[\gindom] For $a\inn V$, if for all $b\inn\suczz(a)$, $\Bsubb$ is a finite genetic bar on the fann $\Vsubb$ , then $\bigcupp_{b\in\suczs(a)}\Bsubb$ is a finite genetic bar on $\Vsuba$ since $\suczz(a)$ is finite since \Vnatt\ is a fann.
\ee
Therefore all genetic bars on basic subspraids of \Vnatt\ are finite.
\eprf

\crl (Inductive Heine-Borel for fanns)
Every inductive cover of a fann \Vnatt\ has a finite subcover, and (from the pointwise perspective:) every inductive open cover of \Vnatt\ has a finite open subcover.
\ecrl

\prf
By the theorem, every inductive cover $C$ of $V$ descends from a finite genetic bar $B$ (meaning: for all $b\inn B$ there is a $c\inn C$ with $b\leqc c$), so we find a finite subset $C'$ of $C$ which  is already a basic cover of $V$. From the pointwise perspective, let $\ucov\iz\{\hattr{U_n}\midd\ninn\}$ be an inductive open cover of \Vcalt, then $C\iz\bigcupp_{n\in\N} U_n$ is an inductive cover of $V$ (see def.\,\ref{indcov}), therefore we find a finite $C'\iz\{c_i\midd i\leqq N\}\subseteqq C$ which is a basic cover of $V$. For each $i\leqq N$ we can determine an $n_i$ such that $c_i\inn U_{n_i}$, and so $\{\hattr{U_{n_i}}\midd i\leqq N\}$ is a finite open cover of \Vcalt.
\eprf

\sbsc{Inductive Heine-Borel for subfanns (including the real interval $[\alpha, \beta]$)}\label{HBforsubfanns}\hspace*{-8.8pt}\newlyne
Mostly however, we are interested in fanns as subspraids of larger spaces. We obtain an inductive Heine-Borel property for subfanns as the corollary of a basic proposition about genetic bars on subspraids:

\clearpage
\prp
Let \Wnatt\ be a subspraid (derived from \Wprenatt) of a spraid \Vnatt\ derived from \Vprenatt. Let $a\inn W$ and let $B$ be a genetic bar on $\Vsuba$. Then $B$ contains a genetic bar on the basic subspraid $\Wsuba$ of \Wnatt.
\eprp

\prf
By genetic induction:
\be
\item[\gindz]  $B\iz \{\maxdot_a\}$, then we are done.
\item[\gindom] Else, $B\iz \bigcupp_{b\in\suczs(a)} \Bsubb$ where for each $b\inn\suczz(a)$ the genetic bar $\Bsubb$ on $\Vsubb$ is such that $b\inn W$ implies that $\Bsubb$ contains a genetic bar $C_b$ on $\Wsubb$. This means that $C\iz\bigcupp_{b\in\suczs(a)\cap W}C_b$ is a genetic bar on $\Wsuba$ contained in $B$.
\ee
\eprf

\crl 
\be
\item[(i)] If $E\subseteqq W$ and $F\subseteqq V$ and $E\oicovbv F$ then $E\oicovbw (F_{\!\!\leqc}\!\capp W)$.
\item[(ii)]
(Inductive Heine-Borel for subfanns, \HBind)
If \Wnatt\ is a subfann (derived from \Wprenatt) of a spraid \Vnatt, and $C$ is an inductive cover of $W$ in \Vnatt, then $C$ contains a finite subcover of $W$ in \Vnatt. From the pointwise perspective, if \ucovt\ is an inductive open cover of \Wcalt\ in \Vnatt, then \ucovt\ contains a finite open cover of \Wcalt\ in \Vnatt.  
\ee
\ecrl

\prf
Ad (i): Let $E\oicovbv F$, then for $e\inn E$ there is a genetic bar $B$ on $\Vsube$ such that $F$ descends from $B$, which means that $B\subseteqq F_{\!\!\leqc}$. By the proposition $B$ contains a genetic bar $C$ on $\Wsube$, and trivially $C\subseteqq (F_{\!\!\leqc}\!\capp W)$. Therefore $\{e\}\oicovbw (F_{\!\!\leqc}\!\capp W)$ and since $e$ is arbitrary we see that $E\oicovbw (F_{\!\!\leqc}\!\capp W)$.
\parr
Ad (ii): Determine $c\iz\maxdotW\inn V$. Under conditions as stated, there is a genetic bar $B$ on $\Vsubc$ such that $C$ descends from $B$. By the proposition, $B$ contains a genetic bar $B'$ on the fann $\Wsubc\iz W$ which forms \Wnatt. By theorem \ref{HBforfanns}\ $B'$ is finite. Since $B'\subseteqq B$ and $C$ descends from $B$, we find that for all $b\inn B'$ there is a $c\inn C$ with $b\leqc c$. So we find a finite $C'\subseteqq C$ such that $C'$ is a cover of $W$ (in \Vnatt). 

From the pointwise perspective, let $\ucov\iz\{\hattr{U_n}\midd\ninn\}$ be an inductive open cover of $\Wcal$, then $C\iz\bigcupp_{\ninn} U_n$ is an inductive cover of $W$ (see def.\,\ref{indcov}), therefore we find a finite $C'\iz\{c_i\midd i\leqq N\}\subseteqq C$ which is a basic cover of $W$. For each $i\leqq N$ we can determine an $n_i$ such that $c_i\inn U_{n_i}$, and so $\{\hattr{U_{n_i}}\midd i\leqq N\}$ is a finite open cover of \Wcalt\ in \Vnatt.
\eprf

\rem
From this we conclude that for $\alpha\smlr\beta\inn\R$ the real interval $[\alpha, \beta]\subsett\Rnat$ has the inductive Heine-Borel property, also from the pointwise perspective.\footnote{A similar pointfree result is proved in \NegCed\ for the formal interval $[\alpha, \beta]$. The author hasn't come across a general proof of the Heine-Borel property for `fanlike' formal-topological subspaces, but his knowledge of formal topology is limited. A related monograph in many aspects is \MLof.} 
This because it is fairly easy to indicate a finitely branching full subtrea $W$ of \sigleqcRt\ such that $[\alpha, \beta]\iz\closr{W}$. Therefore, any inductive (open) cover of $[\alpha, \beta]$ is an inductive (open) cover of the subfann \Wcalt\ derived from $W$, and so contains a finite (open) cover of $[\alpha, \beta]$.
\erem

\sectionn{Inductive morphisms}

\sbsc{Inductive morphisms: definition}\label{defindmorf}
Having established some basic properties of inductive covers, we turn to inductive morphisms. Inductive morphisms are those morphisms which respect inductively acquired covers, inversely (looking at the pre-image). Therefore they inversely preserve inductive Heine-Borel properties. In metric spaces this implies that inductive morphisms are uniformly continuous on compact subspaces. 

\parr
In order to deal with trail morphisms elegantly, we turn to the unglueing of \Vnatt. But the reader can safely concentrate on refinement morphisms, as we show later on. 

\defi
Let $f$ be a \leqc-morphism between the spraids \Vnatt\ and \Wnattot, derived from \Vprenatt\ and \Wprenattot. Let $g$ be a \pthh-morphism from \Vnatt\ to \Wnattot\ (so $g$ is a \leqc-morphism from \Vunglnatt\ to \Wnattot)\footnote{Any \leqc-morphism from \Vpathnatt\ to \Wnattot\ is determined by its restriction to \Vunglnatt.}. Recall the definition of $\idstr$ in the proof of theorem \ref{pathmorphisms}.
\be
\item[(i)] For $c\inn W$ put $\finv(c)\iz\{b\inn V\midd f(b)\leqcto c\}$, for $C\subseteqq W$ put $\finv(C)\isdef\bigcupp_{c\inn C}\finv(c)$. 
\item[(ii)] For $c\inn W$ put $\ginvstr(c)=\idstr(\ginv(c))=\!\hspace*{-1.2pt}\{b\inn V\midd\!\driss b'\inn\Vungl [\idstr(b')\iz b \weddge g(b')\leqcto c\}$, for $C\subseteqq W$ put $\ginvstr(C)\isdef\bigcupp_{c\inn C\,}\ginvstr(c)$.
\item[(iii)] We call $f$ resp. $g$ an \deff{inductive morphism}\ iff for any genetic bar $G$ on $W$ we have that $\finv(G)$ resp. $\ginvstr(G)$ contains a genetic bar $H$ on $V$.
\item[(iv)] If the context is clear, we also simply write $\ginv$ for $\ginvstr$.
\ee
\edefi

\prp
With notation as above, we have that $g$ is an inductive trail morphism iff $g$ is inductive as a \leqc-morphism from \Vunglnatt\ to \Wnattot.
\eprp

\prf
In the appendix \ref{prfdefindmorf}\ we show more, namely that genetic bars on $V$ correspond to genetic bars on \Vunglt\ in a precise way. We could call this the unglueing of genetic bars on $V$.
\eprf 

The proposition together with its proof in \ref{prfdefindmorf}\ shows that we can safely restrict ourselves from now on to refinement morphisms. If at any time we need to use (inductive) trail morphisms, then we know that there is a direct correspondence between (inductive) trail morphisms on \Vnatt\ and (inductive) refinement morphisms on \Vunglnat. This explains our next:

\conv
From now on, unless stated otherwise explicitly, our morphisms are refinement morphisms.
\econv

\parr
 
Next we intend to show that if $f$ is an inductive morphism between the spraids \Vnatot\ and \Wnattot, then $A\oicovbw B$ implies $\finv(A)\oicovbv \finv(B)$, for subsets $A,B\subseteqq W$. This is relatively straightforward, if we first tackle a few extra details about genetic bars (which are nice enough in their own right).\footnote{One reason to to proceed like this is once more to keep our concepts intuitively manageable. Instead of using `for all subsets $A, B$ such that $A\oicovb B$' we only use decidable genetic bars. But it also provides for shorter proofs.}

\sbsc{Genetic bars are extendable and reducible}\label{genbarsext}
Another basic property of genetic bars is that for $c\leqc a\inn V$ we can reduce a genetic bar on $\Vsuba$ to a genetic bar on $\Vsubc$, and vice versa if we have a genetic bar on $\Vsubc$ then we can expand it to a genetic bar on $\Vsuba$. We need some minor technicalities for this.

\lem
Let \Vnatt\ be a spraid with corresponding \Vprenatt. Let $a\inn V$, then for all \ninnt\ the set $\Vsuban\iz\{b\inn \Vsuba\midd \grdd(b)\iz \grdd(a)\pluz n\}$ is a genetic bar on $\Vsuba$.
\elem

\prf
By induction on \ninnt. If $n\iz 0$ then $\Vsuban\iz \{a\}\iz\{\maxdot_a\}$ so we are done. Suppose the lemma holds for given $n$ (and for all $a\inn V$), then we show it holds for $n\pluz 1$ as well. For then we know that for each $c\inn\suczz(a)$ the set $\Vsubcn$ is a genetic bar on $\Vsubc$. And so by \gindomb\ (see def.\,\ref{indcov}) $\Vsubano\iz \bigcupp_{c\inn\suczz(a)}\Vsubcn$ is a genetic bar on $\Vsuba$.
\eprf

For $c\precc a\inn V$, the basic subspraid $\Vsubc$ is of course covered elementarily by $\{a\}$ which in this respect fulfills the same maximal role as $\{\maxdot_c\}\iz\{c\}$. But (for reasons of elegance) $\{a\}$ is not a genetic bar on $\Vsubc$. We resolve this by introducing the `reduction' of a genetic bar on $\Vsuba$ to $\Vsubc$.

\defi
Let \Vnatt\ be a spraid, and $c\leqc a\inn V$. If $B$ is a genetic bar on $\Vsuba$, then $B^{\uparrow c}\iz B\cupp\{\maxdot_c\midd\driss b\inn B[c\leqc b]\}$ is called the \deff{reduction} of $B$ to $\Vsubc$. If $B$ is a genetic bar on $\Vsubc$ then $B^{\downarrow a}\iz \{b\inn \Vsuba\midd \grdd(b)\iz\grdd(c) \weddge b\notiz c\}\cupp B$ is called the \deff{expansion} of $B$ to $\Vsuba$.
\edefi

\prp
Let \Vnatt\ be a spraid derived from \Vprenatt. Let $c\leqc a\inn V$. 
\be
\item[(i)] If $B$ is a genetic bar on $\Vsuba$ then $B^{\uparrow c}$, the reduction of $B$ to $\Vsubc$, contains a genetic bar on $\Vsubc$.
\item[(ii)] If $B$ is a genetic bar on $\Vsubc$ then $B^{\downarrow a}$, the expansion of $B$ to $\Vsuba$, is a genetic bar on $\Vsuba$
\ee
\eprp

\prf
Ad (i): If $a\iz c$ then we are done. Else, $c\precc a$. Then by genetic induction:
\be
\item[\gindz] If $B\iz \{\maxdot_a\}$, then $B^{\uparrow c}\iz \{\maxdot_a, \maxdot_c\}$ and so contains the genetic bar $\maxdot_c$ on $\Vsubc$.
\item[\gindom] Else $B\iz \bigcupp_{b\in\suczs(a)}\Bsubb$ where for all $b\inn\suczz(a)$, $\Bsubb$ is a genetic bar on $\Vsubb$ such that if $c\leqc b$, then $\Bsubb^{\hspace*{-0.3em}\uparrow c}$ contains a genetic bar on $\Vsubc$. There has to be at least one $b\inn\suczz(a)$ such that $c\leqc b$, so for such $b$ we find that $\Bsubb^{\hspace*{-0.25em}\uparrow c}$ contains a genetic bar on $\Vsubc$. And so $B^{\uparrow c}$ also contains a genetic bar on $\Vsubc$.
\ee 

Ad (ii): By induction on $n\iz\grdd(c)\minuz \grdd(a)$. If $n\iz 0$ then $c\iz a$ and we are done. Suppose the lemma holds for given $n$ (and for all $a, c\inn V$) then we show it holds for $n\pluz 1$ as well. So then let $\grdd(c)\iz \grdd(a) \pluz (n\pluz 1)$. We can locate a $b\inn\suczz(a)$ such that $c\leqc b$. Then  
by induction $B^{\downarrow b}$ is a genetic bar on $\Vsubb$. By the above lemma, for all $d\inn \suczz(a), d\notiz b$ we know that $\Vsubdn$ is a genetic bar on $\Vsubd$. From this we conclude that $V^{\downarrow a} \iz\bigcupp_{d\inn\suczz(a), d\neq b}\Vsubdn\cupp B^{\downarrow b}$ is a genetic bar on $\Vsuba$.
\eprf

For spreads, one can define $B^{\uparrow c}$ in such a way that it is itself a genetic bar. 

\sbsc{Inductive morphisms respect inductive covers (and Heine-Borel)}\label{indmorfresp}

We need a preparatory lemma, after which we can prove the basic theorem with respect to inductive morphisms. The lemma shows that for an inductive morphism, we can transfer the inductive pre-image property of the whole space to the basic subspraids. (We could have taken that as definition too, but we believe that in practice it leads to longer proofs.)

\lem
Let $f$ be an inductive morphism between the two spraids \Vnatot\ and \Wnattot, with corresponding pre-natural spaces \Vprenatot\ and \Wprenattot. Let $a\inn W$, and let $G$ be a genetic bar on $\Wsuba$. Then for all $d\inn V$: if $f(d)\inn \Wsuba$, then $\finv(G)$ contains a genetic bar on $\Vsubd$.
\elem

\prf
The proof is surprisingly involved, we give it in \ref{prfindmorfresp}. The reader is welcome to try her/himself, it should be a good exercise.
\eprf

\thm
Let $f$ be an inductive morphism between the spraids \Vnatot\ and \Wnattot, and let $A, B\subseteqq W$ where $A\oicovbw B$. Then $\finv(A)\oicovbv \finv(B)$. 
\ethm

\prf
Let $a\inn A$, then since $A\oicovbw B$ there is a genetic bar $G$ on $\Wsuba$ such that $B$ descends from $G$. Now consider $d\inn\finv(a)\subseteqq V$, which is equivalent to $f(d)\inn \Wsuba$. By the previous lemma $\finv(G)$ contains a genetic bar on $\Vsubd$. Since $\finv(G)\subseteqq\finv(B)$ we see that $\finv(B)$ contains a genetic bar on $\Vsubd$, so $\{d\}\oicovbv\finv(B)$. Since $a, d$ are arbitrary, we find that $\finv(A)\oicovbv\finv(B)$.
\eprf

\crl
\be
\item[(i)]
If $Z$ forms a subspraid \Znatot\ of \Vnatot, then the restriction $f_Z$ of $f$ to $Z$ is an inductive morphism from \Znatot\ to \Wnattot.
\item[(ii)]
If $K$ forms a subfann \Knatot\ of \Vnatot, then $f(K)$ is contained in the trea $E$ of a subfann \Ecalt\ of \Wnatto, where $f(\Kcal)$ is dense in \Ecalt\ (so $f(\Kcal)\subseteqq\Ecal$ but equality is not always the case).
\ee
\ecrl

\prf
Ad (i): it suffices to show that $\finvz(A)\oicovbz \finvz(B)$ for $A\oicovbw B$. We already know by the theorem that $\finv(A)\oicovbv \finv(B)$. Therefore trivially $\finvz(A)\oicovbv \finv(B)$. By corollary \ref{HBforsubfanns}(i) we find that $\finvz(A)\oicovbz ((\finv(B))_{\leqc}\!\capp Z)\iz\finvz(B)$, and we are done.  
\parr
Ad (ii): by (i) above we can take (the restriction of) $f$ to be an inductive morphism from \Knatot\ to \Wnattot. Determine $d\iz f(\maxdot_K)$. To see that $f(K)$ is a contained in the trea $E$ of a subfann \Ecalt\ of \Wnattot, consider for each \ninnt\ the genetic bar $G_{_{\!}n}\iz\{b\inn \Wsubd\midd \grdd(b)\iz \grdd(d)\pluz n\}$ on $\Wsubd$ (see lemma \ref{genbarsext}). For each \ninnt\ we have $\{d\}\oicovbw G_{_{\!}n}$, so by the theorem $\{\maxdot_K\}\oicovbk\finv(G_{_{\!}n})$. Since $K$ is a fann, there is a finite genetic bar $D_{_{\!}n}$ on $K$ such that $\finv(G_{_{\!}n})$ descends from $D_{_{\!}n}$. Therefore we can finitely determine the subset $E_{_{\!}n}\iz\{e\inn G_{_{\!}n}\midd \driss a\inn K\,[f(a)\leqc e]\}$ of $G_{_{\!}n}$, since $E_{_{\!}n}$ equals $\{e\inn G_{_{\!}n}\midd \driss d\inn D_{_{\!}n}\,[f(d)\leqc e]\}$. Now we can simply take $E\iz 
\bigcupp_{n\in\N}E_{_{\!}n}$ to fulfill the corollary.
\eprf

To develop constructive mathematics without Brouwer's Thesis \BT\ (and/or the weaker fan theorem \FT), one hopes that inductive covers and inductive morphisms enable reproducing much of classical compactness. We indicate some problems with this approach later on, but first we turn to the question which continuous functions between topological spaces are easily seen to be representable by an inductive morphism. 

\sbsc{Bishop-continuous \R-to-\R-functions are inductively representable}\label{biscontind}
In formal topology (see \PalCon) a \bisco\ function from \Rt\ to \Rt\ is representable by a formal mapping from the formal reals to the formal reals, and vice versa each such mapping represents a \bisco\ function. We repeat this insight in our setting, but first we give the definition of `\bisco' for real-valued functions on \Rt.

\defi (in the pointwise setting of \bish)
Let $f$ be a function from \Rt\ to \Rt, then $f$ is \deff{\bisco}\ iff $f$ is uniformly continuous on every \bish-compact (meaning complete and totally bounded) subspace of \Rt. (Equivalently, iff $f$ is uniformly continuous on every closed interval $[-n, n]$ for \ninn).
\edefi 

\prp
Let \ft\ be a \bisco\ function from \Rt\ to \Rt. Then there is an inductive morphism \fstart\ from \sigRt\ to \sigRt\ such that for all $x\inn\sigR$ we have $f(x)\equivv\fstar(x)$ (where we identify \Rt\ and \sigRt\ for convenience). Conversely, if \gt\ is an inductive morphism from \sigRt\ to \sigRt, then as a function \gt\ is uniformly continuous on each compact subspace of \Rt.
\eprp

\prf 
See appendix \ref{prfbiscontind}.
\eprf

\rem
If we replace the image space by \Rplust, then the situation is quite different. The statement that every uniformly continuous function from \zort\ to \Rplust\ is representable by an inductive morphism from \zort\ to \Rplust\ (as a natural space) is equivalent to the fan theorem \FT.
\erem

For the above remark to be precise, we first need to define the positive reals \Rplust\ as a a natural space. We likewise define the apart-from-zero reals \Raprtz. 

\defi
Let $\Rratp\iz\{[a, b]\inn\Rrat\midd a\bygr 0\}\cupp\{\maxdotR\}$. The space of the natural \deff{positive real numbers}\ \Rplusnatt\ is the natural subspace of \Rnatt\ derived from the pre-natural space $(\Rratp, \leqcR, \aprtR)$. Put $\Rratnp\isdef\{[a, b]\inn\Rrat\midd b\smlr 0 \vee a\bygr 0\}\cupp\{\maxdotR\}$. The space of the natural \deff{apart-from-zero real numbers}\ \Raprtzt\ is the natural subspace of \Rnatt\ derived from the pre-natural space $(\Rratnp, \leqcR, \aprtR)$.
\parr
Next, put $\sigRplus\isdef \{[a, b]\inn\sigR\midd a\bygr 0\}\cupp\{\maxdotR\}$ and $\sigRaprtz\isdef 
\{[a, b]\inn\sigR\midd b\smlr 0 \vee a\bygr 0\}\cupp\{\maxdotR\}$, and we see that \sigRplust\ and \sigRaprtzt\ are spraids representing \Rplusnatt\ and \Raprtzt\ respectively.
\edefi

\lem
The statement that every uniformly continuous function from \zort\ to \Rplust\ is representable by an inductive morphism from \sigzort\ to \sigRplust\ is equivalent to the fan theorem \FT.
\elem

\prf
This follows from the well-known result in \JulRic\ that \FT\ is equivalent to the statement that each uniformly continuous $f$ from \zort\ to \Rplust\ is bounded away from $0$. For completeness we detail this easy consequence in the appendix \ref{prfbiscontind}.
\eprf

This lemma foreshadows paragraph \ref{metaruss}, where we discuss pointwise problems for \bish\ which arise from the absence of \BT.

\sbsc{Inductive Baire morphisms are constructible}
As another positive example, we can construct inductive Baire morphisms in the following way. In steps, we build an inductive morphism \ft\ by determining first for all $\alpha\inn\Bnat$ a nontrivial first value $f(\alpha)(1)$ of the inductive morphism $f$. For this we need a genetic bar $B_{\maxdots}$ on \Bnat, then for each $a\succ b\inn B_{\maxdots}$ we put $f(a)\iz\maxdot$, and to each $b\inn B_{\maxdots}$ we assign an $n_b\inn\N$ and put $f(b)\iz n_b$. Next, for all $\alpha\inn\Bnat$ we determine a nontrivial refinement $f(\alpha)(2)$ of $f(\alpha)(1)$. We do so by constructing, for each $b\inn B_{\maxdots}$, a genetic bar $\Bsubb$. Then we assign to each $c\inn \Bsubb$ an $n_c\inn\N$ and put $f(c)\iz f(b)\starr n_c$ (and for all $c\leqc d\leqc b$ we put $f(d)\iz f(b)$). And so on...
\parr
Although laborious, the process above really is a construction. We can show that 
it yields an inductive morphism and that all inductive Baire morphisms are equivalent to a morphism which is constructed in the above way. (This is another partial explanation of our preference for genetic bars. Genetic bars can in our eyes be constructed in an intuitively manageable inductive way.

\sbsc{Kleene's realizability, \BTg\ and inductive morphisms}\label{kleimpind}
The situation for other spraids is complicated by the apartness relation. But if we are willing to accept Brouwer's Thesis (\BT) then any morphism is inductive. The reader can ponder on the question whether it is more elegant to inductivize all the definitions (see also the rest of this section, because we are still not done yet) or to accept \BTg. With \BTg, most of the inductive machinery that we developed here becomes superfluous.
\parr
A different indication of \BT's constructive content comes from Kleene's realizability results on intuitionistic mathematics. 
Kleene's formalization of intuitionistic mathematics is usually denoted \FIM. After a remarkable effort, Kleene proved in \KleRea\ that if we can prove the existence of a Baire function \ft\ in \FIM\,\footnote{Equivalently we can say: `if we can define \ft\ in \FIM,...'}, then this function is representable by a general recursive Baire morphism \fslangt. This is sometimes called Church's Rule (\CRo) for \FIM. Notice that \CRo\ gives a true construction for the general recursive \fslangt\ which is derived canonically from the existence proof of \ft\ in the formal system \FIM.
\parr
The same holds for a decidable thin bar $B$: if we can prove its existence in \FIM, then there is a recursive representative of $B$. In \WaaArt\ it is shown that \BT\ is equivalent to the combination of two axioms \BID\ and \BDD, where \BID\ is Kleene's decidable-Bar Induction (see \ref{defbarinduction}, and $^x$26.3 in \KleVes), and \BDD\ follows from \acoz\ (see \ref{axiomscont}, \ref{axiombdd}, and 27.1 in \KleVes). It is also shown that the concepts of `genetic bar' and `decidable thin bar' coincide under assumption of \BT.

\parr
We wish to turn this result to our advantage. Our strategy is clear: we wish to derive from a \FIM\ existence proof of a Baire function \ft, a general recursive Baire morphism \fslangt\ representing \ft\ such that in addition \fslangt\ is inductive. We are straightaway confident that Kleene's realizability fulfills this extra aspect, by the equivalence of \BT\ with the combination of the \FIM-valid axioms \BID\ and \BDD. Therefore we think that it should be possible to prove that the general recursive Kleene-realizing morphisms are also inductive.\footnote{Probably the proof should be on a meta-level, since we do not see how to formalize the genetic property within \FIM. Perhaps it can be done in an easy extension of \FIM, ... but the author is not knowledgeable enough in these matters.}

\parr
However, we are no experts on this subject, and we can only kindly invite those who are to take this issue under consideration. As stated above, we are confident that it is possible to prove this desirable extra inductive property. So we return to our strategy. We are looking for a way to construct inductive morphisms, and if our conjecture above is true, then we have found an elegant route.
\parr
For then we can use intuitionistic theory to derive existence of a morphism \ft\ between spraids, and -if the inductivity of \ft\ is not immediately apparent already- use Kleene's realizability to construct from this `abstract' \FIM-existence proof a general recursive \fslangt\ representing \ft\ which is also inductive.
\parr
Notice that we are still working entirely within \bish. It might by the above reasoning seem that we can avoid endorsing \BT\ and the continuity principle \CP\ and just incorporate the relevant conditions into our definitions. But we cannot expect this to work as easy as all that. This because in \FIM, the pointwise setting plays an integral part. Information such as `$\all x\inn\zor[f(x)\bygr 0]$' has to be seen as having been acquired inductively, and in the absence of an endorsement of \BT, it is not so easy to incorporate this type of information into the definitions. There is always a simple litmus test: whether the definitions work in \russ.

\sectionn{Pointwise problems in {BISH} and final inductivization}\label{pointprobfinind}

\sbsc{Pointwise problems in the absence of Brouwer's Thesis}\label{metaruss}
We return for a moment to the discussion started in the introduction of this section. Outside of the restricted class of Bishop-locally-compact spaces\footnote{Metric completeness is required, \Rplust\ and $(0,1)$ are not locally compact in \bishf.}, the property of being uniformly continuous on compact subspaces is a consequence of inductivity, but it doesn't by itself imply inductivity. In hindsight, it seems as if Bishop underestimated the necessity of an inductive machinery in order to build a smooth theory of compactness related to uniform continuity.
\parr
In \SchuCon, there seems to be a feeling that the inductive approach of formal topology solves these issues. Notwithstanding our own inductive treatment, we are not convinced that these issues are satisfactorily solved for pointwise settings such as \bish. And, our own pointfree machinery notwithstanding, we are not convinced that the pointwise setting should be abandoned in favour of the pointfree setting. 
Below we list some pointwise problems for \bish\ regarding our inductive approach, which as far as we can tell also hold for formal topology. It therefore seems to us that the conclusions in \SchuCon\ are too optimistic.

\mthm

In \russ\ (and by implication \bish) we have the following problems regarding pointwise use of inductive definitions:
\be
\item[P$_1$] Uniform continuity of a function $f$ does not imply that there is an inductive morphism representing $f$. Counterexamples can be given even for uniformly continuous functions from \zort\ to \Rplus. Shortly put: uniform continuity does not imply inductive representability.
\item[P$_2$] Weak completeness\footnote{The property for a located subset $A$ of a metric space \xdt\ that for all $x\inn X$: if $x\aprt a$ for all $a \inn A$, then $d(x, A)\iz\inf(\{d(x,a)\midd a\inn A\})\bygr 0$.}\ of a compact space is not preserved under inductivity. In \russ, even for an inductive morphism from \zort\ to \zor, the image of a compact subspace may be strongly incomplete.
\item[P$_3$] Inductive representability is not preserved under the restriction of a function to its pointwise image space. This follows from the counterexamples for P$_1$, since every uniformly continuous function from \zort\ to \Rt\ is inductively representable by proposition \ref{biscontind}. Therefore we can expect problems with the reciprocal function $x\rightarrow\frac{1}{x}$, and must continually address these problems by adapting our definitions.
\ee

Therefore in \bish, the desirable properties associated with the problems above cannot be shown to hold without further assumptions. In fact, assertion of any of these properties implies the fan theorem \FT. 
\emthm

\prf
The proof is given in the appendix \ref{prfmetaruss}. It is derived from the construction of a compact subspace \contrCant\ of \zort\ such that if we write \Canzort\ for the standard embedding of Cantor space in \zort, we see: $\dR(\contrCan, \Canzor)\iz 0$ and yet in \russ\ we also have $\dR(x, \Canzor)\bygr 0$ for $\all x\inn\contrCan$.  The `ContraCantor space' \contrCant\  is defined using the Kleene Tree.
\eprf

\sbsc{Inductivization of the basic definitions in the absence of \BTg}\label{indbasdef}
A final important issue for inductivity concerns our basic definitions. We started out by defining natural spaces and pre-natural spaces, using basic dots and a pre-apartness on the basic dots. If we accept \BTg, then these definitions suffice for building a classically valid theory in which compactness and inductivity are naturally incorporated, and which largely resembles \intu. 
\parr
However, if one wishes to build an inductive theory for natural spaces without accepting \BTg, then one has to `inductivize' the basic definitions. (Without \BTg\ pointwise problems persist though, by thm.\,\ref{metaruss}.). 

\rem To do this thoroughly it seems attractive to abandon pointwise notions, since they generally require a pointfree translation to use inductive information. This starts already with the definition of topology itself. The `arbitrary union' requirement (\Topth) is pointwise in our constructive setting, see \ref{defnattop}(iii). It is possible to remedy this with a pointfree notion `toipology', and to thus develop a completely pointfree version of natural topology (say `natural toipology'). But we believe one should remember that the concept of `point' is actually the same as the concept of a countable sequence and the concept of countable infinity. These concepts are already heavily involved in the very definition of $V$. Therefore the foundational gain to the author seems smaller than one might think at first glance, and there is a price to pay in terms of readability. For this reason, we will continue here with the basic pointwise notions.
\erem

For inductivization of the basic definitions, it suffices to exact that covers which are required to exist by the basic definitions are inductive. These covers originate from the definition of points and \aprt-open sets. 

\parr
From the definition of points, for a spraid \Vnatt\ derived from \Vprenatt\ and basic dots $a\aprt b$ in $V$ we know that $\all x\inn\Vcal\driss \ninn\,[x_n\aprt a \vee x_n\aprt b]$. This means that the subset $C\iz\{c\inn V\midd c\aprt a\vee c\aprt b\}$ is a bar on \Vnat. For a really smooth working of inductivity, we should know that $C$ is an inductive bar, and adapt our definition accordingly.
\parr
Other covers which arise directly from our basic definitions are the ones associated with \aprt-open subsets $\Ucal\subseteqq\Vcal$. Recall that \Ucalt\ is \aprt-open iff for each $ x\inn\Ucal$ and $y\inn\Vcal$ we can determine (non-exclusively) $x\aprt y$ and/or  there is an \minnt\ with $\hattr{y_m}\subseteqq\Ucal$. But for given $x\inn\Ucal$ this is equivalent to saying that the set $B\iz\{b\inn V\midd b\aprt x_{\grds(b)}\vee\hattr{b}\subseteq\Ucal\}$ is a bar on $V$. For a really smooth working of inductivity, we should know that this bar is inductive, and we should add this to our definition of \aprt-open accordingly.
\parr
Our solution is to add the word `inductive(ly)' to the original definitions:

\defi
Let \Vnatt\ be a spraid derived from \Vprenatt. 
An \aprt-open subset $\Ucal\subseteqq\Vcal$ is called \deff{inductively open}\ in \Vnatt\ iff for any $x\inn\Ucal$ the set $\BUx\iz\{b\inn V\mid b\aprt x_{\grds(b)}\vee\hattr{b}\subseteq\Ucal\}$ is an inductive bar on $V$. We then write: \Ucalt\ is \sucz-open. 
The collection of \sucz-open subsets of \Vcalt\ is called the \deff{inductive apartness topology}\ on \Vnatt, notation \Topsuczt.
\parr
We call \Vnatt\ an \deff{inductive}\ spraid (we write: `a \sucz-spraid') iff:
\be
\item[(i)] For all $a,b$ in $V$ with $a\aprt b$, the set $C\iz\{c\inn V\midd c\aprt a\vee c\aprt b\}$ is an inductive bar on $V$. 
\item[(ii)] The inductive apartness topology coincides with the apartness topology, that is $\Topaprt\iz\Topsucz$.
\ee
Finally, some abbreviating notation will be also useful.
For subsets $A, B$ of $V$ we write $A\aprt B$ iff $a\aprt b$ for all $a\inn A, b\inn B$. We write $A\touch B$ iff $a\touch b$ for some $a\inn A, b\inn B$. We shortly write $a\aprt B$, $a\touch B$ for $\{a\}\aprt B$, $\{a\}\touch B$ respectively.
We define: $^nV\isdef \{a\inn V\midd \grdd(a)\iz n\}$, for \ninn.
\edefi

\rem It follows from \BT\ that every spraid is inductive. We believe that any spraid which is \FIM-definable\ will be inductive, by our remarks in \ref{kleimpind}. One easily sees for instance that Baire space and \sigRt\ (the spraid representing \Rt) are inductive. We can prove that a metric spraid is inductive when its trea is given by shrinking metric balls as basic dots, where each dot of $\grdd(n)$ is a metric ball of diameter less than $2^{-n}$ and apartness of dots implies a positive distance between the dots. This shows that complete metric spaces (by a standard completion procedure) are homeomorphic to an inductive spread, see paragraph \ref{compmetind}.
In the absence of \BT\ it seems practical to demand inductiveness by definition.
\erem

For completeness we still need to show that \Topsuczt\ is indeed a topology. We add to this a simple extension of the first property of inductive spraids (we need this extended property later on in our metrization theorem): 

\prp
\be
\item[(i)] For a spraid \Vnatt\ with corresponding \Vprenat, the collection \Topsuczt\ is a topology which is refined by \Topaprt.
\item[(ii)] Let \Vnatt\ be an inductive spraid derived from \Vprenatt. 
Then for finite subsets $A\aprt B$ of $V$, the subset $C\iz\{c\inn V\midd c\aprt A\vee c\aprt B\}$ is an inductive bar on \Vnatt.
\ee
\eprp

\prf
The proof is a bit involved, we give it in the appendix see \ref{prftopsucz}. 
\eprf

\sbsc{Inductivization of other definitions}\label{indothdef}\hspace*{-1.3pt}
In the absence of \BTg, for completeness one should also inductivize some other definitions given in earlier sections. The most important definition in this respect is the definition of `fanlike' since in \russ\ (where \BTg\ fails) Baire space is fanlike under the non-inductivized definition, by proposition \ref{canbairus}\ . The isomorphism that we constructed in the proof is an example of a non-inductive morphism.

To define `\sucz-fanlike' for spraids we could use the inductive morphisms already defined,
but we wish to generalize this definition to natural spaces, and to inductivize the definition of `spreadlike' also. In general, natural spaces need not be given by an apartness on a trea, so we must first expand our definition of `inductive morphism' to include such natural spaces, and for this we use the formal inductive covering relation \icovbt\ (see def.\,\ref{formindcov}) in the obvious way.

\defi
Let $f$ be a morphism from a natural space \Vnatot\ to a natural space \Wnattot. We say that $f$ is \deff{inductive}\ iff for all $A, B\subseteqq W$ where $A\icovbw B$ we have that $\finv(A)\icovbv \finv(B)$ (this agrees on spraids with the earlier definition by theorem \ref{indmorfresp}\ and proposition \ref{indcov}). We say that \Vnatt\ is \deff{\sucz-fanlike}\ resp. \deff{\sucz-spreadlike}\ iff there is an inductive isomorphism from \Vnatt\ to a \sucz-fan resp. a \sucz-spread\ (with an inductive inverse).
\edefi

We can now show: Baire space is not \sucz-fanlike. In fact any spraid which is \sucz-fanlike contains a subfann on which the identity is an isomorphism with the whole space, see corollary \ref{indmorfresp}(ii) and its proof.

\rem
In our \bish\ framework therefore, Heine-Borel compactness is best characterized as `$\!$\sucz-fanlike'. We will use this to phrase a natural-topology version of Tychonoff's theorem in the next section.
\erem

\sbsc{Complete metric spaces have an inductive representation}\label{compmetind}\hspace*{-.3pt}
We aim to show the viability of the concept `inductive spraid', by proving that a complete metric space has a representation as a \sucz-spraid. (In fact it suffices to look closely at the proof of theorem \ref{sepmetnat}, and adapt the natural space constructed there).

\thm
Every complete metric space \xdt\ is homeomorphic to a \sucz-spraid.
\ethm

\prf
We give the proof in the appendix \ref{prfcompmetind}. The idea is not difficult: in the proof of theorem \ref{sepmetnat}\ we constructed, for a complete metric space \xdt, a natural space \Vnatt\ homeomorphic to \xdt. If we look more carefully, we see that its trail space \Vpathnatt\ contains a (homeomorphic) \sucz-spread.
\eprf

\sectionnb{(In)finite products and Tychonoff's theorem}\label{prodtych}

\sbsc{Products, lazy convergence and isolated points}\label{prodlazyiso}
For a (\bish) natural-topological version of Tychonoff's theorem,
we must define (in)finite-product spaces. The idea is straightforward, but there are technical issues related to our `lazy convergence' of points, and isolated points, see also \ref{disclazycon}. 

\parr

Isolated points are points which as a set are open in the topology, for instance in a one-point natural space. Generally, consider a natural space \Vnatt\ derived from \Vprenatt, where $a\inn V$ is 'isolated': for all $b\aprt c\inn V$ one has $a\aprt b \vee a\aprt c$. Then both $a$ and the infinite sequence $\underline{a}\iz a, a, a,\ldots$ have the \aprt-characteristics of a point, except (and this is crucial) that their image under a morphism generally does not share those point-characteristics.

\defi
Let \Vnatt\ be a natural space derived from \Vprenatt. Put $\Viso\izdef$ $\{a\inn V \midd \all b,c\inn V [b\aprt c\rightarrow (a\aprt b \vee a\aprt c)]\}$, then elements of \Visot\ are \deff{isolated} basic dots. We call \Vnatt\ a \deff{decidable-isolation} space iff \Visot\ is a decidable subset of $V$, and a \deff{perfect}\ space iff $\Viso\iz\emptyy$.
\edefi

An obvious try for a definition of the product of two natural spaces \Vnatot\ and \Wnattot\ (derived from \Vprenatot, \Wprenattot) is to take $V\Timez W$ as the set of basic dots, and to define $(c,d)\leqc (a,b)$ iff $c\leqco a \weddge d\leqcto b$ and $(c,d)\aprt (a,b)$ iff $c\aprto a \vee d\aprtto b$. This works fine for perfect spaces \Vnatot\ and \Wnattot. But for isolated $a\inn V$ and $q\inn\Wcal$ the sequence $(a,q_0), (a,q_1), (a,q_2),\ldots\iz (\underline{a},q)$ becomes a point in $\Vcal\Timez\Wcal$. This is unwanted, since we need the coordinate projections $\pi_0: (x,y)\rightarrow x$ and $\pi_1:(x,y)\rightarrow y$ to be morphisms.

\sbsc{Definition of (in)finite-product spaces}\label{defprod}
We look to retain the elegance of the simple approach when possible, and yet ensure at all times that the coordinate projections are morphisms. For this we concentrate first on the most important class of spreads/spraids, since we wish the product to be a spread/spraid as well.

\parr
We could restrict ourselves to spraids w.l.o.g., but other representations interest us also.
The simple approach works for perfect spaces. They form a subclass of the decidable-isolation spaces where a slightly adjusted approach works in all but trivial cases. Finally we give a general definition which works for all spaces, but is somewhat less elegant.

\defi
Let $((\Vcal_n, \Topaprtn))_{n\in\N}$ be natural spaces derived from the corresponding pre-natural $((V_n, \aprtn, \leqcn))_{n\in\N}$ with maximal dots $(\maxdot_n)_{n\in\N}$.
We form $\Vpin=\Pi_{i \leq n}V_i=\{(a_0, \ldots, a_{n})\midd \all i\leqq n [a_i\inn V_i]\}$ for \ninnt, and also 
$\Vpi\iz\bigcupp_{n\in\N} \Vpin$.

\parr 

For $n\leqq\minn, a\inn\Vpin, b\inn\Vpim$ put $a\aprtpi b$ iff there is $i\leqq n$ with $a_i\aprti  b_i$. 
Also:

\parr

$b\leqcPi a$ iff $\all i\leqq n [b_i\leqci a_i]$

$b\precPio a$ iff $b\precPi a\weddge \all i\leqq n [a_i\inn {\Viso}_i\rightarrow b_i\precci a_i]$ ; $b\precPioo a$ iff $b\precPio a\weddge m >n$

$b\precPic a$ iff $\all i\leqq n [b_i\precci a_i]$ ; $b\precPicc a$ iff $b\precPic a\weddge m >n$

\parr
Then \leqcPiot\ is decidable when $((\Vcal_n, \Topaprtn))_{n\in\N}$ are all decidable-isolation spaces. If $((\Vcal_n, \Topaprtn))_{n\in\N}$ are perfect spaces then \leqcPiot\ equals \leqcPit. 

\parr
For $((\Vcal_n, \Topaprtn))_{n\in\N}$ as above 
the \deff{simple finite product}\ $\Pi_{i\leq n}(\Vcal_i, \Topaprti)$ is the natural space derived from $(\Vpin, \aprtpi, \leqcPi)$ with maximal dot $\maxdotPin=\maxdot_0, \ldots, \maxdot_{n}$. We also write $\Vcal_0\Timez\Vcal_1\ldots\Timez\Vcal_n$ for $\Pi_{i\leq n}(\Vcal_i, \Topaprti)$.

The \deff{simple infinite product}\ $\Pi_{n\in\N}(\Vcal_n, \Topaprtn)$ is the natural space derived from  $(\Vpi, \aprtpi, \leqcPi)$ with maximal dot $\maxdotPi\isdef\,_{{\,}}\maxdot_0$. We often omit the word `simple'. For perfect spaces, the simple (in)finite products suffice.

\parr

When $((\Vcal_n, \Topaprtn))_{n\in\N}$ are spraids, then $(\Vpin,\leqcPi)$ and $(\Vpi, \leqcPi)$ are treas. Yet $\Pi_{i\leq n}(\Vcal_i, \Topaprti)$ and $\Pi_{n\in\N}(\Vcal_n, \Topaprtn)$ usually fail to be a spraid (also see \ref{sprddefcom}) since for a spraid, infinite $\precc$-trails have to define a point. So we put $\Vpisign \izdef$ $\{a\inn\Vpin\midd \all i,j\leqq n [\grdd(a_i)\iz \grdd(a_j)]\}$, $\Vpisig\isdef\bigcupp_{n\in\N}\{a\inn \Vpisign\midd \all i\leqq n [\grdd(a_i)\iz n]\}$. 

\parr

The \deff{finite-product spraid}\ \PisigVnt\ is the spraid (or spread) derived from $(\Vpisign, \aprtpi, \leqcPi)$ with maximal dot $\maxdotPin$. We frequently write $\Vcal_0\Timsig\Vcal_1\ldots\Timsig\Vcal_n$ for $\PisigVn$, and for $\Vpisign$ we also write $V_0\Timsig V_1\ldots\Timsig V_n$.

The \deff{infinite-product spraid}\ \PisigVNt\ is the spraid (or spread) derived from $(\Vpisig, \aprtpi, \leqcPi)$ with maximal dot $\maxdotPi\iz\maxdot_0$.

\parr

When $((\Vcal_n, \Topaprtn))_{n\in\N}$ are decidable-isolation spaces, then
the \deff{finite $\circ$-product}\ $\PiO_{i\leq n}(\Vcal_i, \Topaprti)$ is the natural space derived from $(\Vpin, \aprtpi, \leqcPio)$ with maximal dot $\maxdotPin$. 
We also write $\Vcal_0\Timiso\Vcal_1\ldots\Timiso\Vcal_n$ for $\Pi_{i\leq n}(\Vcal_i, \Topaprti)$.

The \deff{infinite $\circ$-product}\ $\PiO_{n\in\N}(\Vcal_n, \Topaprtn)$ is derived from $(\Vpi, \aprtpi, \leqcPioo)$ with maximal dot $\maxdotPi$. We can replace $\leqcPioo$ with $\leqcPio$ when $\all\ninn\driss m\bygr n[\maxdot_m\notinn{\Viso}_m]$.

\parr For $((\Vcal_n, \Topaprtn))_{n\in\N}$ as above 
the \deff{strict finite product}\ $\PiC_{i\leq n}(\Vcal_i, \Topaprti)$ is the natural space derived from $(\Vpin, \aprtpi, \leqcPic)$ with maximal dot $\maxdotPin$.  
We also frequently write $\Vcal_0\Timstr\Vcal_1\ldots\Timstr\Vcal_n$ for $\Pi_{i\leq n}(\Vcal_i, \Topaprti)$.

The \deff{strict infinite product}\ $\PiC_{n\in\N}(\Vcal_n, \Topaprtn)$ is derived from $(\Vpi, \aprtpi, \leqcPicc)$ with maximal dot $\maxdotPi$. We can replace $\leqcPicc$ with $\leqcPic$ when $\all\ninn\driss m\bygr n[\maxdot_m\notinn{\Viso}_m]$.
We often omit the word `strict'.
\edefi

The reader may verify that these definitions are valid. Some key properties are stated in the next paragraphs, one of which contains a constructive form of Tychonoff's theorem. 

\sbsc{Projections and the Tychonoff topology}\label{projtych}
(Spaces and notation continued from the previous paragraph:) we wish to throw some light on the relation between the product-apartness topology and the Tychonoff (product) topology. Classically, the Tychonoff topology is the 'natural' product topology, it being the coarsest topology rendering the coordinate projections continuous. For natural spaces however, these projections moreover need to be morphisms.  
For (in)finite products of (weak) basic neighborhood spaces (by far the most important class) the two topologies coincide, which is not surprising by theorem \ref{basicneighbor}. In \bish, for other spaces this remains elusive to us, although the product-apartness topology refines the Tychonoff topology. In \class\ and \intu\ we can show that also for the (in)finite product of star-finitary spaces (see def.\,\ref{starfinsprd}) the two topologies coincide.

\parr

Therefore, calling an (in)finite product `faithful' if the apartness topology and the Tychonoff topology coincide, there will certainly be no easy examples of unfaithful products. Faithfulness is an important property for products, which we now define along with `weak basic neighborhood space' and the relevant coordinate projections $\pi_i$, using the spaces and (in)finite products from the previous paragraph.

\defi (notation from \ref{defprod})
Let $a\iz a_0, \ldots, a_s\inn\Vpi$, for $i\leqq s$ put $\pi_i(a)\isdef a_i$ and for $i\bygr s$ put 
$\pi_i(a)\isdef \maxdot_i$. For $i\leqq m$ and points $x\inn\Pi_{j\leq m}\Vcal_j$ and $y\inn\Pi_{n\in\N}\Vcal_n$, put $x_{[i]}\isdef\pi_i(x)=\pi_i(x_0), \pi_i(x_1),\ldots$ and $y_{[i]}\isdef\pi_i(y)\iz\pi_i(y_0), \pi_i(y_1),\ldots$.

\parr
Let \Wnatt\ be any of the products defined in \ref{defprod}, then \Wnatt\ is \deff{faithful} iff the product-apartness topology \Topaprtt\ and the Tychonoff topology coincide.
\parr
A natural space \Vnatt\ is a \deff{weak basic neighborhood space} iff for every $x\inn\Vcal$ there is an $x'\equivv x$ such that ${x'}\!_n$ is a basic neighborhood of $x$ for every \ninnt\ (so every basic neighborhood space is a weak basic neighborhood space by prp.\,\ref{basicneighbor}).
\edefi

\prp (about def.\,\ref{defprod})
\be
\item[(i)] For $((\Vcal_n, \Topaprtn))_{n\in\N}$ and the defined simple (in)finite products, the appropriate projections $\pi_i$ are morphisms iff all the $\Vcal_j$'s involved in the product are perfect spaces.
\item[(ii)] For the other defined (in)finite-product spaces, the appropriate projections $\pi_i$ are always morphisms.
\item[(iii)] For (weak) basic-neighborhood spaces $((\Vcal_n, \Topaprtn))_{n\in\N}$, the (in)finite products are faithful (weak) basic neighborhood spaces. 
\item[(iv)] In \class\ and \intu, also the (in)finite products of star-finitary spaces are faithful.
\item[(v)] For other spaces, the natural topology of the (in)finite products refines the Tychonoff topology.\footnote{It remains elusive to us whether unfaithful products exist in \classf, \russf, or \intuf.} 
\ee
\eprp

\prf
We leave (i) and (ii) as exercise for the reader, and prove (iv) in \ref{prfinttych}. For (iii) first let $((\Vcal_n, \Topaprtn))_{n\in\N}$ be basic-open spaces. Let \ninnt, and consider $a\iz a_0,\ldots, a_n\inn\Vpin$. We hold that $\hattr{a}$ is open in the relevant defined products (see def.\,\ref{defprod}), since the set $\hattr{a}$ equals $\bigcapp_{i\leq n}{\pi_i}^{-1}(\hattr{a_i})$ where the $\pi_i$'s are morphisms by (ii) (and so continuous), and the $\hattr{a_i}$'s are open. This also shows that the basic open set $\hattr{a}$ is open in the Tychonoff topology, therefore the natural topology and the Tychonoff topology coincide.

\parr

So for basic-open $((\Vcal_n, \Topaprtn))_{n\in\N}$ the relevant (in)finite products are again basic-open. For basic neighborhood spaces $((\Vcal_n, \Topaprtn))_{n\in\N}$ we can now use the coordinate-wise isomorphisms with a basic-open space to see that the relevant (in)finite products are again basic neighborhood spaces.

\parr For weak basic neighborhood spaces $((\Vcal_n, \Topaprtn))_{n\in\N}$, a very similar argument using \aczo\ can be given. Let $x$ be in an (in)finite product, then for each relevant $i$ there is $y_{(i)}\equivvi x_{[i]}$ with $(y_{(i)})_n$ a basic neighborhood of $x_{[i]}$ for all \ninnt. Using \aczo\ we can determine a $y\equivvPi x$ in the same product such that $(y_{[i]})_n$ is a basic neighborhood of $x_{[i]}$ for all \ninnt. This implies that $y_n$ is a basic neighborhood of $x$ for all \ninnt. 

If our $x$ is in some open \Ucalt, then we find $\hattr{y_n}\subseteqq\Ucal$ for some \ninnt. However, $\closr{y_n}$ equals $\bigcapp_{i\leq M}{\pi_i}^{-1}(\closr{\pi_i(y_n)})$ for some \Minnt, and so is a neighborhood of $x$ in the Tychonoff topology. This shows that the product is faithful.

\parr

Finally (v): for other spaces, the natural topology of the relevant (in)finite products renders the projections as morphisms and therefore continuous, so it trivially refines the Tychonoff topology. 
\eprf

\sbsc{Spraid products and Tychonoff's theorem}\label{sprdtych}
Our prime purpose with this section was to define (in)finite products for spraids. Having done so, we phrase some (expected) nice properties of these products as a theorem. There is one important property which we can so far only prove in \bish\ if the product is faithful, namely that if $((\Vcal_n, \Topaprtn))_{n\in\N}$ are \sucz-spraids, then a faithful (in)finite product is also a \sucz-spraid. 

\parr

Specified to \sucz-fanns, this yields a natural-topology version of Tychonoff's theorem (since we characterized Heine-Borel compactness as `\sucz-fanlike' in \ref{indothdef}). In keeping with the rest of the chapter, to show in \bish\ that the faithful product of inductive spraids is again inductive requires quite some work. We start with the following definition and lemma:

\defi
Let $((\Vcal_n, \Topaprtn))_{n\in\N}$ be spraids derived from the corresponding pre-natural $((V_n, \aprtn, \leqcn))_{n\in\N}$. 
Given subsets $A_i\subseteqq V_i$ for $i\leqq \ninn$, we put 

\parr

$\Pisig{n}\hspace*{.05em} A_i\iz A_0\Timsig \ldots\Timsig A_n\isdef \{a\inn\Vpisign\midd \allz i\leqq n\,[\pi_i(a)\inn (A_i)_{\leqc}] \weddge \dris j\leqq n\, [\pi_j(a_j)\inn A_j]\}$

$\overline{\Pisig{n}\hspace*{.05em} A_i}\isdef \{a\inn\Vpisig\midd \allz i\inn\N\,[\pi_i(a)\inn (A_i)_{\leqc}]\}$

\parr

(this aligns with the definition of $V\Timsig W$ in def.\,\ref{defprod}).
\edefi

\lem
(Notations as above) Let $a\inn V_0, b\inn V_1$ and let $G, H$ be genetic bars on $(V_0)_a, (V_1)_b$ respectively. Then $G\Timsig H$ is an inductive bar on $(V_0)_a\Timsig (V_1)_b$.
\elem

\prf
By double genetic induction, see \ref{prfsprdtych}. 
\eprf

\crl
For $i\leqq n$ let $B_i$ be an inductive bar on $V_i$. Then $\Pisig{n}B_i$ is an inductive bar on $\Vpisign$ and $\overline{\Pisig{n}B_i}$ is an inductive bar on $\Vpisig$. 
\ecrl

\thm 
 For (star-finite) spreads (spraids, fanns) $((\Vcal_n, \Topaprtn))_{n\in\N}$ the products \PisigVnt\ and \PisigVNt\ are in turn (star-finite) spreads (spraids, fanns). For \sucz-spraids the products \PisigVnt\ and \PisigVNt, if faithful, are in turn inductive. 
\ethm

\prf
We prove that for \sucz-spraids $((\Vcal_n, \Topaprtn))_{n\in\N}$ the products \PisigVmt\ and \PisigVNt, if faithful, are in turn \sucz-spraids. The other (combinations of) properties are left as an exercise for the reader.

So let \PisigVmt\ be a faithful finite product of the \sucz-spraids $((\Vcal_n, \Topaprtn))_{n\in\N}$. To show that \PisigVmt\ is inductive, first let $a\aprtpi b$ in $\Vpisigm$. We must show that the bar $B\iz\{c\inn\Vpisigm\midd c\aprtpi a\vee c\aprtpi b\}$ is inductive. There is $i\leqq m$ with $a_i\aprti b_i$, so $C_i\iz\{c\inn V_i\midd c\aprti a_i\vee c\aprti b_i\}$ is an inductive bar on $V_i$. For $j\notiz i$ let $C_j\iz\{\maxdot_j\}$, then by the above lemma and corollary $C\iz C_0\Timsig\ldots\Timsig C_n$ is an inductive bar on $\Vpisigm$. It is easy to see that $C\subseteqq B$ so $B$ is inductive as well.

\parr
Second, let \Ucalt\ be open in \PisigVmt. We must show that \Ucalt\ is \sucz-open. For this let $x\inn\Ucal$. To see that the bar $\BUx\iz\{b\inn\Vpisigm\mid b\aprtpi x_{\grds(b)}\vee\hattr{b}\subseteq\Ucal\}$ is an inductive bar on $\Vpisigm$, determine \sucz-open $\Ucal_i\nni x_{[i]}$ in each $\Vcal_i$ such that $\bigcapp_i {\pi_i}^{-1}(\Ucal_i)\subseteq \Ucal$ (remember \PisigVmt\ is faithful). Then for $i\leqq m$ the bars $\BUxi\iz\{b\inn V_i\mid b\aprti (x_{[i]})_{\grds(b)}\vee\hattr{b}\subseteq\Ucal_i\}$ are inductive and monotone. 

\parr
For $i\notiz j\leqq m$ let $D_{i,j}\iz\{\maxdot_j\}$ and $D_{i,i}\iz\BUxi$, then by the above lemma and corollary $D_i\iz D_{i,0}\Timsig\ldots\Timsig D_{i,m}$ is a monotone inductive bar on $\Vpisigm$. By lemma \ref{prftopsucz}, the bar $D\iz\bigcapp_i D_i$ is a monotone inductive bar on $\Vpisigm$. Clearly $D\subseteqq\BUx$, so $\BUx$ is inductive as well and we are done for \PisigVmt.

\parr

Finally, let \PisigVNt\ be faithful. We can copy the arguments above for \PisigVmt, replacing $C, D$ with $\overline{C}, \overline{D}$, to see that \PisigVNt\ is inductive.
\eprf

\rem
Our \bish\ countable version of Tychonoff's theorem is that for \sucz-fanns $((\Vcal_n, \Topaprtn))_{n\in\N}$ the products \PisigVnt\ and \PisigVNt, if faithful, are in turn \sucz-fanns. In \class\ and \intu\ the inductive properties follow from \BT, and the countable Tychonoff's theorem is easy.
\erem

\chepter{Chapter four}{Metrizability, and natural topology in physics}{Metric spaces are perhaps the most important topological spaces, and the question of metrizability has been fundamental for the development of classical topology. In constructive topology, the concept of `located in' is usually tied to a metric. We discuss the alternative `(topologically) strongly halflocated in', which is transitive. We give the basic definitions, and show first that there are interesting non-metrizable natural spaces. We can define Silva spaces also in the context of natural topology, and infinite-dimensional Silva spaces are non-metrizable.
\parr
We translate the intuitionistic metrizability of star-finitary apartness spreads to \bish, using the inductive definitions of the previous chapter. The resulting metrization theorem for star-finitary natural spaces closely resembles classical metrization results for strongly paracompact spaces. Star-finite metric developments are also of interest for efficient computing with complete metric spaces.
\parr
In the last section, we again discuss which mathematics are suited for physics. We also present an informal two-player model of \intu\ in \russ, called `Limited Information for Earthlings' (\LIfE). 
}

\sectionnb{Metric spaces and (non-)metrizability}\label{metspamet}

\sbsc{Metric spaces and (non-)topological notions}\label{metnat}\hspace*{-.25pt}
The best-studied topological spaces in constructive mathematics are metric spaces, since a metric gives a lot of constructive traction in defining compactness, continuity and open covers. Still, different metrics can give rise to the same topology.\footnote{In \raisebox{0pt}[0ex][0ex]{\WaaThe}\ these metrics are then called `equivalent' and `strongly equivalent' if they give rise to the same \raisebox{0pt}[0ex][0ex]{Cauchy-sequences}.}\ Conversely many notions defined for a metric space are not `topological' -- that is they are not necessarily preserved under homeomorphisms. Metrical completeness is the standard example of a non-topological notion, which leads us to the definition of `topologically complete' (namely `homeomorphic to a complete metric space').\footnote{E.g. `locally compact' in \bishf\ is a non-topological notion relying on metrical completeness. We believe this to be unwieldy for topology.}

\parr
A more important example of a non-topological constructive notion is that of a subset $A$ being `located' in a metric space \xdt, meaning that the distance to $A$ can be calculated for any $x\inn X$. Although the notion `located in' is extensively used, it has substantial drawbacks especially in the context of topology. We discuss some alternatives given in \WaaThe\ which resolve these drawbacks. A form of locatedness is important also for natural topology. We briefly discuss this here, to expand on in the appendix.

\sbsc{Various concepts of locatedness}\label{varconloc}
What are these drawbacks of the concept `located in'?
First of all, the notion is not transitive, which is unpractical when working with extensions and subspaces of \xdt. Second, even for a closed located $A\subsett X$, the notion gives little handhold for $x\inn X$ to find $a\inn A$ such that $x\aprt a$ implies $x\aprt A$, which is an important prerequisite for many constructions involving $A$. Thirdly, as mentioned, the notion is non-topological and this means we cannot use it easily in the context of topology.

\parr

In \WaaThe\ several alternative notions are given in \bish, of which `strongly halflocated in' seems more fruitful than `located in'. The notion is transitive, unlike `located in'. To see that it gives results, note its use in the \bish-proof of the Dugundji extension theorem in \WaaThe. Another result is that every complete metric space can be isometrically embedded in a normed linear extension such that it becomes strongly halflocated in this extension -- and where we know of no general proof that it is located. 
\parr
`Topologically strongly halflocated' is seen in \intu\ to be equivalent on complete metric spaces to a topological locatedness property called `strongly sublocated in'. This notion can also be defined for the apartness topology of natural spaces, and seems relevant. We repeat the definitions for these properties in the appendix \ref{appvarconloc}, in order not to lose reading pace.

\sbsc{Natural metric spaces}\label{defnatmet}
Let us start by defining natural metric spaces, as natural spaces with a metric (respecting the apartness). The induced metric topology need not coincide with the apartness topology however. For (a suitable natural representation of) a complete metric space we can show that the apartness topology coincides with the metric topology (see thm.\,\ref{sepmetnat}, which is analogous to \Rnattopt\ being homeomorphic to \Rdr). 

Weakly metrizable natural spaces are those spaces on which we can construct a metric respecting the apartness, and they are metrizable if the metric topology coincides with the apartness topology. Complete metric spaces are thus metrizable. We will see that there are weakly metrizable spaces which are non-metrizable. 

\defi
For a natural space \Wnatt\ a morphism $d$ from $\Wcal\timez\Wcal$ to \Rnatt\ is called a \deff{metric} on \Wnatt\ iff for each $\alpha,\beta,\gamma\inn\Wcal$ we have 

\be
\item[(i)]$d(\alpha,\beta)\equivv d(\beta,\alpha)\geqq 0$
\item[(ii)]$d(\alpha,\beta)\bygr 0$ iff $\alpha\aprt\beta$
\item[(iii)]$d(\alpha,\gamma)\leqq d(\alpha,\beta)\pluz d(\beta,\gamma)$
\ee

If \Tcalt\ is a topology on \Wcal, then we say that $d$ \deff{metrizes}\ \Tcalt\ iff the metric topology determined by $d$ is the same as \Tcal.
\Wnatt\ is called \deff{weakly metrizable} if there is a metric $d$ on \Wnatt, and \deff{metrizable}\ if this metric metrizes \Topaprt.
\edefi

In this section we show that a vast class of natural spaces is metrizable. This class is not limited to (natural representations of) locally compact spaces, in fact the `star-finitary' property that we use resembles (strong) paracompactness very closely. Star-finite open covers are a useful tool in general topology, and accessible to constructive topology in the context of metric spaces
(see \WaaThe).
\parr
Before turning to metrizable spaces, let us consider the question whether we can find interesting and important non-metrizable natural spaces. 

\sbsc{Are there interesting non-metrizable spaces?}\label{nonmetspaces}\hspace*{-1.8pt}
Apartness topology (and possible metrizability) is extensively studied in \WaaThe, in the context of \intu. Still, only simple examples of non-metrizable natural spaces are given. They seem to us of little value other than for demonstration purposes. 

\parr

In this monograph we do not delve into non-metrizability either, but we show that in infinite-dimensional topology the construction of direct limits leads us to interesting non-metrizable natural spaces. These spaces have already been studied in a \TTE\ setting by others, see \mbox{e.g.}\ \KunSch\ for a nice article discussing also practical computational aspects.\footnote{As stated earlier \TTE\ is not an area in which the author is knowledgeable, it would seem that results from \classf\ are taken for granted, and then translated to a computability setting.}\ The examples given in \KunSch\ are: the space of polynomials with real coefficients, the space of real analytic functions on \zort, and the space of distributions with compact support. It should be possible to translate these spaces to the setting of natural spaces, since the \TTE\ treatment resembles our development, but there might be some work involved in constructivizing the classical results which are used.

\parr

Our simple example will consist of taking the direct limit of all the Euclidean spaces $(\R^n)_{n\in\N}$, which is equivalent to the space of `eventually vanishing real sequences', for which limit we show that the resulting classical limit-topology is in fact the apartness topology, and that this topology is non-metrizable. The non-metrizability is well-known, but the equivalence to the apartness topology is a nice illustration, we hope, that the concepts in this paper are also of interest for classical mathematicians.  

\parr

The space of eventually vanishing real sequences is one way of representing the space of polynomials with real coefficients. Therefore our example also illustrates how to translate the spaces in \KunSch\ to our setting. 

\sbsc{Direct limits of natural spaces and Silva spaces}\label{directlimits}
We define the direct limit of an increasing sequence of natural spaces along classical lines, where the apartness topology provides a very easy way to do so.  

\defi
Let \Vnatt\ be an natural space derived from \Vprenatt. Suppose there is a sequence of subsets $(V_n)_{n\in\N}$ of $V\iz\bigcupp_{n\in\N}V_n$, where for each \ninnt\ the pre-natural space $(V_n, \aprt, \leqc)$ gives rise to a natural subspace $(\Vcal_n, \Topaprt)$ of \Vnatt, and where in addition $V_n\subseteqq V_{n+1}$.

Then \Vnatt\ is called the \deff{direct limit}\ of $((\Vcal_n, \Topaprt))_{n\in\N}$. 
\edefi

This very simple definition suffices for our purposes here, but it also is, we believe, a faithful mirror of the traditional definition for \mbox{e.g.}\ Silva spaces arising as direct limits of Banach spaces, since on a separable complete metric space the metric topology coincides with the apartness topology (see theorem \ref{sepmetnat}).

\sbsc{The direct limit of the Euclidean spaces is non-metrizable}\label{dirlimnonmet}
Firstly, let $\Rom=\{x\inn\RN\midd \driss\ninn \all \minn, m\geqq n [x_{[m]}\equivv 0]\}$ be the subset of \RNt\ of eventually vanishing real sequences.\footnote{Remember that we write $x_{[m]}$ for the real number which is the $m$-th coordinate of $x$ as an infinite sequence of real numbers, and $x_m$ for the $m$-th basic dot of $x$ as a sequence of basic dots. Also, we drop the subscript `nat' in \Rnatt\ for notational simplicity.}\ As a subset of \RNt, it is at least weakly metrizable, but a different topology arises when we consider \Romt\ as the direct limit of the Euclidean spaces $(\R^n)_{n\in\N}$. This direct limit is a very simple example of an infinite-dimensional Silva space (meaning the direct limit of an increasing sequence of Banach spaces, with an extra compactness condition on the closed inclusion images of the respective unit spheres). 
We show that this direct limit arises naturally in our context, when building \Romt\ as the point set of a natural space. Notice that \Romt\ is not a subspace of \RNt\ under our definition of `natural subspace'. 

\parr

So to build \Romt\ as a natural space, we cannot simply use the infinite-product definitions used in defining \RNt\ (see definition \ref{defprod}), although we can use the same collection of basic dots. We  should distinguish between \Romt\ as a point-subset of \RNt, and \Romt\ as the point set of the natural space \Romnatt, but to avoid tedious notation we will not do so.

\defi
Let ${\Rrat}\!\!\!\!^*\iz(\Rrat)_{_{\Pi}}\iz\bigcupp_{n\in\N}{\Rrat}\!\!\!\!^n$ be the set of finite sequences of closed rational intervals. Denote the empty sequence as \maxdotRNt. We define \RNnatt\ as the infinite product $\Pi_{n\in\N} \Rnat$, see definition \ref{defprod}, but we write \aprtRNt\ for the product apartness
and $\leqcRN$ for the product refinement (and often \RNt\ for \RNnatt)

\parr

To define \Romt\ as the point set of a natural space, we need to sharpen the product apartness $\aprtRN$ and the refinement $\leqcRN$. To understand this, note that for \ninnt\ we now consider the basic dots in ${\Rrat}\!\!\!\!^{n+1}$ as belonging to points $x\inn\RN$ with $x_{[m]}\equivv 0$ for all $m\bygr n\pluz 1$. 

So let
$a\iz a_0, \ldots, a_n\inn{\Rrat}\!\!\!\!^{n+1}$ and $b\iz b_0, \ldots b_m\inn{\Rrat}\!\!\!\!^{m+1}$ with $m\geqq n$, then we put $a\aprtRlim b$ iff: there is $i\leqq n$ with $a_i\aprtR b_i$ and/or there is $n\smlr i\leqq m$ with $0\notinn b_i$ (which is decidable since $b_i$ is a rational interval for all $i\leqq m$). Similarly we put $b\leqcRlim a$ iff: $b_i\leqcR a_i$ for all $i\leqq n$ and $0\inn b_i$ for all $n\smlr i\leqq m$.  
\parr
The \deff{natural space of eventually vanishing real sequences}\ is \Romnatt, derived from the pre-natural space \Rompret. We also for each \ninnplt\ define: $\Rnstar=\{x\inn\Rom\midd \all \minn\,[x_m\inn\bigcupp_{1\leq i\leq n}{\Rrat}\!\!\!\!^i\,\,]\}$. We see: if $x\inn\Rnstar$, then $x_{[m]}\equivv 0$ for all $m\bygr n$. Notice that there is only a trivial difference between \Rnstart\ and $\R^n$ defined as finite product space in \ref{defprod}.
\edefi

\prp
\Romnatt\ is the non-metrizable direct limit of $(\Rnstar, \TopRlim)_{n\in\N}$,  where $(\Rnstar, \TopRlim)$ is \leqc-isomorphic to the Euclidean space $(\R^n, \TopRaprt)$ for \ninnt. There is a continuous injective surjection (which does not have a continuous inverse) from \Romnatt\ to $(\Rom, \TopRNaprt)$ as subspace of $(\RN, \TopRNaprt)$.
\eprp

\prf
We prove this in the appendix \ref{prfdirlimnonmet}, but the reader should have no problem with it as an exercise.
\eprf

\rem
\be
\item[(i)] More generally one can prove in this fashion that infinite-dimensional Silva spaces are non-metrizable natural spaces. We do not go into this, but it is an easy generalization of the above constructions for $(\R^n)_{n\in\N}$ and \RNt. 
\item[(ii)] Perhaps interesting: \Rnatt\ `equals' the direct limit of $([-n, n]_{\scriptscriptstyle\rm rat}^*, \TopRaprt)_{n\in\Nplus}$, where $[-n, n]_{\scriptscriptstyle\rm rat}^*=\{[a, b]\inn\Q\timez\Q\midd -n\leqq a\smlr b\leqq n\}\cupp\{\maxdotR\}$.
\ee
\erem

For non-metrizable natural spaces resembling examples in \Urya\ we refer the reader to \WaaThe. It should be possible to transfer the separation properties T$_1$ through T$_4$ given there to our setting of natural spaces. 

\sbsc{Metrizability results for paracompact-like spaces}\label{metrizparacomp}\hspace*{-1.35pt}
Metrizability of topological spaces has played a very important part in the development of classical general topology. In constructive mathematics, the focus until recent years has been primarily on metric spaces. In \Freu, an intuitionistic metrization result for compact spaces was obtained. This result was extended in \TroThe\ to locally compact spaces. In \WaaThe, using a different technique, these results were reobtained for apartness spreads, and metrizability was extended to `star-finitary' spaces. The latter resembles classical metrization results for strongly paracompact spaces.
\parr  
For now, we turn to the metrizability results in \WaaThe. For completeness, let us repeat the classical definitions of `paracompact' and `star-finite cover'. In \class, a topological space \xTopt\ is called \emph{paracompact} iff every open cover \Ucalt\ of the space has a locally finite refinement, meaning a refining open cover \Wcalt\ such that each $x\inn X$ is contained in only finitely many open sets in \Wcalt. An open cover \Ucalt\ of \xTopt\ is called \emph{star-finite} iff for each $U\inn\Ucal$ there are only finitely many elements $W\inn\Ucal$ such that $U\capp W\notiz\emptyy$. So we see that a star-finite open cover is always locally finite, but the converse is not true. Therefore a space is called `strongly paracompact' iff every open cover of the space has a star-finite refinement. 

\parr
There are many classical results on metrizability in regard to paracompactness and strong paracompactness, also specifically for separable Lindel\"of spaces (which resembles our setting). We profess not to have a clear assessment of the relation between our result and these classical results, although we expect the below constructive metrization theorem to follow in \class\ from an existing classical metrization theorem. The classical theorem most resembling our constructive treatment seems the one stating that a T$_0$ space is metrizable iff it has a development $(\Wcal_n)_{n\in\N}$ where $\Wcal_{n+1}$ is a star-finite refinement of $\Wcal_n$ for each \ninnt.

\sbsc{Star-finitary spaces}\label{starfin}
From the definition of spraids, it might seem at first glance that for a given spraid \Vnatt\ derived from the pre-natural space \Vprenat\ and $a\inn V$, there are only finitely many $b\inn V$ with $\grdd(b)\iz\grdd(a)$ and $\Vsubb\capp \Vsuba\notiz\emptyy$ (one could call this a $\grdd(a)$ star-finite basic intersection property, but we will not use this). However, this does not follow from the definition if one looks carefully, and counterexamples are easy to construct.    
\parr
Even if a $\grdd(n)$ star-finite basic intersection property holds for each \ninnt, then because of the complicating apartness relation, it tells us little about the natural topology of the space \Vnatt. Therefore we will define our star-finite property in terms of the apartness relation \aprt\ (or rather its complementing touch-relation \touch) on basic dots. For metrizability results we then look at `star-finitary' natural spaces,  which are \sucz-isomorphic to a star-finite \sucz-spraid \Vnatt\ (where the inductivity of \Vnatt\ ensures nice behaviour of the apartness relation).

\defi
Let \Vnatt\ be a spraid derived from \Vprenat, where $\touch\subsett V\timez V$ is the complement of \aprt.  We say that \toucht\ as well as \Vnatt\ and \Vprenatt\ are \deff{star-finite}\ iff for each $a\inn V$ the subset $\{b\inn V\midd\grdd(b)\iz\grdd(a) \weddge b\touch a\}$ is finite. A natural space \Wnatt\ is called \deff{star-finitary}\ iff \Wnatt\ is \sucz-isomorphic to a star-finite \sucz-spraid. 
\edefi

The class of star-finitary spaces is large since countable products of star-finitary spaces are again star-finitary. Baire space is a star-finite \sucz-spraid, and \RNt\ is star-finitary (a star-finite \sucz-spraid representing \RNt\ can be obtained by forming the countable product of copies of the star-finite \sucz-spraid \sigRt).

\sbsc{Star-finitary spaces are metrizable}\label{starfinmet}
In this paragraph we translate an intuitionistic result regarding metrizability of apartness spreads. This translation (especially its proof) takes some time, since the original result makes multiple use of the intuitionistic axioms \acoz\ and \FT. Our translation below (and especially its proof) is therefore also a good example how to circumvent the use of these axioms if one is so inclined, and simultaneously shows that -as should be expected- intuitionistic theorems carry `effective' content. \footnote{In \PalInt\ the usage of intuitionistic axioms is called non-effective, which indicates that their constructive content might not be immediately apparent, a reason to elucidate this constructive content here.}      

\thm
Every star-finitary natural space is metrizable.
\ethm 

\prf
It suffices to prove the theorem for a star-finite \sucz-spraid\ (remember that `\sucz-spraid' is shorthand for `inductive spraid'). For this we translate the proof of the corresponding intuitionistic theorem in \WaaThe\ for star-finite apartness spreads. This translation is in principle unproblematic, since the use of \FT\ in that proof can now be replaced by the use of the inductive Heine-Borel property for subfanns (\HBind), and the use of \acoz\ becomes superfluous by our definition of inductive spraid.
 
Still, since the original result is not trivial, our translation also involves quite some work. We give this translation in the appendix, see \ref{prfstarfinmet}. We also sketch the main idea of the proof below, when discussing the corollary.
\eprf

\crl
Every \sucz-fann is metrizable (`every compact space is metrizable'), and every space with a one-point \sucz-fanlike extension is metrizable (`every locally compact space is metrizable').
\ecrl

\prf For didactical reasons we will first prove the corollary directly in \ref{prfcrlstarfinmet}, before proving all of the theorem above. Let us sketch the strategy for a \sucz-fan \Wnatt\ derived from \Wprenatt\ already here. To each $a\aprt b\inn ^nW$ (for certain \ninn) we can construct a canonical morphism $f_{a,b}$ from \Wnatt\ to  \zortert\ such that $f_{a,b}\equivvR 0$ on $\Wsuba$ and $f_{a,b}\equivvR 1$ on $\Wsubb$. Then let $h:\N\rightarrow W\timez W$ be an enumeration of all apart pairs $a\aprt b\inn ^nW$ (for all \ninn), and put $d(x, y)\isdef\sum_{n\in\N}2^{-n}\cdott\midd f_{h(n)}(x)-f_{h(n)}(y)\midd$ to obtain the desired metric $d$.

\parr

For $a\aprt b\inn ^nW$ (for certain \ninn), in order to construct a canonical morphism $f_{a,b}$ from \Wnatt\ to  \zortert\ such that $f_{a,b}\equivvR 0$ on $\Wsuba$ and $f_{a,b}\equivvR 1$ on $\Wsubb$, we split $W$ first in $W_0$, $W_1$ and $W_2$ where $a\touch W_0$, $a\aprt W_1\aprt b$ and $W_2\touch b$ and moreover $W_0\aprt W_2$.

\parr

For this we use the splitting lemma \ref{prfcrlstarfinmet}\,(\bol{B}), which tells us that finite apart $A\aprt B\subsett W$ always lead (for big enough \Ninn) to a partition of ${^N}\!W$ in three sets $C, D, E$ where $A\touch C$, $A\aprt D\aprt B$, $E\touch B$ and moreover $C\aprt E$.

\parr

Then we iterate this procedure for each $W_i$ (where $i\inn\{0, 1,2\}$), splitting \mbox{e.g.}\ $W_0$ in $W_{00}, W_{01}, W_{02}$ where $a\touch W_{00}$, $a\aprt W_{01}\aprt W_1$ and $W_{02}\touch W_1$ and moreover $W_{00}\aprt W_{02}$. 

In this way, we actually construct a function from $W$ to \zotstart, such that this function represents a morphism from \Wnatt\ to \sigthrRt\ (see example \ref{morfbinreals}), where the latter is isomorphic to \zortert.

\parr

To generalize this strategy for \sucz-fans to a star-finite \sucz-spread involves a lot of extra work, but the basic idea remains the same. For the complete proofs, we refer the reader to \ref{prfcrlstarfinmet}\ and \ref{prfstarfinmet}.
\eprf 

\rem
\be
\item[(i)] The theorem seems to cover much of what is possible, metrization-wise. One can prove in \intu\ (and likely \class, \russ) that complete metric spaces are weakly star-finitary (which in \class\ is equivalent to being star-finitary). We discuss this in the next paragraph.

\item[(ii)] The example in \ref{nonarc}\ of a possible spherical completion of \Cpt\ might pose an interesting application for the metrization theorem above, but we confess to not having studied this issue in any detail.
\ee
\erem

\sbsc{Complete metric spaces are weakly star-finitary in \intu}\label{convstarfinmet}
It is shown in \WaaThe\ for \intu\ that every (separable) complete metric space is weakly star-finitary. If we look at the proof, we see that it essentially uses only countable choice and the axiom \BDD\ (see the appendix \ref{axiombdd}), which is true in \class, \intu\ and \russ. This means that the theorem likely holds also in \class\ and \russ. In \class\ `weakly star-finitary' is equivalent to `star-finitary'. 

\defi
Let \Vnatt\ be a spraid derived from \Vprenat, where $\touch\subsett V\timez V$ is the complement of \aprt.  
We call \toucht\ and \Vnatt\ and \Vprenatt\ \deff{weakly star-finite}\ iff for each $a\inn V$ the subset $\{b\inn V\midd\grdd(b)\iz\grdd(a) \weddge b\touch a\}$ is subfinite, meaning there is \Ninnt\ such that $\{b\inn V\midd\grdd(b)\iz\grdd(a) \weddge b\touch a\}$ contains at most $N$ elements. A natural space \Wnatt\ is called \deff{weakly star-finitary}\ iff \Wnatt\ is \sucz-isomorphic to a weakly star-finite \sucz-spraid. 
\edefi

\thm (in \intu\ and likely \class\ and \russ)
Every complete metric space is homeomorphic to a weakly star-finite \sucz-spread.
\ethm

\prf (detailed sketch)
The proof for \intu\ is given in \WaaThe\ (thm.\,3.1.8) using \acoz\ in several lemmas. However, upon reflection one sees that each use of \acoz\ can be substituted by a use of \BDD\ (see the appendix \ref{axiombdd}), which is true in \class, \intu\ and \russ. 

\parr

We sketch the translation to our setting in some detail. Let \xdt\ be a complete metric space. We saw in (the proof of) theorem \ref{compmetind}\ that there is a \sucz-spread \Vdiamcalt\ which is homeomorphic to \xdt. This spread \Vdiamcalt\ is derived from the neighborhood development of \xdt\ formed by the collection $\{B(a_n, 2^{-s})\midd n, \sinn\}$ where \anninnt\ is a dense sequence in \xdt. The strategy of the proof is to sharpen this neighborhood development into a strongly star-finite neighborhood development as follows.

\parr

By \WaaThe, thm.\,3.1.1 (valid in \bish) a per-enumerable open cover of \xdt\ has a strongly star-finite refinement.

We construct a strongly star-finite refinement $\Ucal_0$ of $\{B(a_n,2^{-0})\midd\ninn\}$. Then we construct a strongly star-finite refinement $\Ucal_1$ of $\{B(a_n,2^{-1})\midd\ninn\}$ , and then of $\,\{B(a_n,2^{-2})\midd\ninn\}\,$, etc. We obtain a sequence $(\Ucal_n)_{n\in\N}$ of strongly star-finite covers such that for each $U$ in $\Ucal_n$: $\,diam(U)\smlr\twominn\,$. Now we use the elements of $\Ucal_n$ as the $\grdd(n)$ basic dots of our weakly star-finite spread \Wnat, where we need to be precise in defining \leqct, since these $\grdd(n)$ dots belong to particular $\grdd(n\minuz 1)$ dots. Here some real work needs to be done, also to define a weakly star-finite touch-relation \toucht\ on $W$ such that for the corresponding apartness \aprt\ we have: \Wnatt\ coincides with \xdt.

\parr

The entire proof, in order to be precise and correct, already for \intu\ involves some more work than one might expect. We refer the reader to \WaaThe, and leave this proof-sketch for what it is: a sketch.
\eprf

\sbsc{Star-finite developments and `efficient' metric spraids}\label{starfinsprd}
We turn once more to the \appl\ perspective. In paragraph \ref{refvstrail2}\ we showed that our representation of the real numbers has computational advantages regarding continuous \R-to-\R-functions. The general situation for complete metric spaces is more complicated. Yet very often we can find similar efficient representations for complete metric spaces. We do not delve into this in all detail, but we illustrate the basic idea.

\defi
Given a metric space \xdt, let $(\Ucal_n)_{n\in\N}$ be a family of open (or neighborhood) covers where (i) $\Ucal_n$ consists of inhabited sets of diameter less than \twominnt\ (ii) $\Ucal_n$ is star-finite (iii) $\Ucal_{n+1}$ is a (star-finite) refinement of $\Ucal_n$ for each \ninnt. Then we call $(\Ucal_n)_{n\in\N}$ a \deff{regular star-finite development}\ of \xdt. 
\edefi

Suppose we have a complete metric space \xdt\ and a regular star-finite development $(\Ucal_n)_{n\in\N}$ of \xdt. By the previous paragraph, this should not be a rare occurrence at all. Then we can build a spraid \Vnatt\ (with corresponding \Vprenatt) representing \xdt, where $V$ represents $\bigcupp_{n\in\N}\Ucal_n$ (each basic dot corresponds to an element of some $\Ucal_n$, except for \maxdott). This is much more direct and efficient than the spread-representations of \xdt\ given through theorems \ref{sepmetnat}\ and \ref{compmetind}. 

\parr

In this situation we would like a general way to replace a \pthh-morphism from \Vnatt\ to another natural space \Wnattot\ by an equivalent \leqc-morphism. We are confident that most important cases will resemble proposition \ref{refvstrail2}. One direct example of this are of course the Euclidean spaces $\R^n$ and \RNt. But we follow remark \ref{refvstrail2}\ a bit further to say something in general as well:

\parr
If \Wnattot\ is \leqc-isomorphic to its `finite-intersection-of-touching-basic-dots space', then we can follow the proof of proposition \ref{refvstrail2}\ closely to see that indeed a \pthh-morphism from \Vnatt\ to \Wnattot\ is represented by an equivalent \leqc-morphism.

\parr

We leave it to the interested reader to further ponder on the efficiency of different representations of complete metric spaces.

\sectionnb{From natural to general topology}\label{natgentop}

\sbsc{Invitation to general natural topology}\label{genmettop}
We think that from here one can develop a substantial part of general topology constructively, using pointwise concepts. We will however not tackle this significant undertaking in this monograph (which is lengthy enough as it is), preferring to invite others to do so. A very incomplete list of interesting subjects to cover/prove (some of which are also partly covered in \WaaThe): 

\be
\item Metric topology and (in)finite-dimensional topology (see \Mill\ for an excellent exposition of the latter)
\be
\item Star-finite covers and partitions of unity
\item Topological and transitive notions of locatedness
\item Universal metric spaces and isometric (normed linear) extensions
\item Absolute retracts, Dugundji Theorem, Michael Selection Theorem
\item Compactification
\item Dimension theory and fractals
\item Topological manifolds
\item $\ldots$
\ee
\item Algebraic topology, homotopy theory, $\ldots$
\item Non-metric topology, topological lattice theory, separation axioms, $\ldots$
\item Topological groups
\item Topological model theory
\item $\ldots$
\ee
 
The list is probably easily expandable, since so little of classical general topology has been charted constructively. Apart from mathematics, this undertaking will have importance also for physics, we believe.
\parr
The relation to physics is one aspect that we wish to examine just a little bit more in this monograph, in our next and final section.

\sectionnb{Natural topology and physics}\label{intuinrussdet}

\sbsc{A philosophical apology}
This section is essentially philosophical in nature, although it employs mathematics as a vehicle. Its purpose is more to raise questions, than to give answers. The reader should ideally look at these questions with an open mind, to consider whether they merit more attention than usually given. Although some issues in this section are speculative to a certain degree, the author feels it would be a pity to ignore them.

\sbsc{The strained relation between finiteness and infinity}\label{finvsinf}\hspace*{-.91pt}
One could argue that all constructive troubles regarding the foundations of mathematics stem from one strained relation: the one between finiteness and infinity. \hspace*{-.15pt}Brouwer's critique of \class\ can be summarized as: `the notion in \class\ of infinity as a completed entity is fundamentally unsound'. According to Brouwer and other (pre-)intuitionists, we can only construct \Nt\ as a `potentially infinite' set, that is: growing in time beyond any fixed boundary, as time progresses (potentially) infinitely.
\parr
From the past century it becomes clear that the mathematical problem for `infinitism' (denoted here as \infi) is to build a foundationally sound and satisfying theory. This has proved to be difficult mathematically, philosophically and not in the least socially because of the added difficulty to get mathematicians to agree on what is `foundationally sound' and what is `satisfying'.

\parr
Potential infinity is also contested by some, a position called (ultra)finitism (denoted here as \fini).\footnote{Well-known advocates are Alexander Esenin-Vol'pin and Edward Nelson.}\ The mathematical problem for \fini\ is the same: to build a foundationally sound and satisfying theory. This has proved to be difficult also -apart from the social factor mentioned above- because of the mathematical and philosophical troubles that arise when trying to fix a good upper bound on \Nt. 

\parr
It should be clear that the strain between \infi\ and \fini\ also has direct roots in physics. The `simple' questions: `What is time?' and `Is the universe (fundamentally) finite?' have not been answered to any satisfactory degree, as far as we are aware.

\parr
It is no surprise therefore, that a lot of `action' in the foundational debate takes place in the arena where we go from finiteness to infinity. For example, each axiom used or discussed in this paper (see the axioms' section \ref{axioms} in the appendix) is tailored to describe certain transitions from finiteness to infinity or vice versa.\footnote{`From here to eternity and back again'.}

\parr
And this transition is also precisely the thrust of `natural spaces'.

\conv
In the following (as in the previous...) we silently assume the world to be potentially infinite. The main reason to do so is that if the world is fundamentally finite, then philosophically speaking we can limit ourselves to discrete combinatorics, we think. So it seems more interesting for the debate to assume \infi, even though \fini\ makes as much sense.
\econv

\sbsc{Which mathematics are suited for physics?}\label{whichmathphys}
Historically, it bears little surprise that classical mathematics (\class) is the mathematics of preference for physicists, which they generally use without question to develop physics.    However, the author believes that such an unquestioning use of \class\ is difficult to support in the light of the developments in the foundations of mathematics.
\parr
As we have described in this monograph, there are several good reasons to doubt the inherent suitability of \class\ to describe our physical world. We do not mean to say that it is impossible to do physics well in \class, since any other approach can always be translated to a model in \class. But we are concerned about suitability, and possible blind spots arising from the exclusive use of \class. \footnote{As an example of what we believe to be unawareness of constructive foundational issues, one can take Hawking's recent book `God created the Integers' (\HawGod) in which he discusses 31 great mathematicians of all times which according to him have been of fundamental importance to mathematics and physics. Brouwer is not among them, nor is he even referenced, even though both his topological work and the foundational crisis that he discovered have had a major impact.}   
\parr
As stated, we believe that from our current state of knowledge in the foundations of mathematics one cannot simply choose one of the main alternatives \class, \intu, or \russ\ as the preferred mathematics for physics. One of the motivations for this paper is to show that natural topology yields a simple model in \class\ of basic principles in \intu, and to couple this model with an explanation why we think that natural topology is better suited for physics than the standard practice.
\parr
Yet as we have seen, this leaves unresolved another fundamental aspect of the relation between physics and mathematics, namely the question what reason we have to assume that reality is modeled better by a compact unit interval than by a non-compact unit interval. Or by paradigm leap: what reason do we have to assume that Nature is capable of producing non-recursive real numbers? A strong case can be made, we believe, for a form of Laplacian determinism in which Nature can only produce recursive sequences.\footnote{A similar belief seems to underly Stephen Wolfram's `A New Kind of Science' (\WolNew), but without a sharp foundational motivation that we know of.}\ This belief is often written as \CTphys, where \CT\ stands for Church's Thesis: `every infinite sequence of natural numbers is given by a recursive algorithm'.  
\parr
This question is debated in some circles, yet in the author's eyes the debate mostly lacks mathematical precision. In \WaaArt, a precise mathematical setting is given in which \CTphys\ can even be experimentally tested, although the latter involves some serious computational effort and has not yet been done. In addition there is no guarantee of result for this experiment. We will comment on this further on, because even as a thought-experiment alone, doubt is cast on the intrinsic validity of the standard probability approach in science. Laplace himself already raised similar doubts but resolved them using a model of `limited information' (\LapDif, also see \LapPhi):

\begin{quote}
\emph{"Before going further, it is important to pin down the sense of the words chance and probability. We look upon a thing as the effect of chance when we see nothing regular in it, nothing that manifests design, and when furthermore we are ignorant of the causes that brought it about. Thus, chance has no reality in itself. It is nothing but a term for expressing our ignorance of the way in which the various aspects of a phenomenon are interconnected and related to the rest of nature."}
\end{quote}

A different model of `limited information' will also serve us well.

\parr

Apart from foundational probability issues, in \WaaArt, the open question is raised whether \russ\ could be a better model for physics than \intu\ or \class.
What remains lacking in the literature (as far as we are aware) is a sharp analysis of the possible consequences of \CTphys\ for physics.

\rem In \WaaArt\ a possible partial reconciliation between \russ\ and \intu\ was left uninvestigated. We will describe such a partial reconciliation in this section, using a model of limited information which resembles `type two effectivity' (\TTE). This partial reconciliation is relevant for physics, we believe, in the following sense. Suppose that \CTphys\ is one day seen to hold.\footnote{For instance by a positive outcome to the experiment described in \WaaArt.} Then the reconciliation that we propose below offers an explanation why the physical world would still (appear to) conform to \BT.  
\erem 

\sbsc{Kleene's function realizability and type two effectivity}\label{KleeneTTE}
One of the developments with some similarity to natural spaces is Weihrauch's type two effectivity.\footnote{\TTE, to which of course many others have also contributed, see the literature.}\ \TTE\ uses classical logic, and can therefore in our eyes not be considered a constructive theory, yet there is a certain conceptual overlap we believe.

\parr

Kleene's function realizability used in a setting similar to \TTE\ may give us a constructive way to formally interpret \intu\ as part of \russ. We leave this as an open question.\footnote{Which most likely has been answered already outside the author's awareness, see also paragraph \ref{closapo}.}

\parr
The basic idea of \TTE\ is not too difficult. One considers Turing machines which can handle infinite input tapes, and which produce infinite output tapes step-by-step (like Brouwer's choice sequences). In other words, a given machine $T_e$ is fed an infinite input sequence $\alpha$ (not necessarily recursive!) and it is possible that after some finite time $T_e$ computes some first result from an initial segment of $\alpha$, which is then the first output value of a (possibly) infinite output sequence 
$\beta\iz T_e(\alpha)$. We refer the reader to the literature, notably \WeiCom.
\parr
Since we wish to limit the scope of this article to our own limited level of knowledge, we will not delve into \TTE\ and Kleene realizability. Instead we present a simple two-player game which we believe gives a similar (but informal) model of \intu\ in \russ.

\sbsc{An informal two-player game representing \intu\ in \russ}\label{intuinruss}
We describe a game called `Limited Information for Earthlings' (\LIfE). We start out by introducing our two players, and describing their rules of conduct. 

\parr

Player I is called `Giver of Digits' (\GoD), and player II is called `Humble Mathematician accepting Numbers' (\HuMAN).

\parr

The basic idea of the game \LIfE\ is that \GoD\ gives infinite sequences of numbers to \HuMAN\, where all sequences are in fact given by a recursive algorithm, but \HuMAN\ does not know this. We list the rules of \LIfE:

\be

\item[\lone] \GoD\ has at her disposal all infinite recursive sequences (all $\alpha\inn\NN$ for which there is a Turing algorithm which computes $\alpha$).\,\footnote{This is equivalent to: $\dris e\inn\N\all\ninn\driss k\inn\N\ [\, T(e, n, k)\weddge\outc(k)\iz\alpha(n)]$, if we use Kleene's $T$-predicate.}\ \GoD\ hands out such  sequences $\alpha$ to \HuMAN, but is allowed to do so step-by-step without disclosing any information about the algorithm computing $\alpha$. 

\item[\lto] \GoD\ may lie in two different ways (and does). First by telling \HuMAN\ that not all sequences are recursive. Secondly, \GoD\ is allowed to cheat with information on the sequences and algorithms involved, if \HuMAN\ cannot know the difference. For example, if at stage $n$ \GoD\ has given out only the finitely many values $\alfstr(n)$ for a certain computable $\alpha$, then \GoD\ may switch to any computable $\beta$ with $\betastr(n)\iz\alfstr(n)$. However, during the game and when playing a sequence, \GoD\ may only switch a finite number of times. Compliance to this can be determined post-game, \HuMAN\ remains ignorant.

On the other hand, \GoD\ also can disclose part or all of the algorithm computing $\alpha$. This must of course be done consistently, so \GoD\ can only switch to other sequences conforming to earlier disclosures.

\item[\lthr]
The player \HuMAN\ has to build mathematics from the sequences given by \GoD, but \HuMAN\ can also at any time use his own recursively defined sequences. \HuMAN\ believes that \GoD\ also has non-recursive sequences at her disposal (although this is in fact false). \GoD\ is obliged to give to \HuMAN\ as much initial values of a sequence $\alpha$ as \HuMAN\ wants. (But \GoD\ is allowed to cheat by \lto). 

\item[\lfo]
\HuMAN\ can only use constructive logic, implying for instance that existential statements and disjunction (A or B) must have constructive content. 

\ee

\parr

So what are the truths of \LIfE\ for \HuMAN? We think that for \HuMAN, from our perspective, \LIfE\ is a model of Brouwer's intuitionism. The sequences that form \HuMAN's mathematical universe correspond to the universal spread \sigomt, or natural Baire space through Brouwer's eyes. Brouwer's motivation for the Continuity Principle \CP\ transfers directly to our game \LIfE\ (we can even prove it, see below). The same holds for \BT, since we see no other way for \HuMAN\ to know that a certain $B\subset\Nstar$ is a bar on \Bnatt\ other than by having constructed a genetic bar on \Bnatt\ from which $B$ descends. (We can prove that if \GoD\ is omniscient, then $\negg\dris B\subsett\Nstar [B$ is a non-inductive bar on $\Bnat]$ holds in \LIfE.)
 
\parr

Therefore, we think that the truths of \LIfE\ for \HuMAN\ (seen from our perspective) are precisely the theorems of \intu. 
But from \GoD's perspective, everything in \LIfE\ is algorithmically determined, in other words \LIfE\ is part of \russ. 

\rem 
The situation in which from an omniscient perspective sequences are recursively determined, whereas from the limited-information perspective the mathematical truths are intuitionistic, might provide -should the need arise\footnote{Which probably occurs, we believe, if one day \CTphys\ should be validated.}- a way to explain that nature `behaves inductively', so to speak. 

\parr
Or to put it differently: even if one day \CTphys\ is believed to hold, then we will be unlikely to encounter a Kleene Tree in nature, because nature itself is most likely built inductively. The experiment in \WaaArt\ however indicates a (perhaps practically infeasible) way to expose the statistical anomalies which would result from \CTphys. This corresponds to a strange truth in \LIfE: if \HuMAN\ is given enough time, \HuMAN\ can discover that \GoD\ is lying about having non-recursive sequences, by using the Kleene Tree and our standard model of probability for the testing of physics' hypotheses.

\parr
Yet the real rub if one believes \CTphys, is that our standard model of probability used to prove physics' theories would seem to need serious revision. Ironically, this model is largely due to Laplace. He could however point out, in accordance to the quote above (\ref{whichmathphys}), that knowledge of the truth of \CTphys\ alters the state of our ignorance, and therefore by necessity also our probability models. 
\erem

We now formulate our last meta-theorem concerning the game \LIfE. It should be seen as an alternative illustration of our attempts to further the `reunion of the antipodes' \footnote{Derived from Schuster's (et al.) terminology `reuniting the antipodes', see \SchuReu.}\ in constructive mathematics, and not as an important result in itself.

\mthm
\be
\item[(i)]In \LIfE, we can prove \CP. 
\item[(ii)]
Suppose \GoD\ is omniscient. Then we can prove $\negg\dris B\subsett\Nstar [B$ is a non-inductive bar on $\Bnat]$ for \LIfE.
\item[(iii)] Given enough time, \HuMAN\ can discover that by an overwhelming odds ratio, \GoD\ plays only recursive sequences. 
\ee
\emthm

\prf
We prove this metatheorem in the appendix \ref{prfcplife}.
\eprf

\crl
In \class, we can prove \CP\ and \BT\ for the game \LIfE. (This would seem to express the situation in \TTE).
\ecrl

\xam
We show that, assuming \GoD's omniscience, the Kleene Tree ($K_{\rm bar}$, see \ref{canbairus}) is not a bar on \sigtot\ in \LIfE\ even though all \sigto-sequences in \LIfE\ are in fact recursive. For suppose \HuMAN\ proves that a subset $B\subsett \zostar$ is a bar such that $K_{\rm bar}\subseteqq B$. Then \GoD, with omniscience, can choose a recursive sequence $\alpha$ such that $\{\alfstr(i)\}_{\leqc}\capp K_{\rm bar}$ is infinite for any $i\leqq N$, where $N$ is the index such that $\alfstr(N)\inn B$. (\GoD\ plays a $\beta$ step-by-step and at each step determines whether $\betastr(n)\star 0 \capp K_{\rm bar}$ is infinite, if so then the next choice for $\beta$ is a $0$, if not then the next choice is a $1$. At a certain point in time $M$, \HuMAN\ must produce $N$ such that $\betastr(N)\inn B$. At this point, \GoD\ fixes $\alpha\iz\betastr(N)\starr\zero$, which is a recursive sequence. We see that \GoD\ has played by the rules.).

\parr

By the properties of $\alpha$, we see that $\alfstr(N)\inn B, \alfstr(N)\notinn K_{\rm bar}$, in other words $B\notiz K_{\rm bar}$.
\exam

\appendix\parbox{0pt}{\chapter{}}
\chepterap{Appendix}{Additionals, Examples, Proofs} {The appendix starts with closing remarks and personal acknowledgements. 
Then we give some motivation and historical background of this paper, and a small recap. We add a brief history of apartness topology.
\parr The appendix is predominantly taken up by examples and proofs. The examples are worked out in some (but not all) detail, leaving some aspects as exercises to the reader. The proofs are worked out in almost all details, yet occasionally we also leave some exercise to the reader.
\parr
To be foundationally precise (although more precision is certainly possible) we present and explain the axioms used and discussed in this monograph.
\parr
We end with the bibliography.
}

\sectionnb{Closing remarks and acknowledgements}

\sbsc{Some closing remarks}\label{closapo}
Our introductory exposition of natural topology ends here (bar the appendix). Although the author would have liked to be able to present more, it proved to be a truly energy-consuming undertaking to work out, in an elegant way,
the precise machinery needed for natural topology.

Not having been employed as a mathematician for the past 15 years, there is a limit to how much energy one can reasonably spend on such undertakings. However, this particular project is worthwhile for the author personally, since it has answered some questions which kept on resurfacing after writing his PhD-thesis in 1996.

\parr

We hope that others will also see some merit in this project. Experience teaches that new concepts take time and effort to get used to, and the author is no exception. The effort which went into writing this study has precluded him from really studying various promising other concepts, such as Abstract Stone Duality (see \TayBau) and formal topology (see \mbox{e.g.}\ \SamFor). Also, any real comparative study should be done from a deeper understanding of various related fields than the author possesses. He apologizes for his omissions-by-ignorance and inaccuracies on this count.

\parr

Reactions to this study are most welcome. Please do not hesitate to point out errors, omissions, other points of view, interesting sources for an update of the bibliography, solutions to questions, results of this paper which were already proved elsewhere, etc.

On this last item, we should state that it is likely that some of our results resemble or even equal results proved elsewhere. Not mentioning a source in this respect only shows that the author is unaware of such a source, which he will be happy to name when pointed out.

Finally, the appendix contains a `background and motivation' section with a short historical overview. This overview is likely to contain inaccuracies, and suggestions for improvement are most welcome.

\sbsc{Acknowledgements of the author}
First and foremost, I would like to thank Wim Couwenberg for his continuing interest and active participation is this seemingly unending project. But what a beautiful project it is, in my humble opinion, and it was sparked off by Wim's insight and insistence that a simple presentation of `vlekjesruimten' (`dotty spaces'...which for marketing reasons we have renamed natural spaces...) should be possible. Wim's basic idea was perhaps not so simple to develop as we had hoped, but now allows for a nice and faithful classical interpretation of many results in intuitionistic topology. 

So the credit for the concept of `natural space' goes to Wim, but more importantly his support and friendship have kept this project afloat. There were several times when I was convinced that practically nobody is interested in intuitionism, and that in natural spaces there was little worth mentioning when compared to formal topology and other disciplines.  

Wim convinced me otherwise, by repeating his questions about intuitionism and other constructive disciplines and by becoming captivated by the constructive approach and foundations in general. Wim underlined the need for a simple exposition of intuitionistic ideas for classical mathematicians, and convinced me that this monograph fulfills a worthwhile role in that respect. He also actively participated in thrashing out the first and most illuminating versions of `natural space' and `natural morphism'. His sparring role in many other discussions has simply been invaluable.

\parr

It turned out however to be still a daunting task to get down to the nitty-gritty and have everything in precise mathematical working order. Since this nitty-gritty lies close to my PhD thesis, and the comparative study of constructive foundations done in \WaaArt, it turned out that I had better write this monograph in solo fashion.\footnote{(Final) definitions, results, proofs, examples, and errors therefore are the author's.} We decided this after a joint talk last January, which revealed that the foundational complexities surrounding `natural spaces' could not be disregarded. I can only hope that with this text I have done justice to all our discussions (especially the fun of discovering yet another snag!) and working together. 

\parr

Second, I would like to dedicate this paper to Wim Veldman, who introduced me to foundations, intuitionism and constructive mathematics, and was my doctoral advisor during my PhD research. Wim Veldman's precise and elegant style reflects on his insistence that a structural framework for constructive mathematics should be both elegant and foundationally precise. I hope that this paper passes muster in that respect. Wim's active work in intuitionism is also an inspiration, I took just one of his results as an important illustration in this text.

\parr

Furthermore there are many people, too many to list individually, who have contributed in some way or other to the epigenesis of this paper. I would like to mention Bas Spitters who is always willing to discuss and explain formal topology beyond my limited knowledge and understanding. Some twelve years ago, Peter Schuster (et al.) already contributed his apt terminology `reuniting the antipodes', which has been a continuing inspiration. Giovanni Sambin wrote me some very kind words when I felt miserable for mistakenly having claimed to have spotted an error in a formal-topology paper. The tireless work of Douglas Bridges to carry on the program started by Bishop has also played an important part. It is my hope that this monograph will help to continue investigating `constructive mathematics' in the spirit of Bishop, by showing how to inductivize pointwise notions. 

\parr

Perhaps one day the difference between Brouwer's and Bishop's approach will be felt to be far less important than their 
correspondence. Studying constructive mathematics in the pointwise style of classical mathematics seems to me to remain the most attractive choice. I do however feel that a transparent axiomatization similar to \FIM\ is always called for, even though it shouldn't have to dominate the presentation. This is a matter of taste also, but recent work in constructive reverse mathematics by various authors (notably Ishihara, see \mbox{e.g.}\ \IshRev, and Veldman, see \mbox{e.g.}\ \VelFan)  has shown the benefits of further explicitizing relative axiomatic dependencies.\footnote{Veldman has also developed `Basic Intuitionistic Mathematics' (\texttt{BIM}) as a formalization comparable to \FIM\ in which reverse intuitionistic mathematics can be carried out.}

\parr

The support from my family and my friends, although mostly in other areas than mathematics, has played a major role in this project. Finally, the love of my wife Suzan and my daughters Nora and Femke (not to mention their patience with this project), has been quite indispensable.

\parr

I would like to thank all these people mentioned above very sincerely. 

\parr
\texttt{(the author, july 2011)}

\sbsc{Hommage to Brouwer, Kleene, Bishop and...}\hspace*{-2.65pt}
One may also see this paper as a hommage to Brouwer, Kleene and Bishop, but let us not forget all those other mathematicians who have worked hard to both expand and simplify mathematics. This is an ongoing collective endeavour, notwithstanding differences of style, character and opinion. The bibliography, incomplete as it may be, should serve to illustrate just how many people are dedicated to developing constructive mathematics. Behind each name and reference lies a body of related work which is left unmentioned, but with the help of the internet should be easily findable.

\sectionnb{Background, motivation, and recap}\label{background}

\sbsc{Background}\label{backback}
Much of theoretical mathematics is built on idealizations which fail in real life. A simple example is that of a floating point representation of a real number $x$ very close to 0, where the decision whether $x=0$ or $x\neq 0$ can be needed for further computations. Theoretically this decision is trivial, but in practice we cannot always determine whether $x=0$ or $x\neq 0$. That is one situation where `applied mathematics' comes in, with mathematicians and computer scientists working to translate theoretical mathematical ideas to `real life' situations.

 \parr

One often comes across excellent practical solutions, which are yet ad hoc in character. People working in the field of applied mathematics do not seem to consider this a problem, but we find it interesting to note that a more coherent and unified approach is possible. This might hopefully shed some new light on theoretical science as well, and help explain some important ideas of intuitionism and constructive mathematics to the `working' mathematician who is used to classical mathematics.

\parr

This paper is concerned with such an approach, from the standpoint of topology. A lot of work towards bridging the gap between theoretical and practical mathematics has been done in what is known as constructive mathematics. Although constructive mathematics in its essence is as old as mathematics itself, one can still consider Brouwer to be its founding father. Brouwer was the first to critically analyze the body of classical mathematics to come to the conclusion that the principle of the excluded middle (\PEM) could only be maintained for infinite sets at the cost of constructivity.  Brouwer sharply demonstrated this mathematically, with various clever examples, showing \mbox{e.g.}\ that for an arbitrary real-valued continuous function $f$ on the real interval \zort\ one cannot always construct an $x$ in \zort\ where \ft\ assumes a maximum. By the very nature of his critique, the foundations of mathematics in general, and especially the then popular and newly evolving mathematical discipline of set theory (started by Cantor) were shaken badly. It must be said that already Kronecker and Poincar\'{e} had serious reservations about Cantor's treatment of infinity, and also among his contemporaries Brouwer was not the only one with constructive views.
\parr

However Brouwer, the master topologist, did not content himself with criticism alone. From 1912 to roughly 1927 he developed a new constructive framework for mathematics, called intuitionism. One of his insights was that an everywhere-defined real function (to the reals) has to be continuous (we rediscover the reason for this in this paper). But he also encountered some difficult obstacles in building this constructive framework. He presented his solutions to these obstacles as theorems, albeit with rather unorthodox proofs, in none too easy language to the mathematical community of his time. This mathematical community was mostly unreceptive to Brouwer's critique of the classical foundations, and unwilling to change its comfortable views on classical mathematics as being the only viable framework for doing math. Hilbert (who in 1928 ousted Brouwer from the board of the Mathematische Annalen, showing how deep his resentment of Brouwer's views had become) is famously quoted to have said: `No one shall expel us from the paradise that Cantor created for us'. 

\parr

Deeply disappointed, Brouwer more or less retreated in his own activities and never regained his former prolific productivity. Still, Brouwer had followers, notably his student Heyting who simplified and formalized his mentor's intuitionism to make it more accessible. Later on, the computability expert Kleene became a supporter of Brouwer's ideas. Kleene managed to axiomatize intuitionism in a very clear way in \KleVes, to prove relative consistency (through realizability) which further opened the door for interested mathematicians. Kleene discovered that under these axioms Brouwer's notion of `choice sequence' could not coincide with the notion of `computable sequence', a result which will also concern us in this paper since the `true' reason for this is topological in nature. Namely in recursive mathematics the Cantor space \zoNt\ is not compact. Kleene thus also showed that Brouwer's Fan Theorem (\FT) (stating compactness of the Cantor space \zoN) was truly an axiom, in the sense that it could not be proved from the other axioms.

\parr 

Around that time, Bishop also became convinced of the intrinsic worth of the constructive framework. However, Bishop was not attracted to the foundational discussions, involving a great deal of logic, logicism and axiomatics. Bishop wanted to build a solid body of constructive mathematics, picking up where Brouwer had left around 1930, and also in such a way that the results would be acceptable to classical mathematicians as well. By this time, computers had already entered the scene, and mathematical awareness of computability and computational issues had increased the reception of constructivism.  Bishop-style mathematics (\bish) has increasingly become popular, not in the least because of its refocussing on `plain mathematics' instead of logic and foundations, and its down-to-earth approach. This approach started with a comprehensive treatise of constructive analysis, in the context of metric spaces (\BisFou, \BisBri). Bishop stated that there was little need for general topology, and that `mystic' axioms like the Fan Theorem (although classically true) were unnecessary if one chose the right definitions. However, it was later shown (\WaaThe, \WaaArt) that Bishop's definitions practically imply the Fan Theorem. Bishop's work is carried on by many, including notably Bridges. Brouwer's intuitionism seems to attract less mathematical attention, yet is carried on notably by Veldman. 

\parr

From the 1960's on, a parallel development of many authors (Scott, Martin-L\"{o}f, Fourman, Sambin, et al.) led to the field of domain theory, pointfree topology and formal topology. Actually already Freudenthal started with (intuitionistic) pointfree topology in \Freu, reacting to Brouwer and using ideas of Alexandroff and Hurewicz. In recent developments in formal topology, a number of the above issues have been dealt with in such a way that one can view this as a reunion of parts of the different approaches. This topological setting is no coincidence, since Brouwer was a brilliant topologist, and Brouwer's intuitionism was built with the backing of his topological expertise. Bishop seems to have underestimated the topological necessity of Brouwer's Fan Theorem in order to build a constructive model of analysis in which continuous functions on compact spaces are uniformly continuous.

\sbsc{Motivation and some results}\label{motivres}
The accomplishments of formal topology notwithstanding, we feel that formal topology in many of its current presentations (a growing number of papers and tutorials in a formal categorical style) lacks the intuitive appeal of both Brouwer and especially Bishop. Much of this is due to the fact that both Brouwer and Bishop concentrated on `points' and `spaces' in the usual mathematical way, and limited themselves to some form of separability (enumerable bases, enumerable dense subsets) in order to achieve constructivity. Also, there are still some foundational concerns surrounding compactness which we feel merit attention.\footnote{An alternative approach called Abstract Stone Duality (ASD) was recently developed by Taylor and Bauer, see also \TayBau. This development seems partly motivated by a search similar to ours for a simpler and more directly appealing approach to constructive topology. However, ASD draws on many areas in which the author is not proficient, and therefore comparison here is not possible. Kalantari and Welch (see \mbox{e.g.}\ \KalWel) have also been developing related concepts in a recursive-computable framework.} 

\parr

Interestingly enough, such a separability approach is ideally suited for a simple and elegant version of pointfree topology which we think deserves the name `natural topology' for three reasons. First of all, we think that natural topology is ideally suited for dealing with the study of nature (in other words physics...), since natural phenomena can only be observed and measured by scientists in a manner corresponding to these definitions. Secondly, in natural topology one builds a space of points in a natural way, and immediately sees a topology on this space arising from the construction, matching the space. Thirdly, we believe this allows for a natural pointwise style of mathematics, now that the most foundational aspects have been established.\footnote{This monograph needs to focus on foundations also, which makes for less easy reading.} 

\parr

Another advantage of this approach we hope to show in this paper: the easier translation of existing (intuitionistic) results to this setting. One such result is that all complete separable metric spaces are representable as a natural space. Another is that all natural spaces arise as a quotient space of Baire space, something which Brouwer already incorporated into his theory of spreads. Still another important translated result is the Heine-Borel property for inductive covers of subfanns (\HBind, see \ref{HBforsubfanns}). This relates to (and generalizes we believe) the result in \NegCed\ concerning the Heine-Borel property for (formal-topological) inductive covers of the formal real interval $[\alpha, \beta]$. As an icing to the cake we prove (in \bish) the metrizability of star-finitary natural spaces. The corresponding intuitionistic result from \WaaThe\ has not been transposed to a formal-topological setting (as far as we are aware). The constructive theorem seems comparable to the classical metrization results of strongly paracompact spaces.    

\parr
\enlargethispage{2mm}
Since we can define natural spreads and fans (and even spraids and fanns), we are led to the meta-theorem that natural Baire spreads (fans) correspond precisely to Brouwer's spreads (fans). Moreover, natural morphisms correspond precisely to Brouwer's spread-functions. Natural morphisms are suited, we think, for computational purposes also.\footnote{In formal topology mappings are usually defined as multivalued relations, which seems computationally more complicated, but the author's knowledge on this is again limited.} The definition of spraids and fanns is made to facilitate both computational practice and (topological) lattice theory. For nice papers on computation and constructive topology see \GNSWcom, which also contains a large list of references and detailed historical background, and \BauKav. In \BauKav, recommendations are made for efficient exact real arithmetic, which seem to match our definition of \sigRt\ and \leqc-morphisms quite closely.

\sbsc{To recap}\label{recap}
In this monograph we develop a simple topological framework, which can serve as a theoretical yet applicable basis for dealing with real-world phenomena. The paper is self-contained, but some familiarity with basic topology is probably necessary for understanding its build-up. It turns out that this approach covers not only the real numbers, but in fact all `separable spaces', meaning topological spaces having an enumerable dense subset (for \Rt\ consider \mbox{e.g.}\ \Qt).

\parr

In the first half of the paper we present the framework and formulate the basic theorems and properties. We discuss the strong connection with applied mathematics and physics, and give examples of computational and topological issues in applied mathematics. When developing the theory further, we discover some foundational issues surrounding compactness, leading also to questions on the topological character of our physical universe. 

\parr

The second half of the paper is therefore more foundational in nature, exploring possible avenues to resolve these issues. Links with existing frameworks and theories are discussed, especially classical mathematics (\class), recursive mathematics (\russ), intuitionism (\intu), Bishop-style mathematics (\bish), and formal topology. Natural topology can also be seen as a simplified version of formal topology, or of domain theory. The simplification fits in the historical line of simplification efforts of Heyting, Freudenthal, Kleene, Scott, Bishop, Martin-L\"{o}f, Bridges, Veldman\footnote{Wim Veldman's lecture notes \VelLec\ (in Dutch) are a very nice exposition of intuitionistic mathematics, one hopes for an English translation some day.} and many others (see also Troelstra's and van Dalen's standard treatise on constructive mathematics \TroDal).

\parr

In the final section, we return to discuss the relation with physics. We think that the question of which mathematics to choose for physics deserves more attention. We believe that natural spaces provide a strong conceptual reference frame for physics. Compactness issues also play a role here, since the question `can Nature produce a non-recursive sequence?' finds a negative answer in \CTphys. \CTphys, if true, would seem at first glance to point to \russ\ as the mathematics of choice for physics.  To discuss this issue, we wax more philosophical. We present a simple informal model of \intu\ within \russ, in a two-player game called `Limited Information for Earthlings' (\LIfE) with players `Giver of Digits' (\GoD) and `Humble Mathematician accepting Numbers' (\HuMAN). We also point to \WaaArt\ for a physical experiment which could cast light on \CTphys.

\parr

In the appendix we work out more interesting details regarding the examples given in the first half of the paper. One of these examples is the line-calling decision-support system \hwki, for which we recommend a \textscc{let}\ feature. The appendix also contains most of the proofs. 

\parr

We hope that this monograph will serve as a welcoming introduction to the varieties of (constructive) topology. Our aim is to compare these varieties in such a way that their presentation is simplified and their mutual differences are reduced to their essence. One result that we are happy to mention is that the presented framework of natural spaces gives a faithful classical representation of basic intuitionistic results.

\parr

Looking at the results in this monograph, the author comes to the conclusion that Brouwer's concepts and axioms are still of exceptional elegance and relevance for constructive mathematics. In fact \intu\ is the only constructive theory we know with an elegant pointwise approach which solves all the compactness issues that we studied in chapter three. Perhaps even more relevant: the axioms of \intu\ are precise and appeal directly to the author's mathematical intuition. This monograph should therefore also be seen as strongly supportive of further development of \intu.

\sbsc{Brief historical note on apartness topology}\label{apartnesshistory}
We can say that apartness topology started with Brouwer, and was given a pointfree flavour by Freudenthal in \Freu. This since intuitionistic topology practically entails all the phenomena of apartness topology. Troelstra also studied intuitionistic topology in \TroThe. Martin-L\"{o}f provided a strong germ for constructive pointfree topology in \MLof\ (also see \ref{backback}). In \CoqFor, Coquand even asks the question whether starting formal topology from an inequality relation (apartness) would be worthwhile. Kalantari and Welch (see \mbox{e.g.}\ \KalWel) have been developing related concepts in a recursive-computable framework.

\parr

As far as we can tell, the first definition of `apartness topology' and `apartness space' was given in \WaaThe, a study of modern intuitionistic topology which was set up in such a way as to attract attention also from people in Bishop's school. Apartness topology plays a central role in \WaaThe, just as in this account. The definition in \WaaThe\ was given in \bish, but with the use of \CP\ in mind to arrive at basically the same topology as the natural topology of the current paper.\footnote{The idea of `apartness topology' came from an earlier study of intuitionistic model theory (see \VelWaa), where the current author discovered that first-order sentences using only apartness can describe topological properties, due to the presence of \CP. This gives a correspondence between intuitionistic model theory and classical topological model theory.}

\parr

Some years later, Bridges and V\^{i}\c{t}\u{a} developed a related but different notion, also called apartness topology (or apartness spaces), which does not depend on the use of \CP\ (see \BriVitthr). There have been a number of articles by different authors on this different notion of apartness topology since then. In \BriVitto\ there is also a treatment of lattices and pointfree machinery, but it seems different from the treatment here.

\parr

These articles do not, as far as we are aware, address the question whether we can find interesting apartness spaces which are not already equivalent to a metric space. As we show in this paper, non-metrizable apartness spaces arise naturally in the context of infinite-dimensional topology, an area where constructive methods should be fruitful but which has been somewhat lagging behind in constructive topological investigations, we believe.

\parr

The above should illustrate that the main ideas in this monograph can already be found in older sources. What natural topology has to offer, is a new combination of these ideas (with a high level of detail). Much of classical separable topology (and mathematics) is still uncharted from a constructive perspective. There is in other words yet a nice long road ahead of us.

\sectionnb{Examples}\label{examples}

\sbsc{Hawk-Eye}\label{exhawkeye}
(See \ref{hawkeye}:)
The line-calling decision-support system \hwki\ was critically analyzed by Collins and Evans in \ColEva, in which they raise a concern which resembles our claim (repeated below) and many other concerns about error margins and measuring. However, the precise nature of the problem associated to what they call `digitizing' (making a decision based on continuous data) is left undiscussed. We therefore state explicitly: 
 
\clmm
\hwki, irrespective of the precision of the cameras, will systematically call \textscc{out}\ certain balls which are measurably \textscc{in}\ or vice versa.
\eclmm

\parr
To see why this is so, it is enough to notice that any real measurement and calculation derived from this measurement, resulting in either \textscc{in}\ or \textscc{out}, correspond to a function from \Rratstart\ to $\{\mbox{\textscc{in}}, \mbox{\textscc{out}}\}$. For simplicity's sake let's put the border of a line at the natural real number $0\inn\Rnat$, where a given trajectory end $x$ being \textscc{in}\ corresponds model-wise to $x\geqq 0$. To give \hwki\ credit, we will assume that trajectories are calculated mathematically correctly from data entered. Now the situation for a ball to just touch or just not touch the line can be translated by looking at a shrinking sequence of rational intervals hovering around $0$.

\parr 
\enlargethispage{2mm}
This means that \hwki, for each such shrinking sequence, must yield either \textscc{in}\ or \textscc{out}, after only a finite number of intervals in the sequence (a Wimbledon match must be finished before August, say). Taking the shrinking sequence $x=[-1,1], [-\half,\half], [-\foth,\foth]\ldots$ we determine \hwki's decision on $x$, say \textscc{in}, which is arrived on say at interval $[-2^{-m},2^{-m}]$. Clearly then there are balls whose translated trajectory starts out with $[-1,1], [-\half,\half], [-\foth,\foth],$ {$\ldots, [-2^{-m},2^{-m}]$}, which are nonetheless \textscc{out}\ by a margin of $2^{-m-1}$ (surely measurable, if cameras of \hwki\ can measure up to $[-2^{-m},2^{-m}]$).

\parr

The claim is not per se important. \hwki\ admits to an inaccuracy of around 3 mm.\footnote{Still, one sees `sure' decisions being pronounced by the system where the margin is smaller. A Nadal-Federer match in which this occurred for a margin of 1 mm attracted some media attention to \hwki's (in)accuracy.}\ This is usually blamed on inaccuracy of the camera system. But regardless of camera precision we cannot expect to solve the topological problem that there is no natural morphism from the real numbers to a two-point space $\{\mbox{\textscc{in}}, \mbox{\textscc{out}}\}$ which takes both values \textscc{in}\ and \textscc{out}. To make this precise we define:

\defi
For $m, \ninn$ we write $\overline{\underline{m}}(n)$ for the sequence $m, \ldots, m$ of length $n$. Let $T_{\twop} = \{\zeron\midd\ninn\}\cupp\{\onen\midd\ninn\} = \{0\}^{\!*}\cupp\{1\}^{\!*}$. Similarly, let $T_{\thrp} = \{0\}^{\!*}\cupp\{1\}^{\!*}\cupp\{2\}^{\!*}$. Put $\twos\aprt\zerom\aprt\onen\aprt\twos$ for all $n,m,s>0$. Likewise put $\zerom\leqc\zeron$, $\onem\leqc\onen$ and $\twom\leqc\twon$ for all $m\geq n$. Then $(T_{\twop}, \aprt, \leqc)$ and $(T_{\thrp}, \aprt, \leqc)$ are pre-natural spaces with as maximal dot the empty sequence of length $0$. Write \twotopt\ and \thrtopt\ for the corresponding natural spaces, which up to equivalence contain precisely two (resp. three) points $\underline{0}=0, 00, 000, \ldots $ and $\underline{1}=1, 11, 111, \ldots $ (and resp. $\underline{2}=2, 22, 222, \ldots $).
\edefi

\parr

From theorem \ref{naturalmorphisms}\ it follows directly that any natural morphism $f$ from \Rnatt\ to \twopt\ is constant (meaning either $f(x)\equivv\underline{0}$ for all $x$ in \Rnat\ or $f(x)\equivv\underline{1}$ for all $x$ in \Rnat). So there is no surjective natural morphism from \Rnatt\ to \twopt. 

\rem
This is not the end of the line though for \hwki-like applications. One restriction on morphisms can and should be relaxed in cases like \hwki, namely the restriction that morphisms respect the apartness relation \aprtR. By this we mean that we should turn to morphisms on the unglueing \sigRunglt\ of \sigRt, where \sigRunglt\ is equipped with the finer apartness \aprtomt\ given by $a\aprtom b$ iff $(a\nleqc b \weddge b\nleqc a)$. 

\parr

If we turn to the natural space \Raprtomt\ derived from $(\sigRungl, \aprtom, \leqc)$, then we see that there are many surjective morphisms from \Raprtomt\ to \twopt. This means that we can for example represent the situation $\all x\inn \R [ x\bygr 0 \vee x\smlr 1]$ by a morphism $h$ from \Raprtomt\ to \twopt\ such that $h(x)\iz 0$ implies $x\bygr 0$ and $h(x)\iz 1$ implies $x\smlr 1$.
\erem

So finally, how should \hwki\ be amended? Clearly our recommendation to \hwki\ is to introduce a \textscc{let}-feature. Suppose for simplicity that \hwki's camera-cum-software margin of error can safely be taken to be $2.5$ mm $\approx 2^{-8}$ m.\footnote{There have been concerns about \hwki's accuracy in this respect, raised by Collins and Evans in \ColEva.}\ Since \hwki\ calculates the trajectory end $x$ of a ball from several camera measurements, one should conclude these calculations up to reaching an interval $x_n\iz [a, b]$ of width $2^{-8}$ m (where $[a, b]\iz[\frac{s}{2^9}, \frac{s+2}{2^9}]$ for some $s\inn\Z$). Then one checks whether $0\inn x_n$, and if so, \hwki\ should return the value \textscc{let}, which in turn should lead to a replaying of the point in the tennis match. If $0\smlr a$ (or likewise $b\smlr 0$), then \hwki\ can safely return the value \textscc{in} (or likewise \textscc{out}), with the same consequences for play as are currently in use. 

\parr

This illustrates our remark above, since adapting \hwki\ in this way corresponds to 
creating a morphism from \Raprtomt\ to \thrtopt. This type of problem occurs extremely frequently of course in applied mathematics, and one may think we are merely going over well-trodden grounds. But the topological cadre of reference is seldom explicitized, and systems like \hwki\ illustrate that awareness of this type of problem can still be improved.

\parr

As for visualizing our recommendation for \hwki, why not introduce a narrow gray line (representing \textscc{let}) separating the outer part of the white line (representing \textscc{in}) from the green (representing \textscc{out}). The gray line should be half on the former white band, and half on the former outer green expanse. Then the \textscc{let}-situation occurs only if the ball is calculated to be largely in the outer green, but with a small overlap of gray and no overlap of white. This should be easily understandable to public and players.

\sbsc{Non-archimedean intermezzo: \Cp}\label{nonarc}
For $p$ a prime number, the complex $p$-adic numbers \Cpt\ can be defined constructively by first defining a valuation-norm on the algebraic $p$-adic numbers $\Ap=\{a_m\midd\minn\}$, and then defining \Cpt\ to be the metrical completion of \Apt\ (see \MiRiRu). \Cpt\ has an interesting feature which relates to our study of natural representations of metric spaces:  it is a complete metric space which is not \deff{spherically complete}. There is a series of shrinking closed spheres $(\overline{B}(c_m,q_m))_{\minn}$ where $c_m\inn \Ap$, $q_m\inn \Q$, $\overline{B}(c_{m+1},q_{m+1})\subseteqq \overline{B}(c_m,q_m)$ for all $m$ and yet $\bigcapp_{\minn}\overline{B}(c_m,q_m)=\emptyy$. 

\parr

This means, that if for \Cpt\ we proceed creating a natural space \Vnatt\ as in remark \ref{sepmetnat}, taking $V=\{\overline{B}(a_m, q)\midd q\inn\Q,\minn\}$, then the `shriveling' sequence $(\overline{B}(c_m,q_m))_{\minn}$ becomes a `new' point in \Vcalt\ which has no corresponding point in \Cpt. In fact we do not know whether this \Vnatt\ is metrizable. Classically, \Vnatt\ can be weakly metrized as an extension of \Cp:
\parr
(in \class:) For $p,q\inn\Vcal$, put $d(p,q)=\lim_{n \rightarrow \infty}d(p_n,q_n)$, where for \ninnt\ we let $d(p_n, q_n)\isdef\inf(\{d(x,y)\midd x\inn p_n, y\inn q_n\})$. ($p_n, q_n$ are closed spheres in $(\Cp, d)$.).
\parr
In \class, the limit exists for elements of \Vcalt, and defines a metric which extends the non-archimedean metric $d$ on \Cp.\footnote{Consider that shriveling sphere-sequences $p, q$ converge downwards to a limit radius $r_1, r_2$. If $r_1\smlr r_2$ then $d(p,q)\iz d(p_n, q_n)\bygr r_2$ for some $n$. If $r_1\iz r_2\bygr 0$ then $d(p,q)\iz 0$ or $d(p,q)\iz d(p_n, q_n)\bygr r_2$ for some $n$. If $r_1\iz r_2\iz 0$ then $p, q\inn\Cp$.}\ We think $(\Vcal,d)$ is spherically complete. However, $(\Vcal,d)$ is not separable (so $d$ does not metrize \Vnatt), and we do not know whether an alternative (constructive, separable) metric exists which might enable us to work with some form of spherical completion of \Cpt\ constructively as well. Our metrization theorem for star-finitary spaces (thm.\,\ref{starfinmet}) suggests that one should find out whether \Vnatt\ is star-finitary.

Whether all this leads to new topological spaces, we do not know. \Cpt\ with the usual metric topology is homeomorphic to Baire space, this might also hold for the \Vnatt\ indicated here.\footnote{Perhaps a nice subject for a Master's thesis?}

\sbsc{The Cantor function and other morphisms to the binary reals}\label{cantorfun}
\hspace*{-.9pt}(Con\-ti\-nued from paragraph \ref{morfbinreals}). To show the equivalence of the binary reals with the reals allowing a binary expansion, let \zostart\ be the set of finite sequences of elements of \zot. Now natural Cantor space is the natural subspace \Cnattopt\ of Baire space \Bnatt\ formed by the pre-natural space \Cprenatt\ and its set of points \Cnatt\ (see paragraph \ref{defcantor}\ for the precise definition). 
We also write \sigtot\ for the Cantor space.

\parr

Let us inductively define a surjective \leqc-morphism \fevlbint\ from \Cnatt\ to \zorbint. First we put $\fevlbin(\maxdotB)\iz \maxdot_{\zor}\iz [0,1]\inn \Rratbin$. Now let $a\iz a_0, \ldots, a_{n-1} \inn\zo$ where $\fevlbin(a)$ has been defined and equals $[d,e]\inn\Rratbin$. Then, with $a\starr 0 = a_0, \ldots, a_{n-1}, 0$ and $a\starr 1\iz a_0, \ldots, a_{n-1}, 1$ we define $\fevlbin(a\starr 0)\isdef [d, \frac{d+e}{2}]$ and $\fevlbin(a\starr 1)\isdef [\frac{d+e}{2}, e]$. This inductively defines \fevlbint\ on all of \zostart.

\parr

The reader can simply verify that \fevlbint\ is a surjective morphism from \Cnatt\ to \zorbint. Next, we
pull back the apartness \aprtRt\ on \zorbint\ to \Cprenatt\ by stipulating, for $a, b\inn\zostar$ that $a\aprtR b$ iff $\fevlbin(a)\aprtR\fevlbin(b)$. Then it is easy to see that \fevlbint\ is a \leqc-isomorphism from the natural space derived from $(\zostar, \aprtR, \leqcom)$ to \zorbint.

\parr

Completely similar, we define $\sigthr=(\zotstar, \aprtom, \leqcom)$ and a surjective morphism \fevltert\ from \sigthrt\ to \zortert, such that pulling back \aprtRt\ using \fevltert\ we see that \fevltert\ is a \leqc-isomorphism from $\sigthrR\isdef (\zotstar, \aprtR, \leqcom)$ to \zortert.
We have use for an inverse of \fevltert\ given explicitly, we therefore define the appropriate morphism \fevlterinvt\ inductively. Put $\fevlterinv([0,1])\iz\maxdotB$. Next, suppose $\fevlterinv(a)$ has been defined for a given interval $a\iz[\frac{n}{3^m}, \frac{n+1}{3^m}]\inn\Rratter$. Then we define: $\fevlterinv([\frac{n+0}{3^{m+1}}, \frac{n+1}{3^{m+1}}])=\fevlterinv(a)\starr 0$, $\fevlterinv([\frac{n+1}{3^{m+1}}, \frac{n+2}{3^{m+1}}])=\fevlterinv(a)\starr 1$ and $\fevlterinv([\frac{n+2}{3^{m+1}}, \frac{n+3}{3^{m+1}}])=\fevlterinv(a)\starr 2$. 

\parr

The Cantor function \fcant\ is most easily defined as a morphism from \sigthrt\ to \sigtot.  First take $a\iz a_0, \ldots, a_{n-1}\inn \{0,2\}^*$. For all $i\smlr n\iz\grdd(a)$ let $b_i\iz\min(a_i, 1)$. Then put $\fcan(a)\isdef b_0, \ldots, b_{n-1}$. Next let $b\inn\zotstar, b\notinn\{0,2\}^*$. 
Determine $a\inn \{0,2\}^*$ and $c\inn\zotstar$ such that $b\iz a\starr 1\starr c$. Let $m\iz\grdd(c)$, and let $\zerom$ be the sequence $0, \ldots, 0$ of length $m$. 
Put $\fcan(b)\isdef \fcan(a)\starr 1\starr\zerom$.

\parr

Using \fevlterinv\ and \fevlbint, this completely defines the Cantor function as a \leqc-morphism from \zortert\ to \zorbint. This morphism cannot however be extended to a \leqc-morphism from \zort\ to \zorbint, for the reasons described in chapter two. Proposition \ref{refvstrail2}\ shows how we can represent the Cantor function as a \leqc-morphism from \zort\ to \zort. 

\parr

We can however also define the Cantor function as a \pthh-morphism from \zort\ to \zorbint. We leave this as a non-trivial exercise to the reader interested in applied mathematics and representation issues. For our narrative we turn to the promised property that any morphism $f$ from \zort\ to \zorbint\ is `locally constant' around the $f$-original of a binary rational $\frac{n}{2^m}$ ($n\leqq 2^m$).

\prp
Let $f$ be a \pthh-morphism from \zort\ to \zorbint\ such that $x\leqqR y$ implies $f(x)\leqqR f(y)$, for $x,y\inn\zor$, and such that $f(0)\equivvR 0$ and $f(1)\equivvR 1$. Suppose $z\inn\zor$ is such that $f(z)\equivvR \half$. Then there is a rational interval $[a, b]\inn\Rrat$ such that $f(x)\equivvR \half$ for all $x\inn [a,b]$.    
\eprp

\prf
Clearly we can find a sequence $z'\iz ([a_n, b_n])_{n\in\N}$ of strictly shrinking rational intervals such that $z\equivvR z'$ ($a_n\smlr a_{n+1}\smlr b_{n+1}\smlr b_n$ for all \ninn). We determine a value of $f(z')$ which is not equal to the maximal dot $\maxdot$, say $f(\overline{z'}(m))\precc\maxdot$. Then we have either \kase{0}: $f(z')(m)\leqcR [0, \half]$, then since $f(z')\equivvR\half$ and $f$ is \leqqR-preserving, we see that $f(x)\equivvR\half$ for all $x\inn [b_{m}, b_{m-1}]$ or \kase{1}: $f(z')(m)\leqcR [\half, 1]$, then since $f(z')\equivvR\half$ and $f$ is \leqqR-preserving, we see that $f(x)\equivvR\half$ for all $x\inn [a_{m-1}, a_{m}]$.
\eprf

\crl
\zort\ and \zorbint\ are not isomorphic.
\ecrl

\parr

We end with a number of statements which we do not prove (exercise):

\be
\item[(i)] The $n$-ary reals are isomorphic to the $m$-ary reals, for $n, \minn$. 
\item[(ii)]
The $n$-ary reals can be identically embedded (meaning $f(x)\equivvR x$ for all $x$) in the $m$-ary reals iff there is a $b\geqq 1$ in \Nt\ such that $m$ divides $n^b$.
\item[(iii)] Addition and multiplication cannot be represented as morphisms from \zorbint\ to itself, in other words the $n$-ary reals are not closed under addition and multiplication. This makes the $n$-ary reals unsuited for computational purposes, in our eyes.
\item[(iv)] Adding the extra digit $-1$ to the digits $0,\ldots, n\mino$ solves the current problems with $n$-ary digital representation.  
\ee
 
We believe the above to be of interest for representation issues of real numbers. One can also entertain an independent topological interest, see theorem \ref{patharc}\ on \pathnat\ connectedness and the next example of ContraCantor space. 
Finally, we note that the ternary reals play a very nice role in the proof of our metrization theorem for star-finitary spaces, see \ref{prfstarfinmet}.

\sbsc{ContraCantor space and \zor-embedded Cantor space}\label{contraCan}

From \KleVes\ we can directly define a decidable countable subset (derived from what Andrej Bauer in \BauKle\ aptly calls the Kleene Tree) $K_{\rm bar}=\{k_n\midd\ninn\}$ of \zostart\ such that $k_n\leqcom k_m$ implies $n=m$ for all $n,\minn$ and in addition such that in \russ\ $\{\hattr{k}\midd k\inn K_{\rm bar}\}$ is an open cover of \Cnattop\ which has no finite subcover.

\parr

We use $K_{\rm bar}$ to define ContraCantor space, which is a compact subspace \contrCant\ of \zort\ such that if we write \Canzort\ for the Cantor set (which is the standard embedding of Cantor space in \zort), we see: $\dR(\contrCan, \Canzor)\iz 0$ and yet in \russ\ we also have $\dR(x, \Canzor)\bygr 0$ for $\all x\inn\contrCan$.

\parr

First let us define \Canzort. For this we first embed $\Cnat\iz\sigto$ in \sigthrt, using the `doubling' morphism $f_{\cdot 2}$ defined by: $f_{\cdot 2}(a_0, \ldots, a_{n-1})\iz 2\cdott a_0, \ldots, 2\cdott a_{n-1}$, for $a_0, \ldots, a_{n-1}$ in \zostar. Clearly $f_{\cdot 2}$ is an embedding morphism from \sigtot\ to \sigthrt. We put $\Cnat_{\cdot 2}\iz f_{\cdot 2}(\Cnat)$. Next, we combine $f_{\cdot 2}$ with the morphism \fevltert\ from \sigthrt\ to \zortert, which we defined in the previous example.

We now put: $\Canzor\isdef \fevlter\circc f_{\cdot 2}(\Cnat)$.

\parr

Next we define the ContraCantor set \contrCant, by first considering the subspace $\sigthrc=\{f_{\cdot 2}(k_n)\starr 1\starr\alpha\mid k_n\inn K_{\rm bar}, \alpha\inn\sigthr\}$ of \sigthrt. 
Notice that \sigthrct\ lies apart from $f_{\cdot 2}(\Cnat)$, but at distance $0$. In \russ, \sigthrct\ is (the point set of) a fan. In \intu\ and \class, \sigthrct\ is not closed, and to obtain a fan we must move to the closure of \sigthrc.

We transfer this situation to \zort\ by letting \contrCant\ be the metrical closure of $\fevlter(\sigthrc)$.

\prp
\contrCant\ is a \bish-compact subspace of \zort\ with the property that 
$\dR(\contrCan, \Canzor)\iz 0$ and yet in \russ\ also $\dR(x, \Canzor)\bygr 0$ for all $x\inn\contrCan$.
\eprp

\prf
By definition \contrCant\ is complete and totally bounded, therefore \bish-compact, just like the Cantor set \Canzort. Clearly $\dR(\contrCan, \Canzor)\iz 0$. In \russ, we have $\contrCan=\fevlter(\sigthrc)$, since \sigthrct\ is already complete. This we see by considering in \russ\ a convergent sequence \xnninnt\ in \sigthrct. We construct a `shadow' sequence \ynninnt\ in $\Cnat_{\cdot 2}\iz f_{\cdot 2}(\Cnat)$ thus: for each \ninnt, $x_n$ equals $f_{\cdot 2}(k_m)\starr 1\starr\alpha$ for certain $\minn, \alpha\inn\sigthr$. We now put $y_n\iz f_{\cdot 2}(k_m)\starr\zero$ (where \zerot\ is the infinite sequence $0, 0, \ldots\inn\sigthr$).
\parr
Clearly \ynninnt\ is convergent in $\Cnat_{\cdot 2}$, we consider the limit $y\inn\Cnat_{\cdot 2}$. There is $z\inn\Cnat$ such that $y\iz f_{\cdot 2}(z)$, and since $K_{\rm bar}$ is a bar on \Cnat, there is \ninnt\ such that $z\leqc k_n$. This however implies that there is \Ninnt\ such that $x_m\leqc f_{\cdot 2}(k_n)\starr\one$ for all $m\geqq N$, showing that \xnninnt\ converges to a limit in \sigthrc.
\parr
Therefore we obtain in \russ: $\dR(x, \Canzor)\bygr 0$ for all $x\inn\contrCan$.
\eprf

This situation in \russ, where two compact spaces have distance $0$ and yet are apart, is well known. Perhaps less known is how it bears on our discussion of inductive morphisms and the pointwise problems arising in \bish\ even if we inductivize our definitions, see paragraph \ref{metaruss}. 

\sbsc{The spraid of uniformly continuous real-valued functions on \zort}\label{contfunasspraid}

Brouwer already showed how to build the space $C^{\rm unif}(\zor, \R)_{\rm nat}$ of uniformly continuous real-valued functions on \zort\ as a spread (see \BroCol, \DalInt). The basic idea is to consider the graph of such functions, one then sees that this graph is a compact subspace of the real plane which is homeomorphic to the line segment \zort. For a uniformly continuous real-valued $f$ on \zort, the graph $G_f$ can thus be built as a line segment in the $x, y$-plane twisting from the vertical line $x=0$ to the vertical line $x=1$ without `doubling back' in the horizontal sense. We can approximate the graph $G_f$ with step-by-step growing precision, by forming at stage $n$ a `tape' $T^f_n$ of rectangles, each with height $2^{-n}$, which encloses $G_f$ and which runs from the line $x=0$ to the line $x=1$. The rectangles all have equal width $2^{-m}$ where $m$ is determined by the uniform-continuity modulus of $f$. And we specify that the corner-coordinates of these rectangles are taken from the set $\{(a,b)\midd [2^{n+1}\cdott a\inn\N\weddge 2^{m+1}\cdott b\inn\N]\}$. We can do this in such a way that each rectangle of $T^f_{n+1}$ lies completely within a rectangle of $T^f_n$, for each \ninnt.
\parr
We turn to the whole space $C^{\rm unif}(\zor, \R)_{\rm nat}$, which we wish to represent as a spraid. Abstracting from the specified $f$ above, the properties of these `tapes' $(T^f_n)_{n\in\N}$ can now be formulated in such a way that  we can take the countable collection of all such tapes as the set of basic dots of our desired spraid. Then the definition of the refinement relation and the apartness relation is a direct consequence of our plan thus far. If we stick to this plan, the resulting points will be seen to represent elements of 
$C^{\rm unif}(\zor, \R)_{\rm nat}$, and vice versa, each element of $C^{\rm unif}(\zor, \R)_{\rm nat}$ will correspond to a point in this spraid.

\parr
This indicates how to build $C^{\rm unif}(\zor, \R)_{\rm nat}$ as a spraid. We do not go into this further here, but leave this as a challenge to the reader. 

\sbsc{A counterexample in \class\ illustrating \acoo}\label{excontmorf}
We elaborate on our remarks in \ref{basicneighbor}, by giving an example in \class\ of a natural space \Vnatt\ where for all $x\inn\Vcal$ there is a $y\equivv x$ in \Vcalt\ such that $\closr{y_n}$ is open for all \ninn, and yet \Vnatt\ is not a basic neighborhood space.

For this, we use the `open' refinement relation \leqcRo, and we turn again to \zort\ and \zorbin\ (see defs. \ref{defrealnat}\ and \ref{binarydecimal}). For each $a\inn\zorrat$ we introduce a copy $a^*\!$ which we call \deff{starred}\ and we put:

\btab
$V$ \= $\isdef$ $\{$ \= $[\half\minuz 2^{-n-2}, \half\pluz 2^{-n-2}]\midd\ninn\}\cupp\zorratbin\cupp\{[p, q]^*\midd [p, q]\inn\zorrat\midd p\geqq \half$\\
\> \> $\vee q\leqq\half\}$.
\etab

For starred $a^*\!\inn V$, where $a\inn \zorrat$, put $i\,(a^*\!)\iz a$. And for unstarred $c\inn V$ put $i\,(c)\iz c$. For $c\inn V$ put $c\leqc [0,1]\iz\maxdotV$. Then for $c, d\inn V, d\notiz [0,1]$  we put $c\leqc d$ iff ($c$ is starred $\leftrightarrow$ $d$ is starred) $\weddge\ i\,(c)\leqcRo i\,(d)$. Finally put $c\aprt d$ iff $i\,(c)\aprtR i\,(d)$.

Now let \Vnatt\ be the natural space derived from \Vprenat. We identify the elements of \Vcalt\ with the real numbers that they obviously represent.

\clmm (in \class)
For all $x\inn\Vcal$ there is $\Vcal\nni y\equivv x$ where $\closr{y_n}$ is open for all \ninn.
\eclmm

\clmmprf 
If $x\aprt \half$, then we can trivially find $y\equivv x$ in \Vcalt\ such that $y_n$ is starred and moreover $\closr{y_n}$ is open for all \ninnt\ (by our remarks on \leqcRo, in \ref{defrealnat}). If $x\equivv\half$, then take $y$ given by $y_n\iz[\half\minuz 2^{-n-2}, \half\pluz 2^{-n-2}]$ for \ninnt, and so we are done.
\eclmmprf 

\bclmm (in \bish)
\Vnatt\ is not isomorphic to a basic-open space.
\ebclmm

\clmmprf 
Let \Wnattot\ be a basic-open space such that \Vnatt\ is isomorphic to \Wnattot\ under isomorphism $f$ with inverse $g$. Put $h\iz g\circc f$, then $h$ is an identical automorphism on \Vnatt\ where for all $x\inn\Vcal$ in addition $h(x)_n\iz h(\xstr(n))$ is a basic neighborhood of $h(x)$ for all \ninnt\ (w.l.o.g.\ $h$ is a \pthh-morphism). 
Consider the point $x\equivv\half$ given by $x_n\iz [\half\minuz 2^{-n-2}, \half]^*$ for \ninn. By the above and since $h(x)\equivv\half$ is a point, there must be \ninnt\ such that $h(x)_n$ is of the form $[\half\minuz 2^{-m-2}, \half\pluz 2^{-m-2}]$ for some \minn. 
Now look at the point $y\equivv \half\minuz 2^{-n-2}$ in \Vcalt\ given by $\ystr(n)\iz \xstr(n)$ and $y_{n+s}\iz [\half\minuz 2^{-n-2}, \half\minuz 2^{-n-2}\pluz 2^{-n-2-s}]^*$, for all \sinn. 
Since $\ystr(n)\iz \xstr(n)$, we have
$h(y)_n\iz [\half\minuz 2^{-m-2}, \half\pluz 2^{-m-2}]$. Now the only way to \leqc-refine $h(y)_n$ to a point equivalent to $y$ is by using basic dots in \zorratbin. But these dots do not form a neighborhood of $y$, contradiction. 
\eclmmprf 

Thus our example \Vnatt\ is not a basic neighborhood space, and yet in \class: for all $x\inn\Vcal$ there is a $y\equivv x$ in \Vcalt\ such that $\closr{y_n}$ is open for all \ninn.

\rem
The `reason' that \Vnatt\ is not a basic neighborhood space lies in the fact that in the statement: `$\all x\inn\Vcal\driss y\equivv x\all\ninn\,[\closr{y_n}\inn\Topaprt]$' the information is not given `continuously', that is by a morphism. In fact one can question the statement precisely on this account, since by taking a $z$ resembling the $x$ in the claim-proof above, we see that we have no method to assign to $z$ an appropriate $y$ unless we already know in advance the infinite behaviour of $z$ (of which in general we are ignorant).
\parr
This is precisely the gist of the intuitionistic axiom \acoo\ which states that if we really know $\all x\driss y\,[P(x, y)]$, then this information must be given continuously (by a morphism), otherwise we will always have sequences $x$ for which we have no method to produce $y$ with $P(x,y)$.
\erem

\sectionnb{Proofs and additional definitions}\label{proofs}
\sbsc{Proof of theorem \ref{defrealnat}}\label{realshomeo}

\thm (repeated from \ref{defrealnat})
\Rnattopt\ is a natural space which is homeomorphic to the topological space of the real numbers \Rt\ equipped with the usual metric topology.
\ethm  

\prf
We are a bit free here, since for a classical theorist we should first move to the quotient space of equivalence classes. We consider this a cumbersome practice, and prefer to give a `direct' proof. In this proof we consider \Rt\ to be given as the collection of all Cauchy-sequences in $(\Q, \dR)$. 
\parr
To any $x\iz ([a_n, b_n])_{n\in\N}\inn\Rnat$ we assign the Cauchy-sequence $f(x)\iz\anninn$ in \Rt. We leave it to the reader to verify that (i) for any $y\inn\R$ there is a $z\inn\Rnat$ with $f(z)\iz y$, (ii) for all $x,y\inn\Rnat$ we have $x\aprtR y$ iff $d(f(x), f(y))\bygr 0$. Therefore $f$ is surjective and injective.
\parr
It is easy to see that $f$ is continuous, so to see that $f$ is a homeomorphism we must show that $f$ is open. For this let $U\subseteqq\Rnat$ be \aprt-open. We must show that $f(U)$ is open in \Rdrt. For this let $z\inn f(U)$, determine $x\iz ([a_n, b_n])_{n\in\N}\inn U$ such that $f(x)\iz z$. It is not hard to construct a $y\iz ([c_n, d_n])_{n\in\N}\inn\Rnat$ such that $y\equivvR x$ and for all \ninnt: $c_{n+1}\minuz c_n\bygr \foth\cdott(d_n\minuz c_n)$ and $d_{n}\minuz d_{n+1}\bygr \foth\cdott(d_n\minuz c_n)$.
\parr
Since $U$ is open and $y\equivvR x\inn U$, we can find \ninnt\ such that $\hattr{[c_n, d_n]}\subseteqq U$. So $f(U)$ contains the interval $[c_n, d_n]$, and also $z\iz f(x)\equivv f(y)\inn [c_{n+1}, d_{n+1}]$. From this we conclude that $f(U)$ contains the metric ball $B(z,\foth\cdott(d_n\minuz c_n))$, showing that $f(U)$ is open.
\eprf

\sbsc{Proof of theorem \ref{basicneighbor}}\label{prfcontmorf}

For this proof we use theory from later sections. Especially simplifying is theorem \ref{Baireuni}\ which states that every natural space is spreadlike, and even isomorphic to a spread whose tree is $(\Nstar,\leqcom)$. The corollary from its proof in \ref{Baireuniprf}\ shows that a basic-open space is isomorphic to a basic-open spread whose tree is $(\Nstar,\leqcom)$. We also use the terminology of later chapters, notably chapter three.

\defi  
Let \Vnatt\ be a spread derived from \Vprenatt, and let $B$ be a bar on $V$ (see def.\,\ref{indcov}). Then we say that $B$ is a \deff{thin bar}\ on $V$ iff for all $a\inn B$ and $b\precc a$ we have that $b\notinn B$. (Then for a successor point $x\inn\Vcal$ there is a unique \ninnt\ with $x_n\inn B$, see \ref{prfcrlstarfinmet} (\bol{a})).
\edefi

\thm (in \class, \intu\ and \russ, repeated from \ref{basicneighbor}):

Let $f$ be a continuous function from a natural space \Vnatot\ to a basic neighborhood space \Wnattot.  Then there is a natural morphism $g$ from \Vnatot\ to \Wnattot\ such that for all $x$ in \Vcalt: $f(x)\equivvto g(x)$.

\ethm 

\prf
We give a unified proof for \class, \intu\ and \russ\ derived from the common Lindel\"{o}f axiom \BDD$^*$ defined in \ref{axiombdd}, which states that every bar on \Bnatt\ (or equivalently \Nstar) descends from a thin decidable bar.

\parr

By theorem \ref{Baireuni}\ and its corollary (proved in \ref{Baireuniprf}) it suffices to prove the theorem for the case where \Vnatot\ is a spread derived from $(\Nstar,\aprto,\leqcom)$ and \Wnattot\ is a basic-open spread  derived from $(\Nstar,\aprtto, \leqcom)$. So in the following keep in mind that $\closr{b}$ is $\aprtto$-open for every $b\inn W\iz\,\Nstar$.

\parr

We will inductively define a sequence of thin bars $(C_n)_{n\in\N}$ on $V\iz\,\Nstar$ and simultaneously construct the desired morphism $g$, as follows. First let $n\iz 0$, put $C_0\iz\{\maxdotV\}$ and put $g(\maxdotV)\iz\maxdotW$. Let $c\inn C_0$, then since $f$ is a continuous function, for all $x\inn\hattr{c}\iz\Vcal$ there are \sinnt\ with $x_s\precc c$ and $b\inn \suczz(g(c))\iz\suczz(\maxdotW)$ such that $f(\hattr{x_s})\subseteqq \closr{b}$.

\parr
Therefore the set $B_1\iz\{a\inn\Nstar\midd \driss c\inn C_0[a\precc c \weddge \driss b\inn\suczz(g(c)) [f(\hattr{a})\subseteqq \closr{b}]\}$ is a bar on $V\iz\,\Nstar$. By \BDD$^*$ we find a decidable thin bar $C_1$ on $V\iz\,\Nstar$ from which $B_1$ descends. 
Then for all $a\inn C_1$ we have: $\dris b\inn\suczz(g(\maxdotV)) [f(\hattr{a})\subseteqq \closr{b}]$.

\parr

Using countable choice (\aczz) we can now assign to each $a\inn C_1$ a value $g(a)$ in $\suczz(g(\maxdotV))=\suczz(\maxdotW)$ such that $f(\hattr{a})\subseteqq\closr{g(a)}$. Also, to each $b\inn V$ for which there is $d\precc b\precc c$ with $d\inn C_1, c\inn C_0$ we assign: $g(b)\iz g(c)\iz\maxdotW$.

\parr

We can repeat this process for $n\iz 1$ and $C_1$. For let $c\inn C_1$, then since $f$ is continuous, for all $x\inn\hattr{c}$ there are \sinnt\ with $x_s\precc c$ and $b\inn\suczz(g(c))$ such that $f(\hattr{x_s})\subseteqq \closr{b}$.

\parr

Therefore the set $B_2\iz\{a\inn\Nstar\midd \driss c\inn C_1[a\precc c \weddge \driss b\inn\suczz(g(c)) [f(\hattr{a})\subseteqq \closr{b}]\}$ is a bar on $V\iz\,\Nstar$. By \BDD$^*$ we find a decidable thin bar $C_2$ on $V$ from which $B_2$ descends. 
Then for all $a\inn C_2$ we have: $\dris c\inn C_1\driss b\inn\suczz(g(c)) [f(\hattr{a})\subseteqq \closr{b}]$

\parr

Using countable choice (\aczz) we can now assign to each $a\inn C_2$ a value $g(a)$ in $\suczz(g(c))$ (where $c\inn C_1$ is such that $a\precc c$) such that $f(\hattr{a}\subseteqq\closr{g(a)}$. Also, to each $b\inn V$ for which there is $d\precc b\precc c'$ with $d\inn C_2, c'\inn C_1$ we assign: $g(b)\iz g(c')$.

\parr

We can repeat this process for $n\iz 2$ and $C_2$, etc. In this way, by using dependent countable choice (\dco) we can construct $(C_n)_{n\in\N}$ and simultaneously define $g$ on all of $V\iz\Nstar$. It is easy to see that $g(x)\equivvto f(x)$ for all $x\inn\Vcal$, which by continuity of $f$ shows that $g$ is a morphism. 
\eprf

\sbsc{Proof of theorem \ref{sepmetnat}}\label{sepmetnatprf}

\thm (repeated from \ref{sepmetnat})
Every complete separable metric space \xdt\ is homeomorphic to a basic-open space \Vnatt. 
\ethm

\prf
The rough idea is simple: for a separable metric space \xdt\ with dense subset \anninnt, let for each $n,\sinn$ a basic dot be the open sphere $B(a_n,2^{-s})=\{x\inn X\midd d(x, a_n)\smlr 2^{-s}\}$.  Then we can take as set of dots $V=\{B(a_n, 2^{-s})\midd n,\sinn\}\cupp\{\maxdotV\}$. The technical trouble now is to define \aprt\ and \leqct\ constructively, since in general even for $s\bygr t$ the containment relation $B(a_n, 2^{-s})\subseteqq B(a_m,2^{-t})$ is not decidable. However, this containment relation has an enumerable subrelation which also does the trick. This because for all $(a_n, s)$ and $(a_m, t)$ with $s\bygr t$ there is $i\inn\{0,1\}$ such that:
\parr
$(i\iz 0 \weddge d(a_n, a_m)\smlr 2^{-t}\minuz 2^{-s})$ or $(i\iz 1\weddge d(a_n, a_m)\bygr 2^{-t}\minuz 2^{-s}\minuz 2^{-2s})$
\parr
Using \aczz\ (countable choice) we can define a function $h$ fulfilling the above statement. Now we put $B(a_n, 2^{-s})\precc B(a_m, 2^{-t})$ iff $h((a_n, s), (a_m, t))\iz 0$.
Likewise we define \aprt, since for all $(a_n, s)$ and $(a_m, t)$ there is $j\inn\{0,1\}$ such that:\parr
$({^{\,}}j\iz 0 \weddge d(a_n, a_m)\smlr 2^{-s\!}\pluz 2^{-t\!}\pluz 2^{-s-t} )$ or $({^{\,}}j\iz 1\weddge d(a_n, a_m)\bygr 2^{-s\!}\pluz 2^{-t\!}\pluz 2^{-s-t-1} )$
\parr
Using \aczz\ we can define a function $g$ fulfilling the above statement. Now we simply put $B(a_n, 2^{-s})\aprt B(a_m, 2^{-t})$ iff $g((a_n, s), (a_m, t))\iz 1$.

It is not difficult to see that \Vprenatt\ generates a basic-open natural space which is homeomorphic to \xdt. 
\eprf

\rem
Notice that by our definition of $\aprt$, if $B(a_n, 2^{-s})\aprt B(a_m, 2^{-t})$, then $d(a_n, a_m)\bygr 2^{-s}\pluz 2^{-t}\pluz 2^{-s-t-1} )$. This is an important detail for proving theorem \ref{compmetind}.
\erem

\crl 
In \class, \intu\ and \russ\ the following holds:
\be
\item[(i)]
A continuous function $f$ from a natural space \Wnatt\ to a complete metric space \xdt\ can be represented by a morphism from \Wnatt\ to a basic neighborhood space \Vnatt\ homeomorphic to \xdt, by theorem \ref{basicneighbor}. 
\item[(ii)] A representation of a complete metric space as a basic neighborhood space is unique up to isomorphism.
\ee

In \bish\ the following holds:
\be
\item[(iii)] If \xdt\ and \Vnatt\ are as above in the theorem, then we can define a metric $d'$ on \Vnatt\ (see def.\,\ref{metnat}) by defining $d'(x, y)\iz d(h(x), h(y))$ for $x,y\inn\Vcal$ and $h$ a homeomorphism from \Vnatt\ to \xdt. This metric can be obtained as a morphism from $(\Vcal\timez\Vcal, \Topaprt)$ to \Rnatt\ by the construction of \Vnatt. We then see that the apartness topology and the metric $d'$-topology coincide, in other words \Vnatt\ is metrizable.
We conclude: on (this natural representation of) a complete metric space, the metric topology coincides with the apartness topology.
\ee
\ecrl

\sbsc{Proof of theorem \ref{patharc}}\label{binrealpatharc}

\thm (repeated from \ref{patharc})
\Rbint\ (as well as \Rtert, \Rdect) is a \pathnat\ connected space which is not \arcnat\ connected.
\ethm

\prf
A detailed constructive proof for \Rtert\ is given in \WaaThe, this construction can be literally transposed to our setting to show that \Rtert\ is \pathnat\ connected. That \Rtert\ is not \arcnat\ connected follows from our work in example \ref{cantorfun}. We sketch an alternative proof using the Cantor function (see example \ref{cantorfun}\ and paragraph \ref{morfbinreals}). 

\parr

Suppose $x\smlrR y\inn\zorbin$, we want to show that there is a morphism $f$ from \zort\ to \zorbint\ such that $f(0)\equivvR x$ and $f(1)\equivvR y$. It is not so difficult to see that \zorbint\ is isomorphic to $\{z\inn\zorbin\midd x\leqqR z\leqqR y\}$ under an isomorphism $g$ with $g(0)\equivv x$ and $g(1)\equivv y$. This means that we can take
$f\isdef g\circc\fcan$. For $y\smlrR x$ we can mirror this argument.

\parr

Now if $x, y\inn\zorbin$ such that at stage $n$ we still cannot determine $x\aprtR y$, then we can still start constructing $f$, sending initial values of $0\inn\zor$ to initial values of $x$ and intial values of $1\inn\zor$ to initial values of $y$, in such a way that if at any later stage $m$ we see $x(m)\aprtR y(m)$, then we can continue as above in the case where $x\smlrR y$ or $y\smlrR x$.

\parr

Of course, to complete the proof one must show that \fcant\ can indeed be given as a \pthh-morphism from \zort\ to \zorbint\ (we left this as a non-trivial exercise in example \ref{cantorfun}). One also needs to extend the proof for \zorbint\ to all of \Rbint, which involves some extra work since \Rbint\ is not closed under addition and multiplication. Finally, one can use the exercise that all the $n$-ary reals are isomorphic to transfer the \pathnat\ connectedness of \Rbint\ to all the $n$-ary reals.
\eprf 

\sbsc{Proof of theorem \ref{Baireuni}}\label{Baireuniprf}
We need a preparatory definition.

\defi
Let \Vnatt\ be a natural space derived from \Vprenatt. We introduce a formal element \mindott\ not contained in $V$, and put $\mindot\aprt a$ for all $\maxdotV\notiz a\inn V$, and $\negg(\mindot\aprt\mindot)$. Also put $\Vmindot\iz V\cupp\{\mindot\}$. Then the set $\aprtVmin\iz\{(c,d)\inn \Vmindot\timez \Vmindot\midd c\aprt d\}$ is countable.\footnote{This is the reason for introducing \mindott, since for a space containing just one point up to equivalence, $\aprtV$ is empty, and there are many spaces of which we don't know whether they contain more than one point.}\ 
\parr
We say that $e: \N\rightarrow\aprtVmin$ is a \deff{pregrade}\ on $V$ iff $e$ is an enumeration of either $\aprtVmin$ or of $\aprtV\iz \{(c,d)\inn V\timez V\midd c\aprt d\}$ \footnote{This last addition is to avoid cumbersome notation if we have a space containing at least two points.}. Let $e$ be a pregrade on $V$, then for \ninnt\ and a basic dot $a$ we say that $a$ chooses on $e_n\iz (c,d)$ iff $a\aprt c$ and/or $a\aprt d$.  We now define a decidable gradation on basic dots as follows: for every \ninnt\ a basic dot $a$ is of \deff{$e$-grade $n$}\ (notation $\grad_{e}(a)\geqq n$) iff $a$ chooses on $e_i$ for every $i\smlr n$. 
\parr
In addition, let $v: \N\rightarrow V$ be an injective enumeration of $V$. 
For each \ninnt\ we define the decidable set $\Bven\iz\{v_m\inn V\midd m\geqq n \weddge \grad_{e}(v_m)\geqq n\}$. 

Finally, for \ninnt\ a basic dot $a$ is of \deff{$(v,e)$-grade $n$}\ iff $a$ is in $\Bven$ and not in $B^{v,e}_{n+1}$. This is also decidable, so every basic dot has a unique decidable $(v,e)$-grade. Notice that $B^{v,e}_0\iz V$, and $\maxdotV\notinn B^{v,e}_1$.
\edefi

\lem
Let \Vnatt\ be a natural space derived from \Vprenatt. Let $v: \N\rightarrow V$ be an injective enumeration of $V$, and let $e: \N\rightarrow\aprtVmin$ be a pregrade on $V$.
Then a sequence $x\iz x_0\geqc x_1\geqc x_2 \ldots$ of basic dots in $V$ forms a point in \Vcalt\ iff for every \ninnt, there is an \minnt\ such that $x_m\inn \Bven$.
\elem

\prf
The proof is a simple checking of the definitions, which we leave to the reader as an exercise.
\eprf

We are now ready to prove the main theorem:

\thm (repeated from \ref{Baireuni})
Every natural space is spreadlike. In fact, every natural space \Vnatt\ is isomorphic to a spread \Wnatt\ whose tree is $(\Nstar, \leqcom)$.
\ethm

\prf
Let \Vnatt\ be a natural space derived from \Vprenatt. Let $v: \N\rightarrow V$ be an injective enumeration of $V$, and let $e: \N\rightarrow\aprtVmin$ be a pregrade on $V$.
For $d\inn V$ we put $\{d\}_{\!\leqc}\iz\{c\inn V\midd c\leqc d\}$.
\parr
We inductively define a sequence of functions $(h_n)_{n\in\N}$ from $^{n}\Nstar = \{a\inn\Nstar \mid \grdd(a)\iz n\}$ to $\Bven$ such that for $b$ in $^{n+1}\Nstar$ and $a$ in $^{n}\Nstar$: if $b\inn\suczz(a)$ then $h_{n+1}(b)\precc h_n(a)$. 
In addition, the functions $(h_n)_{n\in\N}$ will be `surjective enough', meaning that every point in \Vcalt\ will be represented in the end, when we join all the $h_n$'s to a single morphism $h$ from Baire space to \Vnatt.
\parr
First, let $h_0$ be the function from $\{\maxdotB\}$ to $B^{v, e}_0$ given by $h_0(\maxdotB)\iz\maxdotV$.
\parr
We turn to $n\iz 1$. We can determine a bijection $g_{\maxdotBs}\iz h_1$ from $\sucz(\maxdotB)\iz^1\Nstar$ to $B^{v, e}_1$. Trivially for $b$ in $^{1}\Nstar$ and $a$ in $^{0}\Nstar\iz\{\maxdotB\}$ we have: if $b\inn\suczz(a)$ then $h_{1}(b)\precc h_0(a)$. 
\parr
Next we turn to $n\iz 2$. For $a\inn^{1}\Nstar$ we look at $h_1(a)\inn V$. Remember that $\{h_1(a)\}_{\!\precc}=\{b\inn V\midd b\precc h_1(a)\}$. We hold: $\Bsubba\isdef B^{v, e}_2\capp \{h_1(a)\}_{\!\precc}$ is a countable set (since it is infinite, and all the relevant relations are decidable). Therefore we can determine a bijection $g_a$ from $\suczz(a)$ to $\Bsubba$. This means that we can take $h_2\iz\bigcupp_{a\in {^{1}\Nstar}}g_a$, and see that for $b$ in $^{2}\Nstar$ and $a$ in $^{1}\Nstar$ we have: if $b\inn\suczz(a)$ then $h_{2}(b)\precc h_1(a)$. 
\parr
Now we are in business, since for $n\iz 3$ and so on the above process can be continued verbatim, changing only the index $n$. This yields two sequences of functions, which arise intertwinedly. The first sequence is $(h_n)_{n\in\N}$ from $^{n}\Nstar\iz\{a\inn\Nstar\midd \grdd(a)\iz n\}$ to $\Bven$ such that for $b$ in $^{n+1}\Nstar$ and $a$ in $^{n}\Nstar$: if $b\inn\suczz(a)$ then $h_{n+1}(b)\precc h_n(a)$. The second sequence is $(g_a)_{a\in\Nstar}$, where for $a\inn ^{n}\Nstar$, $g_a$ is a bijection from $\suczz(a)$ to $\Bsubba\iz B^{v,e}_{n+1}\capp \{h_n(a)\}_{\precc}$.  

(We do not really use \dco\ since all these functions can be found canonically once we have fixed our enumerations $v$ and $e$).

\clmm
If for $a\inn ^n\Nstar$ we put $h(a)\iz h_n(a)$, then $h$ is a surjective morphism from Baire space \Bnatt\ to \Vnatt. If we introduce an apartness $\aprtW$ on \Nstart\ by putting $a\aprtW b$ iff $h(a)\aprtV h(b)$, then $h$ is an injective morphism from the spread \Wnatt\ derived from $(\Nstar, \aprtW, \leqcom)$ to \Vnatt.
\eclmm

\clmmprf 
By the above lemma, $h$ is surjective since $g_a$ is a bijection from $\suczz(a)$ to $\Bsubba$ for $a\inn \Nstar$. 
Now \Wnatt\ is a spread, since for $a\inn ^{n}\Nstar$ we have $h(a)\inn\Bven$, so infinite $\preccom$-trails define points. 
\eclmmprf 

To complete the proof, we need to show that $h$ is an isomorphism from \Wnatt\ to \Vnatt.
Therefore we need to construct an inverse $h^{\rm inv}$ for $h$. In general this inverse $h^{\rm inv}$ can only be constructed as a trail morphism (see def.\,\ref{pathmorphisms}), in other words a refinement morphism on the trail space of \Vnatt.  So we turn to the trail space \Vpathnatt, derived from \Vpathpret. We also use the above definition and lemma, and the sequences of functions $(h_n)_{n\in\N}$ and $(g_a)_{\raisebox{.2ex}{$\maxdots$}\neq a\in\Nstar}$ above.

\parr

To start put $h^{\rm inv}(\maxdotp)\iz\maxdotB$. Next, let $a\iz a_0, \ldots, a_n$ be in \Vpatht\ (so $a_0\succ\ldots\succ a_n$). We associate to $a$ a unique `minimal grade sequence' $p_a\iz p_0, \ldots, p_{j-1}$ with $j\leqq n\pluz 1$, where $p_a$ is a subsequence of $a$. First we take $i_0\iz \mu k\leqq n [a_k\inn B^{v,e}_1]$ if such $k$ exists.  If such $k$ does not exist, we are done and the minimal grade sequence associated to $a$ is the empty sequence and $j\iz 0$. Next, if $j\notiz 0$ we take $i_1\iz \mu k\leqq n [k\bygr i_0 \weddge a_k\inn B^{v,e}_2]$ if such $k$ exists (else $j\iz 1$ and we are done). And generally let $i_s\iz \mu k\leqq n [k\bygr i_{s-1}\weddge a_k\inn B^{v,e}_{s+1}]$, until we have exhausted $a$. Then $p_a\iz p_0, \ldots, p_{j-1}\iz a_{i_0}, \ldots, a_{i_{j-1}}$.

\parr

To define $h^{\rm inv}(a)$ we use the minimal grade sequence $p_0, \ldots, p_{j-1}$. If $j\iz 0$ then we put $h^{\rm inv}(a)\iz \maxdotB$. Else, we can first turn to $p_0$. Since $h_1$ is a bijection, we can determine $b_0\iz h_1^{-1}(p_0)\inn ^1\Nstar$. If $j\iz 1$ we are done, else we know that $g_{b_0}$ is a bijection from $\suczz(b_0)$ to $B_{p_0}\iz B^{v,e}_2\capp \{p_0\}_{\!\precc}$. Then we can put $b_1\iz g_{b_0}^{-1}(p_1)$. In this way we can continue, putting  $b_{s+1}\iz g_{b_s}^{-1}(p_{s+1})$ for all $s\smlr j\minuz 1$. Finally we put $h^{\rm inv}(a)\iz b_0, \ldots, b_{j-1}$. Since $a$ is arbitrary, this suffices to define $h^{\rm inv}$ on \Vpatht, and also on \Vnatt\ by putting $h^{\rm inv}(x)\iz h^{\rm inv}(\xsstr(0)), h^{\rm inv}(\xsstr(1)), \ldots$ for $x\inn\Vcal$.

\clmm
$h^{\rm inv}$ is a trail morphism from \Vnatt\ to \Wnatt\ such that $h^{\rm inv}\circ h(x)\equivv x$ for all $x\inn \Wcal$ and $h\circ h^{\rm inv}(y)\equivv y$ for all $y\inn\Vcal$.
\eclmm

\clmmprf 
The only real concern for showing that $h^{\rm inv}$ is actually a \pthh-morphism with the correct properties, is that $h^{\rm inv}$ sends points in \Vcalt\ to points in \Wcalt. So let $x\inn\Vcal$, we need to show that $h^{\rm inv}(x)$ is in \Wcalt. Let $n\inn\N, n\geqq 1$ be arbitrary. By the above lemma, we know that there is $m_0\inn\N$ such that $x_{m_0}\inn B^{v,e}_1$, then $m_1\inn\N, m_1\bygr m_0$ with $x_{m_0}\succ x_{m_1}\inn B^{v,e}_2$, etc., until we find $m_{n-1}\bygr m_{n-2}$ with $x_{m_{n-2}}\succ x_{m_{n-1}}\inn \Bven$. Then the minimal grade sequence associated to $\xsstr(m_{n-1}\pluz 1)$ has at least length $n$. Therefore $h^{\rm inv}(\xsstr(m_{n-1}\pluz 1))$ is a sequence of length at least $n$. Thus we see that $h^{\rm inv}(x)$ is an infinite sequence of ever-shrinking basic dots of Baire space, and so a point in Baire space. 
\parr
To show that $h^{\rm inv}(x)$ is also a point in \Wcalt, consider $a, c\inn\Nstar$ with $a\aprtW c$. This by definition of $\aprtW$ means that $h(a)\aprtV h(c)$. Therefore $h(a)\aprtV h(c)$ is one of the $e$-enumerated pairs of apart dots, say that $(h(a), h(c))\iz e_M$. By the above reasoning which demonstrated that $h^{\rm inv}(x)$ is a point of Baire space, we know that there is an \minnt\ such that the minimal grade sequence $p_0, \ldots, p_{j-1}$ associated to $\xsstr(m)$ has length $j\iz M\pluz 1$. By carefully looking at the construction of $h^{\rm inv}$ we see that $h^{\rm inv}(\xsstr(m))$ is a sequence $b_0, \ldots, b_{M}$ where $b_M\inn B^{v,e}_{M+1}$ and so a fortiori $b_{M}\aprtW a$ and/or $b_{M}\aprtW c$.
\parr
Therefore $h^{\rm inv}(x)$ is also a point in \Wcalt. We leave it to the reader to verify that $h^{\rm inv}\circ h(x)\equivv x$ for all $x\inn \Wcal$ and $h\circ h^{\rm inv}(y)\equivv y$ for all $y\inn\Vcal$.
\eclmmprf 
\minnewlyne
\eprf

\crl
\be
\item[(i)] Let \Vnatt\ be a natural space, then there is a surjective \leqc-morphism from Baire space to \Vnatt. (`Baire space is a universal spread', `every natural space is the natural image of Baire space', `every natural space is a quotient topology of Baire space'). 
\item[(ii)] If \Vnatt\ is a basic-open space (see definition \ref{basicneighbor}) then \Vnatt\ is isomorphic to a basic-open spread \Wnatt\ whose tree is $(\Nstar, \leqcom)$. 
\ee
\ecrl

\sbsc{Proof of theorem \ref{Cantuni}}\label{Cantuniprf}

That Cantor space is a universal fan is an easy and -in different terminology- well-known result, but we prove it anyway.

\thm (repeated from \ref{Cantuni})
Let \Vnatt\ be a fann, then there is a surjective morphism from Cantor space to \Vnatt.
\ethm

\prf
For each $a\inn V$ we need to fix an element $x^{a}\inn\Vcal, x^{a}\leqc a$. For this we use a bijection $v: \N\rightarrow V$ with inverse $v^{-1}$. Now for $a$ in $V$, let $x^{a}$ be the unique point in \Vcalt\ such that $x^{a}_0\iz a$ and for each \ninnt\ we have: $v^{-1}(x^{a}_{n+1})=\mu\sinn [ v(s)\precc x^{a}_{n}]$.\,\footnote{This construction also works for spraids, not just for fanns.}

\parr

We define the surjective morphism $f$ inductively. Determine the finite subset $\suczz(\maxdotV)\iz\{a_i\midd i\leqq m\}\subsett V$ (in some enumeration, for some \minn). Let \ninnt\ be such that $2^n\smlr m\leqq 2^{n+1}$. Putting $^{n+1}\zostar\iz\{b\inn\zostar\midd \grdd(b)\iz n\pluz 1\}$, we see that $^{n+1}\zostar$ contains sufficient elements $b_i$ (numbered in the obvious lexicographical/binary way) to put $f(b_i)=a_i$ for all $i\leqq m$. For $m\smlr i\leqq 2^{n+1}$, we fix $f$ by sending all of $\{y\inn\Cnat\midd y\leqc b_i\}$ to the fixed $x^{a_0}\inn \Vcal, x^{a_0}\leqc a_0$.

\parr

Now notice that $V_{{\!}a_i}$ again determines a fann for all $i\leqq m$, and so we can repeat this process for $V_{{\!}a_i}$ and $\{y\inn\Cnat\midd y\leqc b_i\}$ which is isomorphic to \Cnat.
This inductively defines $f$. 
We leave it to the reader to verify that $f$ is the desired surjective morphism.
\eprf

\sbsc{Proof of theorem \ref{indcov}}\label{prfindcov}
We prove that for spraids the formal inductive covering relation \icovbt\ equals the genetic inductive covering relation \oicovbt.

\thm (from \ref{indcov})
Let \Vnatt\ be a spraid derived from \Vprenatt, and let $E, F\subseteqq V$. Then $F\icovb E$ iff $F\oicovb E$.
\ethm

\prf 
We necessarily use both \PFI\ and \PGIg. First let $F\oicovb E$, which means by definition that $\{a\}\oicovb E$ for all $a\inn F$. Let $a\inn F$, then $\{a\}\oicovb E$ which means $E$ descends from a genetic bar $G$ on $\Vsuba$. We now show by genetic induction on $G$ that $\{a\}\icovb E$. 

\be
\item[\gindz] If $G\iz\{\maxdot_a\}\iz\{a\}$, then $a\inn E_{\!\leqc}$, so by \indob\ and \indthrb\ we see that  $\{a\}\icovb E$.
\item[\gindom] Let $G\iz\bigcupp_{b\in\suczs(a)} \Bsubb$ where for all $b\inn\suczz(a)$ we have that $\Bsubb$ is a genetic bar on $\Bsubb$ such that if $E$ descends from $\Bsubb$, then $\{b\}\icovb E$. However, we do indeed know that $E$ descends from $\Bsubb$ for all $b\inn\suczz(a)$ by the assumption on $G$. Therefore we find: $\{b\}\icovb E$ for all $b\inn\suczz(a)$. This means that $\suczz(a)\icovb E$ by \indtob. On the other hand, one sees by \indob\ and \indtob\ that $\{c\midd c\precc a\}\icovb\suczz(a)$. Combining this with \indfib\ we obtain $\{a\}\icovb\{c\midd c\precc a\}\icovb\suczz(a)\icovb E$, and so by \indfob, we see that $\{a\}\icovb E$.
\ee

Since $a$ is arbitrary, this means that by \indtob\ we can conclude $F\icovb E$.

\parr

Now for the implication in the other direction, let $P(F, E)$ be the property: $F\oicovb E$. We prove that $P$ satisfies \indo\ through \indfi\ (for subsets $F, E, D$ of $V$):

\bei
\itemmz{\indo} if $b\leqc c$ in $V$, then $\{b\}\oicovb\{c\}$ since $\{c\}$ descends from the genetic bar $\{\maxdot_b\}$ on $\Vsubb$.

\itemmz{\indto} if for all $a\inn F$ we have $\{a\}\oicovb E$, then by definition of \oicovbt, we have $F\oicovb E$.

\itemmz{\indthr} if $F\oicovb D$ and $D\subseteqq E$, this trivially implies that $F\oicovb E$.

\itemmz{\indfo} suppose $F\oicovb D\oicovb E$. Let $a$ in $F$, then $\{a\}\oicovb D\oicovb E$. We must show that $\{a\}\oicovb E$. Let $G$ be a genetic bar on $\Vsuba$ such that $D$ descends from $G$. To show that $\{a\}\oicovb E$ we use genetic induction and propositions \ref{genbarsext}\ and \ref{prfindmorfresp}\ to prove the statement: `Let $B$ be a genetic bar on $\Vsuba$ such that $D$ descends from $B$, then $\{a\}\oicovb E$'.

\be
\item[\gindz] If $B\iz\{\maxdot_a\}\iz\{a\}$, then since $D$ descends from $B$ there is a $d\inn D$ with $a\leqc d$. Since $\{d\}\oicovb E$ we find a genetic bar $H$ on $\Vsubd$ such that $E$ descends from $H$. Now by proposition \ref{genbarsext}, the reduction $H^{\uparrow a}$ of $H$ to $\Vsuba$ contains a genetic bar $H'$ on $\Vsuba$. Clearly $E$ descends from $H'$ as well, showing that $\{a\}\oicovb E$.
\item[\gindom] Let $B\iz\bigcupp_{b\in\suczs(a)}\Bsubb$ where for each $b\inn\suczz(a)$ we have that $\Bsubb$ is a genetic bar on $\Vsubb$ satisfying: if $D$ descends from $\Bsubb$ then $\{b\}\oicovb E$. Still, we already know that $D$ descends from $\Bsubb$ for each $b\inn\suczz(a)$. Therefore for each $b\inn\suczz(a)$ we know $\{b\}\oicovb E$. This means that for each $b\inn\suczz(a)$ we have a genetic bar $H_b$ on $\Vsubb$ such that $E$ descends from $H_b$. Now by proposition \ref{prfindmorfresp}\ we see that $H\iz\bigcupp_{b\in\suczs(a)}H_b$ is a genetic bar on $\Vsuba$ such that $E$ descends from $H$, showing that $\{a\}\oicovb E$.
\ee

By \PGIg\ our statement is proven, and we conclude (since $G$ is a genetic bar on $\Vsuba$ such that $D$ descends from $G$) that $\{a\}\oicovb E$. Since $a\inn F$ is arbitrary, this shows that $F\oicovb E$.

\itemmz{\indfi} $b\oicovb\{d\midd d\precc b\}$, since $\{d\midd d\precc b\}$ contains the genetic bar $\suczz(b)$ on $\Vsubb$.

\eei 

By \PFI\ we now conclude that $F\icovb E$ implies $F\oicovb E$.
\eprf

\sbsc{Proof of proposition \ref{defindmorf}}\label{prfdefindmorf}

We wish to show that for a spraid \Vnatt\ derived from \Vprenat,  genetic bars on $V$ correspond to genetic bars on \Vunglt\ in a precise way. We could call this the unglueing of genetic bars on $V$. From this correspondence it follows that trail morphisms are inductive iff they are inductive as refinement morphism.

\lem
Let \Vnatt\ be a spraid derived from \Vprenat. Let $c\inn V, d\inn\Vungl$ and let $G\subseteqq V$ be a genetic bar on $\Vsubc$ and $H$ a genetic bar on $\Vungl_d$. Then:

\be
\item[(i)] $\idstr(H)$ is a genetic bar on $V_{\idstrs(d)}$.
\item[(ii)] for all $c'\inn\Vungl$ with $\idstr(c')\iz c$ there is a genetic bar $G'$ on $\Vungl_{c'}$ such that $\idstr(G')\iz G$.
\ee
\elem

\crl
There is a direct correspondence between genetic bars on $V$ and genetic bars on \Vungl.
\ecrl

\prf
Ad (i): by genetic induction:
\be
\item[\gindz] If $H\iz\{\maxdot_d\}\iz\{d\}$, then $\idstr(H)\iz\{\idstr(d)\}\iz\{\maxdot_{\idstrs(d)}\}$, and we are done.
\item[\gindom] Let $H\iz\bigcupp_{b\in\suczs(d)}\Bsubb$ where for each $b\inn\suczz(d)$ we know that $\idstr(\Bsubb)$ is a genetic bar on $V_{\idstrs(b)}$. Notice that $\idstr(\suczz(d))\iz\suczz(\idstr(d))$ in a trivial bijective correspondence. Put $b'\iz\idstr(b)$ for $b\inn\suczz(d)$ and $B'_{b'}\iz\idstr(\Bsubb)$, then $\idstr(H)\iz\bigcupp_{b'\in\suczs(\idstrs(d))}B'_{b'}$. This shows $\idstr(H)$ is a genetic bar on $V_{\idstrs(d)}$, and we are done.
\ee

Ad (ii): also by genetic induction. Take any $c'\inn\Vungl$ such that $\idstr(c')\iz c$. We again use the  bijective correspondence between $\suczz(c')\subsett\Vungl$ and $\suczz(c)\subsett V$, but now for $b\inn\suczz(c)\subsett V$ we let $b'\inn\suczz(c')\subsett\Vungl$ be such that $\idstr(b')\iz b$.

\be
\item[\gindz] If $G\iz\{\maxdot_c\}$, then take $G'\iz\{\maxdot_{c'}\}$.
\item[\gindom] Let $G\iz\bigcupp_{b\in\suczs(c)}\Bsubb$ where for each $b\inn\suczz(c)$ we have a genetic bar $B'_{b'}$ on $\Vungl_{b'}$ such that $\idstr(B'_{b'})$ is $\Bsubb$. Now take $G'\iz\bigcupp_{b'\in\suczs(c')}B'_{b'}$, and we are done.
\ee
\eprf

\prp (from \ref{defindmorf})
Let $g$ be a \pthh-morphism from \Vnatt\ to \Wnattot. Then $g$ is an inductive \pthh-morphism iff $g$ is inductive as a \leqc-morphism from \Vunglnatt\ to \Wnattot.
\eprp

\prf
With the above lemma, the proof now is trivial.
\eprf

\sbsc{Proof of lemma \ref{indmorfresp}}\label{prfindmorfresp}
The proof of lemma \ref{indmorfresp}\ is quite involved. We even need an extra proposition, which has intrinsic value since it can be put to use often in genetic induction proofs:
 
\prp
Let \Vnatt\ be a spraid derived from the pre-natural space \Vprenatt. Let $a\inn V$ and suppose $G$ is a genetic bar on $\Vsuba$ where for all $e\inn G$ we have a genetic bar $D_e$ on $\Vsube$. Then $D\iz\bigcupp_{e\in G}D_e$ is a genetic bar on $\Vsuba$.
\eprp

\prf
By genetic induction on $G$:
\be
\item[\gindz] If $G\iz\{\maxdot_a\}\iz\{a\}$ then we are trivially done.
\item[\gindom] Let $G\iz\bigcupp_{b\in\suczs(a)}\Bsubb$ where for all $b\inn\suczz(a)$: if for all $e\inn \Bsubb$ there is a genetic bar $D_e$ on $\Vsube$, then $D^b\iz \bigcupp_{e\in \Bsubb}D_e$ is a genetic bar on $\Vsubb$. However, we already know for all $b\inn \suczz(a)$ that for all $e\inn \Bsubb$ there is a genetic bar $D_e$ on $\Vsube$ and so $D^b$ is a genetic bar on $\Vsubb$. Therefore $D\iz\bigcupp_{b\in\suczs(a)}D^b$ is a genetic bar on $\Vsuba$.
\ee
\eprf

\lem (repeated from \ref{indmorfresp})
Let $f$ be an inductive morphism between the two spraids \Vnatot\ and \Wnattot, derived from \Vprenatot\ and \Wprenattot\ respectively. Let $a\inn W$, and let $G$ be a genetic bar on $\Wsuba$. Then for all $d\inn V$: if $f(d)\inn \Wsuba$, then $\finv(G)$ contains a genetic bar on $\Vsubd$.
\elem

\prf
By genetic induction on the bar $G$ on $\Wsuba$. 
\be
\item[\gindz] $G\iz\{\maxdot_a\}\iz\{a\}$. Let $d\inn V$ with $f(d)\inn \Wsuba$. Then $d\inn\finv(a)\iz\finv(G)$ so $\{\maxdot_d\}$ is a genetic bar on $\Vsubd$ contained in $\finv(G)$.
\item[\gindom] $G\iz\bigcupp_{b\in\suczs(a)}\Bsubb$ and for all $b\inn\suczz(a)$ we know that for all $d\inn V$: if $f(d)\inn \Wsubb$, then $\finv(\Bsubb)$ contains a genetic bar on $\Vsubd$ (the induction property).

Now let $d\inn V$ with $f(d)\inn \Wsuba$. If $f(d)\precc a$, then there is $b\inn\suczz(a)$ with $f(d)\inn \Wsubb$ which by the induction premise implies that there is a genetic bar on $\Vsubd$ which is contained in $\finv(\Bsubb)\subseteqq\finv(G)$ and we are done. 

Else, $f(d)\iz a$. Then we consider the genetic bar $H\iz\suczz(a)$ on $\Wsuba$. By proposition \ref{genbarsext}\ the expansion $H^{\downarrow\maxdots}$ of $H$ to $W$ is a genetic bar on $W$, so $\finv(H^{\downarrow\maxdots})$ contains a genetic bar $E$ on $V$. Again by proposition \ref{genbarsext}, the reduction $E^{\uparrow d}$ of $E$ to $\Vsubd$ contains a genetic bar $D$ on $\Vsubd$.

By construction of $E$, we see that for $e\inn D$ there is a $b\inn\suczz(a)$ such that $f(e)\inn \Wsubb$. Therefore by the induction assumption, $\finv(\Bsubb)$ contains a genetic bar $K_e$ on $\Vsube$. Now by the proposition above, $K\iz\bigcupp_{e\in D}K_e$ is a genetic bar on $\Vsubd$ contained in $\finv(G)$.
\ee
\eprf 

\sbsc{Proof of proposition and lemma \ref{biscontind}}\label{prfbiscontind}

In formal topology (see \PalCon) a \bisco\ function from \Rt\ to \Rt\ is representable by a formal mapping from the formal reals to the formal reals, and vice versa each such mapping represents a \bisco\ function. We repeat this insight in our setting (recall for the proof that $^s\sigRungl\iz\{a\inn\sigRungl\midd \grdd(a)\iz s\}$):

\prp
Let \ft\ be a \bisco\ function from \Rt\ to \Rt. Then there is an inductive morphism \fstart\ from \sigRt\ to \sigRt\ such that for all $x\inn\sigR$ we have $f(x)\equivv\fstar(x)$ (where we identify \Rt\ and \sigRt\ for convenience). Conversely, if \gt\ is an inductive morphism from \sigRt\ to \sigRt, then as a function \gt\ is uniformly continuous on each compact subspace of \Rt.
\eprp

\prf 
We first create \fstart\ as a trail morphism, since then we need to keep only little track of the intersection properties of intervals. Therefore we use the unglueing \sigRunglt\ of \sigRt, which was defined in \ref{Baireuni}.

\parr

We turn to the uniform-continuity properties of \ft\ on compact subsets of \Rt. We divide \Rt\ in the pairwise overlapping compact spaces $([m, m\pluz 2])_{m\in\Z}$, and by the uniform-continuity properties of \ft\ (using \aczz) we can determine a sequence $(\delta(m,n))_{m\in\Z, n\in\N}$ of elements of \Nt\ (with $\delta(m,n)\smlr\delta(m, n\pluz 1)$) such that if $a\iz [\frac{k}{2^{\delta(m,n)}}, \frac{k+2}{2^{\delta(m,n)}}]$ is a subinterval of $[m, m\pluz 2]$, then $f(a)\subsett [\frac{s}{2^{n}}, \frac{s+2}{2^{n}}]$ for some $s\inn\Z$.
\parr
Again using \aczz, we can turn this into a function \fslangt\ on all pairs $(m, a)$ where $a\iz [\frac{k}{2^{\delta(m,n)}}, \frac{k+2}{2^{\delta(m,n)}}]$ is a subinterval of $[m, m+2]$. We hereby define \fslangt\ in such a way that $\fslang((m,a))\iz [\frac{s}{2^{n}}, \frac{s+2}{2^{n}}]$ for some $s\inn\Z$ and $f(a)\subsett \fslang((m,a))$ for all such pairs $(m,a)$.

\parr

This will yield the definition of the desired morphism \fstart\ on basic dots $a$, by induction on $t\iz\grdd(a)$. For $t\iz 0$ put $\fstar(\maxdot^*)\isdef\maxdotR$. Next, suppose \fstart\ has been defined on $^s\sigRungl$ for $s\smlr t\bygr 0$, and let $a\iz a_0, \ldots, a_{t-1}\inn ^t\sigRungl$. Determine $m\inn\Z$ with $a_0\iz [m, m+2]$ and the smallest \ninnt\ such that $t\iz\grdd(a)\smlr\delta(m,n)\pluz 1$. 

\parr

If $n\iz 0$, this means that for any $c\geqcstar a_{t-1}$, in the process of defining \fslangt\ we did not assign to $(m,c)$ an interval $b$ such that $f(c)\subsett b$, because even $a_{t-1}$ is still too large. Therefore we put $\fstar(a)\isdef\maxdotR$.

If on the other hand $n\bygr 0$, then $t\minuz 1\geqq\delta(m, n\minuz 1)$. Putting $s\iz\delta(m, n\minuz 1)$, we define: $\fstar(a)\isdef\fslang(m, a_s)\capp \fstar(a_0, \ldots, a_{s-1})$, interpreted as intervals.

\parr

We leave it to the reader to verify that \fstart\ is a \leqc-morphism from \sigRunglt\ to \sigRt\ (and therefore a \pthh-morphism from \sigRt\ to itself) which represents \ft.

\parr

To finish the proof of the first part of the theorem, we need to show that \fstart\ is an inductive morphism. For this, let $G$ be a genetic bar on \sigRt. We need to show that $\finvstar(G)$ contains a genetic bar on \sigRungl. If $G\iz\{\maxdotR\}$ then we are trivially done. Else consider $\finvstar(G)$ on the subfan $\rho^*_{[m, m+2]}\iz \{a\inn\sigRungl\midd a\leqcstar [m, m+2]\}$ for given \minn. By the uniform continuity of \ft, we find \Ninnt\ such that $f([m, m+2])\subsett [-N+1, N-1]$. Also, $[-N, N]$ determines a subfann $\tau\iz\tau_{[-N, N]}$ of \sigRt, and we know that $G_{\tau}\iz G\capp\tau$ is finite by theorem \ref{HBforfanns}.  

\parr

And so we find a least $M\inn\N$ such that $G_{\tau}$ descends from $\{a\inn\tau\midd\grdd(a)\iz M\pluz 1\}$, in other words such that the intervals of the type $[\frac{k}{2^M}, \frac{k+2}{2^M}]$ form a refinement of $G_{\tau}$. Now on $\rho^*_{[m, m+2]}$ we see that for each element $b$ of the genetic bar $H_{G, [m, m+2]}= \{a\inn\rho^*_{[m, m+2]}\mid\grdd(a)\iz\delta(m,M)\pluz 1\}$, there is a $c\inn G_{\tau}$ such that $\fstar(b)\leqc c$. We conclude: $H_{G, [m, m+2]}$ is contained in $\finvstar(G)$. Thus we canonically obtain a sequence of genetic bars $(H_{G, [m, m+2]})_{m\in\Z}$ (each on the respective $\rho^*_{[m, m+2]}$) such that $H_{G, [m, m+2]}$ is contained in $\finvstar(G)$ for each $m\inn\Z$. 

\parr

It suffices to consider that by definition of genetic bars, $H\iz\bigcupp_{m\in\Z}H_{G, [m, m+2]}$ is a genetic bar on \sigRunglt\ (trivially $H$ is contained in $\finvstar(G)$). This finishes the proof of the first half of the theorem.

\parr

For the second half of the theorem, let \gt\ be an inductive morphism from \sigRt\ to \sigRt. Copying from the reasoning and constructions above, we look at the genetic bar $G_n = \{a\inn\sigR\midd\grdd(a)\iz n\pluz 1\}$, in other words the intervals of the form $[\frac{k}{2^n}, \frac{k+2}{2^n}]$. We see that $\ginv(G_n)$ contains a genetic bar $H_n$, which is finite on each subfan of \sigR. This shows that for a subfan \taut\ and any \ninnt, there is an \minnt\ such that for all $x,y\inn\tau$ with $\dR(x,y)\leqq 2^{-m}$ we have $\dR(g(x), g(y))\leqq 2^{-n}$. A \bish-compact subspace $X$ is always contained in a subfan \taut\ of \sigR, so we see that \gt\ is indeed uniformly continuous on any \bish-compact subspace of \Rt.
\eprf

\rem
Proposition \ref{refvstrail2}\ (and its proof) shows how to represent \fstart\ as a \leqc-morphism\ from \sigRt\ to \sigRt.
\erem

We also need to prove our remark that the situation is very different when we replace the image space with \Rplust.

\lem
The statement that every uniformly continuous function from \zort\ to \Rplust\ is representable by an inductive morphism from \sigzort\ to \sigRplust\ is equivalent to the fan theorem \FT.
\elem

\prf
The proof uses the equivalence of \FT\ to the statement that each uniformly continuous function from \zort\ to \Rplust\ is bounded away from $0$, which was already proved in \JulRic.

\parr

First let $f$ be a uniformly continuous function from \zort\ to \Rplust\ which is representable by an inductive morphism \fstart\ from \sigzort\ to \sigRplus. Notice that $G\iz\suczz(\maxdotR)$ is a genetic bar on \sigRplus, therefore $\finvstar(G)$ contains a (finite) genetic bar $H$ on \sigzort. Looking at the finite set $\fstar(H)\leqc G$ we conclude that there is \ninnt\ such that $f(x)\bygr 2^{-n}$ for all $x\inn\zor$, in other words $f$ is bounded away from $0$.

\parr

Conversely, suppose $f$ is bounded away from $0$, in other words there is \ninnt\ such that $f(x)\bygr 2^{-n}$ for all $x\inn\zor$. By the uniform continuity of $f$, we can now proceed almost exactly as in the proof of the proposition above to construct an inductive morphism \fstart\ from \sigzort\ to \sigRplust\ representing $f$. We leave this to the reader.

\parr

To finish the proof, we now invoke the well-known result that \FT\ is equivalent to the statement that each uniformly continuous function from \zort\ to \Rplust\ is bounded away from $0$ (see \JulRic, or for another elegant proof, see proposition 4.2 in \WaaArt). 
\eprf

\clearpage

\sbsc{Proof of proposition \ref{indbasdef}}\label{prftopsucz}

For this we first need a lemma on genetic and inductive bars:

\lem
Let \Vnatt\ be a spraid derived from the pre-natural space \Vprenatt, and let $D_0, D_1$ be (i) genetic (ii) inductive bars on $V$. Then:

\be

\item[(i)] $\min(D_0, D_1)\iz\{d\inn V\midd \driss c\inn V\driss i\inn\{0,1\}\,[d\inn D_i\weddge c\inn D_{1-i}\weddge d\leqc c]\}$ is a genetic bar on $V$.

\item[(ii)] $(D_0)_{\leqc}\capp (D_1)_{\leqc}$ is an inductive bar on $V$.

\ee

\elem

\prf
For genetic bars $D_0, D_1$ we prove (i) by genetic induction on $D_0$:
\be
\item[\gindz] If $D_0\iz \{\maxdot\}$, then $\min(D_0, D_1)\iz D_1$ and we are done.
\item[\gindom] Else $D_0\iz\bigcupp_{b\inn\suczs(\raisebox{.2ex}{$\maxdots$})}\Bsubb$ where for all $b\inn\suczz(\maxdot)$ we have that $\min(\Bsubb, D_1)$ is a genetic bar on $\Vsubb$. Then $\min(D_0, D_1)\iz\bigcupp_{b\inn\suczs(\raisebox{.2ex}{$\maxdots$})}\min(\Bsubb, D_1)$ is also a genetic bar.
\ee

Ad(ii): for inductive bars $D_0, D_1$ descending from genetic bars $G_0, G_1$ respectively, it suffices to see that $(D_0)_{\leqc}\capp (D_1)_{\leqc}$ contains $\min(G_0, G_1)$, which is a genetic bar by (i). For let $d\inn\min(G_0, G_1)$. Without loss of generality $d\inn G_1$ and there is $d\leqc c\inn G_0$. Then we can determine $d\leqc b\inn D_1$ and $c\leqc a\inn D_0$, whence $d\leqc c\leqc a$ and so $d\inn (D_0)_{\leqc}\capp (D_1)_{\leqc}$. 
\eprf

\crl
By induction on $n$ we can now conclude: if $B_0,\ldots, B_n$ is a finite sequence of (i) genetic (ii) monotone inductive bars, then 
\be
\item[(i)] $\min(B_0,\ldots, B_n)\iz \{d\inn V\midd \drisz i\leqq n\,[d\inn B_i\weddge \allz j\leqq n, j\notiz i\,\drisz c\inn B_{j} [d\leqc c]\}$
is a genetic bar on $V$.
\item[(ii)] $\bigcapp_{i\leq n}B_i$ is a monotone inductive bar.
\ee
\ecrl

\rem
An interesting exercise for the reader is to see why for inductive bars $D_0, D_1$ descending from genetic bars $G_0, G_1$ respectively, the bar $\min(D_0, D_1)$ does not necessarily descend from $\min(G_0, G_1)$. For a spread \Vnatt\ however $\min(D_0, D_1)$ does descend from $\min(G_0, G_1)$.
\erem

\prp (repeated from \ref{indbasdef})
\be
\item[({\small\textscc{i}\normalsize})] For a spraid \Vnatt\ with corresponding pre-natural space \Vprenatt, the collection \Topsuczt\ is a topology which is refined by \Topaprt. 
\item[({\small\textscc{ii}\normalsize})] Let \Vnatt\ be an inductive spraid with corresponding pre-natural space \Vprenatt. 
Then for finite subsets $A\aprt B$ of $V$, the subset $C\iz\{c\inn V\midd c\aprt A\vee c\aprt B\}$ is an inductive bar on \Vnat.
\ee
\eprp

\prf
Ad ({\small\textscc{i}\normalsize}):
That \Topaprtt\ refines \Topsuczt\ is trivial. To show that \Topsuczt\ is a topology, we check the definitions:
\bei
\itemmz{\Topo} Trivially the empty set \emptyyt\ and \Vcalt\ are in \Topsuczt\ (for any $x\inn\Vcal$ the set $B_{\Vcall}^{x}\iz V$ is an inductive bar). 
\itemmz{\Topto} Let $\Ucal, \Wcal\inn\Topsucz$, we must show that $\Ucal\capp\Wcal\inn\Topsucz$.
For this, let $x\inn\Ucal\capp\Wcal$, and consider that $B_{\Ucall\cap\Wcall}^{x}\iz B_{\Ucall}^{x}\capp B_{\Wcall}^{x}$ since for all $b\inn B_{\Ucall}^{x}\capp B_{\Wcall}^{x}$ we can decide: \kase{\Ucalt(1)} $b\aprt x_{\grds(b)}$ (done) or \kase{\Ucalt(2)} $\hattr{b}\subseteqq\Ucal$, in which case we consider \kase{\Wcalt(1)} $b\aprt x_{\grds(b)}$ (done) or \kase{\Wcalt(2)} $\hattr{b}\subseteqq\Wcal$ whence we see that $\hattr{b}\subseteqq\Ucal\capp\Wcal$ and we are done too. The bars $B_{\Ucall}^{x}$ and $B_{\Wcall}^{x}$ are monotone and inductive, and so $B_{\Ucall\cap\Wcall}^{x}\iz B_{\Ucall}^{x}\capp B_{\Wcall}^{x}$
is a monotone inductive bar also by (ii) of the above corollary.
\itemmz{\Topth}
Let $\Ucal\subseteqq\Vcal$ be such that for all $x\inn\Ucal$ there is a $\Wcal\nni x$ in \Topsuczt\ such that $\Wcal\subseteqq\Ucal$. Let $x\inn\Ucal$, determine $\Wcal\nni x$ in \Topsuczt\ such that $\Wcal\subseteqq\Ucal$. Then clearly $B_{\Wcall}^{x}\subseteqq B_{\Ucall}^{x}$ therefore $B_{\Ucall}^{x}$ is an inductive bar, which since $x$ is arbitrary shows that \Ucalt\ is in \Topsuczt.
\eei

Ad ({\small\textscc{ii}\normalsize}): First suppose either $A$ or $B$ is empty, then the conclusion is trivially fulfilled. Else, let $A\iz\{a_0, \ldots, a_n\}$ and  $B\iz\{b_0, \ldots, b_m\}$ with $a_i\aprt b_j$ for all $i\leqq n, j\leqq m$. By definition of `inductive spraid', for any $i\leqq n, j\leqq m$ the set $D_{i,j}\iz\{d\inn V\midd d\aprt a_i \vee d\aprt b_j\}$ is an inductive bar (which is also clearly monotone). Therefore by (ii) of the above corollary, $D\iz\bigcapp_{i\leq n, j\leq m}D_{i,j}$ is a monotone inductive bar. However, 
$D$ actually equals $C\iz\{c\inn V\midd c\aprt A\vee c\aprt B\}$, therefore $C$ is an inductive bar as promised.
\eprf

\sbsc{Proof of meta-theorem \ref{metaruss}}\label{prfmetaruss}

In paragraph \ref{metaruss}\ we indicated pointwise problems regarding inductive morphisms in \bish:

\mthm (repeated from \ref{metaruss})

In \russ\ (and by implication \bish) we have the following problems regarding pointwise use of inductive definitions:
\be
\item[P$_1$] Uniform continuity of a function $f$ does not imply that there is an inductive morphism representing $f$. Counterexamples can be given even for uniformly continuous functions from \zort\ to \Rplus. In other words: uniform continuity does not imply inductive representability.
\item[P$_2$] Weak completeness\footnote{The property for a located subset $A$ of a metric space \xdt\ that for all $x\inn X$: if $x\aprt a$ for all $a \inn A$, then $d(x, A)\iz\inf(\{d(x,a)\midd a\inn A\})\bygr 0$.}\ of a compact space is not preserved under inductivity. In \russ, even for an inductive morphism from \zort\ to \zor, the image of a compact subspace may be strongly incomplete.
\item[P$_3$] Inductive representability is not preserved under the restriction of a function to its pointwise image space. This follows from the counterexamples for P$_1$, since every uniformly continuous function from \zort\ to \Rt\ is inductively representable by proposition \ref{biscontind}. Therefore we can expect problems with the reciprocal function $x\rightarrow\frac{1}{x}$, and must continually address these problems by adapting our definitions.
\ee

Therefore in \bish, the desirable properties associated with the problems above cannot be shown to hold without further assumptions. In fact, assertion of any of these properties implies the fan theorem \FT. 
\emthm

\prf
The proof is derived from the construction of ContraCantor space, which is a compact subspace \contrCant\ of \zort\ such that if we write \Canzort\ for the standard embedding of Cantor space in \zort, we see: $\dR(\contrCan, \Canzor)\iz 0$ and yet in \russ\ we also have $\dR(x, \Canzor)\bygr 0$ for $\all x\inn\contrCan$.  ContraCantor space is defined in the examples' section of the appendix \ref{contraCan}\ using the Kleene Tree, and the properties above are proved in proposition \ref{contraCan}.

\parr

To show P$_1$, consider the uniformly continuous function $d_{\mathcal{C}}$ from \zort\ to $[0, \frac{1}{6}]$ given by $d_{\mathcal{C}}(x)\iz \dR(x, \Canzor)$.\footnote{The distance of a point to a compact subspace can always be constructively calculated.}\ The restriction of $d_{\mathcal{C}}$ to \contrCant\ is an example in \russ\ of a uniformly continuous function (from \contrCant\ to \Rplus) which cannot be represented as an inductive morphism from the subfann \contrCant\ to \Rplusnatt. Notice that $d_{\mathcal{C}}$ is a beautiful function, and that \contrCant\ is a beautiful compact space, so this type of problem can be expected to crop up at any time.

\parr

The above example also shows what we mean with P$_3$, since $d_{\mathcal{C}}$ seen as a uniformly continuous funtion from \contrCant\ to \Rt\ is representable as an inductive morphism from the subfann \contrCant\ to \Rnatt\ (see prp.\,\ref{biscontind}\ and \ref{prfbiscontind}).

\parr

To show P$_2$, we consider the subset $d_{\mathcal{C}}(\contrCan)$ of $[0, \frac{1}{6}]$, which is obviously strongly incomplete, since $0\aprt d_{\mathcal{C}}(\contrCan)$ and yet $\dR(\contrCan, \Canzor)\iz 0$.

\parr

Another example in \russ\ of a strongly incomplete fann-image under an inductive morphism was already given in \WaaArt. This example consists of $\zorbin^{\N}$ and a recursive \betat\ in \zorNt\ such that $\dRN(\beta,\zorbin^{\N})\iz 0)$ and yet $\beta\aprtRN x$ for all recursive $x\inn\zorbin^{\N}$.
\eprf

\clearpage

\sbsc{Proof of theorem \ref{compmetind}}\label{prfcompmetind}

We prove that every complete metric space is homeomorphic to a \sucz-spread, by adapting our proof in \ref{sepmetnat}\ that every complete metric space is homeomorphic to a natural space.

\thm (from \ref{compmetind})
Every complete metric space \xdt\ is homeomorphic to a \sucz-spraid.
\ethm

\prf
For a complete metric space \xdt\ with dense subset \anninnt, we constructed a homeomorphic \Vprenatt\ where $V=\{B(a_n, 2^{-s})\midd n,\sinn\}$. If we look more carefully, we see that the trail space of this \Vnatt\ contains a (homeomorphic) \sucz-spread. 

\parr
We copy the notation in \ref{sepmetnatprf}. 
We defined \leqct\ and \aprt\ on pairs of basic dots in $V$ such that $B(a_n, 2^{-s})\precc B(a_m, 2^{-t})$ implies $s\bygr t$ and $d(a_n, a_m)\smlr 2^{-t}\minuz 2^{-s}$, and secondly $B(a_n, 2^{-s})\aprt B(a_m, 2^{-t})$ implies $d(a_n, a_m)\bygr 2^{-s}\pluz 2^{-t}\pluz 2^{-s-t-1}$.

\parr

We now define new basic dots, similar to forming the trail space (see \ref{pathmorphisms}).
Firstly, for a basic dot $a\iz B(a_n, 2^{-s})\inn V$ we put $\diam(a)\isdef s$.

Let $p$ be a point in \Vcalt, then remember we write $\pstr(n)$ for the finite sequence $p_0, \ldots, p_{n-1}$ of basic dots in $V$. Notice that by definition $p_0\geqc\ldots\geqc p_{n-1}$. A finite sequence $a\iz a_0\geqc\ldots\geqc a_{n-1}$ of basic dots in $V$ is called a \deff{graded trail}\ in \Vleqct\ (of length $n$) iff $\diam(a_i)\iz i$ for $0\leqq i\smlr n$. The empty sequence is the unique graded trail of length $0$, and denoted \maxdotpt. The countable set of  graded trails in \Vleqct\ is denoted \Vdiam, notice that $\Vdiam\subsett\Vpath\iz\{\pstr(n)\midd\ninn, p\inn\Vcal\}$. 

\parr
The pre-natural space $(\Vdiam, \leqcstar, \aprtstar)$ (see def.\,\ref{pathmorphisms}) induces a spread, which we call \Vdiamcalt. We show that \Vdiamcalt\ is an inductive spread by two claims:

\clmm
For $c\aprtstar d\inn\Vdiam$ there is a genetic bar $B$ on $\Vdiam$ such that for all $b\inn B$ we have: ($b\aprtstar c \vee b\aprtstar d)$. 
\eclmm

\clmmprf 
Let $c\iz c_0, \ldots, c_s$ and $d\iz d_0, \ldots, d_t$, with $c_s\iz B(a_n, 2^{-s})$, $d_t\iz B(a_m, 2^{-t})$ for certain $n, m, s, t\inn\N$. By definition $c\aprtstar d$ means 
$B(a_n, 2^{-s})\aprt B(a_m, 2^{-t})$.

\parr
By the properties of $\aprt$ for \Vnatt\ (see above), if $B(a_n, 2^{-s})\aprt B(a_m, 2^{-t})$, then $d(a_n, a_m)\bygr 2^{-s}\pluz 2^{-t}\pluz 2^{-s-t-1}$. This
means that for any basic dot $e$ in $V$ with $\diam(e)\geqq s\pluz t\pluz 3$ we can decide: $e\aprt B(a_n, 2^{-s})$ or $e\aprt B(a_m, 2^{-t})$, by looking at $e\iz B(a_k, 2^{-\diam(e)})$ to see that either $d(a_k, a_n)$ or $d(a_k, a_m)$ is big enough.  We conclude that the genetic bar $B\iz \{b\inn\Vdiam\midd \grdd(b)\iz s\pluz t\pluz 3\}$ satisfies the requirements.
\eclmmprf 

\bclmm
Let \Ucalt\ be open in \Vdiamcalt, and $x\inn \Ucal$. 
Then $B_{\Ucall}^{x}\iz\{b\inn \Vdiam\midd b\aprtstar x_{\grds(b)}\vee\hattr{b}\subseteq\Ucal\}$ is an inductive bar on $\Vdiam$. 
\ebclmm

\clmmprf 
(Remember the definition of $\idstr$ in \ref{pathmorphisms}.) Looking at $x$, we see that there are $a\inn V$ and $t, n, i\inn\N, i\geqq 1$ such that $\hattr{x_{t}}\subseteqq\Ucal$ and $\idstr(x_t)\iz B(a_n, 2^{-i+1})$. Determine $j, m\inn\N$ such that $\idstr(x_j)\leqc B(a_m, 2^{-i})\leqc B(a_n, 2^{-i+1})$. 
By definition of \leqct\ (see \ref{sepmetnatprf}), we have:
$d(a_n, a_m)\smlr 2^{-i+1}\minuz 2^{-i}\iz 2^{-i}$. Determine \sinnt\ such that $d(a_n, a_m)\smlr 2^{-i}\minuz 2^{-s+2}$. Determine $u\inn\N$ such that $\diam(\idstr(x_u))\geqq s$.

\parr

Now let $e\inn V$ with $\diam(e)\iz s$, say $e\iz B(a_k, 2^{-s})$. By the properties above, we can determine: $e\aprt \idstr(x_u)$ or $e\iz B(a_k, 2^{-s})\subseteqq B(a_n, 2^{-i+1})\iz \idstr(x_t)$. But then in turn for any $b\inn \Vdiam, \grdd(b)\geqq \max(s, u)$ we can determine: $b\aprtstar x_{\grds(b)}$ or $\hattr{b}\subseteq\Ucal$. So $B_{\Ucall}^{x}$ contains the genetic bar $\{b\inn\Vdiam\midd\grdd(b)\iz \max(s,u)\}$, and is therefore inductive.
\eclmmprf 

By the two above claims, we see that \Vdiamcalt\ is a \sucz-spread. We leave it to the reader to verify that $(\Vdiamcal, {\Topaprt}^{\raisebox{.25ex}{$\scriptstyle\hspace*{-.4em}\wr$}})$ is homeomorphic to \xdt. 
\eprf

\rem
The question which representation to choose for complete metric spaces can probably not be answered in just one way. We believe that detailed study of this question, for various spaces, yields both theoretical and practical advantages (from the \appl\ perspective).
\erem

\sbsc{Proof of proposition \ref{projtych} (iv)}\label{prfinttych}
We prove that an (in)finite product of star-finite spreads (see def.\,\ref{starfinsprd}) is faithful in \class\ and \intu. For the proof we will use some theory developed in the proof of theorem \ref{starfinmet}, in paragraph \ref{prfstarfinmet}.

\prp (\class, \intu, from \ref{projtych})
Let $((\Vcal_n, \Topaprtn))_{n\in\N}$ be star-finite spraids derived from the corresponding pre-natural $((V_n, \aprtn, \leqcn))_{n\in\N}$. Then the finite spraid-products \PisigVnt\ (for \ninnt) and the infinite spraid-product \PisigVNt\ are faithful.
\eprp

\prf We only prove the infinite-product case, the finite-product case is completely similar and easier. Let \Ucalt\ be open in \PisigVNt, determine $x\inn\Ucal$. 
Since \PisigVNt\ is star-finite, following paragraph \ref{prfstarfinmet} (def.\,(\bol{d})) we can define the subfann \Wxt\ of \Vpisigt. By \BTg\ (see \ref{brother}), \Ucalt\ is inductively open, therefore the bar $\BUx\iz\{b\inn\Vpisig\midd b\aprt x_{\grds(b)}\vee\hattr{b}\subseteqq\Ucal\}$ is inductive. By \HBind\ (prp.\,\ref{HBforsubfanns}) $\BUx$ contains a finite subbar $D$ on \Wxt, therefore we can find \Ninnt\ with $^{N}\Wx\subseteqq D_{\leqc}$. 
\parr
Determine $a\iz a_0,\ldots, a_N\inn\Vpisig$ such that $x\precpi a$. Then $a\inn ^{N}\Wx$ and $\hattr{a}\subseteqq\Ucal$. For $b\iz b_0,\ldots, b_N\inn\Vpisig$, by the nature of $D$, if $b\touchpi a$, then $\hattr{b}\subseteqq\Ucal$. This implies that $\bigcapp_{i\leq N}{\pi_i}^{-1}(\hattr{C_i})\subseteqq\Ucal$, where $C_i\iz\{b\inn ^{N}V_i\midd b\touchi a_i\}$, for $i\leqq N$. However, $\hattr{C_i}$ is a neighborhood of $x_{[i]}$ in $(\Vcal_i, \Topaprti)$ for each $i\leqq N$, which shows that \Ucalt\ is open in the Tychonoff topology. 
\eprf

\sbsc{Proof of lemma \ref{sprdtych}}\label{prfsprdtych}
We prove the technical lemma which is needed for theorem \ref{sprdtych} (entailing a \bish\ version  of Tychonoff's theorem). 

\lem
(Notations as in \ref{defprod}) Let $a\inn V_0, b\inn V_1$ and let $G, H$ be genetic bars on $(V_0)_a, (V_1)_b$ respectively. Then $G\Timsig H$ is an inductive bar on $(V_0)_a\Timsig (V_1)_b$.
\elem

\crl
For $i\leqq n$ let $B_i$ be an inductive bar on $V_i$. Then $\Pisig{n}B_i$ is an inductive bar on $\Vpisign$ and $\overline{\Pisig{n}B_i}$ is an inductive bar on $\Vpisig$.
\ecrl

\prf
By (double) genetic induction. For notational simplicity put $V_0\iz V$ and $V_1\iz W$. First let $P_W$ be the following property of genetic bars $B$ on basic subspraids $W_b$ of $W$: `for all $a\inn V$: $\{\maxdot_a\}\Timsig B$ is an inductive bar on $V_a\Timsig W_b$'. We show by genetic induction that all genetic bars $B$ on basic subspraids $W_b$ have property $P_W$. Let $a\inn V, b\inn W$.

\be
\item[\gindz] If $B\iz\{\maxdot_b\}$, then $\{\maxdot_a\}\Timsig B\iz\{(\maxdot_a,\maxdot_b)\}$ and we are done trivially.
\item[\gindom] Else $B\iz\bigcupp_{b'\inn\suczs(b)}B_{b'}$ where $B_{b'}$ has property $P_W$ for all $b'\inn\sucz(b)$. Then $\{a'\}\Timsig B_{b'}$ is an inductive bar on $V_{a'}\Timsig W_{b'}$ for all $a'\inn\sucz(a)$, which implies that $D\iz\bigcupp_{a'\in\suczs(a), b'\in\suczs(b)}\{a'\}\Timsig B_{b'}$ 
is an inductive bar on $V_a\Timsig W_b$. Clearly $\{\maxdot_a\}\Timsig B$ descends from $D$ and so is inductive also.
\ee

By symmetry, for all genetic bars $B$ on basic subspraids $V_a$ of $V$ we also obtain $P_V$: `for all $b\inn W$: $B\Timsig\{\maxdot_b\}$ is an inductive bar on $V_a\Timsig W_b$'.

Next we show that all genetic bars $B$ on basic subspraids $V_a$ of $V$ have property $Q$: `for all $b\inn W$ and every genetic bar $C$ on $W_b$: $B\Timsig C$ is an inductive bar on $V_a\Timsig W_b$'. Let $a\inn V, b\inn W$.

\be
\item[\gindz] If $B\iz\{\maxdot_a\}$, then we are done since every genetic bar on $W_b$ has property $P_W$.
\item[\gindom] Else $B\iz\bigcupp_{a'\inn\suczs(a)}B_{a'}$ where $B_{a'}$ has property $Q$ for all $a'\inn\sucz(a)$. Let $C$ be a genetic bar on $W_b$, we proceed by genetic induction on $C$.

\be
\item[\gindz] If $C\iz\{\maxdot_b\}$, then we are done since every genetic bar on $V_a$ has property $P_V$.
\item[\gindom] Else $C\iz\bigcupp_{b'\inn\suczs(b)}C_{b'}$ where by induction $B_{a'}\Timsig C_{b'}$ is an inductive bar on $V_{a'}\Timsig W_{b'}$ for all $a'\inn\sucz(a), b'\inn\suczs(b)$. Therefore we see that $D\iz\bigcupp_{a'\in\suczs(a)\\b'\in\suczs(b)}\{a'\}\Timsig B_{b'}$ is an inductive bar on $V_a\Timsig W_b$. Clearly $B\Timsig C$ descends from $D$ and so is inductive also.
\ee
\ee
This proves the lemma. The first part of the corollary follows by induction. For the second part we show by genetic induction that for every genetic bar $B$ on a basic subspraid $(\Vpisign)_a$ of $\Vpisign$ we have: $\overline{\{a\}}\oicovb\overline{B}$ in \Vpisigt. Let $a\inn \Vpisign$.

\be
\item[\gindz] If $B\iz\{\maxdot_a\}$, then trivially $\overline{\{a\}}\oicovb\overline{B}$.
\item[\gindom] Else $B\iz\bigcupp_{b\inn\suczs(a)}B_{b}$, where by induction $\overline{\{b\}}\oicovb \overline{B_{b}}$ for all $b\inn\sucz(a)$. We have $\overline{\{a\}}\oicovb\bigcupp_{b\inn\suczs(a)}\overline{\{b\}}$ so $\overline{\{a\}}\oicovb\bigcupp_{b\inn\suczs(a)}\overline{B_{b}}\iz\overline{B}$ and we are done.
\ee

Finally, $\overline{\{\maxdotPin\}}$ equals $\overline{^n\Vpisig}$ which is an inductive bar on \Vpisigt\ by lemma \ref{genbarsext}. So any inductive bar on $\overline{\{\maxdotPin\}}$ is an inductive bar on \Vpisigt.
\eprf

\sbsc{Defining various concepts of locatedness}\label{appvarconloc}
\hfill\mbox{(Partly repeating \ref{varconloc}:)}
What are the drawbacks of the concept `located in'?
First of all, the notion is not transitive, which is unpractical when working with extensions and subspaces of \xdt. Second, even for a closed located $A\subsett X$, the notion gives little handhold for $x\inn X$ to find $a\inn A$ such that $x\aprt a$ implies $x\aprt A$, which is an important prerequisite for many constructions involving $A$. Thirdly, as mentioned, the notion is non-topological and this means we cannot use it easily in the context of topology.

\parr

In \WaaThe\ several alternatives 
are given in \bish, of which `strongly halflocated in' (transitive) seems the most fruitful.\footnote{By transitive we mean: if $(B,d)$ is (strongly) halflocated in $(A,d)$ which is (strongly) halflocated in $(X,d)$ then $(B,d)$ is (strongly) halflocated in $(X,d)$.}\ It gives results such as in the \bish-proof of the Dugundji extension theorem in \WaaThe. Another result is that every complete metric space can be isometrically embedded in a normed linear extension such that it becomes strongly halflocated in this extension -- and where we know of no general proof that it is located. 
\parr
`Topologically strongly halflocated' in \intu\ is equivalent on complete metric spaces to a topological locatedness property called `strongly sublocated in'. This notion can also be defined for the apartness topology of general natural spaces, and seems to us important. Our definition of `located in' is easily seen to be equivalent to the traditional definition, and in this form it opens the door for adaptations.

\defi
Let \adt\ be a subspace of \xdt, a metric space. Then \adt\ is (i)\deff{located}, (ii)\deff{halflocated},
(iii)\deff{sublocated} in \xdt\ iff: \adt\ is inhabited and

\be
\item[(i)]$\allz{^{\!}} D\inn\R_{>1} \allz x\inn X \all m\inn\Z \,[\,\drisz a\inn A\,[d(x,a)\smlr D^{m+1}]\orr
          \allz a\inn A\,[d(x,a)\bygr D^m]]$
\item[(ii)]$\drisz D\inn\R_{>1} \allz x\inn X \all m\inn\Z \,[\,\drisz a\inn A\,[d(x,a)\smlr D^{m+1}]\orr
          \allz a\inn A\,[d(x,a)\bygr D^m]]$
\item[(iii)]$\all x\inn X \all m\inn\Z \,[\,\dris a\inn A\,[d(x,a)\smlr 2^{m}]\orr\driss\ninn\
          \all a\inn A\,[d(x,a)\bygr\twominn]]$
\ee
If $D\inn\R_{>1}$ realizes (ii), we say that \adt\ is halflocated in \xdt\
\deff{with parameter} $D$. We now strengthen the definition: 

\parr

\adt\ is (i)$^*$\deff{strongly located}, (ii)$^*$\deff{strongly halflocated}, (iii)$^*$\deff{strongly sublocated} in \xdt\ iff:

\be
\item[(i)$^*$]$\all D\inn\R_{>1} \all x\inn X  \driss y\inn A\
                  \all a\inn A \,[\,d(x,y)\leqq D\cdott d(x,a)\,]$
\item[(ii)$^*$]$\dris D\inn\R_{\geq 1} \all x\inn X  \driss y\inn A\
                  \all a\inn A \,[\,d(x,y)\leqq D\cdott d(x,a)\,]$
\item[(iii)$^*$] $\all x\inn X \driss y\inn A \,[\,x\aprt y\,\rightarrow\,\driss\ninn\
          \all a\inn A\,[d(x,a)\bygr\twominn]]$
\ee

If $D\inn\R_{>1}$ realizes (ii), then \adt\ is strongly halflocated in \xdt\
          \deff{with parameter}\ $D$. If (ii) is realized by $D\iz 1$, then \adt\ is 
          \deff{best approximable}\ in \xdt.

\parr

Now let \Vnatt\ be a natural space derived from \Vprenatt, and let \Wnatt\ be a subspace. We say that \Wnatt\ is (i)$^\star$\deff{\aprt-sublocated}, (ii)$^\star$\deff{strongly \aprt-sublocated} in \Vnatt\ iff: 

\be
\item[(i)$^\star$] $\all x\inn\Vcal\all\Ucal\inn\Topaprt, x\inn\Ucal\,[\,\dris y\inn\Wcal [y\inn\Ucal]\vee \all y\inn\Wcal [ x\aprt y]]$.
\item[(ii)$^\star$] $\all x\inn\Vcal\driss y\inn\Wcal\,[\,x\aprt y\,\rightarrow\,\all z\inn\Wcal\,[x\aprt z]]$.
\ee
\edefi

\rem
The word `strongly' for (i)$^*$-(iv)$^*$ and (ii)$^\star$ is justified.
For $x\inn X$, if \yyt\  realizes (ii)$^*$ with parameter $D$, we can decide:
$d(x,y)\smlr D^{m+1}$ or $d(x,y)\bygr D^{m}$, for \minzt. But $d(x,y)\bygr D^{m}$ implies that 
for all $a\inn A$: $d(x,a)\bygr D^{m-1}$. This shows that if \adt\ is strongly halflocated in
\xdt, with parameter $D$, then \adt\ is halflocated in \xdt\ with parameter $D^2$. Then 
`strongly located in' implies `located in' since $\{D^2\midd D\inn\Rbigro\}\iz\Rbigro$. The other implications can be obtained in a similar but easier fashion. Notice 
that if \adt\ is strongly (half,sub)located in \xdt, then \adt\ is closed in \xdt. 
For a more extensive treatment of these properties, see \WaaThe.
\erem

\sbsc{Proof of proposition \ref{dirlimnonmet}}\label{prfdirlimnonmet}

Our example of the eventually vanishing real sequences still needs to be proven non-metrizable. We repeat from \ref{dirlimnonmet}:

\prp
\Romnatt\ is a natural space which is the non-metrizable direct limit of $(\Rnstar, \TopRlim)_{n\in\N}$,  where $(\Rnstar, \TopRlim)$ is \leqc-isomorphic to the Euclidean space $(\R^n, \TopRaprt)$ for \ninnplt. There is a continuous injective surjection (which does not have a continuous inverse) from \Romnatt\ to $(\Rom, \TopRNaprt)$ as subspace of $(\RN, \TopRNaprt)$.
\eprp

\prf
That \Romnatt\ is non-metrizable is seen thus: consider a metric $d$ on \Romnat, then $d$ is also a metric on the respective $(\R^n, \TopRaprt)$ for \ninnplt\ which respects the canonical inclusion relation $i_n: \R^n\rightarrow\R^{n+1}$. For \ninnplt\ let $0_{\R^n}$ be the origin in $\R^n$ and let $B_d^n(0_{\R^n}, 2^{-n})$  be the open $d$-sphere around this origin with radius $2^{-n}$. We can construct a series of sets $(U_n)_{n\in\Nplus}$ where each $U_n$ is a $d$-open neighborhood in $\R^n$ of $0_{\R^n}$ and $i_n(U_n)\subsett U_{n+1}$, and where in addition $B_d^n(0_{\R^n}, 2^{-n})$ contains a point which is not contained in $U_n$. 

\parr
Remember that for $x\inn\Rom$ the $m$-th basic dot $x_m$ (if not equal to the maximal dot) is a finite sequence of closed rational intervals $([a_j, b_j])_{j\leq s}$ for some \sinnt\ and $a_j\smlr b_j$ for all $j\leqq s$. We identify $x_m$ with the cartesian product set $\Pi(x_m)\iz \Pi_{j\leqq s}[a_j, b_j]$ in the Euclidean space $\R^{s+1}$.  Now it is not difficult to see that the subset $\{y\in\Rom\mid \driss x\equivv y [\,\all\ninn [\Pi(x_n)\subsett U_n]\,]\}$ is open in \TopRlimt, but cannot be open in the metric topology generated by $d$. 
 
\parr

Further, the identity on \Romt\ is a continuous injective surjection (which does not have a continuous inverse) from \Romnatt\ to $(\Rom, \TopRNaprt)$ as subspace of $(\RN, \TopRNaprt)$. One sees this by considering that on \RNt\ the apartness topology \TopRNaprt\ coincides with the metric $\dRN$-topology.

\parr
Finally, that $(\Rnstar, \TopRlim)$ is \leqc-isomorphic to the Euclidean space $(\R^n, \TopRaprt)$ for \ninnplt\ is left to the reader as an exercise.
\eprf

\sbsc{Proof of corollary \ref{starfinmet}}\label{prfcrlstarfinmet}

For didactical reasons we prove the metrizability of \sucz-fans first (the corollary of theorem \ref{starfinmet}), and metrizability of star-finite \sucz-spreads (theorem \ref{starfinmet}\ itself)
 in the next paragraph. The proof of the theorem for star-finite \sucz-spreads employs the same main strategy, but is lengthy and involves some hard work on details, which tends to obscure this strategy.

\parr

The following two paragraphs contain the longest proof in the monograph. In order to prove that every star-finitary space is metrizable, we will need several lemmas, to which we assign uppercase letters \bol{A}, \bol{B}, $\ldots$. In a way we follow the classical strategy used by Urysohn (see \Uryb) to show that (classically) a normal space \xtopt\ with a countable base is metrizable. Urysohn embeds such an \xtopt\ in the Hilbert cube, and this is what we do also for a star-finite \sucz-spread, although our construction of such an embedding is different.\footnote{For one thing, it is a real construction. The key idea can also be seen as an onion strategy,but different from the classical one.}

\parr

For this construction we use the ternary real numbers \Rtert, and specifically \zortert. We refer the reader to the examples' section \ref{cantorfun}\ for the relevant definitions.

\parr

We rearrange the ingredients and the route followed in \WaaThe, partly in order to avoid the use of \acoz.\footnote{We comment on the relation between the intuitionistic proof in \WaaThe\ and the proof here, in the comments' section, paragraph \ref{starfinmetcom}.}\ This calls for a short description of our route beforehand. It is possible to prove metrization for \sucz-fans, using the Urysohn metrization lemma (\bol{A}), a splitting lemma (\bol{B}), and the Urysohn function lemma (\bol{C}) below.\footnote{The term `Urysohn's lemma' is usually reserved for the related classical result that a space is `normal' iff two disjoint closed subsets can be separated by a continuous function to \zort. See our lemma (\bol{C}) and its generalization (\bol{E}) which we call `Urysohn function lemma'.}\ To generalize this to a star-finite \sucz-spread \Vnatt, we associate to points $x$ in \Vcalt\ a subfan \Wxt\ of \Vnatt\ which almost acts as a neighborhood of $x$ in \Vnatt. Using an elaborate adaptation of lemma (\bol{C}) we can again apply the Urysohn metrization lemma to conclude metrizability of \Vnatt. All in all this will take up numerous pages.

\parr

In this and the following paragraph we need some definitions, to which we assign lowercase letters \bol{a}, \bol{b}, etc. 

\defi (\bol{a})

Let \Vnatt\ be a spraid derived from \Vprenatt. Let $A\subseteqq V$, then we write $^nA$ for $\{a\inn A\midd \grdd(a)\iz n\}$, for \ninn, and $A_{\leqc}$ for $\{b\inn V\midd \driss a\inn A [b\leqc a]\,\}$. Also remember that we write \toucht\ for the touch-relation on basic dots which is the complement of the pre-apartness relation \aprt.

\parr

Now let $x\inn\Vcal$. We say that $x$ is a \deff{successor point}\ iff $\grdd(x_n)\iz n$ for all \ninnt\ (which implies $x_0\iz\maxdot$ and $x_{n+1}\sucz x_n$ for \ninnt). One easily sees that any $y\inn\Vcal$ contains a unique subsequence $y^{\!\sucz}\equivv y$ such that $y^{\!\sucz}$ is a successor point.
This means that without loss of generality we can conveniently concentrate on successor points.

\parr

We wish to expand \Vnatt\ with a single isolated point, which w.l.o.g. can be taken to be $\myno=\bullet, \bullet, \ldots$. Strictly speaking we put $\Vex\isdef V\cupp \{\bullet\}^*$ where $\{\bullet\}^*\iz\{\mynon\midd\ninn, n\geqq 1\}$ is the set of all finite sequences of the symbol $\bullet$, and specify that $\maxdot\geqc\mynon\geqc \overline{\myno}(n\pluz 1)$ and $\mynon\aprt a$ for all $\maxdot\notiz a\inn V$, $n\geqq 1$. The resulting spraid is denoted $(\Vexcal,\Topaprt)$ or simply \Vexcalt.
\edefi

\lem (\bol{A}) (Urysohn metrization lemma)

Let \Vnatt\ be a spraid derived from \Vprenatt, where we write \toucht\ for the complement of the pre-apartness relation \aprt. Suppose that:

\be
\item[(i)] For all \ninn, for all $a\aprt b\inn ^n\Vex$ there is a given morphism $f_{a,b}$ from \Vexcalt\ to \zort, such that $f_{a,b}\restrct{\Vexa}\equivvR 0\,$ and $f_{a,b}\restrct{\Vexb}\equivvR 1$. 

\item[(ii)] For all \ninnt\ and $a\aprt b\inn ^n\Vex$: if $c\inn ^n\Vex$ with $a\aprt c \aprt b$ then $f_{a,b}\restrct{\Vexc} \subseteqq [\frac{1}{3}, \frac{2}{3}]$
\item[(iii)] For any \aprt-open $\Ucal\subseteqq\Vcal$ and successor point $x\inn\Ucal$ there is an \ninnt\ such that for all $a\inn ^nV$ we have:  $a\touch x_n$ implies $\hattr{a}\subseteqq\Ucal$. 
\ee

Then \Vnatt\ is metrizable.

\elem
 
\prf 
There is of course a canonical morphism $f$ from \Vexcalt\ to \zort\ such that $f(x)\equivvR 0$ for all $x\inn\Vexcal$. So by the assumption (i), we can define for all \ninn\ and all $a, b\inn ^n\Vex$ a morphism $f_{a,b}$
such that $a\touch b$ implies $f_{a,b}(x)\equivvR 0$ for all $x\inn\Vexcal$, and $a\aprt b$ implies $f_{a,b}\restrct{\Vexa}\equivvR 0\,$ and $f_{a,b}\restrct{\Vexb}\equivvR 1$.

\parr

Now let $h$ be an enumeration of $\bigcupp_{n\in\N}{^n\Vex}\timez ^n\Vex$. Define a metric $d$ on \Vexcalt\ by putting $d(x,y)\iz \sum_{m\in\N} 2^{-m}\cdott\midd f_{h(m)}(x)\minuz f_{h(m)}(y)\midd$. Then we see that $d(x,y)\bygr 0$ iff $x\aprt y$, and that $d$ is a metric on \Vexcalt\ (this shows weak metrizability).

\parr

To show that $d$ metrizes \Vnatt, let \Ucalt\ be \aprt-open in\Vnatt. We show that \Ucalt\ is $d$-open as well. For this let $x\inn\Ucal$ be a successor point. By assumption (iii) there is an \ninnt\ such that for all $a\inn ^nV$ we have:  $a\touch x_n$ implies $\hattr{a}\subseteqq\Ucal$. Now suppose we have a successor point $y\inn\Vcal$ such that $y_n\aprt x_n$. Determine \minnt\ such that $h(m)\iz (\mynon, x_n)$. Clearly $\mynon\aprt y_n\aprt x_n\aprt\mynon$, so by assumption (ii) we see that $f_{\mynon,x_n}(y)\inn [\frac{1}{3}, \frac{2}{3}]$, whereas $f_{\mynon,x_n}(x)\equivvR 1$. 

This implies that $d(x, y)\bygr 2^{-m-1}\cdott\frac{1}{3}$ and so in turn that $B(x, 2^{-m-1}\cdott\frac{1}{3})\subseteqq\Ucal$. Since $x\inn\Ucal$ is arbitrary, this shows that \Ucalt\ is $d$-open as well. (In fact $d$ metrizes \Vexcal, but this is unimportant).
\eprf

We turn to our splitting lemma. Remember that
for a spraid \Vnatt\ derived from \Vprenatt, and subsets $A, B$ of $V$ we write $A\aprt B$ iff $a\aprt b$ for all $a\inn A, b\inn B$. We write $A\touch B$ iff $a\touch b$ for some $a\inn A, b\inn B$. We shortly write $a\aprt B$, $a\touch B$ for $\{a\}\aprt B$, $\{a\}\touch B$ respectively.

\lem (\bol{B}) (splitting lemma)

Let \Wnatt\ be a \sucz-fan derived from \Wprenatt. Suppose $A, B$ are finite subsets of $W$ such that $A\aprt B$. Then there is an \Ninnt\ such that for all $c, d \inn {^N}^{\!}W$ we have: $(c\touch A \weddge d\touch B)$ implies $c\aprt d$ (and $(c\touch {^N}^{\!}{\!}A_{\leqc} \weddge d\touch {^N}^{\!}{\!}B_{\leqc})$ implies $c\aprt d$).

\elem

\prf
By proposition \ref{indbasdef}(ii), $C\iz\{c\inn W\midd c\aprt A\vee c\aprt B\}$ is an inductive bar on \Wcalt. By \HBind\ (crl.\,\ref{HBforsubfanns}) $C$ contains a finite bar $C'$. Let $A'\iz \{c\inn C'\midd c\touch A\}$, then $A'$ is finite and $A'\aprt B$. So by repeating our argument we find a finite bar $C''$ on \Wcalt\ such that for all $c\inn C''$ we have $c\aprt A'$ or $c\aprt B$. Put $N\iz\max(\{\grdd(c)\mid c\inn C''\cupp A\cupp B\})$. Then for all $c, d \inn {^N}^{\!}W$ we have: $(c\touch A \weddge d\touch B)$ implies $c\aprt d$ (and $(c\touch {^N}^{\!}{\!}A_{\leqc} \weddge d\touch {^N}^{\!}{\!}B_{\leqc})$ implies $c\aprt d$).
\eprf

The splitting lemma tells us that for a \sucz-fan \Wnatt, finite apart $A\aprt B\subsett W$ lead (for big enough \Ninn) to a partition of ${^N}\!W$ in three sets $C, D, E$ where $A\touch C$, $A\aprt D\aprt B$, $E\touch B$ and moreover $C\aprt E$. We will use such partitions obtained by lemma (\bol{B}) to construct morphisms from \Wnatt\ to \sigthrRt. In these constructions, to such $C$ we assign a $0$, to such $D$ a $1$ and to such $E$ a $2$.

\parr

For this, remember our definition in \ref{cantorfun}\ of $\sigthr=(\zotstar, \aprtom, \leqcom)$ and the surjective morphism \fevltert\ from \sigthrt\ to \zortert. Pulling back \aprtRt\ using \fevltert\ we saw that \fevltert\ is a \leqc-isomorphism from $\sigthrR\iz (\zotstar, \aprtR, \leqcom)$ to \zortert.

\defi (\bol{b})

We define a lexicographical ordering \lexomt\ on $\Nstar$ putting, for $a=a_0,\ldots a_{n-1}$ and $b\iz b_0, \ldots, b_{m-1}$ in \Nstart:

\parr

$a\lexom b$ iff ($b\precc a \vee \driss i\smlr n, m [ a_0, \ldots, a_{i-1}= b_0, \ldots, b_{i-1} \weddge a_i\smlr b_i] $).

\parr

By slight abuse of notation we write $^n\sigthr$ for $\{a\inn\zotstar\midd\grdd(a)\iz n\}$.
For each \ninnt, \lexomt\ induces a finite linear ordering on $^n\sigthr$. For $a\inn ^n\sigthr, a\notiz \zeron$ we write $\prdd(a)$ for the immediate predecessor of $a$ in this ordering. For $a\inn ^n\sigthr, a\notiz \twon$ we write $\sczz(a)$ for the immediate successor of $a$ in this ordering. Additionally we put $\prdd(\zeron)\isdef -1$ and $\sczz(\twon)\isdef 3$.
\edefi

\lem (\bol{C}) (Urysohn function lemma for fans)

Let \Wnatt\ be a \sucz-fan derived from \Wprenat. Let $a\aprt b\inn ^mW$ for certain \minn. Then there is a canonical morphism $f_{a,b}$ from \Wcalt\ to \zort\ such that:
\be
\item[(i)] $f_{a,b}\restrct{\Wsuba}\equivvR 0\,$ and $f_{a,b}\restrct{\Wsubb}\equivvR 1$. 

\item[(ii)] If $c\inn ^mW$ with $a\aprt c \aprt b$ then $f_{a,b}\restrct{\Wsubc} \subseteqq [\frac{1}{3}, \frac{2}{3}]$
\ee

\elem

\prf
We inductively define, for all $i\inn \zotstar\cupp\{-1, 3\}$, a subset $W_i\subsett W$ such that for all \ninnt:

\be
\item[\azt ]  For $i, j\inn ^n\sigthr\cupp\{-1, 3\}$ the set $W_i$ is decidable. Moreover $W_i\aprt W_j$ whenever  $j\notinn\{\prdd(i), i, \sczz(i)\}$, and if $i\notiz j$ are in $^n\sigthr$ then $W_i\capp W_j\iz\emptyy$.
  
\item[\aztazt ] For $i\inn\zotstar$ and $s\inn\{0, 1, 2\}$ we have $W_{i\strcomp\star s}\subsett W_i$. 

\item[\aztaztazt ] There is $t\inn\N$ such that all $c\inn ^tW$ are in some $W_i$ where $i\inn  ^n\sigthr$. 
\ee

Basis: for $n\iz 0$, put $W_{\maxdots}\iz W$, $W_{-1}\iz \Wsuba$ and $W_3\iz \Wsubb$, then $W_{-1}\aprt W_{3}$, and \aztazt\ and \aztaztazt\ are trivially fulfilled so we are done. 

\parr

Induction: suppose for \ninnt\ and $i\inn ^n\sigthr$ the subsets $W_i$ have been defined such that \azt, \aztazt, and \aztaztazt\ above hold. 

Let $i\inn ^n\sigthr$. Determine the least $t\inn\N$ such that all $c\inn ^tW$ are in some $W_j$ for $j\inn  ^n\sigthr$. Now put $A\iz W_{\prdd(i)}\capp ^tW$ and $B\iz W_{\sczz(i)}\capp ^tW$. Then by \azt\ we see that $A\aprt B\subsett W$ are finite subsets of $W$, so by lemma (\bol{B}) we can find the smallest $N\geqq t\pluz 1$ such that for all $c, d \inn {^N}^{\!}W$ we have: $(c\touch A\weddge d\touch B)$ implies $c\aprt d$.

Let $C\iz\{c\inn {^N}^{\!}W\midd c\touch A\}$, $D\iz \{c\inn {^N}^{\!}W\midd A\aprt c\aprt B\}$, and $E\iz \{c\inn {^N}^{\!}W\midd c\touch B\}$. Then $A\aprt D\aprt B$ and $C\aprt E$.
Put $W_{i_{\,}\star 0}\iz C_{\leqc}\capp W_i$, $W_{i_{\,}\star 1}\iz D_{\leqc}\capp W_i$, and $W_{i_{\,}\star 2}\iz E_{\leqc}\capp W_i$.

\parr

This defines $W_i$ for all $i\inn^{n+1}\sigthr$. It is straightforward to see that \azt\ and \aztazt\ hold for $n\pluz 1$. To see that \aztaztazt\ holds as well, it suffices to consider that $^n\sigthr$ is finite, and that therefore in our procedure above there was $j\inn ^n\sigthr$ with (copying notations) a maximal $N_j\geqq t\pluz 1$ compared to other $i\inn ^n\sigthr$. This means that all $c\inn{^{N_j}}{\!}W$ are in some $W_i$ for $i\inn ^{n+1}\sigthr$.

\parr

We define a \leqc-morphism $h_{a,b}$ from \Wcalt\ to \sigthrRt\ by specifying $h_{a,b}$ on $W$. Let $c\inn W$, then there is a maximal $\sinn, s\leqq\grdd(c)$ such that $c$ is in $W_{i}$ for a unique $i\inn ^s\sigthr$. We now simply put $h_{a,b}(c)\iz i$.

\clmm
$h_{a,b}$ is a \leqc-morphism from \Wcalt\ to \sigthrRt.
\eclmm

\clmmprf 
It follows from \azt, \aztazt\ and \aztaztazt\ above that for $c\leqc d\inn W$ we have $h_{a,b}(c)\leqc h_{a,b}(d)$. Let $x\inn\Wcal$ be a successor point, we need to show that $h_{a,b}(x)$ is a point in \sigthr.
This follows easily however from \aztaztazt\ above, since we see that for all \ninnt\ there is a $t\inn\N$ such that $h_{a,b}(x_t)\inn ^n\sigthr$. Finally, suppose $y\inn\Wcal$ is a successor point such that $h_{a,b}(y)\aprtR h_{a,b}(x)$, then we must show $y\aprt x$. Since $h_{a,b}(y)\aprtR h_{a,b}(x)$, we can determine \ninnt\ and $i, j\inn ^n\sigthr$ such that $i\aprtR j$ and $h_{a,b}(y)\precc i$ and $h_{a,b}(x)\precc j$. However, $i\aprtR j$ equals $j\notinn\{\prdd(i), i, \sczz(i)\}$. There is \sinnt\ such that $y_s\inn W_i$ and $x_s\inn W_j$, where moreover $W_i\aprt W_j$ since $j\notinn\{\prdd(i), i, \sczz(i)\}$. So we see that $y_s\aprt x_s$ and so $y\aprt x$.
\eclmmprf 

We can now define $f_{a,b}\isdef \fevlter\circ h_{a,b}$ (which is canonical since $h_{a,b}$ and \fevltert\ are constructed canonically). Clearly $f_{a,b}\restrct{\Wsuba}\equivvR 0\,$ and $f_{a,b}\restrct{\Wsubb}\equivvR 1$ (proving (i) of the lemma). So the only thing left to prove is that if $c\inn ^mW$ with $a\aprt c \aprt b$ then $f_{a,b}\restrct{\Wsubc} \subseteqq [\frac{1}{3}, \frac{2}{3}]$. And this follows trivially from our construction, since there is \sinnt\ such that $\{d\inn \Wsubc\midd \grdd(d)\geqq s\} \subsett W_{1}$ and $f_{a,b}\restrct{W_1}\subseteqq [\frac{1}{3},\frac{2}{3}]$.
\eprf

\prp 
Every \sucz-fan is metrizable.
\eprp

\prf
Let \Wnatt\ be a \sucz-fan. Then trivially \Wexcalt\ is also a \sucz-fan. Therefore by lemma (\bol{C}) above, we find: 

\be
\item[(i)] For all \ninn, for all $a\aprt b\inn ^n\Wex$ there is a given morphism $f_{a,b}$ from \Wexcalt\ to \zort, such that $f_{a,b}\restrct{\Wexa}\equivvR 0\,$ and $f_{a,b}\restrct{\Wexb}\equivvR 1$. 

\item[(ii)] For all \ninnt\ and $a\aprt b\inn ^n\Wex$: if $c\inn ^n\Wex$ with $a\aprt c \aprt b$ then $f_{a,b}\restrct{\Wexc} \subseteqq [\frac{1}{3}, \frac{2}{3}]$.
\ee

\clmm 
Let $\Ucal\subseteqq\Wcal$ be \aprt-open and let $x\inn\Ucal$ be a successor point. Then there is an \ninnt\ such that for all $a\inn ^nW$ we have:  $a\touch x_n$ implies $\hattr{a}\subseteqq\Ucal$. 
\eclmm

\clmmprf 
By definition \ref{indbasdef}\ of `\sucz-spread', \Ucalt\ is \sucz-open, which means that $\BUx\iz\{b\inn W\midd b\aprt x_{\grds(b)}\vee\hattr{b}\subseteq\Ucal\}$ is an inductive bar on $W$. 
Now since \Wcalt\ is a fan, by \HBind\ (crl.\,\ref{HBforsubfanns}) we find a finite subbar $B'\subseteqq \BUx$ on $W$. Let $n$ be the maximum of $\{\grdd(b)\midd b\inn B'\}$. Let $a\inn ^nW$. We know by the properties of $B'$ that $a\aprt x_n$ or $\hattr{a}\subseteqq\Ucal$. Therefore $a\touch x_n$ implies $\hattr{a}\subseteqq\Ucal$.
\eclmmprf 

By the claim, we see that \Wnatt\ satisfies the conditions (i), (ii) and (iii) of the Urysohn lemma (\bol{A}). Therefore \Wnatt\ is metrizable.
\eprf

To prove corollary \ref{starfinmet}, we define `one-point \sucz-fanlike extension' (see also \ref{indothdef}) to represent `locally compact' spaces,  just as \sucz-fanlike spaces represent `compact'  spaces.\footnote{The difference with the \classf\ and \bishf\ notions is that our analogons need not be locally metrically complete. This can be addressed satisfactorily, in the sense that local completeness turns out to be a natural-topological property (invariant under isomorphisms). But we will leave this for subsequent expositions on natural topology, hopefully written by others.}\ 

\defi (\bol{c})

Let \Vnatt\ be a spread, then \Vnatt\ has a \deff{one-point \sucz-fanlike extension}\ iff there is a \sucz-fan \Wnattot\ (derived from \Wprenatto) and a function $f$ from \Vext\ to $W$ such that putting $a\aprto b$ iff $f(a){\aprtto}{_{_{\,}}} f(b)$, we have that $f$ is an inductive isomorphism from $(\Vexcal, \Topaprto)$ to \Wnattot\ where in addition $f(\myno){\aprtto}{_{_{\,}}} f(x)$ for all $x\inn\Vcal$ and \Vnatot\ is identically automorphic to \Vnatt. More generally, a natural space \Vnatt\ is said to have a one-point \sucz-fanlike extension iff \Vnatt\ is \sucz-isomorphic to a spread with a one-point \sucz-fanlike extension.
\edefi

\crl (which is the same as corollary \ref{starfinmet})
Every \sucz-fanlike space is metrizable (`every compact space is metrizable'), and every space with a one-point \sucz-fanlike extension is metrizable (`every locally compact space is metrizable').
\ecrl

\prf
We only need to prove that a spread \Vnatt\ (derived from \Vprenatt) with a one-point \sucz-fanlike extension is metrizable. For this, let \Wnattot\ be a \sucz-fan (derived from \Wprenatto) and $f$ a function from \Vext\ to $W$ such that putting $a\aprto b$ iff $f(a){\aprtto}{_{_{\,}}} f(b)$, we have that $f$ is an inductive isomorphism from $(\Vexcal, \Topaprto)$ to \Wnattot\ where in addition $f(\myno){\aprtto}{_{_{\,}}} f(x)$ for all $x\inn\Vcal$. 
\parr
By the above proposition \Wnattot\ is metrizable. This means that 
$(\Vexcal, \Topaprto)$ is metrizable, say by metric $d$. Then the restriction of $d$ to \Vcalt\ metrizes \Vnatt. This follows trivially from the fact that $\Topaprt\restrct{\textstyle\Vcall}=\Topaprto\restrct{\textstyle\Vcall}$, which in turn follows trivially from the fact that $\myno\aprto x$ for all $x\inn\Vcal$. 
\eprf

\sbsc{Proof of theorem \ref{starfinmet}}\label{prfstarfinmet}

The key to the generalization of our results in paragraph \ref{prfcrlstarfinmet}\ is the following simple observation. In a star-finite \Vnatt, for each $x$ in \Vcalt\ the equivalence class of $x$, that is $\{y\inn\Vcal\midd y\equivv x\}$ is contained in a subfan \Wxt\ of \Vnatt. The spine of this subfan \Wxt\ is formed by $\{a\inn V\midd a\touch x_m\midd \grdd(x_m)\iz\grdd(a)\}$. But we need to do some extra work to turn this spine into a subfan. Specifically, if $a$ is in the spine but all continuations $b\precc a$ are seen to be not in this spine (which is decidable since \Vnatt\ is star-finite), we need to ensure that $a$ still contains a point in \Wxt. 

\parr

To avoid cumbersome repetition of the prevailing conditions, we state the following: 

\conv From now on, without loss of generality, we assume \Vnatt\ to be a star-finite \sucz-spread, derived from \Vprenatt\ where $v:\N\rightarrow V$ is an enumeration of $V\iz\{v_n\midd \ninn\}$ such that for each $n, \minn$ we have that $v_n\precc v_m$ implies $n\bygr m$. The touch-relation \toucht\ is the complement of \aprt\ on $V\timez V$. We concentrate on successor points (see def.\,(\bol{a}) above). Also we write `$m\iz\mu\sinn [{\rm P}(s)]$' as abbreviation for `$m$ is the smallest natural number for which P$(m)$ holds'.
\econv

\defi (\bol{d})

Let $x\inn \Vcal$ be a successor point. Using the enumeration $v$ we define a subfan \Wxt\ of \Vnatt\ by inductively describing $^n\Wx$ for each \ninnt. Of course $^0\Wx\isdef \{\maxdot\}$. Now suppose that $^n\Wx$ has been defined for given \ninnt, then put

\btab

$^{n+1}\Wx\isdef\{v_m\inn ^{n+1}V\mid $ \= $v_m\touch x_{n+1} \vee $ \\
\> $\dris a\inn ^n\Wx \all b\inn\suczz(a) [b\aprt x_{n+1}\weddge m\iz\mu\sinn [v_s\precc a]]\}$.

\etab

Then \Wxt\ is a subfan of \Vnatt\ such that $\{y\inn\Vcal\midd y\equivv x\}\subseteqq\Wx$. For \ninnt, we are mostly interested in the subset $\Wnxtch\isdef\{a\inn ^n\Wx\midd a\touch x_n\}$ of $^n\Wx$. Finally, let $y\inn\Vcal$ be arbitrary, then $y$ has as subsequence the successor point $y^{\!\sucz}\equivv y$ (see def.\,(\bol{a}) above), and we put $W_{y}\iz W_{{y^{\!\sucz}}}$.
\edefi

Next we show that \Wxt\ acts almost like a neighborhood of $x$ in \Vnatt. We need this later on, as prerequisite (iii) of the Urysohn metrization lemma. 

\lem (\bol{D})

Let \Ucalt\ be open in \Vnatt, and let $x\inn\Ucal$. Then there is \Ninnt\ such that for all  $a\inn ^NV$ we have: $a\touch \xN$ implies $\hattr{a}\subseteqq\Ucal$.
\elem

\prf
By definition \ref{indbasdef}\ of `\sucz-spread', \Ucalt\ is \sucz-open, which means that $\BUx\iz\{b\inn V\midd b\aprt x_{\grds(b)}\vee\hattr{b}\subseteq\Ucal\}$ is an inductive bar on $V$. Now since \Wxt\ is a subfan and $\BUx$ is monotone, by \HBind\ (crl.\,\ref{HBforsubfanns}) we find a finite subbar $B'\subseteqq \BUx$ on \Wxt. Let $N$ be the maximum of $\{\grdd(b)\midd b\inn B'\}$. Let $a\inn ^NV$. We know that if $a \notinn \Wx$ then $a\aprt \xN$. Else, if $a\inn \Wx$ we know by the properties of $B'$ that $a\aprt \xN$ or $\hattr{a}\subseteqq\Ucal$. We combine this to conclude for all $a\inn ^NV$ that $a\touch \xN$ implies $\hattr{a}\subseteqq\Ucal$.
\eprf

We need a sequence $(\Wxn)_{n\in\N}$ of very similar but slightly larger fans than \Wxt\ for our purposes, where $\Wxn\subseteqq\Wxnpl$ for each \ninn. For this we expand our touch relation \toucht\ inductively to equivalent touch-relations $(\!\touchn\!)_{n\in\N}$ (`equivalent' meaning that they induce the same apartness relation on points).  

\defi (\bol{e})

Let $\touch\!'\subsett V\timez V$, and write $\aprt'$ for the complement of $\touch\!'$. We say that $\touch\!'$ is a \deff{$\touch$-equivalent touch-relation}\ iff $\touch\!'$ is decidable and for all points $x, y$ in \Vcalt, we have:  $x\aprt y$ iff there is \ninnt\ with $x_n\aprt' y_n$. We say that $\touch\!'$ is star-finite
iff for each $a\inn V$ the subset $\{b\inn V\midd\grdd(b)\iz\grdd(a) \weddge b\touch\!' a\}$ is finite. 

\parr

We inductively define decidable relations $(\!\touchn\!)_{n\in\N}$ which are \touch-equivalent for $n\geqq 1$ as follows. Let $\touchze\iz\{(a,b)\midd a,b\inn V\midd a\leqc b \vee b\leqc a\}$. Suppose for \ninnt\ that $\touchn$ has been defined. Let $a, b\inn V$ with $m\iz\grdd(a)\leqq\grdd(b)$. Then $a\touchnpl b$ and $b\touchnpl a$ iff there is $c\inn {^m}V$ such that
         $a\touchn c$ and $c\touch b$. (The idea is that $a\touchn b$ iff there is a \touch-trail of length $n$ from $a$ to $b$ which does not use basic dots $d$ with $\grdd(d)\smlr\grdd(a)\leqq\grdd(b)$. Notice that $\touchone$ equals \toucht.).\footnote{The touch-relation $\touchze$ could also be named $\touch_{\!\omega}$ since it corresponds to the `naked' spread $(\Vcalt, \Topaprtom)$ where we have stripped \Vcalt\ of the apartness relation \aprt. Compare this to our discussion of \Raprtomt\ derived from $(\sigRungl, \aprtom, \leqc)$ and \hwki\ in example \ref{exhawkeye}.}

\parr

Let $a\inn ^mV$ for certain \minn, then we let $\starn(a)\iz\{b\inn^mV\midd a\touchn b\}$, which is a finite set. For $\staro(a)$ we also simply write ${\rm St}(a)$. For the complement of \touchnt\ we write $\touchnotn$ and for subsets $A, B$ of $V$ we write $A\touchnotn B$ iff $a\touchnotn b$ for all $a\inn A, b\inn B$. We write $A\touchn B$ iff $a\touchn b$ for some $a\inn A, b\inn B$. We shortly write $a\touchnotn B$, $a\touchn B$ for $\{a\}\touchnotn B$, $\{a\}\touchn B$ respectively.

\parr

For $x\inn\Vcal$ we define \Wxnt\ as the subfan of \Vnatt\ which we obtain by substituting $\touchn$ for \toucht\ and $\touchnotn$ for $\aprt$ in the definition (\bol{d}) of \Wxt\ above.  

\parr

Finally, for $n, \minn$ put $^n\Wxmtch\isdef\{a\inn ^n\Wxm\midd a\touchm {x^{\propto}}\!\!\!\!_n\}$, where $x^{\propto}\equivv x$ is the relevant successor point, see def.\,(\bol{a}).
\edefi 

To see that the definition of \Wxnt\ is valid, it suffices to check that $\touchn$ is again star-finite (for $a\inn V$, the set $\{b\inn V\midd \grdd(b)\iz\grdd(a) \weddge b\touchn a\}$ is finite). 

\parr

Lemma (\bol{B}) above can now be abbreviated thus: let \Wnatt\ be a \sucz-fan derived from \Wprenatt. Suppose $A, B$ are finite subsets of $W$ such that $A\aprt B$. Then there is an \Ninnt\ such that ${^N}^{\!}{\!}A_{\leqc}\touchnotthr {^N}^{\!}{\!}B_{\leqc}$.\ \footnote{This shows that $\touchn$ is \touch-equivalent for $n\geqq 1$.}

\parr

We use lemma (\bol{B}) to generalize the Urysohn function lemma (\bol{C}) to \Vnatt. This is an arduous task, which we try to make palatable by dividing the proof in two parts. In the first part we detail the construction (for $a\aprt b\inn {^M}{^{\!}}V$ and for $i\inn \zotstar\cupp\{-1, 3\}$) of subsets $V_i\subsett V$ such that in analogy to lemma (\bol{C}) for $i, j\inn ^n\sigthr\cupp\{-1, 3\}$ the set $V_i$ is decidable; if $j\notinn\{\prdd(i), i, \sczz(i)\}$ then $^mV_i\aprt ^mV_j$ for all $m\geqq M$, and for $i,j\inn ^n\sigthr$: if $i\notiz j$ then $V_i\capp V_j\iz\emptyy$. Also for $i\inn\zotstar$ and $s\inn\{0, 1, 2\}$ we have $V_{i\strcomp\star s}\subsett V_i$.

In the second part of the proof these sets will help us to define a morphism $h_{a,b}$ to \sigthrRt, which using \fevltert\ can be turned into the desired morphism $f_{a,b}$. 

\lem (\bol{E}) (Urysohn function lemma, generalizing lemma (\bol{C}))

Let the basic dots $a\aprt b$ be in ${^M}{^{\!}}V$ for certain \Minn. Then there is a canonical morphism $f_{a,b}$ from \Vcalt\ to \zort\ such that:
\be
\item[(i)] $f_{a,b}\restrct{\Vsuba}\equivvR 0\,$ and $f_{a,b}\restrct{\Vsubb}\equivvR 1$. 

\item[(ii)] If $c\inn {^M}{^{\!}}V$ with $a\aprt c \aprt b$ then $f_{a,b}\restrct{\Vsubc} \subseteqq [\frac{1}{3}, \frac{2}{3}]$
\ee

\elem

\prf

\parr

\underline{Part one}

We first put $V_{-1}\iz \Vsuba$ and $V_3\iz \Vsubb$. Then $V_{-1}\aprt V_3$.

We then inductively define, for \ninnt\ and $i\inn ^n\sigthr$, subsets $V_i\subsett V$ as well as a subset ${\rm Sec}_n\subsett V$ such that (except for $i\iz \maxdotB\inn ^0\sigthr$) a basic dot $c$ is member of any of these sets iff $\grdd(c)\geqq M$,  and in addition we have:

\be
\item[\azt ]  For $i \inn ^n\sigthr$ the set $V_i$ is decidable and $V_i\iz (V_i)_{\leqc}$. \item[\aztazt] For $i, j\inn ^n\sigthr\cupp\{-1, 3\}$, if $j\notinn\{\prdd(i), i, \sczz(i)\}$ then $V_i\aprt V_j$,  and if $i\notiz j$ are in $^n\sigthr$ then $V_i\capp V_j\iz\emptyy$. 
  
\item[\aztaztazt ] For $j\inn ^{n-1}\sigthr$ and $s\inn\{0, 1, 2\}$ we have $V_{j\strcomp\star s}\subseteqq V_j$. 

\ee

When our construction is done, we find that for all $x\inn\Vcal$ and all \ninnt\ there is $t\inn\N$ and $i\inn ^n\sigthr$ such that $x_t\inn V_i$

\parr

Basis: for $n\iz 0$, the only member of $^n\sigthr$ is $\maxdotB$. Put $V_{\maxdotBs}\iz V$. Then $V_{-1}\aprt V_{3}$, and \azt, \aztazt\ and \aztaztazt\ are trivially fulfilled so we are done. 

\parr

Induction: suppose for \ninnt\ and all $i\inn ^n\sigthr$ the subsets $V_i$ with properties \azt, \aztazt\ and \aztaztazt\ have been defined. Suppose the basic dot $c$ is in some $V_i$ for $i\inn ^n\sigthr$. When trying to classify $c$ on the next level $n\pluz 1$, we can only be sure of making the right choice if all $c$'s neighbors and their neighbors (in other words all members of $\starto(c)$) are classified on level $n$. We therefore first put:

\parr
$\secn\isdef\{c\inn V\midd \grdd(c)\geqq M \weddge \all d\inn\starto(c)\driss j\inn ^n\sigthr\,[d\inn V_j]\,\}$.
\parr

$\secn$ is the decidable set of `$n$-secure' basic dots in our classification scheme. (We take ${\rm Sec}_0\iz ({^M}{^{\!}}V)_{\leqc}$ to fulfill (ii) of the lemma.). Now for $i\inn ^{n}\sigthr$ define: 

\parr

$V_{i_{\,}\star 0}\isdef (\bigcupp_{t\in\N}\{c\inn ^tV_i\midd c\inn\secn\weddge c\touch ^tV_{\prdd(i)}\weddge c\touchnotto ^tV_{\sczz(i)}\,\})_{\leqc}$.

\newlyne

$V_{i_{\,}\star 2}\isdef (\bigcupp_{t\in\N}\{c\inn ^tV_i\midd c\inn\secn \weddge c\touch ^tV_{\sczz(i)}\weddge c\touchnotto ^tV_{\prdd(i)}\,\})_{\leqc}$.

\newlyne

$V_{i_{\,}\star 1}\isdef \bigcupp_{t\in\N}\{c\inn ^tV_i\midd c\inn\secn \weddge c\notinn V_{i_{\,}\star 0} \weddge c\notinn V_{i_{\,}\star 2} \weddge c\aprt ^tV_{\prdd(i)}\weddge c\aprt ^tV_{\sczz(i)}\,\}$.

\parr

Notice that if for $c\inn ^tV$ we know $c\aprt ^tV_{\prdd(i)}\weddge c\aprt ^tV_{\sczz(i)}$, then $c$ can still be in $V_{i_{\,}\star 0}$ or $V_{i_{\,}\star 2}$, simply because there can be a  $d\succ c$ which is in $V_{i_{\,}\star 0}$ or $V_{i_{\,}\star 2}$. 

\parr

The sets $V_{i\strcomp\star s}$ are decidable and monotone, since \toucht\ is decidable and for $c\inn V$ the sets $\starto(c)$ and $\{b\inn V\midd c\precc b\}$ are finite. For $i, j\inn ^{n+1}\sigthr\cupp\{-1,3\}$ we see that $j\notinn\{\prdd(i), i, \sczz(i)\}$ implies $V_i\aprt V_j$, and if $i\notiz j$ are in $^{n+1}\sigthr$ then $V_i\capp V_j=\emptyy$. This establishes \azt\ and \aztazt, and \aztaztazt\ follows straight from the definition.

\parr

We will use the sets $(V_i)_{i\in\zotstar}$ in a way similar to the proof of lemma (\bol{C}). But other than in that proof, for \ninnt\ we cannot determine a uniform level $N$ on $V$ such that for $t\geqq N$ all $c\inn ^tV$ are in some $V_i$ with $i\inn ^n\sigthr$. However, for $x\inn\Vcal$ we can use the fans $(\Wxm)_{m\in\N}$ defined in (\bol{e}) to show that there still is an \Ninnt\ such that $\xN$ is in some $V_i$ with $i\inn ^n\sigthr$. 

This we achieve in the second part of the proof, using a number of claims. We then conclude that the desired morphisms $h_{a,b}$ and $f_{a,b}$ can be derived from our construction above just as in the proof of lemma (\bol{C}). 

\parr

\underline{Part two}

We need to prove that for $x\inn\Vcal$, for all \ninnt\ there is $i\inn ^n\sigthr$ such that $x\inn\hattr{V_i}$ (which is equivalent to there being an \Ninnt\ with $\xN\inn V_i$). Since $V_{\maxdotBs}\iz V$, this is trivial for $n\iz 0$. 

We use the sets $(\secn)_{n\in\N}$ and a form of double induction, on \ninnt\ and on \minnt\ where $m$ is the index of the fans $(\Wxm)_{m\in\N}$ which are needed to make the induction step.

\parr

The induction basis where $n\iz 0$ is trivial: clearly for $m, t\inn\N$ such that $t\geqq M$ we have: $\all c\inn ^t\Wxmtch\ \driss i\inn^n \sigthr\, [c\inn V_i]$, since $V_{\maxdotBs}\iz V$. We now turn to the inductive step going from $n$ to $n\pluz 1$. We have to expand our strategy to \Wxmpltot\ to ensure that $c\inn\secn$ for the relevant $c\inn \Wxm$, so that these $c$ have no unexpected neighbors and can be classified on level $n\pluz 1$.

\clmm
Let $x\inn\Vcal$ and $n, m, t\inn\N$ such that $t\geqq M$ and $\allz_{\!}c\inn ^t\Wxmpltotch\dris i\inn ^n\sigthr{_{\!}}[c\inn V_i]_{\!}$. Then: $\dris\Ninn\,\all c\inn {^N}{^{\!}}\Wxmtch \dris j\inn ^{n+1}\sigthr[c\inn V_j]$. 
\eclmm

\clmmprf  
By the conditions, $c\inn\secn$ for all $c\inn (^t\Wxmtch)_{\leqc}$. For $i\inn ^n\sigthr$ we see by \aztazt\ above: $^tV_{\prdd(i)}\aprt ^tV_{\sczz(i)}$. Put $A\iz ^tV_{\prdd(i)}\capp\Wxmplto$ and $B\iz ^tV_{\sczz(i)}\capp\Wxmplto$, then $A\aprt B$ are finite subsets of \Wxmpltot, so by our splitting lemma (\bol{B}) there is an $N_i\inn\N, N_i\geqq t$ such that for all $c, d \inn {^{N_i}}^{\!}\Wxmplto$ we have: $(c\touch {^{N_i}}^{\!}A_{\leqc} \weddge d\touch {^{N_i}}^{\!}B_{\leqc})$ implies $c\aprt d$. Again, for all $c\inn{^{N_i}}^{\!}\Wxmtch$ we also have $c\inn\secn$.

\parr

By definition (\bol{e}) and the choice of $A, B$ we see that for $c\inn{^{N_i}}^{\!}\Wxmtch$ the conditions $c\touch {^{N_i}}^{\!}A, c\touchto {^{N_i}}^{\!}A$ equal the conditions $c\touch {^{N_i}}^{\!}V_{\prdd(i)}, c\touchto {^{N_i}}^{\!}V_{\prdd(i)}$ and likewise the conditions $c\touch {^{N_i}}^{\!}B, c\touchto {^{N_i}}^{\!}B$ equal the conditions $c\touch {^{N_i}}^{\!}V_{\sczz(i)}, c\touchto {^{N_i}}^{\!}V_{\sczz(i)}$. 

\parr

Therefore for $c\inn{^{N_i}}^{\!}\Wxmtch\capp V_i$ we can decide: \kase{0} $c\touch {^{N_i}}^{\!}V_{\prdd(i)}\weddge c\touchnotto {^{N_i}}^{\!}V_{\sczz(i)}$, so $c\inn V_{i_{\,}\star 0}$ OR \kase{1} $c\aprt {^{N_i}}^{\!}V_{\prdd(i)}\weddge c\aprt {^{N_i}}^{\!}V_{\sczz(i)}$, then we must check \kase{1.0} there is a $d\succ c$ which is already in $V_{i_{\,}\star 0}$, then $c\inn V_{i_{\,}\star 0}$ or \kase{1.2} there is a $d\succ c$ which is already in $V_{i_{\,}\star 2}$, then $c\inn V_{i_{\,}\star 2}$ or else \kase{1.1} $c\inn V_{i_{\,}\star 1}$ OR \kase{2} $c\touch {^{N_i}}^{\!}V_{\sczz(i)}\weddge c\touchnotto {^{N_i}}^{\!}V_{\prdd(i)}$, then $c\inn V_{i_{\,}\star 2}$. 

\parr

This shows that for all $c\inn {^{N_i}}^{\!}\Wxmtch\capp V_i$ we can find $j\inn ^{n+1}\sigthr$ with $c\inn V_j$. Finding for each $i\inn ^n\sigthr$ a likewise $N_i$, put $N\iz\max(\{N_i\midd i\inn ^n\sigthr\})$.
\eclmmprf

\clmm
For all $n,\minn$ there is an \Ninnt\ such that $\all c\inn {^{N}}^{\!}\Wxmtch \driss i\inn^n \sigthr\, [c\inn V_i]$.
\eclmm

\clmmprf 
By double induction. 

Basis: for $n\iz 0$ the statement is trivially true for all \minn. 

Induction: Suppose the statement is true for given \ninnt\ and all \minn. Then for \minn, there is $t\inn\N$ such that $t\geqq M$ and $\all c\inn ^t\Wxmpltotch\ \driss i\inn^n \sigthr\, [c\inn V_i]$. By the previous claim, there is an \Ninnt\ such that $\all c\inn {^{N}}^{\!}\Wxmtch \driss i\inn ^{n+1}\sigthr\, [c\inn V_i]$.

\eclmmprf 

The claim guarantees that our next construction will define a morphism. We use the sets $V_i$ to define a canonical function $h_{a,b}$ from $V$ to \zotstart, quite similar to the proof of lemma (\bol{C}). For \ninnt\ and $c\inn ^nV$ put: 

\btb

$h_{a,b}(c)\isdef$ \= the unique $i$ in \zotstart such that:                                                                      $\grdd(i)\leqq\grdd(c)\iz n$ and $c\inn V_i$ and\\
\>$\all j\inn\zotstar [(\grdd(j)\leqq n \weddge c\inn V_j)\rightarrow i\leqcom j]$.

\etb  

That this definition is valid is easily seen, since there are only finitely many decidable conditions  to check and $c$ is certainly in $V_{\maxdotBs}$.

\clmm
$h_{a,b}$ is a \leqc-morphism from \Vcalt\ to \sigthrRt.
\eclmm

\clmmprf 
It follows from \aztaztazt\ above that for $c\leqc d\inn V$ we have $h_{a,b}(c)\leqc h_{a,b}(d)$. Let $x\inn\Vcal$, we need to show that $h_{a,b}(x)$ is a point in \sigthr.
This follows easily however from the claim above, since we see that for all \ninnt\ there is an \Ninnt\ such that $h_{a,b}(\xN)\inn ^n\sigthr$. Finally, suppose $y\inn\Vcal$ is a point such that $h_{a,b}(y)\aprtR h_{a,b}(x)$, then we must show $y\aprt x$. Since $h_{a,b}(y)\aprtR h_{a,b}(x)$, we can determine \ninnt\ and $i, j\inn ^n\sigthr$ such that $i\aprtR j$ and $h_{a,b}(y)\precc i$ and $h_{a,b}(x)\precc j$. However, $i\aprtR j$ equals $j\notinn\{\prdd(i), i, \sczz(i)\}$. There is \sinnt\ such that $y_s\inn V_i$ and $x_s\inn V_j$ where $V_i\aprt V_j$ since $j\notinn\{\prdd(i), i, \sczz(i)\}$ (by \aztazt), so we see that $y_s\aprt x_s$ and so $y\aprt x$.
\eclmmprf 

We define: $f_{a,b}\isdef \fevlter\circ h_{a,b}$ (canonically since $h_{a,b}$ and \fevltert\ are constructed canonically). Clearly $f_{a,b}\restrct{\Wsuba}\equivvR 0\,$ and $f_{a,b}\restrct{\Wsubb}\equivvR 1$, proving (i) of the lemma. 

\parr

The only thing left to prove is that if $c\inn {^M}^{\!}W$ with $a\aprt c \aprt b$ then $f_{a,b}\restrct{\Wsubc} \subseteqq [\frac{1}{3}, \frac{2}{3}]$. Yet it follows from our construction that $c\inn V_1$. Also trivially $f_{a,b}\restrct{V_1}\subseteqq [\frac{1}{3},\frac{2}{3}]$.
\eprf

\thm (from \ref{starfinmet})
Every star-finitary natural space is metrizable.
\ethm 

\prf
It suffices to prove the theorem for our star-finite \sucz-spread \Vnatt. 
 
\clmm 
\Vnatt\ meets the conditions of the Urysohn metrization lemma (\bol{A}).
\eclmm

\clmmprf  
We need to show:

\be

\item[(i)] For all \ninn, for all $a\aprt b\inn ^n\Vex$ there is a canonical morphism $f_{a,b}$ from \Vexcalt\ to \zort, such that $f_{a,b}\restrct{\Vexa}\equivvR 0\,$ and $f_{a,b}\restrct{\Vexb}\equivvR 1$. 

\item[(ii)] For all \ninnt\ and $a\aprt b\inn ^n\Vex$: if $c\inn ^n\Vex$ with $a\aprt c \aprt b$ then $f_{a,b}\restrct{\Vex_c} \subseteqq [\frac{1}{3}, \frac{2}{3}]$

\item[(iii)] For any \aprt-open $\Ucal\subseteqq\Vcal$ and successor point $x\inn\Ucal$ there is an \ninnt\ such that for all $a\inn ^nV$ we have:  $a\touch x_n$ implies $\hattr{a}\subseteqq\Ucal$. 

\ee

Ad (i) and (ii): Consider that $(\Vexcal, \Topaprt)$ is again a star-finite \sucz-spread. Now apply the Urysohn function lemma (\bol{E}) above.

Ad(iii): this is precisely the content of lemma (\bol{D}) above.
\eclmmprf 

This shows that \Vnatt\ satisfies the requirements of the Urysohn metrization lemma (\bol{A}),
and so \Vnatt\ is metrizable.
\eprf

\sbsc{Proof of meta-theorem \ref{intuinruss}}\label{prfcplife}

The two-player game \LIfE\ serves as an illustration of our belief that \intu\ can be (formally) interpreted in \russ\ also, in an elegant way. We state and prove:

\mthm (repeated from \ref{intuinruss})
\be
\item[(i)]In \LIfE, we can prove \CP. 
\item[(ii)]
Suppose \GoD\ is omniscient. Then we can prove $\negg\dris B\subsett\Nstar [B$ is a non-inductive bar on $\Bnat]$ for \LIfE.
\item[(iii)] Given enough time, \HuMAN\ can discover that by an overwhelming odds ratio, \GoD\ plays only recursive sequences. 
\ee
\emthm

\crl
In \class, we can prove \CP\ and \BT\ for the game \LIfE.
\ecrl

\prf
Ad (i). Suppose that \HuMAN\ has a set $A\subseteq\NN\timez\N$ such that in \LIfE: $\all\alpha\inn\NN\driss\ninn\,[(\alpha, n)\inn A]$. This means: for any sequence $\alpha$ played by \GoD, \HuMAN\ can produce at a finite moment in time an \ninnt\ such that $(\alpha, n)\inn A$.

\parr
Now let $\alpha\inn\NN$. Let \GoD\ play the sequence $\gamma$ which mimicks $\alpha$, starting out as $\gamma\iz\alpha(0), \alpha(1), \ldots$ without revealing any information to \HuMAN. Say at point $m$ in time (when \GoD\ has revealed precisely $\alfstr(m)$, the first $m$ values of $\alpha$), \HuMAN\ produces $n$ such that $(\gamma, n)\inn A$. By the rules of \LIfE, player \GoD\ may switch to any $\beta$ (computable, but \HuMAN\ does not know this) such that $\betastr(m)\iz\gamstr(m)\iz\alfstr(m)$ to conclude that for any such $\beta$ we have $(\beta, n)\inn A$.
Therefore $\all\alpha\inn\NN\driss m,n\inn\N\all\beta\inn\NN\,[ \alfstr(m)=\betastr(m)\rightarrow (\beta, n)\inn A]$, which proves \CP\ for \LIfE.

\parr

Ad (ii). We use the correspondence between natural Baire space \Bnatt\ and \NNt, described in \ref{natBaire}. Suppose \HuMAN\ has a set $B\subsett\Nstar$ which is non-inductive (meaning it does not descend from a genetic bar on \Bnat), and a bar $C\subseteqq\Nstar$ on \Bnatt\ with $B\subseteqq C$. Since $B$ is non-inductive, with omniscience \GoD\ can play a sequence $\alpha\iz p^*$ where $\alfstr\iz p\inn\Bnat$ is given by $p=p_0, p_1,\ldots$ (with $p_0\iz\maxdot$ and $p_{i+1}\sucz p_i$) and where $(p_i)_{\!\leqc}\capp B$ is non-inductive on ${\Nstar}\hspace*{-.6em}_{p_i}\iz (p_i)_{\!\leqc}$ for all $i$ up until and including the $n$ for which $\alfstr(n)\iz p_{n}\inn C$ (since $C$ is claimed to be a bar, \HuMAN\ must produce such $n$ at some finite point in time). 

\parr

Having received from \HuMAN\ the $n$ for which $\alfstr(n)\iz p_n\inn C$, \GoD\ switches to the recursive sequence $\beta\inn\NN$ such that $\beta\iz p_n\starr\zero$. To arrive at $\beta$, \GoD\ can claim to have used the recursive sequences $p_1\starr\zero, p_2\starr\zero,\ldots$ and to have switched $n\minuz 1$ times, therefore \GoD\ has played by the rules.

\parr

We see now, that $\alfstr(n)\iz p_n$ is in $C$ but not in $B$, since $p_n\inn B$ would imply that $(p_n)_{\!\leqc}\capp B$ is inductive on ${\Nstar}\hspace*{-.6em}_{p_n}\iz (p_n)_{\!\leqc}$, contradicting the above. We conclude: $B\notiz C$. This proves (ii).

\parr

Ad (iii). This is basically the same as the physical experiment described in \WaaArt, section 7. \HuMAN\ can build a covering of the recursive unit interval with a sequence of rational intervals $(R_n)_{n\in\N}$ such that the sum of the lengths of these intervals does not exceed $2^{-40}$. Asking \GoD\ for a non-recursive sequence $\alpha\inn\zor$, \HuMAN\ will eventually discover that $\alpha$ falls within one of the intervals $R_n$. (\GoD\ can only switch finitely many times, and then in the end is stuck with a recursive $\alpha$, which will eventually be captured by some $R_n$.) The odds of this happening for a non-recursive $\alpha$ are overwhelmingly small, proving (iii).

\parr

The corollary follows from (i) and (ii), since with classical logic \GoD\ is omniscient and $\negg\dris B\subsett\Nstar [B$ is a non-inductive bar on $\Bnat]$ is equivalent to `every bar on \Bnatt\ is inductive'.
\eprf

\sectionnb{Constructive concepts and axioms used}\label{axioms}

\sbsc{Basic axioms and concepts}
In this section we present some axioms and concepts described in the literature, pertaining to constructive mathematics in particular (\intu, \russ, and \bish). Definitions given in this section may slightly differ from similar earlier definitions given for natural spaces, to conform to standard practice in \intu.

\parr

A relatively short discussion of most of these axioms, and comparisons of interrelative strength, can be found under the same names in \WaaArt. More fundamental discussions on intuitionistic axioms are to be found in \VelThe, \VelBora, \VelBorb, and \VelFan. More fundamental discussions on a large number of constructive axioms (and comparisons of interrelative strength) are to be found in the standard works \BeeFou\ and \TroDal. 

\sbsc{Constructive logic is intuitionistic logic}\label{constrlog}
The common practice in constructive mathematics is to use intuitionistic logic. We explain what we mean by describing the meaning of our quantifiers and expressions involving them. 

Before we continue, let us state that certain mathematical notions will be taken as \deff{primitive}, that is: hopefully understood but not defined in terms of other notions. One of these notions is the notion of a \deff{sequence}, for instance a sequence of natural numbers.

\parr

We call $0,\,1,\,2,\,\ldots $ a sequence of natural numbers. There are many other such sequences of course, for instance the sequence of prime numbers $2,\,3,\,5,\,\ldots $
The set of all sequences of natural numbers is often called \NNt. In intuitionism the tradition is however to call this set \sigomt, for reasons which have hopefully become apparent in our previous narrative.

\parr

Cantor's diagonal argument shows that we cannot produce all sequences of natural numbers one after the other. This exhibits an important difference between \sigomt\ and \Nt. For we do have a way to produce all natural numbers, one after the other, even if we are never done with \Nt\ as a whole. But to produce just one element of \sigomt\ is as much work as producing all of \Nt.
  
\parr

Other primitive notions are those of a `set', a `subset' and an `element' of a subset, along with the notion of a `collection of subsets'. We write \emptyyt\ for the empty subset.

Primitive notions relating to our use of quantifiers are
the notion of `method' and `existence'. We say that a mathematical object such as a natural number, a sequence of natural numbers, or a subset with a certain property P \deff{exists}\ if and only if we have a method to construct it.

Then we write, for example: $\dris\alpha\inn\sigom\,[$P$(\alpha)]$ (\deff{there is}\ an \alfat\ in \sigomt\ with property P). If P is a property applicable to sequences of natural numbers, then we can form the subset $\{\alpha\inn\sigom\middz$P$(\alpha)\,\}$. If P is a property applicable to natural numbers and \nnt\ is in \Nt, then we write $n\iz\mu\sinn \,[\,P(s)\,]$ to mean that \nnt\ is the smallest natural number with property P.

\parr   

We abbreviate `if and only if' with `iff'.

\parr

We assume the reader is familiar with the logical symbols $\all,\dris,\driss!,\weddge,\veee,\negg$ and $\rightarrow$. We often use them to abbreviate otherwise lengthy statements. Let \pt\ be a property applicable to the elements of a set or collection \xt. `$\all x\inn X\,[P(x)]$' means: for all \xxt\ in \xt\ we can prove $P(x)$. `$\dris x\inn X\,[P(x)]$' means: there is an \xxt\ in \xt\ such that $P(x)$, as explained above. `$\driss! x\inn X\,[P(x)]$' means: there is an \xxt\ in \xt\ such that $P(x)$ and for all \yyt\ in \xt: if $P(y)$ then $y\iz x$.
        
If on the other hand \pt\ and \qt\ are statements, then `$P\weddge Q$' means: we can prove both \pt\ and \qt. `$P\vee Q$' means: we can choose either \pt\ or \qt, 
and prove the chosen statement. So in fact  `$P\vee Q$' is equal to:
`$\dris s\inn\zo\,[(s\iz 0\weddge P)$ or $(s\iz 1\weddgez Q)]$'. `$P\wisk{\rightarrow} Q$' means: \pt\ implies \qt\ (we can prove \qt\ from \pt). Finally, `$\negg P$' means that we can prove a contradiction from \pt\ (and our axioms).

\parr

We have to distinguish between \pt\ and $\negg\!\negg P$. Of course $\negg\!\negg P$ follows from \pt, but in general the knowledge that $\negg P$ is impossible does not supply us with a proof of \pt. Similarly we distinguish between $\negg\all x\inn X\,[P(x)]$ and $\dris x\inn X\,[\negg P(x)]$. (There are situations in \intu\ in which we can prove both $\negg\all x\inn X\,[P(x)]$ and $\negg\dris x\inn X\,[\negg P(x)]$.).

\sbsc{Functions are Cartesian subsets}\label{deffunct}

Having taken the notion of `set' and `subset' as primitive, we define a function from an apartness space $(V, \aprto)$ to another apartness space $(W, \aprtto)$ as a subset of the cartesian product $(V\timez W, \aprtpi)$  (see def.\,\ref{defprod}) such that:

\be
\item for all $x\inn V$ there is a $y\inn W$ such that $(x,y)\inn f$ 
\item for all $x, v\inn V$ and $y, z\inn W$: if $(x,y)\inn f$ and $(v,z)\inn f$ and $y\aprtto z$ then $x\aprto v$.
\ee

Then for any pair $(x,y)\inn f$ we write: $f(x)\equivv y$ or $f(x)\iz y$.

\parr
 
The constructive interpretation of the quantifiers `for all' and `there is' ensures in our eyes that this definition nicely captures the connotation of methodicity which the word `function' carries. In this book we almost always work with morphisms anyway, but the definition above is strictly speaking necessary to underpin the theorems on representability of continuous functions by morphisms.

\sbsc{Some intuitionistic definitions}\label{basdefintu}
To define the relevant intuitionistic axioms of (continuous) choice we need a number of straightforward definitions, which closely resemble our earlier definitions regarding natural spaces. The reader should take the definitions below as intuitionistic parallels.

\defi (in \intu)
Let \sigomt\ denote the universal spread of all infinite sequences of natural numbers ($\sigom\iz\N\,^{\N}$). Write \sigombgt\ for the set of finite sequences of natural numbers (often written like this: $\sigombg\iz\N^*$). For \alfat\ in \sigomt\ we write $\alfstr(n)$ for the finite sequence $\alpha(0),\ldots,\alpha(n\minuz 1)$ formed by the first \nnt\ values of \alfat. Then $\alfstr(n)\inn\sigombg$, and vice versa $\sigombg\iz\{\alfstr(n)\midd\alpha\inn\sigom,\,\ninn\}$.
A subset \bt\ of \sigombgt\ is called \deff{decidable} iff for all $a\inn\sigombg$ we have a finite decision procedure to determine whether $a\inn B$ or $a\notinn B$. A subset \bt\ of \sigombgt\ is a \emph{bar}\ on a subset \at\ of \sigomt\ iff $\all\alpha\inn A \driss\ninn \,[\,\alfstr(n)\inn B\,]$, and a \emph{thin bar}\ iff $\all\alpha\inn A \driss!\ninn \,[\,\alfstr(n)\inn B\,]$

\parr
Now let \aat\ be in \sigombgt, then \aat\ is a finite sequence of natural numbers. We write $lg(a)$ for the length of this finite sequence. So if $a\iz a_0,\,\ldots,\,a_{n-1}$ then $lg(a)\iz n$. There is a sequence of length $0$, namely the \deff{empty sequence}\ denoted by $\lcod\rcod$. 
For $i\smlr lg(a)$ we then write $a_i$ for the $i^{\,\mbox{\scriptsize th}\!}$ element of this finite sequence. If $a\iz  a_0, a_1,\ldots, a_{lg(a)-1}$ and $b\iz b_0,b_1,\ldots,b_{lg(b)-1} $ are in \sigombgt\ then we write $a\starr b$ for the \deff{concatenation}\ $ a_0,a_1,\ldots,a_{lg(a)-1},b_0,b_1,\ldots,b_{lg(b)-1}$ of \aat\ and \bbt. We write $a\bgeq b$ iff there is a \cct\ in \sigombgt\ such that $b\iz a\starr c$, and we write $a\bgneq b$ iff in addition $lg(b)\bygr lg(a)$.

\parr
A function \ft\ from \sigomt\ to \Nt\ is called a \deff{spread-function}\ iff there is a function \gt\ from \sigombgt\ to \Nt\ such that for each \alfat\ in \sigomt: $\driss!\ninn\,[g(\alfstr(n))\bygr 0]$ and moreover for all \ninnt: $g(\alfstr(n))\bygr 0\rightarrow f(\alpha)\iz g(\alfstr(n))\mino$.%
\footnote{Notice that $\{a\inn\sigombgf\midd g(a)\bygr 0\}$ is a decidable thin bar. This shows that the concept of spread-function is inherently the same as the concept of a decidable (thin) bar.}\ 
More generally a function \ft\ from \sigomt\ to \sigomt\ is called a \deff{spread-function}\ iff there is a function \gt\ from \sigombgt\ to \sigombgt\ such that for each \alfat\ in \sigomt\ and \ninnt\ there is an \minnt\ such that: $\overline{f(\alpha)}(n)\iz g(\alfstr(m))$, and moreover $g(\alfstr(n))\bgeq g(\alfstr(n\pluz 1))$.
\edefi

\rem
Spread-functions from \sigomt\ to \sigomt\ correspond one-on-one to natural morphisms from \Bnatt\ to \Bnatt. We have declined in this monograph to define natural `morphisms' from \sigomt\ to \Nt, but this is easily remedied.
\erem

\sbsc{Axioms of continuous choice in \intu}\label{axiomscont}
The fundamental intuitionistic axiom of continuous choice \acoo\ can now be formulated as follows:

\xiom{\acoo} Let \at\ be a subset of $\sigom\timez\sigom$ such that:
\stArtabb
$\all\alpha\inn\sigom \driss\beta\inn\sigom \,[\,(\alpha,\beta)\inn A\,]$\etab
Then there is a spread-function \gamt\ from \sigomt\ to \sigomt\ such that for each \alfat\ in \sigomt:
$(\alpha,\gamma(\alpha))$ is in \at. We say that \gamt\ \deff{fulfills}\ \stAr.
\exiom
\parr
We formulate four weaker versions of this axiom: \acoz, \CP, \aczo, and  \aczz. The last two are simple axioms of countable choice, whereas \acoz\ is still an axiom of continuous choice, also known as `Brouwer's principle for numbers'. \acoz\ implies the so-called continuity principle \CP. 
We do not defend the axioms here since they are broadly discussed in the literature (see e.\,g.\ \KleVes, \GieSwaVel, \VelThe\ and \TroDal). We begin with the weaker axioms dealing with continuous choice:

\xiom{\acoz} Let \at\ be a subset of $\sigom\timez\N$ such that:
\stArstArtabb
          $\all\alpha\inn\sigom \driss\ninn \,[\,(\alpha,n)\inn A\,]$
\etab
Then there is a spread-function \gamt\ from \sigomt\ to \Nt\ such that for each \alfat\ in \sigomt:
$(\alpha,\gamma(\alpha))$ is in \at. We say that \gamt\ \deff{fulfills}\ \stArstAr.
\exiom

\xiom{\CP} Let \at\ be a subset of $\sigom\timez\N$ such that:
        \btab
          $\all\alpha\inn\sigom \driss\ninn \,[\,(\alpha,n)\inn A\,]$\etab
         Then:
          $\all\alpha\inn\sigom \driss\ninn \driss\minn \all\beta\inn\sigom\
            [\,\alfstr(m)\iz\betastr(m)\rightarrow (\beta,n)\inn A\,]$.
    \exiom

\sbsc{Axioms of countable choice in \bish}\label{axiomscount}
We present two simple axioms of \deff{countable choice}\ in decreasing order of strength:
\xiom{\aczo} Let \at\ be a subset of $\N\timez\sigom$ such that:
\Astabb
          $\all\ninn \driss\alpha\inn\sigom \,[\,(n,\alpha)\inn A\,]$
\etab
        Then there is a function \hte\ from \Nt\ to \sigomt\ such that for each \ninnt:
        $(n,h(n))$ is in \at. We say that \hte\ \deff{fulfills}\ \azt.
\exiom

\xiom{\aczz} Let \at\ be a subset of $\N\timez\N$ such that:
\AstAstabb
          $\all\ninn \driss\minn \,[\,(n,m)\inn A\,]$
\etab
        Then there is a function \hte\ from \Nt\ to \Nt\ such that for each \ninnt:
        $(n,h(n))$ is in \at. We say that \hte\ \deff{fulfills}\ \aztazt.
\exiom

\sbsc{Axioms of dependent choice in \bish}\label{dependentchoice}%
Likewise we present two axioms of \deff{dependent choice}\ in decreasing order of strength. For an intuitionistic justification of these axioms we refer the reader to \WaaThe.
     
\xiom{\dco} Let \delt\ be in \sigomt, and let \at\ be a subset of \sigomt. 
            Suppose \rt\ is a subset of $A\timez A$ such that:
        \btab
          $\delta\inn A\weddgez\all \alpha\inn A\  \dris\beta\inn A \,[\,(\alpha,\beta)\inn R\,]$\etab
        Then there is a sequence \gamnninnt\  of elements of \sigomt\ such that $\gamma_{0}\iz\delta$ and for each \ninnt:
        $(\gamma_{n},\gamma_{n+1})$ is in \rt. 
     
     \exiom
\xiom{\dcz} Let \sinnt, and let \at\ be a subset of \Nt. Suppose \rt\ is 
            a subset of $A\timez A$ such that:
        \btab
          $s\inn A\weddgez\all n\inn A\  \dris m\inn A \,[\,(n,m)\inn R\,]$\etab
        Then there is an \alfat\ in \sigomt\ such that $\alpha(0)\iz s$ and for each \ninnt:
        $(\alpha(n),\alpha(n\pluz 1))$ is in \rt. 
        
     \exiom

\sbsc{Axiom of extensionality}\label{axiomext}
The axiom of extensionality states that we do not distinguish between infinite sequences which are termwise identical, even though they may have different descriptions/definitions.

\xiom{\bol{Ext}} Let $\alpha, \beta\inn\sigom$ such that $\all\ninn[\alpha(n)\iz\beta(n)]$. Then $\alpha\iz\beta$.
\exiom

\parr
This axiom is often adopted rather silently (like we do also). In \russ, the different algorithms describing infinite sequences play an important role, but still one wishes to see infinite sequences themselves as equal when they are termwise identical.

\sbsc{Bar induction, Brouwer's Thesis and the Fan Theorem}\label{defbarinduction}
To phrase the principle of Bar Induction for Decidable bars (\BID) we need:

\defi
A subset \at\ of \sigombgt\ is called \deff{downwards inductive}\footnote{We must distinguish from the already defined notion `inductive'.}\ iff for all \aat\ in \sigombgt: $\all\ninn \,[a\starr n\inn A]\rightarrow a\inn A$. 
\edefi

\xiom{\BID} Let \bt\ be a decidable bar on \sigomt. Suppose \at\ is a downwards inductive subset of \sigombgt\ such that 
$B\subseteqq A$. Then the empty sequence $\lcod\rcod$ (of length $0$) is in \at.
\exiom 

\rem  In classical mathematics \BID\ can be derived from the principle of the excluded middle. The above version of the bar theorem is therefore classically true. In \WaaArt\   \BID\ is derived from \BT\ (see \ref{brothecom}) which also holds both in \class\ and \intu.
\erem

\parr

One of the results following from \BID\ is the axiom known as the fan theorem \FT. We need a preliminary definition.

\defi
Let \sigtot\ be the binary fan ($\{0,1\}^{\N}$). Write \sigtobgt\ for $\{0,1\}^*$, the set of finite sequences of elements of $\{0,1\}$. Then $\sigtobg\iz\{\alfstr(n)\midd\alpha\inn\sigto,\,\ninn\}$.
A subset \bt\ of \sigtobgt\ is a \deff{bar}\ on \sigtot\ iff $\all\alpha\inn\sigto \driss\ninn\ [\,\alfstr(n)\inn B\,]$.
\edefi

\xiom{\FT} If \bt\ is a decidable bar on \sigtot, then \bt\ contains a finite bar on \sigtot\ (in other words: then $\dris\ninn \all\alpha\inn\sigto \driss m\smlr n \,[\,\alfstr(m)\inn B\,]$).
\exiom

\sbsc{Basic axioms in \russ: Church's Thesis}\label{axiomsrussct}
The basic axiom in \russ\ is of course Church's Thesis: `every sequence of natural numbers is given by a recursive rule' (many results in \russ\ already follow from the weaker statement: `the set of partial functions from \Nt\ to \Nt\ is countable'). A partial recursive function \alfat\ from \Nt\ to \Nt\ will usually be denoted by something like `$\phi_e$' where the natural number \eet\ is the \deff{recursive index}\ of \alfat. This recursive index is nothing but the encoding of the finite algorithm which for each \ninnt\ tries to compute $\alpha(n)$. There is a decidable subset $I(\N,\N)$ of \Nt\ such that each \eet\ in $I(\N,\N)$ is a properly formed recursive index of a partial recursive function from \Nt\ to \Nt, and vice versa for each partial recursive function \alfat\ from \Nt\ to \Nt\ there is an \eet\ in $I(\N,\N)$ such that $\all\ninn \,[\,\alpha(n)\isif\phi_e(n)\,]$, where `$\isif$' stands for: `equal if one of the algorithms terminates, given the input'.  
\parr
It turns out one can canonically encode each finite recursive computation as a natural number, see \KleMet. This is the basis of Kleene's decidable $T$-predicate on triples of natural numbers $(e,n,k)$, given by:
\parr
$\displaystyle T(e,n,k)\iff  e$ is a recursive index and $k$ is the canonical encoding of the computation of $\displaystyle\phi_e(n)$.
\parr
If $T(e,n,k)$, then the algorithm $\phi_e$ terminates on the input \nnt. But we are mostly interested in the result of the computation \kt, and in its length (the number of canonical subcomputations leading to the result). Both can be canonically derived from \kt\ of course, using recursive functions $\outc$ and $\lgth$. So if $T(e,n,k)$, then $\phi_e(n)\iz\outc(k)$ and the length of the computation \kt\ equals $\lgth(k)$.
\parr
In this terminology we formulate the axiom Church's Thesis\footnote{Originally `Church's Thesis' stands for the idea that any `mechanically' obtainable sequence must be computable  by a Turing machine.}\ thus:
\xiom{\CT} $\all\alpha\inn\sigom \driss e\inn I(\N,\N) \all\ninn \,[\,\alpha(n)\iz\phi_e(n)\,]$
\exiom
\parr
If $\allz\ninn\,[\alpha(n)\iz\phi_e(n)]$ for $\alpha\inn\sigom$ and $e\inn I(\N,\N)$, then in particular we have: $\allz_{\!}\ninn_{\,}\drisz_{\!}k\inn\N_{\,}[T(e,n,k)]$. \hspace*{-0.21pt}The set $\TOT\iz\{e\inn I(\N,\N)\midd \allz_{\!}\ninn_{\,}\drisz_{\!} k\inn\N[T(e,n,k)]_{\!}\}$ therefore plays an important role in \russ. 
\parr
The combination of \CT\ with \aczz\ is equivalent to an axiom known as \CTz, which forms a connection between \CT\ and choice axioms:

\xiom{\CTz} Let \at\ be a subset of $\N\timez\N$ such that:
\btab
          $\all\ninn \driss\minn \,[\,(n,m)\inn A\,]$
\etab
        Then there is a \emph{recursive}\ function \hte\ from \Nt\ to \Nt\ such that for each \ninnt:
        $(n,h(n))$ is in \at.
\exiom
\parr
From \CTz\ we can derive a more complex choice axiom \CTzo, which plays a part in the defense of \CToo\ in \WaaArt:

\xiom{\CTzo} Let $A, B$ be subsets of $\N\timez\N$, where \bt\ is \emph{decidable}, such that:
\btab
          $\all\ninn \,[\dris y\inn\N\,[(n,y)\inn B]\rightarrow\driss\minn \,[(n,m)\inn A]\,]$
\etab
        Then there is a \emph{partial recursive}\ function \hte\ from \Nt\ to \Nt\ such that for each \ninnt: if $\dris y\inn\N\,[(n,y)\inn B]$ then $h(n)$ is defined and
        $(n,h(n))$ is in \at.
\exiom

\parr
\CTzo\ follows from \CTz. \CTzo\ is the first step to an even broader choice axiom known as `Extended Church's Thesis' (\ECT), which is widely accepted in \russ. But the phrasing of \ECT\ and its defense (at least in \TroDal) are in logical terms and do not appeal to the author. We present a simpler version \CToo\ for which an intuitive defense is given in \WaaArt.\footnote{We try to adhere to the principle that mathematics needs clear axioms which represent our mathematical intuition (not merely serve mathematical convenience). Some technical axioms however have practical advantages for comparison purposes.} 

\xiom{\CToo}
Let \at\ be a subset of $\N\timez\N$ such that
\btab 
$\all n\inn\TOT \driss\minn \,[\,(n,m)\inn A\,]$. 
\etab
Then there is a partial recursive function \hte\ from \Nt\ to \Nt\ such that for all $n\inn\TOT$ we have: $(n, h(n))\inn A$.
\exiom

\sbsc{Markov's Principle}\label{markov}
A second important axiom in \russ\ is called Markov's Principle: `if it is impossible that a total recursive function \alfat\ does not achieve the value $1$ for some \nnt\ in \Nt, then there is an \nnt\ in \Nt\ with $\alpha(n)\iz 1$'. Formally:

\xiom{\MP} Let \alfat\ be in \sigomt\ such that $\negg\negg\dris\ninn \,[\,\alpha(n)\iz 1\,]$. Then $\dris\ninn \,[\,\alpha(n)\iz 1\,]$.
\exiom

\parr

In \WaaArt\ a plea is given to consider adopting \MP\ in intuitionism as well. Joan Moschovakis (\MosMP) holds a similar view.

\sbsc{Axiom of induction}\label{axiomind}
We present the principle of induction as an axiom.
\xiom{\bol{Ind}} Let \at\ be a subset of \Nt\ such that $0\inn A$ and for all
          \ninnt: $n\inn A$ implies $n\pluz 1\inn A$. Then $A\iz\N$, that is:
          $n\inn A$ for all \ninnt.
\exiom

\sbsc{Axiom of decidable-bar descent for \intu, \russ\ and \class}\label{axiombdd}
We phrase a Lindel\"{o}f-type axiom which holds in \class, \intu\ and \russ.

\defi Let \bt\ and \ct\ be two bars on \sigomt, then \bt\ \deff{descends}\ from \ct\ iff for all \cct\ in \ct\ there is a \bbt\ in \bt\ such that $b\bgeq c$. 
\edefi

\xiom{\BDD}
Every bar on \sigomt\ descends from a decidable bar on \sigomt.
\exiom

\parr

The proof of \BDD\ from \BT\ is immediate (\class\ and \intu), whereas the proof of \BDD\ from \CToo\ in \russ\ is given in \WaaArt, using a result from \IshMar. In \WaaArt, \BT\ is shown to be equivalent to the combination of \BID\ and \BDD. \BDD\ is easily seen to be equivalent to:

\xiom{\BDD$^*$}
Every bar on \sigomt\ descends from a decidable thin bar on \sigomt.
\exiom

\sbsc{More axioms}\label{moreaxioms}
More relevant axioms can be found in the literature, notably \TroDal, \BeeFou, \IshRev\ and \VelFan. In \WaaArt, apart from the above some other interesting axioms and axiomatics relating to constructive mathematics are discussed as well.

\sectionnb{Additional remarks}\label{addrem}
\sbsc{Containment and refinement}\label{containrefine}
Notice that in the case of dots being rational intervals, the refinement relation can be defined completely in terms of a natural containment relation $\subseteq$ derived from the pre-apartness relation as follows: $a \subseteq b$ iff $c \aprt b$ implies $c \aprt a$ for all $c \inn V$. In many of our interest spaces, a similar approach yields a decidable containment relation. But in general this approach involves checking an infinite number of conditions, which leaves the so-defined containment relation non-decidable and therefore unwieldy for practical purposes.

\sbsc{Classical treatment of equivalence}\label{classequiv}
The usual classical approach is to work with equivalence classes as the resulting points (a real number usually is defined as the equivalence class of a Cauchy-sequence of rational numbers). In practice this is cumbersome, since all computations on equivalence classes require working with the representatives of these classes. It is therefore more efficient to work with the original sequences and the original apartness relation directly. For theorists who so desire, the translation to equivalence classes is simple, since all definitions will respect the apartness-induced equivalence relation. 
In topological terms our approach means a both practical and foundational way of dealing with a quotient space of Baire space with respect to a  $\Pi_0^1$-equivalence relation.

\sbsc{Details of proving proposition \ref{basicneighbor}}\label{basneiprf}
In the proof of proposition \ref{basicneighbor}, the isomorphisms $g, h$ between \Vnatt\ and \Wnattot\ are \leqc-morphisms. This is without loss of generality, since by \ref{morfcomp}\ we can always lift arbitrary isomorphisms $g, h$ to \leqc-isomorphisms $g', h'$ between \Vpathnatt\ and \Wpathnattot, where the latter is also basic-open trivially. Thus establishing the proposition for  \Vpathnatt\ with $f'\iz h'\!\circc g'$, we can use the \pthh-automorphism $f\iz\idstr\!\circc f'$ to establish the proposition for \Vnatt.

\sbsc{Lazy convergence and isolated points}\label{disclazycon}
A main theme in defining natural spaces and morphisms is `lazy convergence'. Points may themselves choose, so to say, when a next real step in the refinement takes place. However, if one is not careful this leads to some rather unexpected issues with isolated points (points which as a one-point set are open in \Topaprt). To avoid these issues, we have sharpened the definition of `trail space' and `(in)finite product' given in the first edition of this monograph.

\sbsc{Spaces which cannot be represented as a natural space}\label{notnat}
Our discussion in \ref{Funnotnat}\ shows that the `space' of all Baire morphisms cannot be represented as a natural space. This is a frequently occurring theme for function spaces. It would seem to be highly relevant to study various ways in which we can work with function spaces which do not allow a representation as a natural space. It also seems highly worthwhile to find natural representations for various classes of function spaces. We give one example (already discussed by Brouwer) in the examples' section \ref{contfunasspraid}\ above.

\parr

One way of dealing with a  function space $F$ which cannot be represented as a natural space, is to construct a natural space of which $F$ forms a topological subspace. This is often not really satisfying, but still allows to use the functions as `separably countable' objects.
But we leave this for further research, and what has already been done in the constructive literature on function spaces. The issue is discussed also in chapter zero of \WaaThe, where some (most likely not the best) examples and results are given. 

\sbsc{Definition of spreads and spraids}\label{sprddefcom}
In the definition (\ref{spraids}) of `spread' we exact that each infinite \precc-trail defines a point. This is primarily to obtain a precise match with Brouwer's intuitionistic spreads. Most things seem to work fine without the condition, and we obtain the same inductive spraids with or without. It is however quite instructive to see the effect of this condition on the (in)finite-product definition. For example, if we take the simple finite product $\Pi_{i\leq 1}\sigR$ (which represents $\R^2$, see def.\,\ref{defprod}) then the infinite $\precc_{\!_{\Pi}}$-trail $(([0,\twominn],\maxdotR)_{n\in\N})$ does not define a point. So with the condition, the simple (in)finite-product definition almost never yields a spread even when its constituents are all spreads. Without the condition, the simple (in)finite-product of spreads is again a spread, but now the problem is transferred to the inductiveness of the spreads involved (again use $\Pi_{i\leq 1}\sigR$). 

\parr

For elegance we wish to retain the simple infinite-product for natural spaces in general. This means that for spreads $((\Vcal_n, \Topaprtn))_{n\in\N}$ we need the sharper finite products \PisigVnt\ and infinite product \PisigVNt\ defined in \ref{defprod}. 

\parr

Matching Brouwer's spreads precisely is important to us so we adopt the condition. The above is an illustration of how deeply linked most definitions in this monograph are.

\sbsc{Definition of inductive open covers from \ref{indcov}}\label{addremindcov}
The definition of pointwise inductive open covers in \ref{indcov}\ is trickier than one might expect at first glance (and our first try in the previous edition was inadequate). The approach in \ref{indcov}\ allows the requisite induction to be limited to a subspraid \Wnatt, which is far less restrictive than requiring induction on all of \Vnatt.

\parr

A good counterexample to consider is given by the Kleene Tree $K_{\rm bar}$. If we put $\Ucal\iz\{\hattr{k}\midd k\inn K_{\rm bar}\}$, then we see that \Ucalt\ is an open cover of Cantor space in \russ. But \Ucalt\ is not an inductive open cover of Cantor space in \russ, since $K_{\rm bar}$ is not an inductive cover of \Cnatt.

\parr

We do not delve into this definition since it is little used in this monograph, but one should study pointwise inductive open covers in more detail. This bears directly on the important question of how to `inductivize' results from \BisBri\ and \BriVit\ to natural topology.

\sbsc{Star-finitary metrization in \intu}\label{starfinmetcom}
Our star-finitary metrization theorem in \ref{starfinmet}\ is also intended as an example how intuitionistic results can be translated to our setting by inductivizing the definitions (as in chapter three).

\parr

However, the translation from \WaaThe\ of the corresponding intuitionistic star-finitary metrization theorem was complicated by the discovery that the proof offered in \WaaThe\ contains an error, specifically the proof of lemma 2.4.4. which corresponds to part of our Urysohn function lemma (\bol{E}). 

\parr

This error can be corrected by looking at (genetic) bars which represent the condition $A\aprt B$, as our proof here shows, but a straight translation was rendered impossible. This is the main reason that our proof here is more involved and longer than the proof in \WaaThe.

\sectionnb{Bibliography and further reading references}\label{bibliography}

\sbsc{Further reading, research and researchers}
The bibliography given below contains the necessary research references for reading the monograph. It is also intended as a quick (therefore quite incomplete) overview of current research and researchers related to constructive topology. With current internet access to scientific publications, one should have no difficulty in finding other relevant research and researchers. Suggestions for improvement of the bibliography are welcome, through the website of the author (www.fwaaldijk.nl/mathematics.html\,, where his own few publications can also be found online).

\sbsc{Bibliography}

\blit{\TroDal\hspace*{.5cm}}

\bitem{[AczCur2010]}{P. Aczel and G. Curi}
             {On the T1 axiom and other separation properties in constructive point-free and point-set                topology}
             {Annals of Pure and Applied Logic vol. 161, iss. 4, pp. 560-569}{2010}{\AczCur}

\bitem{[Bau2006]}{A. Bauer}
	{The Kleene Tree}
	{preprint, available from www.andrejbauer.com}{2006}{\BauKle}

\bitem{[BauKav2009]}{A. Bauer and I. Kavkler}
             {A constructive theory of continuous domains suitable for implementation}
{Annals of Pure and Applied Logic vol. 159, iss.1-3, pp. 251-267}{2009}{\BauKav}

\bitem{[Bee1985]}{M. Beeson}
           {Foundations of Constructive Mathematics}
           {Sprin\-ger-Verlag, Berlin Heidelberg}{1985}{\BeeFou}

\bitem{[Bis1967]}{E. Bishop}
           {Foundations of Constructive Analysis}
           {McGraw-Hill, New York}{1967}{\BisFou}

\bitem{[Bis{\&}Bri1985]}{E. Bishop and D.S. Bridges}
                      {Constructive Analysis}
                      {Sprin\-ger-Verlag, Berlin Heidelberg} {1985}{\BisBri}

\bitem{[Bis1985]}{E. Bishop}
           {Schizophrenia in Contemporary Mathematics}
           {Erret Bishop: Reflections on Him and His Research (ed. Rosenblatt), pp. 1-32; AMS, Providence, Rhode Island}{1985}{\BisRos}

\bitem{[Bri1979]}{D.S. Bridges}
 	{Constructive Functional Analysis}
                        {Pitman, London} {1979}{\BriCon}

\bitem{[Bri{\&}Ric1987]}{D.S. Bridges and F. Richman}
	{Varieties of constructive mathematics}
	{London Math. Soc. Lecture Notes no. 93, Cambridge University Press}{1987}{\BriRic}

\bitem{[Bri{\&}V\^{i}\c{t}2006]}{D.S. Bridges and L.S. V\^{i}\c{t}\u{a}}
	{Techniques of constructive analysis}
	{Universitext, Springer Science+Business Media}{2006}{\BriVit}

\bitem{[Bri{\&}V\^{i}\c{t}2009]}{D.S. Bridges and L.S. V\^{i}\c{t}\u{a}}
	{Journey into Apartness Space}
	{in: Logicism, Intuitionism, and Formalism - what has become of them?, Springer Netherlands, pp. 167-187}{2009}{\BriVitto}

\bitem{[Bri{\&}V\^{i}\c{t}2011]}{D.S. Bridges and L.S. V\^{i}\c{t}\u{a}}
	{Apartness and Uniformity}
	{Springer-Verlag Berlin Heidelberg}{2011}{\BriVitthr}

\bitem{[Bro1907]}{L.E.J. Brouwer}
              {Over de Grondslagen der Wiskunde}
              {PhD thesis, Universiteit van Amsterdam}{1907}{\BroThe}

\bitem{[Bro1922]}{L.E.J. Brouwer}
              {Besitzt jede reelle Zahl eine Dezimalbruch-Entwickelung?}
              {Math. Annalen 83, 201-210}{1922}{\BroDez}

\bitem{[Bro1975]}{L.E.J. Brouwer}
              {Collected works}
	{(vol. I, II) North-Holland, Amsterdam}{1975}{\BroCol}

\bitem{[Coq1996]}{T. Coquand}
              {Formal topology with posets}
	{preprint available from www.cse.chalmers.se/~coquand/alt.ps}{1996}{\CoqFor}

\bitem{[Col{\&}Eva2008]}{H. Collins and R. Evans}
	{You cannot be serious! Public Understanding of Technology with special reference                  to `Hawk-Eye'.}
	{Public Understanding of Science, vol. 17, 3}{2008}{\ColEva}

\bitem{[CSSV2003]}{T. Coquand, G. Sambin, J. Smith, S. Valentini}
             {Inductively generated formal topologies}
             {Annals of Pure and Applied Logic, vol. 124, iss. 1-3, pp. 71-106}{2003}{\CSSVind}

\bitem{[vDal1999]}{D. van Dalen}
		  {Mystic, Geometer and Intuitionist --- The Life of L.E.J. Brouwer}
	{vol. 1 The Dawning Revolution, Clarendon Press, Oxford}{1999}{\DalBro}

\bitem{[vDal2003]}{D. van Dalen}
		  {Intuitionism of L.E.J. Brouwer}
	{Clarendon Press, Oxford}{2003}{\DalInt}

\bitem{[DeFin1972]}{B. De Finetti}
             {Probability, Induction and Statistics}
              {Wiley, New York}{1972}{\FinFou}

\bitem{[DeFin1974]}{B. De Finetti}
             {Theory of Probability --- A critical introductory treatment}
              {(vol. I, II) Wiley, New York}{1974}{\FinPro}

\bitem{[Dug1951]}{J. Dugundji}
             {An extension of Tietze's theorem}
              {Pacific Journal of Mathematics 1, pp. 353-367}{1951}{\DugExt}

\bitem{[Ear1986]}{J. Earman}
           {A Primer on Determinism}
           {D. Reidel Publ., Dordrecht}{1986}{\Earman}

\bitem{[FoMforum]}{M. Davis (moderator)}
	{Foundations of Mathematics forum}
	{http://www.cs.nyu.edu/mailman/listinfo/fom/}{1998-2012}{\FomDis}

\bitem{[FouGra1982]}{M.P. Fourman, R.J. Grayson}
             {Formal Spaces}
             {in `The L. E. J. Brouwer Centenary Symposium', North
Holland, Amsterdam, pp. 107-122}{1982}{\FouGra}

\bitem{[Fre1937]}{H. Freudenthal}
              {Zum intuitionistischen Raumbegriff}
               {Compositio Mathematica, vol. 4, pp. 82-111}{1937}{\Freu}

\bitem{[GiSwVe1981]}{W. Gielen, H. de Swart and W.H.M. Veldman}
 	{The Continuum Hypothesis in Intuitionism}
               {Journal of Symbolic Logic 46 (1), 121--136}{1981}{\GieSwaVel}

\bitem{[GNSW2007]}{H. Geuvers, M. Niqui, B. Spitters and F. Wiedijk}
                {Constructive analysis, types and exact real numbers}
               {Mathematical Structures in Computer Science, vol.17, iss. 1, pp 3-36}{2007}{\GNSWcom}

\bitem{[Haw2005]}{S. Hawking}
              {God created the integers}
               {Penguin Books, London (paperback ed. 2006)}{2005}{\HawGod}

\bitem{[Hey1956]}{A. Heyting}
              {Intuitionism --- an Introduction}
               {North-Holland, Amsterdam}{1956 (3rd rev. ed. 1971)}{\HeyInt}

\bitem{[Hil1996]}{T. Hill}
	{A statistical derivation of the significant-digit law}
	{Statistical Science \bol{10}, 354--363}{1996}{\HilSig}

\bitem{[Ish1993]}{H. Ishihara}
	{Markov's Principle, Church's thesis and Lindel\"{o}f's theorem}
	{Indag. Mathem., N.S., 4 (3), 321--325}{1993}{\IshMar}

\bitem{[Ish1994]}{H. Ishihara}
	{A constructive version of Banach's inverse mapping theorem}
	{New Zealand Journal of Mathematics, vol. 23, 71--75}{1994}{\IshBan}

\bitem{[Ish2006]}{H. Ishihara}
	{Reverse Mathematics in Bishop's Constructive Math\-ematics}
	{in: Constructivism: Mathematics, Logic, Philosophy and Linguistics, Philosophia Scientiae, Cahier sp\'ecial 6, pp. 43-59}{2006}{\IshRev}

\bitem{[Joh1987]}{D.M. Johnson}
           {L.E.J. Brouwer's coming of age as a topologist}
           {MAA Studies in Mathematics, vol. 26, Studies in the History of Mathematics (ed. E.R. Phillips), 61--97}{1987}{\Johnso}

\bitem{[Jul{\&}Ric1984]}{W. Julian and F. Richman}
	{A continuous function on \zort\ which is everywhere apart from its infimum}
	{Pacific Journal of Mathematics, vol. 111, iss. 2, 333-340}{1984}{\JulRic}

\bitem{[Kal{\&}Wel2006]}{I. Kalantari and L. Welch}
                     {Larry Specker's theorem, cluster points, and computable quantum functions}
                    {Logic in Tehran, 134-159, Lect. Notes Log., 26, Association for Symbolic Logic}{2006}{\KalWel}

\bitem{[Kle1952]}{S.C. Kleene}
                     {Introduction to metamathematics}
                    {North-Holland, Amsterdam}{1952}{\KleMet}

\bitem{[Kle{\&}\-Ves1965]}{S.C. Kleene and R.E. Vesley}
                     {The Foundations of Intuitionistic Mathematics --- especially in relation to recursive functions}
                    {North-Holland, Amsterdam}{1965}{\KleVes}

\bitem{[Kle1969]}{S.C. Kleene}
                     {Formalized Recursive Functionals and Formalized Realizability}
                    {Memoirs of the American Mathematical Society, vol. 89}{1969}{\KleRea}

\bitem{[Kun{\&}Sch2005]}{D. Kunkl and M. Schr\"{o}der}
                     {Some Examples of Non-Metrizable Spaces Allowing a Simple Type-2 Complexity Theory}
                    {Electronic Notes in Theoretical Computer Science, vol. 120, 111-123}{2005}{\KunSch}

\bitem{[Kus1985]}{B.A. Kushner}
	{Lectures on Constructive Mathematical Analysis}
	{American Mathematical Society, Providence, R.I.}{1985}{\KusLec}

\bitem{[Lap1776]}{P.S. de Laplace}
	{Recherches sur l'int\'egration des \'equations diff\'erentielles aux diff\'erences finies, et sur leur usage dans la th\'eorie des hasards}
	{in: Oeuvres compl\`etes de Laplace, tome 8, this quote: pp. 145, Gauthiers-Villars, Paris (1878-1912)}{1776}{\LapDif}

\bitem{[Lap1814]}{P.S. de Laplace}
	{Essai philosophique sur les probabilit\'{e}s}
	{publ.: Mme. Ve. Courcier, Paris (translated in \HawGod)}{1814}{\LapPhi}

\bitem{[M-L\"{o}f1970]}{P. Martin-L\"{o}f}
              {Notes on Constructive Mathematics}
             {Almqvist {\&} Wiksell, Stockholm}{1970}{\MLof}        

\bitem{[Mic1956]}{E. Michael}
              {Continuous selections (I)}
             {Annals of Mathematics 63, 361-382}{1956}{\MicSel}        

\bitem{[vMil1989]}{J. van Mill}
              {Infinite-dimensional Topology}
             {North-Holland, Amsterdam}{1989}{\Mill}        

\bitem{[MiRiRu1988]}{R. Mines, F. Richman and W. Ruitenburg}
                          {A Course in Constructive Algebra}
                  {Universitext, Springer}{1988}{\MiRiRu} 

\bitem{[Mos1971]}{J.R. Moschovakis}
            {Can there be no nonrecursive sequences?}
            {Journal of Symbolic Logic vol. 36, iss. 2, pp. 309-325}{1971}{\MosRec}

\bitem{[Mos1996]}{J.R. Moschovakis}
            {A classical view of the intuitionistic continuum}
            {Annals of Pure and Applied Logic vol. 81, iss.1-3, pp. 9-24}{1996}{\MosCla}

\bitem{[Mos2002]}{J.R. Moschovakis}
            {Analyzing realizability by Troelstra's methods}
            {Annals of Pure and Applied Logic vol. 114, iss.1-3, pp. 203-225}{2002}{\MosRea}

\bitem{[Mos2012]}{J.R. Moschovakis}
            {Waltzing around Markov's Principle}
            {Talk given in Nijmegen, september}{2012}{\MosMP}

\bitem{[Neg{\&}Ced1996]}{S. Negri and J. Cederquist}
	{A constructive proof of the Heine-Borel covering theorem for formal reals}
	{in: Proceedings TYPES `95, Lecture Notes in Computer Science 1158, Springer Verlag London}{1996}{\NegCed}

\bitem{[Pal2005]}{E. Palmgren}
	{Continuity on the real line and in formal spaces}
	{in: `From Sets and Types to
Topology and Analysis' ed. L. Crosilla and P. Schuster, Oxford University Press}{2005}{\PalCon}

\bitem{[Pal2009]}{E. Palmgren}
	{From intuitionistic to formal topology: some remarks on the foundations of homotopy theory}
	{in: Logicism, Intuitionism, and Formalism - what has become of them?, Springer Netherlands, pp. 237-253}{2009}{\PalInt}

\bitem{[Pal2010]}{E. Palmgren}
	{Open sublocales of localic completions}
	{Journal of Logic \&\ Analysis vol. 2, iss.1, pp. 1-22}{2010}{\PalSub}

\bitem{[Sam2003]}{G. Sambin}
	{Some points in formal topology}
	{Theoretical Computer Science, vol. 305 iss. 1-3, Elsevier}{2003}{\SamFor}

\bitem{[Schu2001]}{P. Schuster, U. Berger, H. Osswald (ed.)}
	{Reuniting the Anti\-podes: constructive and nonstandard views of the continuum}
	{Venice 1999 Symp. proceedings, Kluwer Academic Press}{2001}{\SchuReu}

\bitem{[Schu2005]}{P. Schuster}
	{What is continuity, constructively?}
	{Journal of Universal Computer Science, vol. 11, iss. 12, pp. 2076-2085}{2005}{\SchuCon}

\bitem{[Spi2010]}{B. Spitters}
	{Locatedness and overt sublocales}
	{Annals of Pure and Applied Logic vol. 162, iss.1, pp. 36-54}{2010}{\SpiLoc}

\bitem{[Tay{\&}Bau2009]}{P. Taylor and A. Bauer}
                          {The Dedekind reals in Abstract Stone Duality}
                  {Mathematical Structures in Computer Science
vol. 19, iss. 4}{2009}{\TayBau}

\bitem{[Tro1966]}{A.S. Troelstra}
                          {General Intuitionistic Topology}
                  {PhD thesis, University of Amsterdam}{1966}{\TroThe}

\bitem{[Tro{\&}vDal1988]}{A.S. Troelstra and D. van Dalen}
                          {Constructivism in Mathematics}
                  {(vol. I, II) North-Holland, Amsterdam}{1988}{\TroDal}

\bitem{[Ury1925a]}{P. Urysohn}
                       {\"{U}ber die M\"{a}chtigkeit der Zusammenh\"{a}ngende Mengen}
                         {Mathematische Annalen vol. 94, pp. 262-295}{1925}{\Urya}

\bitem{[Ury1925b]}{P. Urysohn}
                       {Zum Metrisationsproblem}
                         {Mathematische Annalen vol. 94, pp. 309-315}{1925}{\Uryb}

\bitem{[Vel1981]}{W.H.M. Veldman}
              {Investigations in intuitionistic hierarchy theory}
               {PhD thesis, University of Nijmegen}{1981}{\VelThe}

\bitem{[Vel1985]}{W.H.M. Veldman}
              {Intu\"\i tionistische wiskunde}
                {(lecture notes in Dutch) University of Nijmegen}{1985}{\VelLec}

\bitem{[Vel2004]}{W.H.M. Veldman}
              {An intuitionistic proof of Kruskal's theorem}
                {Archive for
Mathematical Logic, vol. 43, pp. 215-264}{2004}{\VelKru}

\bitem{[Vel2008]}{W.H.M. Veldman}
              {The Borel Hierarchy Theorem from Brou\-wer's intuitionistic perspective}
                {Journal of Symbolic Logic, vol. 73, iss. 1, pp. 1-64}{2008}{\VelBora}

\bitem{[Vel2009]}{W.H.M. Veldman}
              {The fine structure of the intuitionistic Borel hierarchy}
               {Review of Symbolic Logic, vol. 2, iss. 1, pp. 30-101}{2009}{\VelBorb}

\bitem{[Vel2011]}{W.H.M. Veldman}
              {Brouwer's Fan Theorem as an axiom and as a
contrast to Kleene's Alternative}
                {Research report, on arxiv.org, arXiv:1106.2738v1}{2011}{\VelFan}

\bitem{[Vel{\&}Waa1996]}{W.H.M. Veldman and F.A. Waaldijk}
	{Some elementary results in intuitionistic model theory}
	{Journal of Symbolic Logic, vol. 62, iss. 3, pp. 745-767}{1996}{\VelWaa}

\bitem{[Waa1996]}{F.A. Waaldijk}
                    {modern intuitionistic topology}
                    {PhD thesis, University of Nijmegen}{1996}{\WaaThe}

\bitem{[Waa2005]}{F.A. Waaldijk}
                    {On the foundations of constructive mathematics --- especially in relation to the                      theory of continuous functions}
                    {Foundations of Science, vol. 3, iss. 10, pp. 249-324}{2005}{\WaaArt}

\bitem{[Wei2000]}{K. Weihrauch}
	{Computable analysis}
	{Springer Verlag, Berlin}{2000}{\WeiCom}

\bitem{[Wol2002]}{S. Wolfram}
	{A New Kind of Science}
	{Wolfram Media, also online at www.wolfram.com}{2002}{\WolNew}

\elit

\end{document}